\documentclass[a4paper, 11pt]{amsart}
\usepackage{amssymb}
\usepackage{amsthm}
\usepackage{fullpage}
\usepackage{graphicx}
\usepackage{subfig}
\usepackage{amsmath,amscd}
\usepackage{xypic}
\usepackage{color}
\usepackage{float}
\usepackage{vmargin}
\usepackage[colorlinks,citecolor=blue,filecolor=black,linkcolor=blue,urlcolor=black]{hyperref}
\usepackage[all,cmtip]{xy}

\newcommand{\Z}{\ensuremath{\mathbb Z}}

\newcommand{\R}{\ensuremath{\mathbb R}}

\theoremstyle{plain}
\newtheorem{thm}{Theorem}[section]
\newtheorem*{thm*}{Theorem}

\newtheorem*{cor*}{Corollary}
\newtheorem*{prop*}{Proposition}

\newtheorem*{lemma*}{Lemma}
\newtheorem*{claim*}{Claim}

\theoremstyle{definition}

\newtheorem*{exmp*}{Example}
\newtheorem*{defn*}{Definition}
\newtheorem*{rem*}{Remark}
\newtheorem*{note*}{Note}

\begin{document}

\title{A four dimensional hyperbolic link complement in a standard $S^4$} 

\author{Hemanth Saratchandran}
\email{hemanth.saratchandran@maths.ox.ac.uk}
\date{\today}

\maketitle 
\parskip=0.2cm
\parindent=0.0cm

\begin{abstract}
Using techniques from the theory of Kirby calculus we give an explicit construction of a four dimensional hyperbolic link complement in a
4-manifold that is diffeomorphic to the standard 4-sphere. 
\end{abstract}

\tableofcontents

\section{Introduction}

In the 1970's W. Thurston initiated a grand study of the topology of 3-manifolds motivated by the realisation that many 3-manifolds 
admitted homogeneous Riemannian metrics. Thurston's insight was that one could use this geometry to study the 
topology of a 3-manifold, the ultimate goal being a unification of the worlds of geometry and topology in dimension three. Two worlds thought
to be completely distinct before Thurston's time. Part of Thurston's study involved a detailed understanding of how various link complements in 3-manifolds
could admit homogeneous geometries. He was able to show that many link complements admitted a hyperbolic geometry thereby throwing the world of hyperbolic
geometry to the forefront (see \cite{thurston_1} chap.3, p.27 and \cite{thurston_2} cor.2.5, p.360). 
This grand vision of Thurston came to be known as the ``geometrisation conjecture'' and was finally
proved by G. Perelman in 2003. 
The geometrisation theorem together with the classification of surfaces leads to a complete understanding of the
worlds of two and three dimensional manifolds. 

The quest for a complete understanding of the world of 4-manifolds, analogous to the two and three dimensional cases, is a futile quest. The main reason
being that, every finitely presented group can arise as the fundamental group of a compact 4-manifold. As there is no algorithm to tell whether two finitely 
presented groups are isomorphic, there is no algorithm to tell if two 4-manifolds have the same fundamental groups. Motivated by Thurston's work on which link complements, in dimension 3, admit hyperbolic geometries we can restrict our attention to the class of hyperbolic 4-manifolds and try to understand which
4-manifolds have link complements that admit a hyperbolic geometry. 
The main problem here is that hyperbolic 4-manifolds are characteristically  
different from hyperbolic 3-manifolds, and the techniques used to show that a 3-manifold admits a link complement that is hyperbolic simply do not have
four dimensional counter parts. Due to this, the quest to understand which 4-manifolds admit link complements that are hyperbolic takes on a very different
feel right from the start. 

In the last decade there have been quite a few constructions of hyperbolic 4-manifolds. One of the simplest constructions was given by
J. Ratcliffe and S. Tschantz in their paper \cite{ratcliffe}. 
Many of the 4-manifolds they construct are non-compact and have ends of the form $E \times [0, \infty)$, with $E$ a closed flat 3-manifold that is a circle bundle over a 
surface.
We can take these ends, chop them off at a certain point, and produce a compact 4-manifold $M_0$ with boundary a certain number of flat closed 3-manifolds 
$E_i$. The region near the boundary is a copy of $E_i \times [0, t]$, for some $t > 0$. These boundary 3-manifolds, being circle bundles, bound a 4-manifold $F_i$, the
associated disk bundle.
We can then glue the manifold $F_i$ to $M_0$ by identifying the boundary components via the identity, we will call this a ``filling'' of $M_0$.

\centerline{\includegraphics[width=8cm, height=8cm]{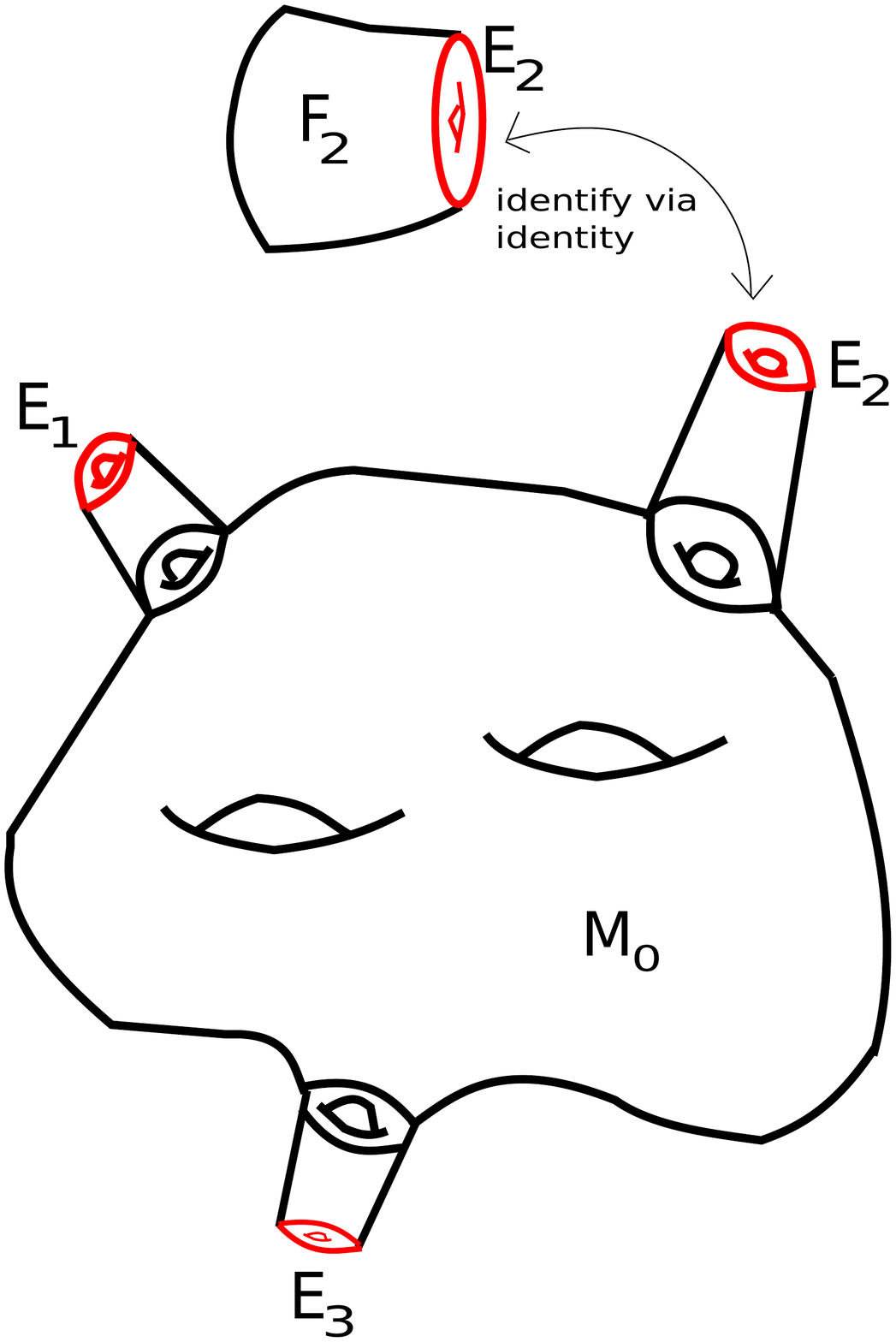}}

It could be the case that the 3-manifold $E$ fibres over a surface in more than one way, hence the filling process will depend on the choice of
$S^1$-fibre.
Carrying out this gluing procedure for each boundary component of $M_0$ produces a closed 4-manifold $\widetilde{M}$. The original hyperbolic 4-manifold
can then be seen to be a codimension two link complement inside of the closed 4-manifold $\widetilde{M}$. At this point one of the most basic questions we can ask is, can we explicitly
identify the manifold $\widetilde{M}$? 
This question is related to the above problem of finding 4-manifolds that have a hyperbolic link complement in that the ability to classify such an $\widetilde{M}$
would then result in an explicit example of a four dimensional hyperbolic link complement. However, one has to be very careful in using the word
``classify''. The main reason being that the
world of 4-manifolds is a truly wild world, there are 4-manifolds out there that do not admit a single smooth structure, others out there
that admit countably many smooth structures, and some even admitting uncountably many smooth structures! This exotic behaviour of 4-manifolds shows the
mathematician that the problem of characterising a 4-manifold up to homeomorphism is a very different problem than to characterise it up to diffeomorphism,
something one does not witness when restricting to manifolds of dimensions two and three. 

In general the identification of the topological/smooth type of the manifold $\widetilde{M}$ can prove to be an impossible task. However, if for example the filled in manifolds
one obtains are simply connected then there is certainly hope. A good way to try and smoothly identify a simply connected 4-manifold 
is to resort to a ``calculus of links'' developed by R. Kirby towards the end of the 70's (see \cite{kirby}). Kirby showed that given a handle decomposition of
a closed 4-manifold, the one and two handles were really what one had to worry about in trying to understand the manifold. The one and two handle
structure of such a manifold can be neatly encoded in a link diagram in $S^3$, which we can view as $\R^3 \cup \{\infty\}$. Therefore, the
one and two-handle structure of such a manifold could be explicitly visualised by a link diagram in $\R^3$, which we now call a Kirby diagram.
Kirby was able to prove that if one applied certain elementary moves to the link diagram, obtaining a new link diagram, the 4-manifold
that corresponded to this new link diagram would be diffeomorphic to the original 4-manifold. In this way Kirby set up a ``calculus'' that one could
appeal to in order to simplify their link diagram but without any compensation being paid on the diffeomorphism type of their closed 4-manifold. 

The aim of this paper is to use the theory of Kirby calculus to construct an explicit example of a four dimensional hyperbolic link complement
in a 4-manifold that is diffeomorphic to the standard 4-sphere.  Our approach begins with the work we started in our paper \cite{sarat}. In that
paper we show how to construct a Kirby diagram for any one of the Ratcliffe-Tschantz hyperbolic 4-manifolds. Using this construction we will show
how to understand a filling on the level of the Kirby diagram. We will then show, using methods from the theory of Kirby calculus, how to reduce
the Kirby diagram of the filling and obtain a closed 4-manifold that has the same Kirby diagram of the standard 4-sphere. This will then tell us that
our original filling is diffeomorphic to the standard 4-sphere, in turn producing a four dimensional hyperbolic link complement in a standard 4-sphere.
The main theorem takes the form:

\begin{thm*}(Theorem \ref{mainthm_1})
There exists a collection $L$ of five linked 2-tori embedded in a smooth 4-manifold $X$ such that $X$ is diffeomorphic to a standard $S^4$, and
$X-L$ admits a finite volume hyperbolic geometry.
\end{thm*}

We would like to mention that D. Ivan$\check{s}$i$\acute{c}$ proves in his paper \cite{ivansic} (see thm.4.3, p.18) that there exists a system of 
five linked tori embedded in a smooth manifold $X$ that is homeomorphic to $S^4$ such that $X-L$ admits a finite volume hyperbolic geometry. The manifold
he uses to construct $X$ is the same manifold we are going to use. However, he does not prove that $X$ is diffeomorphic to $S^4$. Furthermore, 
D. Ivan$\check{s}$i$\acute{c}$, J. Ratcliffe and S. Tschantz in their paper \cite{tschantz}  construct several more examples of hyperbolic
link complements in 4-manifolds that are homeomorphic to $S^4$, they do not prove that any of the fillings they construct are diffeomorphic to $S^4$.
We were unaware of the existence of these two papers during the undertaking of this work,
in fact we were only made aware of these two papers quite recently. In a future paper we will show how to use the methods of this paper and \cite{sarat}
to construct a hyperbolic link complement in 4-manifold that is diffeomorphic to a standard $S^2 \times S^2$.

Unfortunately we have not been able to make this paper as self contained as we would have liked. The main reason being that it would have become
far too long. The constructions and techniques used in this paper heavily depend on those outlined in \cite{sarat} and we recommend the reader look at
that paper to get an idea of what is going on. We also outline, in that paper, the very basic properties of the hyperbolic 24-cell that we will need, and
give some preliminary information on the construction of the Ratcliffe-Tschantz manifolds.

\section*{Acknowledgements}
The author wishes to thank Marc Lackenby for various discussions to do with this work and the several comments/corrections he gave on earlier drafts of this work.
We would also like to thank Andras Juhasz, Panos Papazoglou and John Parker for the comments and corrections they gave on an earlier draft.

\section{Parabolic transformations and the Euclidean structure of a cusp}\label{parabolics}

In this section we are going to explain how to compute parabolic transformations associated to each cusp component of a Ratcliffe-Tschantz hyperbolic
4-manifold. Each of their manifolds has the hyperbolic 24-cell as a fundamental domain.
The 24-cell is a self dual 4-dimensional ideal hyperbolic polyhedron. For the basic construction of the hyperbolic 24-cell that we will use
we refer the reader to section two of \cite{sarat}, for more background information on various properties of the 24-cell we recommend the reader consult
\cite{cox} Ch.4 and \cite{kerckhoff} Sect.3.

Recall that the 24-cell $P$ has twenty four ideal vertices, eight of the form $(\pm 1,0,0,0)$, $(0,\pm 1,0,0)$, $(0,0,\pm 1, 0)$, $(0,0,0,\pm 1)$ and 
sixteen of the form $(\pm 1/2, \pm 1/2, \pm 1/2, \pm 1/2)$. When we apply any group of side pairing transformations defining any one of the
5-cusped Ratcliffe-Tschantz manifolds we find that the ideal vertices split into five equivalence classes. We have four of the form:
\[ \{(1,0,0,0), (-1,0,0,0)\}, \{(0,1,0,0), (0,-1,0,0)\}, \]
\[ \{(0,0,1,0), (0,0,-1,0)\}, \{(0,0,0,1), (0,0,0,-1)\} \]
and one of the form:
\[ \{(\pm 1/2, \pm 1/2, \pm 1/2, \pm 1/2)\}. \]
In order to compute the parabolic isometries associated to these five cusps we need to look at horospherical neighbourhoods about each ideal vertex.
Recall that given any ideal vertex of the form $(\pm 1,0,0,0)$, $(0,\pm 1,0,0)$, $(0,0,\pm 1, 0)$, $(0,0,0,\pm 1)$, such a vertex lies on a given side
$S$ if the centre vector of the sphere defining the side has an equal non-zero entry in the same position as the vertex. For example the ideal vertex 
$(-1,0,0,0)$ lies on the side $S_{(-1,0,0,1)}$ but not on the side $S_{(1,0,0,1)}$ or on the side $S_{(0,1,1,0)}$. Any ideal vertex of the form
$(\pm 1/2, \pm 1/2, \pm 1/2, \pm 1/2)$ lies on a side $S$ if the non-zero entries of the centre of the sphere defining the side have the same sign
as the non-zero entries in the same position of the ideal vertex. For example the ideal vertex $(1/2, 1/2, -1/2, 1/2)$ lies on the side
$S_{(+1,0,-1,0)}$ but not on the side $S_{(-1,+1,0,0)}$.

The above recollection shows us that each ideal vertex lies on precisely six sides. Since each side of the 24-cell $P$ intersects precisely 
four other sides, and does so at right angles, we find that a horospherical neighbourhood of any of the ideal vertices is a cube.
From this fact it follows that for each of the equivalence classes
\[\{(1,0,0,0), (-1,0,0,0)\}, \{(0,1,0,0), (0,-1,0,0)\}, \] 
\[\{(0,0,1,0), (0,0,-1,0)\}, \{(0,0,0,1), (0,0,0,-1)\} \]
the associated boundary cusp will have fundamental domain consisting of two cubes. For example, if we take the vertex class $\{(1,0,0,0), (-1,0,0,0)\}$
then we can think of a fundamental domain for the associated cusp cross section as a horospherical neighbourhood about the ideal vertex $(1,0,0,0)$
together with a horospherical neighbourhood about the ideal vertex $(-1,0,0,0)$. Similarly for the class $\{(\pm 1/2, \pm 1/2, \pm 1/2, \pm 1/2)\}$, we can think
of the boundary cusp cross section as coming from sixteen cubes, each one corresponding to a horospherical neighbourhood about an ideal vertex making up the equivalence class.
In order to understand the structure of these cusp cross sections we need to understand the parabolic isometries associated to them. 
Through the general theory of non-compact hyperbolic manifolds of finite volume we know that the stabiliser subgroup of an ideal vertex induces the
Euclidean structure on the cusp cross section. The way one generally sees this is to observe that if one takes a horospherical neighbourhood about an ideal vertex, then
the stabiliser subgroup acts as a group of Euclidean affine transformations on this neighbourhood. Viewing this neighbourhood as a copy of $\R^3$ one observes that
the quotient via the action of the stabiliser subgroup is a Euclidean 3-manifold, the one corresponding to the cusp cross section
(the reader who is not familiar with this material can consult \cite{benedetti} Thm.D.3.3, p.145).
In order to compute the stabiliser subgroup one can take a horospherical neighbourhood about an ideal vertex and see
how it transforms under the side pairing transformations associated to the sides making up the horospherical neighbourhood. We are going to give explicit examples of
how to do this shortly. Before we do this let us remind the reader that the Euclidean structure associated to each of the five cusps is given
in the Ratcliffe-Tschantz census (see \cite{ratcliffe}) under the column headed \textbf{LT}. A quick glance at the census shows that many of the manifolds in the census have very different
cusp structures, this means that in general one needs to carry out such computations for the particular manifolds they are interested in. That being said, the general
principle of how to compute the cusp structures works for all the Ratcliffe-Tschantz manifolds, and is in fact a standard technique from the theory
of non-compact hyperbolic manifolds of finite volume.

As promised we are going to give a calculation of the cusp structure for a particular manifold using horospherical neighbourhoods as described above. Since we 
have already dealt with aspects of manifold no. 3 we choose to use it again.
For the convenience of the reader we recall its basic construction. The side pairing code for manifold no. 3 is \textbf{1477B8}, decoding this gives the
following side pairing transformations:
\[\xymatrixcolsep{5pc}\xymatrix{ S_{(+1,+1,0,0)}  \ar[r]^a_{k_{_{(-1,+1,+1,+1)}}} & S_{(-1,+1,0,0)} } \hspace{2cm} \xymatrix{S_{(+1,-1,0,0)}  \ar[r]^b_{k_{_{(-1,+1,+1,+1)}}} & S_{(-1,-1,0,0)} } \]

\[\xymatrixcolsep{5pc}\xymatrix{ S_{(+1,0,+1,0)}  \ar[r]^c_{k_{_{(+1,+1,-1,+1)}}} & S_{(+1,0,-1,0)} } \hspace{2cm} \xymatrix{S_{(-1,0,+1,0)}  \ar[r]^d_{k_{_{(+1,+1,-1,+1)}}} & S_{(-1,0,-1,0)} } \]

\[\xymatrixcolsep{5pc} \xymatrix{ S_{(0,+1,+1,0)}  \ar[r]^e_{k_{_{(-1,-1,-1,+1)}}} & S_{(0,-1,-1,0)} } \hspace{2cm} \xymatrix{S_{(0,+1,-1,0)}  \ar[r]^f_{k_{_{(-1,-1,-1,+1)}}} & S_{(0,-1,+1,0)} } \]

\[\xymatrixcolsep{5pc} \xymatrix{ S_{(+1,0,0,+1)}  \ar[r]^g_{k_{_{(-1,-1,-1,+1)}}} & S_{(-1,0,0,+1)} } \hspace{2cm} \xymatrix{S_{(+1,0,0,-1)}  \ar[r]^h_{k_{_{(-1,-1,-1,+1)}}} & S_{(-1,0,0,-1)} } \]

\[\xymatrixcolsep{5pc} \xymatrix{ S_{(0,+1,0,+1)}  \ar[r]^i_{k_{_{(-1,-1,+1,-1)}}} & S_{(0,-1,0,-1)} } \hspace{2cm} \xymatrix{S_{(0,+1,0,-1)}  \ar[r]^j_{k_{_{(-1,-1,+1,-1)}}} & S_{(0,-1,0,+1)} } \]

\[\xymatrixcolsep{5pc} \xymatrix{ S_{(0,0,+1,+1)}  \ar[r]^k_{k_{_{(+1,+1,+1,-1)}}} & S_{(0,0,+1,-1)} } \hspace{2cm} \xymatrix{S_{(0,0,-1,+1)}  \ar[r]^l_{k_{_{(+1,+1,+1,-1)}}} & S_{(0,0,-1,-1)} }. \]

The labelling of the sides is given in the following table: \\

\begin{tabular}{|l|l|l||l|l|l|}
 $A$ &  $S_{(+1,+1,0,0)}$ & $(\frac{1}{\sqrt{2}}, \frac{1}{\sqrt{2}}, 0)$  &  $A'$ & $S_{(-1,+1,0,0)}$ & $(\frac{-1}{\sqrt{2}}, \frac{1}{\sqrt{2}}, 0)$  \\
 $B$ &  $S_{(+1,-1,0,0)}$ & $(\frac{1}{\sqrt{2}}, \frac{-1}{\sqrt{2}}, 0)$ &  $B'$ & $S_{(-1,-1,0,0)}$  & $(\frac{-1}{\sqrt{2}}, \frac{-1}{\sqrt{2}}, 0)$ \\
 $C$ &  $S_{(+1,0,+1,0)}$ &  $(\frac{1}{\sqrt{2}}, 0, \frac{1}{\sqrt{2}})$ &  $C'$ & $S_{(+1,0,-1,0)}$ & $(\frac{1}{\sqrt{2}}, 0, \frac{-1}{\sqrt{2}})$ \\
 $D$ &  $S_{(-1,0,+1,0)}$ &  $(\frac{-1}{\sqrt{2}}, 0, \frac{1}{\sqrt{2}})$ &  $D'$ & $S_{(-1,0,-1,0)}$ & $(\frac{-1}{\sqrt{2}}, 0, \frac{-1}{\sqrt{2}})$ \\ 
 $E$ &  $S_{(0,+1,+1,0)}$ &  $(0, \frac{1}{\sqrt{2}} ,\frac{1}{\sqrt{2}})$ &  $E'$ & $S_{(0,-1,-1,0)}$ &  $(0, \frac{-1}{\sqrt{2}} ,\frac{-1}{\sqrt{2}})$ \\ 
 $F$ &  $S_{(0,+1,-1,0)}$ &   $(0, \frac{1}{\sqrt{2}} ,\frac{-1}{\sqrt{2}})$ &  $F'$ & $S_{(0,-1,+1,0)}$ &  $(0, \frac{-1}{\sqrt{2}} ,\frac{1}{\sqrt{2}})$ \\ 
 $G$ &  $S_{(+1,0,0,+1)}$ &  $(1 + \sqrt{2}, 0, 0)$ &  $G'$ & $S_{(-1,0,0,-1)}$ & $(1 - \sqrt{2}, 0, 0)$ \\ 
 $H$ &  $S_{(+1,0,0,-1)}$ &  $(-1 + \sqrt{2}, 0, 0)$ &  $H'$ & $S_{(-1,0,0,+1)}$ & $(-1 - \sqrt{2}, 0, 0)$ \\
 $I$ &  $S_{(0,+1,0,+1)}$ &  $(0, 1 + \sqrt{2}, 0)$ &  $I'$ & $S_{(0,-1,0,+1)}$ &  $(0, -1 - \sqrt{2}, 0)$ \\
 $J$ &  $S_{(0,+1,0,-1)}$ &   $(0, -1 + \sqrt{2}, 0)$ &  $J'$ & $S_{(0,-1,0,-1)}$ &  $(0, 1 - \sqrt{2}, 0)$  \\
 $K$ &  $S_{(0,0,+1,+1)}$ &   $(0, 0, 1 + \sqrt{2})$ &  $K'$ & $S_{(0,0,+1,-1)}$ &  $(0, 0, -1 + \sqrt{2})$   \\
 $L$ &  $S_{(0,0,-1,+1)}$ &   $(0, 0, -1 - \sqrt{2})$  &  $L'$ & $S_{(0,0,-1,-1)}$ &  $(0, 0, 1 - \sqrt{2})$  
\end{tabular}  \\ \\

Consider the equivalence class of ideal vertices $\{(1,0,0,0), (-1,0,0,0)\}$, if we start with the ideal vertex $(1,0,0,0)$ then by the above discussion it is easy
to see that it lies on the sides $A$, $B$, $C$, $C'$, $G$ and $H$. The vertex $(1,0,0,0)$ can be thought of as the unit basis vector $e_1$ in $\R^4$, a horospherical
neighbourhood represented by pieces of the sides $A$, $B$, $C$, $C'$, $G$ and $H$ can then be thought of as a cube in $\R^3$ 
(viewed as an orthogonal hyperplane to $e_1$). The orientation of the orthogonal copy of $\R^3$ will always be taken to be:

\centerline {\graphicspath{ {manifold_3/parabolics/} } \includegraphics[width=3cm, height=3cm]{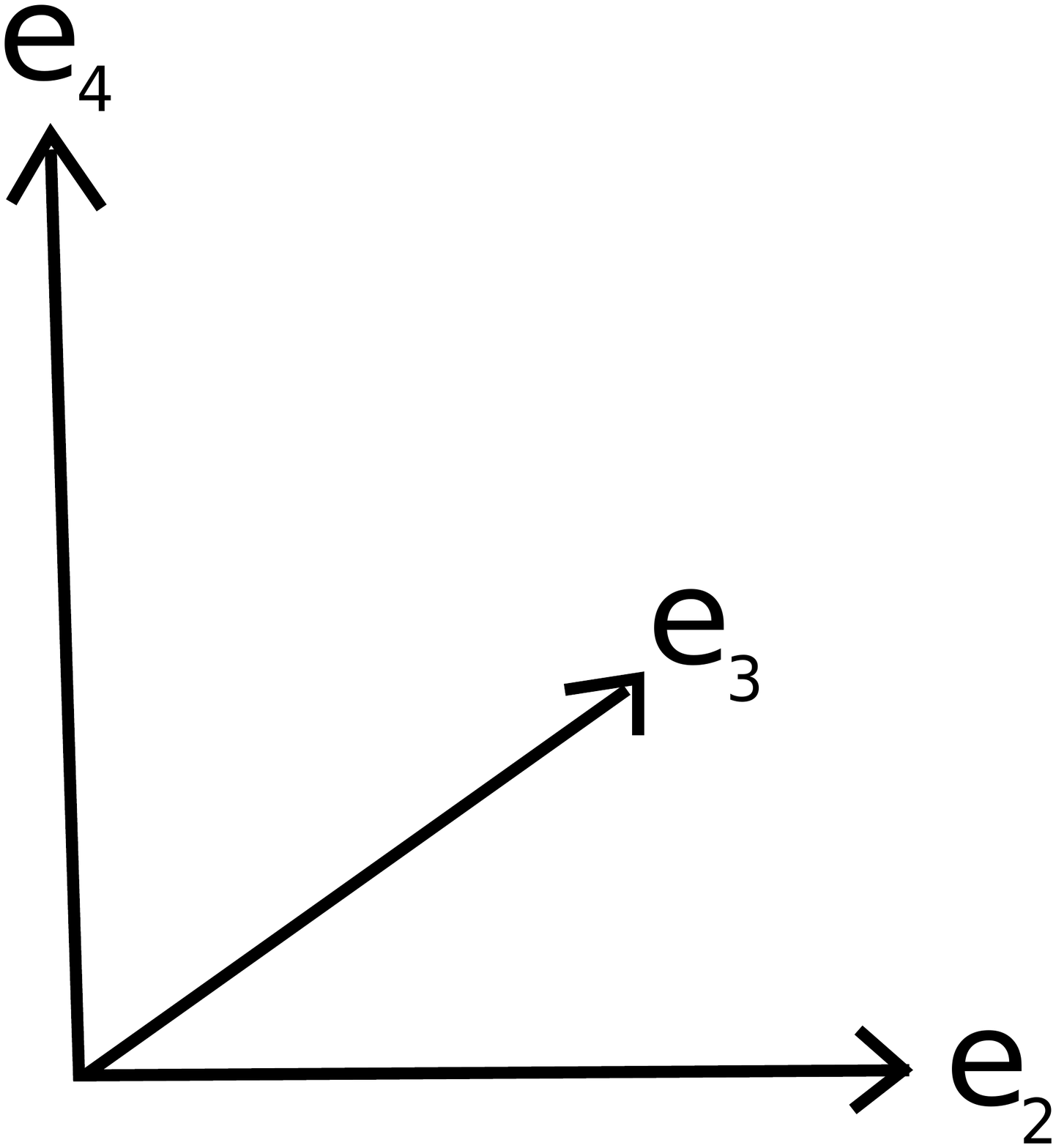}}

The horospherical neighbourhood will then look like:

\centerline {\graphicspath{ {manifold_3/parabolics/} } \includegraphics[width=6cm, height=5cm]{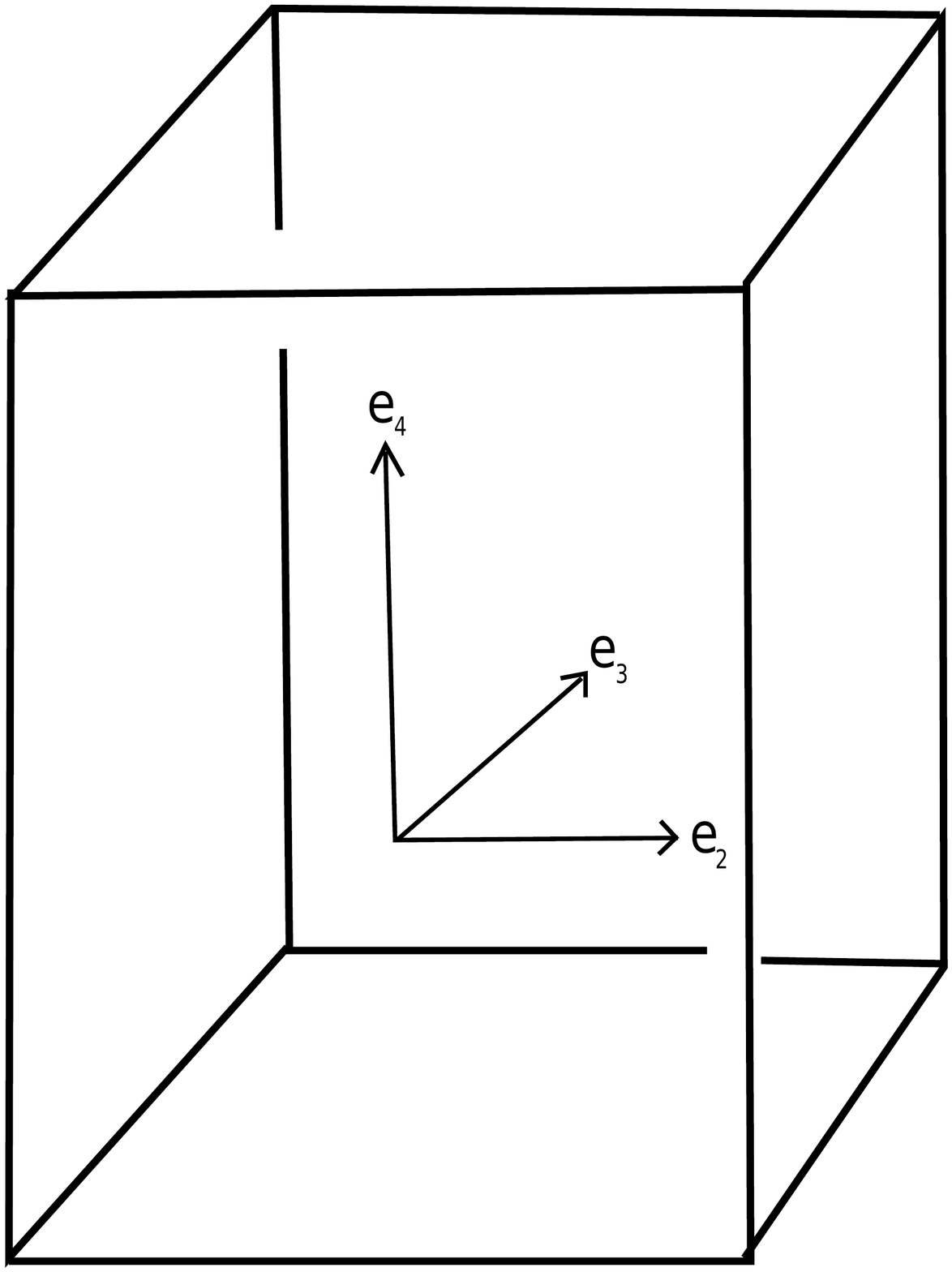}}

As the side $A$ is defined via the sphere $S_{(+1,+1,0,0)}$ with centre vector $(+1,+1,0,0)$ we see that its intersection with a horospherical neighbourhood
about $(1,0,0,0)$ can be thought of as a plane centred at $(x,0,0)$ (here the three co-ordinates are $(e_2, e_3, e_4)$), where $x \in \R$, and parallel to the plane
defined by the equation $e_2 = 0$. The actual value of $x$ is unimportant, it depends on where exactly we form our horospherical neighbourhood about $(1,0,0,0)$, which in turn comes down to how high up into the horoball about $(1,0,0,0)$ we are. The side $B$ is defined via the sphere $S_{(+1,-1,0,0)}$ with centre vector $(+1,-1,0,0)$, its intersection
with a horospherical neighbourhood about $(1,0,0,0)$ can be thought of as a plane centred at $(-x,0,0)$ and parallel to the plane defined by the equation
$e_2 = 0$. 
The side $C$ is defined by the sphere $S_{(+1,0,+1,0)}$ with centre vector $(+1,+1,0,0)$, similar reasoning to the above shows that its intersection with a 
horospherical neighbourhood about the ideal vertex $(1,0,0,0)$ can be viewed as a plane centred at $(0,y,0)$, where $y \in \R$, parallel to the plane $e_3 = 0$. The side
$C'$ can be viewed as a plane centred at $(0,-y,0)$ parallel to the plane $e_3 = 0$.
The side $G$ is defined by the sphere $S_{(+1,0,0,+1)}$ with centre vector $(+1,0,0,+1)$, intersecting it with a horospherical neighbourhood about 
the ideal vertex $(1,0,0,0)$ we obtain a plane centred at $(0,0,z)$, where $z \in \R$ and parallel to the plane $e_4 = 0$. The side $H$ can be viewed as a plane centred
at $(0,0,-z)$ parallel to the plane $e_4 = 0$.
Using this description we can label the above cube, which represents a horospherical neighbourhood based at $(1,0,0,0)$, as follows:

\centerline {\graphicspath{ {manifold_3/parabolics/} } \includegraphics[width=6cm, height=4cm]{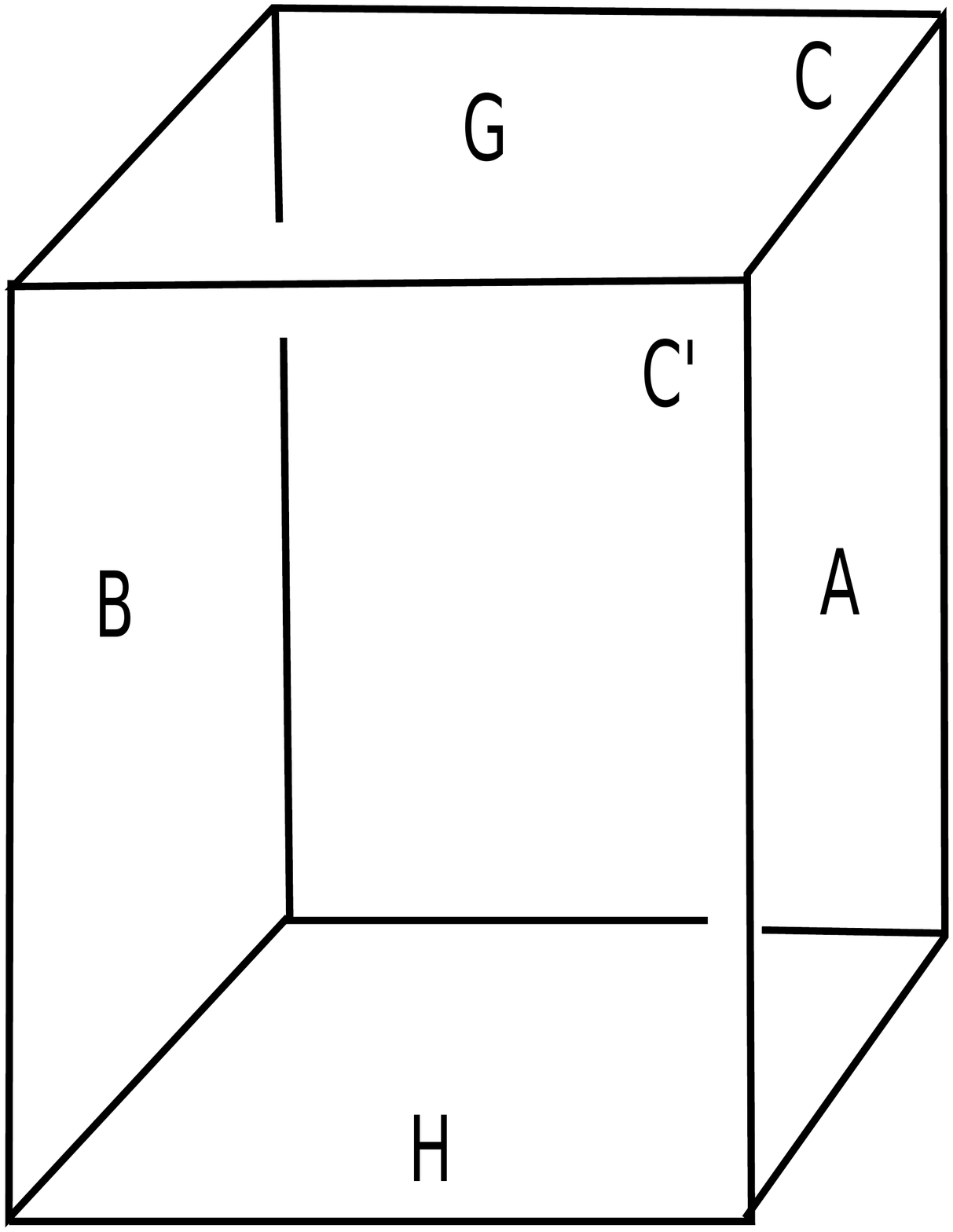}}

The equivalence class corresponding to the ideal vertex $(+1,0,0,0)$ also contains the ideal vertex $(-1,0,0,0)$, therefore in order to
describe the parabolic isometries giving rise to the cusp corresponding to this class we need to describe a horospherical neighbourhood about
$(-1,0,0,0)$. 
The procedure is exactly analogous to what we did above. However, for this vertex we are going to orient the orthogonal copy of $\R^3$ to the vertex 
$(-1,0,0,0)$ as follows:

\centerline {\graphicspath{ {manifold_3/parabolics/} } \includegraphics[width=5cm, height=4cm]{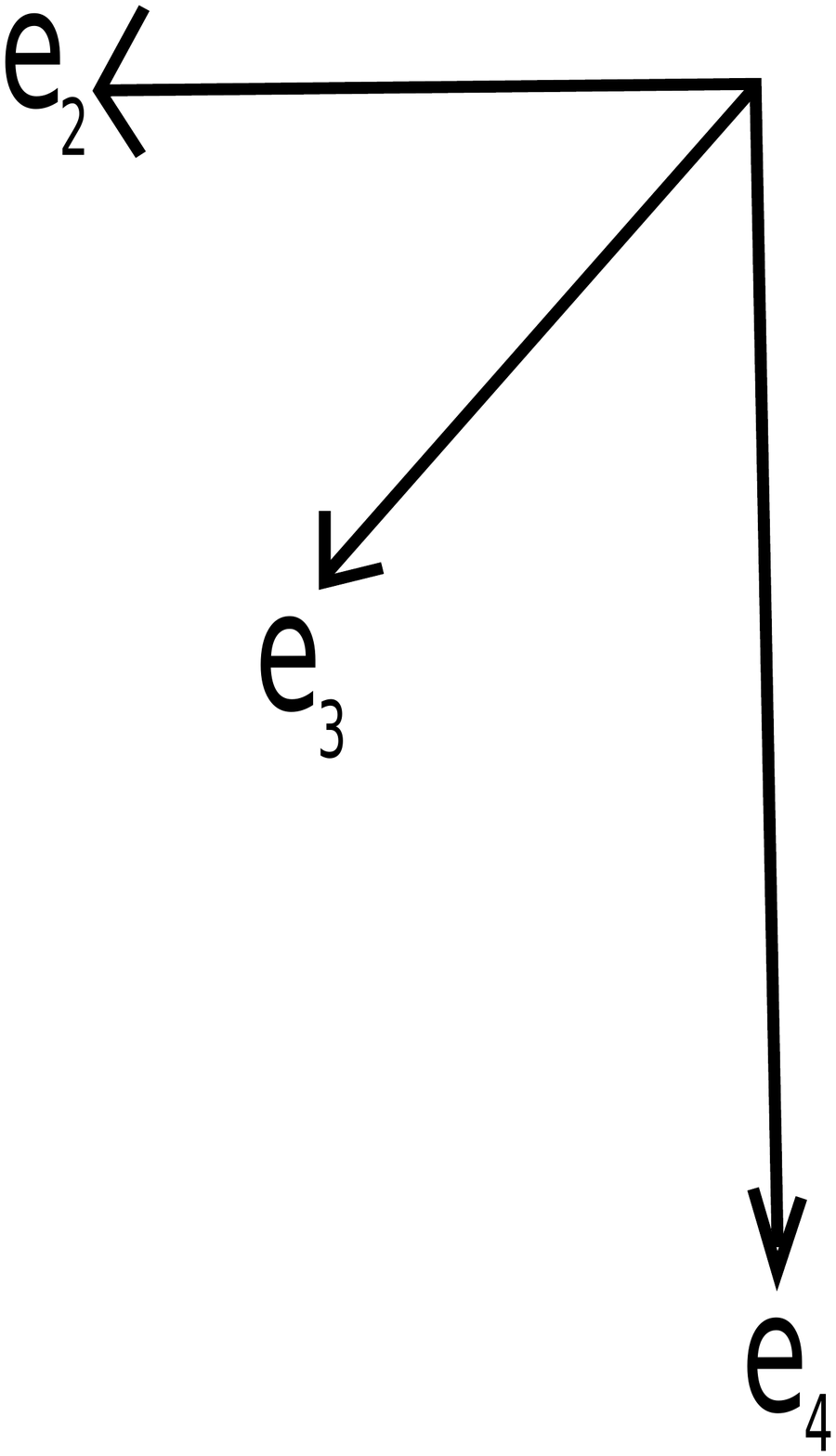}}

The reason for this choice will become apparent shortly, for now let us observe that it leads to the following labelling of a horospherical
neighbourhood:

\centerline {\graphicspath{ {manifold_3/parabolics/} } \includegraphics[width=5cm, height=4cm]{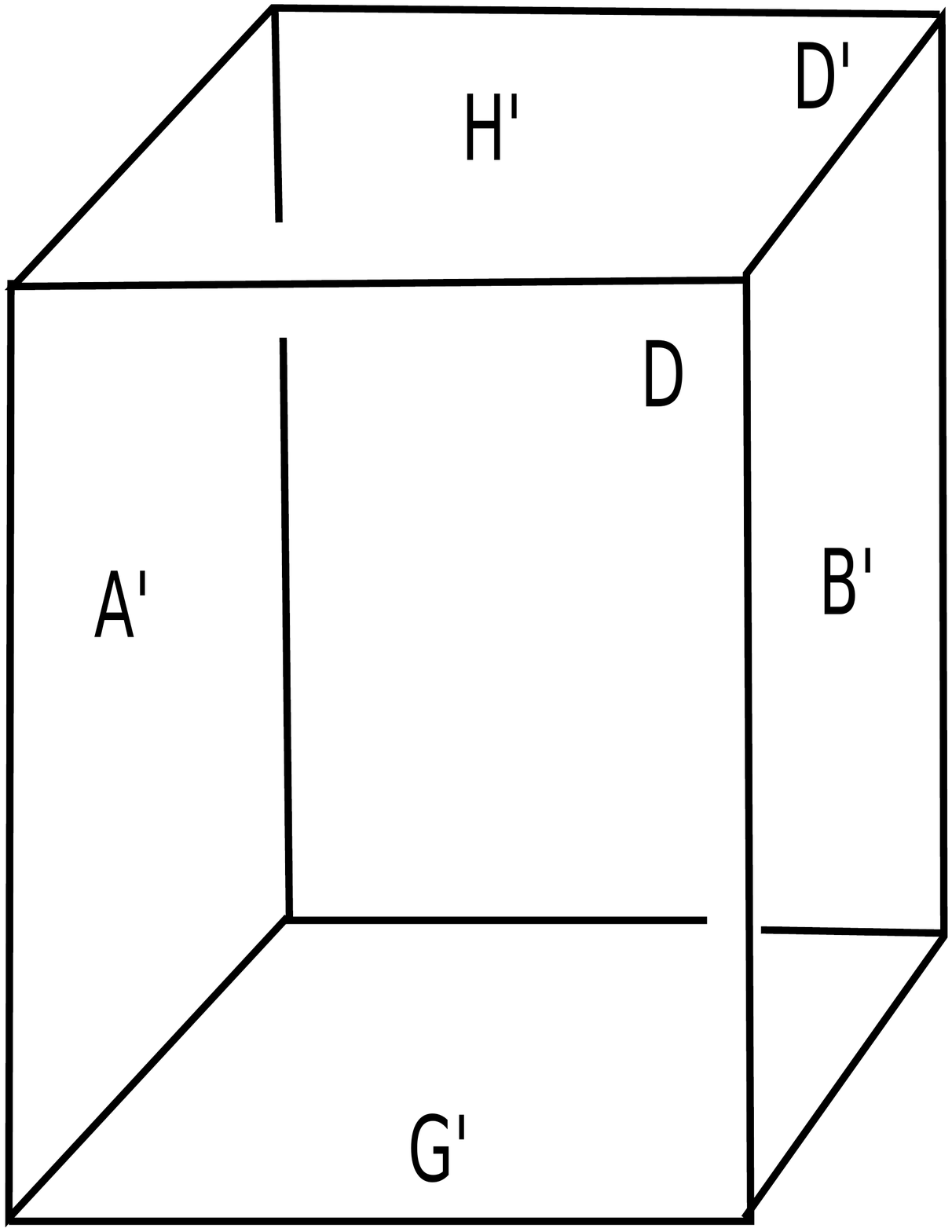}}

A fundamental domain for a cross section of the cusp corresponding to the above ideal vertex equivalence class consists of these two cubes.

\centerline {\graphicspath{ {manifold_3/parabolics/} } \includegraphics[width=6cm, height=5cm]{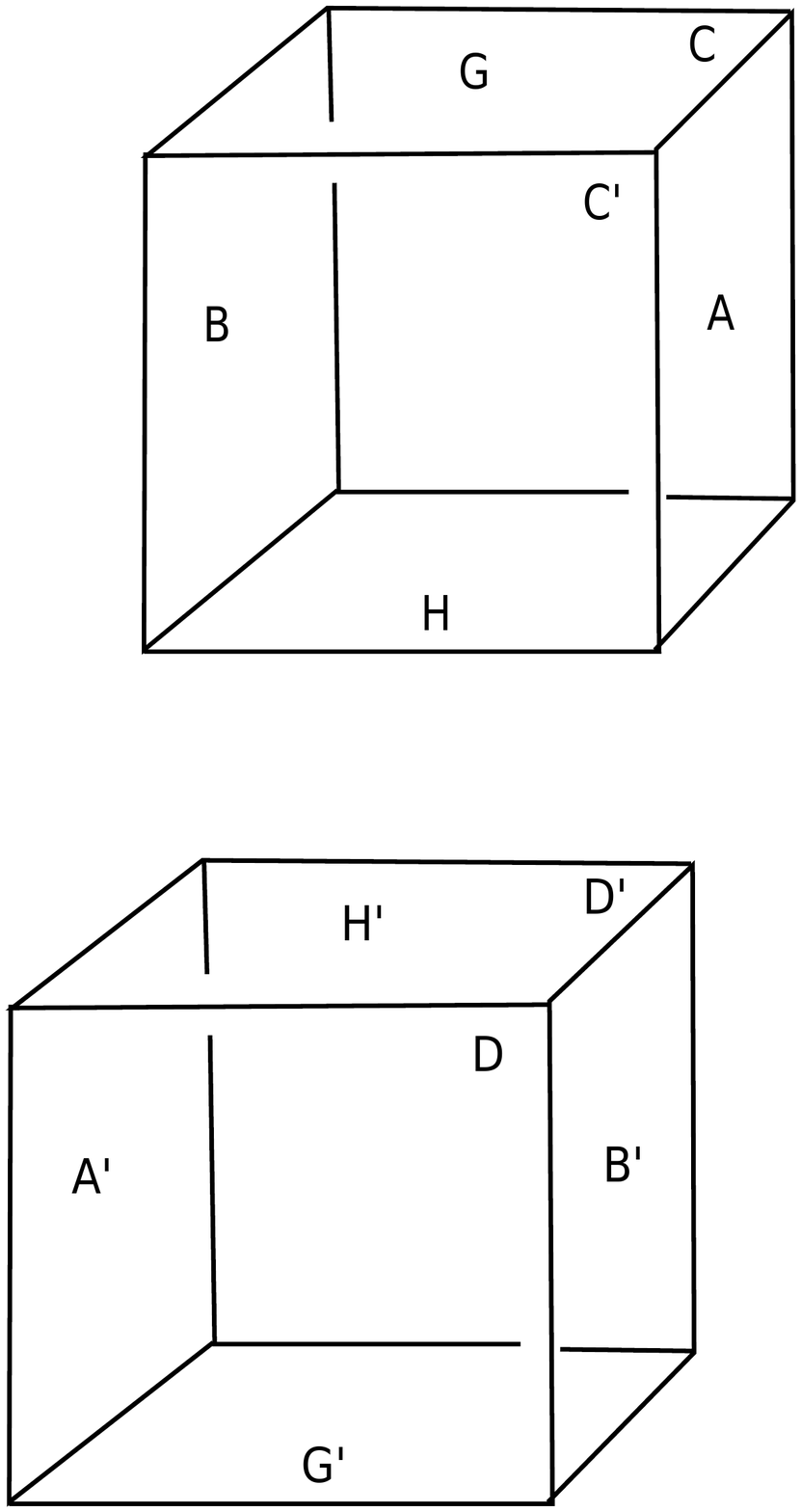}}

We can see that when we apply the transformation $h$ the bottom face of the top box gets joined to the top face of the bottom box. In principle the bottom face
of the top box has four different ways it can be joined to the top face of the bottom box. However, observe that when applying the transformation $h$, the side
$B$ maps to $A'$, side $A$ maps to $B'$, side $C$ maps to $D'$, and side $C'$ maps to $D$.
This means we can
visualise a fundamental domain as a rectangular box with sides labelled as follows:

\centerline {\graphicspath{ {manifold_3/parabolics/} } \includegraphics[width=6cm, height=5cm]{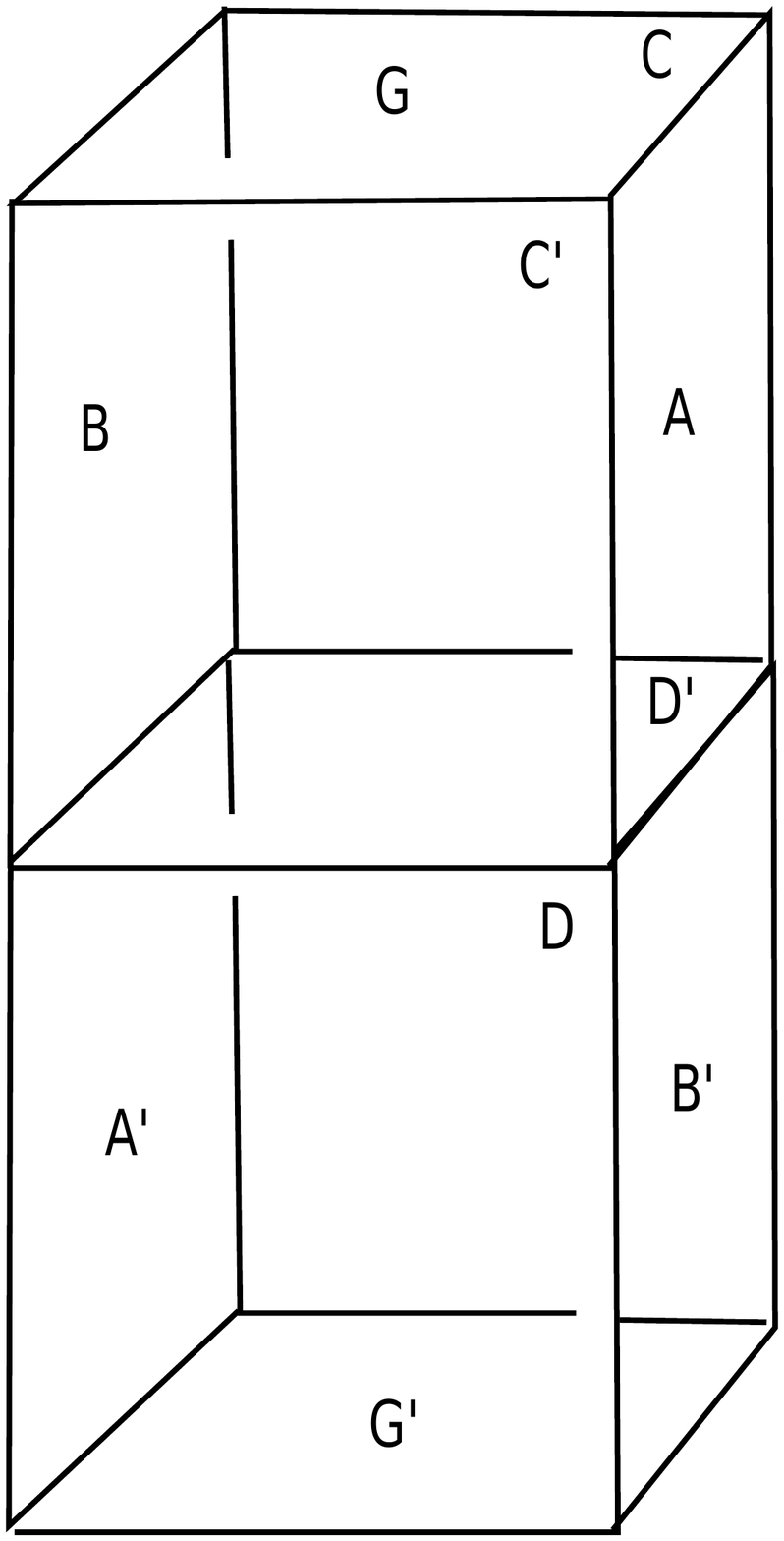}}

This is why we chose to orient the orthogonal copy of $\R^3$ to the ideal vertex $(-1,0,0,0)$ the way we did. We then see that the top face of the rectangular box
will be identified to the bottom face of the rectangular box by the parabolic transformation $g^{-1}h$ (note that $g^{-1}h$ fixes the ideal vertex $(1,0,0,0)$, hence
is clearly an element of the stabiliser subgroup of the ideal vertex $(1,0,0,0)$). Similarly the front face will be identified to the back face via the transformation
$c$ (note one could also take the transformation $d$, they both act in the same way). Finally, the identification of the face on the right side with the left side
is done by the transformation $a^{-1}h$, this is because the top right $A$ get identified to the bottom left $A'$, and the bottom right $B'$ gets identified 
to the top left $B$. Thus what we see that is happening is that the top face is being identified to the bottom face via a translation, the front face is being
identified to the back face via a translation, but the side faces are being identified by a ``twist'' (also called a screw parabolic transformation). This tells us that the stabiliser subgroup
associated to the vertex equivalence class $\{(1,0,0,0), (-1,0,0,0)\}$ is generated by the three transformations $\langle c, g^{-1}h, a^{-1}h \rangle$. Put another way
the parabolic subgroup associated to the vertex class $\{(1,0,0,0), (-1,0,0,0)\}$ is the subgroup $\langle c, g^{-1}h, a^{-1}h \rangle$.

Once one has understood the parabolic subgroup associated to a cusp, one can try to understand how the Euclidean structure on the cusp cross section comes
about. For a general non-compact hyperbolic manifold of finite volume this can be a difficult task, the reason being that in order to work out the
Euclidean structure on the cross section one needs a solid understanding of how elements of the associated parabolic subgroup act as Euclidean transformations on
a horospherical neighbourhood. For the Ratcliffe-Tschantz manifolds the action of the parabolic subgroup on a horospherical neighbourhood can be completely
understood due to some nice symmetry of the 24-cell, and the fact that the side pairing transformations are very easy to describe.

Continuing with the above example, we outline how to describe the Euclidean structure on a cusp cross section associated to the ideal vertex class
$\{(1,0,0,0), (-1,0,0,0)\}$.
We already know that the associated parabolic subgroup is given by $\langle c, g^{-1}h, a^{-1}h \rangle$, we want to understand how these transformations act
on a horospherical neighbourhood centred about $(1,0,0,0)$. A Euclidean transformation on $\R^3$ is just an affine transformation of the form:
\[ x \mapsto \Lambda\cdot x + v \]
where $\Lambda \in O(3)$ and $v \in \R^3$. For each of the transformations $c, g^{-1}h, a^{-1}h$ we want to understand what the matrix $\Lambda$ looks like, and what
the translation vector $v$ looks like.
Each side pairing transformation consists of two components, a k-part and an r-part, the k-part is given by a diagonal matrix, hence its action on the copy
of $\R^3$ is easy to understand. The r-part consists of reflection in the image side of the k-part, in this case its action on $\R^3$ is also easy to understand, this is
because, as was mentioned before, the intersection of the sides that the ideal vertex $(1,0,0,0)$ lies on with a horospherical neighbourhood can be 
thought as the box (remember our orientation convention):

\centerline {\graphicspath{ {manifold_3/parabolics/} } \includegraphics[width=5cm, height=4cm]{horocube_1000}}

Since the sides intersect the horospherical neighbourhood in planes parallel to the $e_2-e_3$, $e_2-e_4$ and $e_3-e_4$ planes we see that the r-part
is an affine transformation with $\Lambda$-matrix corresponding to one of reflections in the $e_2-e_3$, $e_2-e_4$ or $e_3-e_4$ planes. For example suppose
we take the side pairing transformation $c$. From the picture above we see that the part of the side $C'$ that intersects the horospherical neighbourhood  
can be thought of as a plane parallel to the plane $e_3 = 0$ and centred at the point $(0,-y,0)$.

\centerline{\graphicspath{ {manifold_3/parabolics/} } \includegraphics[width=8cm, height=8cm]{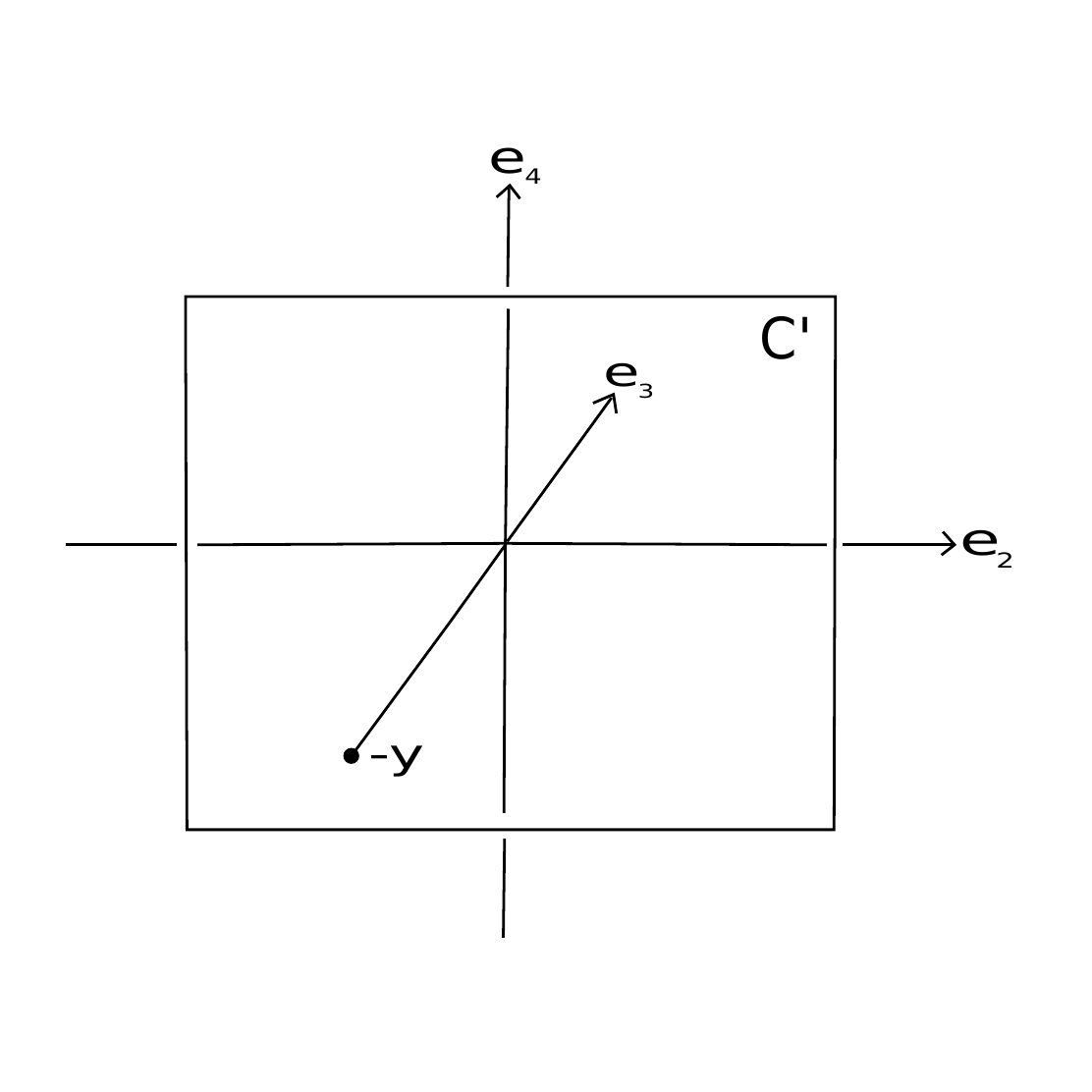}}

This tells us that the r-part of the transformation
$c$, which remember is reflection in the side $C'$, is a reflection in this plane. But it is easy to see that such a reflection is the
affine transformation:
\[ w \mapsto 
\begin{pmatrix}
1 & 0 & 0 \\
0 & -1 & 0 \\
0 & 0 & 1 
\end{pmatrix} 
\cdot w + 
\begin{pmatrix}
0 \\
-2y  \\
0  
\end{pmatrix} \]
The k-part of the transformation $c$ is given by the matrix:
\[
\begin{pmatrix}
1 & 0 & 0 & 0 \\
0 & 1 & 0 & 0 \\
0 & 0 & -1 & 0 \\
0 & 0 & 0 & 1
\end{pmatrix} \]
whose action on the horospherical neighbourhood (which we are viewing as a copy of $\R^3$) is given by the matrix:
\[
\begin{pmatrix}
1 & 0 & 0  \\
0 & -1 & 0 \\
0 & 0 & 1 
\end{pmatrix} \]
Therefore the action of the transformation $c$ on this horospherical neighbourhood being the composition of the k-part and the r-part, 
is given by the transformation:
\[w \mapsto
\begin{pmatrix}
1 & 0 & 0  \\
0 & -1 & 0 \\
0 & 0 & 1 
\end{pmatrix}
\begin{pmatrix}
1 & 0 & 0  \\
0 & -1 & 0 \\
0 & 0 & 1 
\end{pmatrix} \cdot w + 
\begin{pmatrix}
0 \\
-2y  \\
0  
\end{pmatrix} = w + 
\begin{pmatrix}
0 \\
-2y  \\
0  
\end{pmatrix} \]
This tells us that the Euclidean transformation that corresponds to the action of $c$ on a horospherical neighbourhood about $(1,0,0,0)$ is nothing more
than a translation by the vector $(0,-2y,0)$. 

We can apply the same techniques to understand the action of the parabolic transformation $g^{-1}h$ on the horospherical neighbourhood. In this case we need 
to understand the k and r-parts of two transformations, $g^{-1}$ and $h$. 

The k-part of the transformation $h$ is given by the matrix:
\[
\begin{pmatrix}
-1 & 0 & 0 & 0 \\
0 & -1 & 0 & 0 \\
0 & 0 & -1 & 0 \\
0 & 0 & 0 & 1
\end{pmatrix} \]
The side $H'$ is parallel to the plane $e_4 = 0$ and centred at the vector $(0,0,-z)$. Therefore
the r-part of the transformation $h$ is given by the affine transformation:
\[ w \mapsto
\begin{pmatrix}
1 & 0 & 0  \\
0 & 1 & 0  \\
0 & 0 & -1 
\end{pmatrix} \cdot w +  
\begin{pmatrix}
0  \\
0 \\
-2z  
\end{pmatrix} \]
The k-part of the transformation $g^{-1}$ is the same is that of $h$, the r-part is given by reflection in the plane defined by the side $G$. The associated plane
is parallel to the plane $e_4 = 0$ and centred at the point $(0,0,z)$. Therefore the r-part for $g^{-1}$ is given by the affine transformation:
\[ w \mapsto
\begin{pmatrix}
1 & 0 & 0  \\
0 & 1 & 0  \\
0 & 0 & -1 
\end{pmatrix} \cdot w +  
\begin{pmatrix}
0  \\
0 \\
2z  
\end{pmatrix} \]
The composition $g^{-1}h$ consists of the two k-parts one from $g^{-1}$ and one from $h$, as these are both equal they simply give the identity. Therefore
to understand $g^{-1}h$ as an affine transformation we need only understand the composition of the r-parts of $h$ and $g^{-1}$. This is easy to
calculate:
\[ w \mapsto
\begin{pmatrix}
1 & 0 & 0  \\
0 & 1 & 0  \\
0 & 0 & -1 
\end{pmatrix} \left( 
\begin{pmatrix}
1 & 0 & 0  \\
0 & 1 & 0  \\
0 & 0 & -1 
\end{pmatrix} \cdot w +
\begin{pmatrix}
0  \\
0 \\
-2z  
\end{pmatrix} \right) + 
\begin{pmatrix}
0  \\
0 \\
2z  
\end{pmatrix}
\]
which is the affine transformation:
\[ w \mapsto x +
\begin{pmatrix}
0  \\
0 \\
4z  
\end{pmatrix} \]
Finally, let us work out how the parabolic transformation $a^{-1}h$ behaves as a Euclidean transformation on a horospherical neighbourhood.
The k-parts of $h$ and $a^{-1}$ are given by the following matrices:
\[
\begin{pmatrix}
-1 & 0 & 0 & 0 \\
0 & -1 & 0 & 0 \\
0 & 0 & -1 & 0 \\
0 & 0 & 0 & 1
\end{pmatrix}, 
\begin{pmatrix}
-1 & 0 & 0 & 0 \\
0 & 1 & 0 & 0 \\
0 & 0 & 1 & 0 \\
0 & 0 & 0 & 1
\end{pmatrix}
\]
We have already computed the r-part of $h$, as for the r-part of $a^{-1}$ observe that the side $A$ is parallel to the $e_2 = 0$ plane, and centred at
the vector $(x,0,0)$. Reflection in this plane is given by the affine transformation:
\[ w \mapsto
\begin{pmatrix}
-1 & 0 & 0  \\
0 & 1 & 0  \\
0 & 0 & 1 \\
\end{pmatrix} \cdot w +
\begin{pmatrix}
2x \\
0 \\
0 \\
\end{pmatrix} \] 
The composition $a^{-1}h$ is then given by the affine transformation:
\[ w \mapsto
\begin{pmatrix}
-1 & 0 & 0  \\
0 & 1 & 0  \\
0 & 0 & 1 \\
\end{pmatrix} \left( 
\begin{pmatrix}
-1 & 0 & 0  \\
0 & -1 & 0  \\
0 & 0 & -1 \\
\end{pmatrix} \cdot w + 
\begin{pmatrix}
0 \\
0 \\
-2z \\
\end{pmatrix} \right) +
\begin{pmatrix}
2x \\
0 \\
0 \\
\end{pmatrix} \]
which when we expand out we obtain:
\[
\begin{pmatrix}
1 & 0 & 0  \\
0 & -1 & 0  \\
0 & 0 & -1 \\
\end{pmatrix} \cdot w +
\begin{pmatrix}
2x \\
0 \\
-2z \\
\end{pmatrix} \]
The reader will recall that before we started showing how to compute the action of these parabolic isometries on a horospherical neighbourhood, we showed
how a cusp cross section arises through a fundamental domain consisting of a rectangular box. We went on to say that the right and left sides of this fundamental domain
(which were being identified by the parabolic transformation $a^{-1}h$) were identified via a ``twist''. The reader who was not content with the use of the word ``twist'' at that
time should now be at ease, for the ``twisting'' we speak of is given by the matrix
\[
\begin{pmatrix}
1 & 0 & 0  \\
0 & -1 & 0  \\
0 & 0 & -1 \\
\end{pmatrix} \]
in the above formula for the affine transformation associated to the parabolic isometry $a^{-1}h$.

In summary, we have shown that the parabolic subgroup corresponding to the ideal vertex equivalence class $\{(1,0,0,0), (-1,0,0,0)\}$ is given by the group
$\langle c, g^{-1}h, a^{-1}h \rangle$. Furthermore, for each of these parabolic isometries we have shown how they act on a horospherical neighbourhood as
a Euclidean transformation. The following table summarises this information
\begin{center}
\begin{tabular}{|l|l|}
\hline
 Parabolic isometry &  Affine transformation  \\  \hline  

 $c$ &  $\begin{pmatrix}
1 & 0 & 0  \\
0 & 1 & 0  \\
0 & 0 & 1 \\
\end{pmatrix}$,    
$\begin{pmatrix}
0 \\
-2y  \\
0 \\
\end{pmatrix}$ \\ \hline 

 $g^{-1}h$ &   $\begin{pmatrix}
1 & 0 & 0  \\
0 & 1 & 0  \\
0 & 0 & 1 \\
\end{pmatrix}$,    
$\begin{pmatrix}
0 \\
0  \\
4z \\
\end{pmatrix}$ \\ \hline

 $a^{-1}h$ &   $\begin{pmatrix}
1 & 0 & 0  \\
0 & -1 & 0  \\
0 & 0 & -1 \\
\end{pmatrix}$,    
$\begin{pmatrix}
2x \\
0  \\
-2z \\
\end{pmatrix}$ \\  \hline   

\end{tabular}  
\end{center} 
where in the second column the matrix is in $O(3)$ and the vector corresponds to the translation vector.

The information presented above is enough for one to obtain the Euclidean structure of the cusp corresponding to the class $\{(1,0,0,0), (-1,0,0,0)\}$. 
Recall that there are ten classes of closed flat 3-manifolds denoted by $\textbf{A}$, $\textbf{B}$, $\textbf{C}$, $\textbf{D}$, $\textbf{E}$,
$\textbf{F}$, $\textbf{G}$, $\textbf{H}$, $\textbf{I}$, and $\textbf{J}$ in the Hantzsche-Wendt notation see \cite{hantzsche}. These are denoted
by $\mathcal{G}_1$, $\mathcal{G}_2$, $\mathcal{G}_3$, $\mathcal{G}_4$, $\mathcal{G}_5$, $\mathcal{G}_6$, $\mathcal{B}_1$, $\mathcal{B}_2$, 
$\mathcal{B}_3$ and $\mathcal{B}_4$ respectively using the notation of Wolf see \cite{wolf} Thm.3.5.5, p.117. The first six are the orientable ones, with
$\textbf{A}$ being the 3-torus, and $\textbf{B}$ being the orientable $S^1$-fibre bundle over the Klein bottle. The last four are all non-orientable.
Using the classification of orientable compact flat 3-manifolds (see \cite{wolf} Thm.3.5.5, p.117) we see that this cusp has type $\textbf{B}$ 
using the Hantzsche-Wendt notation or type $\mathcal{G}_2$ using Wolf's notation. 

One can carry out an analogous procedure to work out generators for the parabolic subgroups corresponding to the other ideal vertex equivalence classes and their corresponding
affine transformations. The only slight difference is that when dealing with the class $\{(\pm 1/2, \pm 1/2, \pm 1/2, \pm 1/2)\}$ one has to take
sixteen boxes to form a fundamental domain. Also, computing the associated affine transformations in this case is slightly harder as one does not have
such a nice parameterisation of horospherical neighbourhoods as in the case for the other equivalence classes. We will not need explicit formulas for the corresponding
affine transformations, hence will not bother the reader with an explanation of how to find them in this situation.
For the sake of completeness we have included the following table showing the end results of the above computations for all other ideal vertex equivalence classes. \\

\begin{table}[H]
  \resizebox{0.9\textwidth}{10cm}{\begin{minipage}{\textwidth}

\hskip-2.0cm\begin{tabular}{|l | l | l | l |}
\hline

{} & {} & {} & {} \\

Ideal vertex class & Generators for stabiliser subgroup & Affine transformations & Euclidean structure   \\ 

{} & {} & {} & {} \\ \hline

{} & {} & {} & {} \\

$\{(1,0,0,0), (-1,0,0,0)\}$ & $c$ & $\begin{pmatrix}
1 & 0 & 0  \\
0 & 1 & 0  \\
0 & 0 & 1 \\
\end{pmatrix}$,    
$\begin{pmatrix}
0 \\
-2y  \\
0 \\
\end{pmatrix}$ & $\textbf{B}$ ($\mathcal{G}_2$) \\

{} & {} & {} & {} \\

{}         & $g^{-1}h$ & $\begin{pmatrix}
1 & 0 & 0  \\
0 & 1 & 0  \\
0 & 0 & 1 \\
\end{pmatrix}$,    
$\begin{pmatrix}
0 \\
0  \\
4z \\
\end{pmatrix}$ & {} \\

{} & {} & {} & {} \\

{} & $a^{-1}h$ &   $\begin{pmatrix}
1 & 0 & 0  \\
0 & -1 & 0  \\
0 & 0 & -1 \\
\end{pmatrix}$,    
$\begin{pmatrix}
2x \\
0  \\
-2z \\
\end{pmatrix}$ & {} \\ 

{} & {} & {} & {} \\ \hline 

{} & {} & {} & {} \\ 

$\{(0,1,0,0), (0,-1,0,0)\}$ & $a$ & $\begin{pmatrix}
1 & 0 & 0  \\
0 & 1 & 0  \\
0 & 0 & 1 \\
\end{pmatrix}$,    
$\begin{pmatrix}
-2x \\
0  \\
0 \\
\end{pmatrix}$ & $\textbf{A}$ ($\mathcal{G}_1$)  \\ 

{} & {} & {} & {} \\

{} & $e^{-1}f$ & $\begin{pmatrix}
1 & 0 & 0  \\
0 & 1 & 0  \\
0 & 0 & 1 \\
\end{pmatrix}$,    
$\begin{pmatrix}
0 \\
4y  \\
0 \\
\end{pmatrix}$ & {} \\

{} & {} & {} & {} \\

{} & $e^{-1}i$ & $\begin{pmatrix}
1 & 0 & 0  \\
0 & 1 & 0  \\
0 & 0 & 1 \\
\end{pmatrix}$,    
$\begin{pmatrix}
0 \\
2y  \\
-2z \\
\end{pmatrix}$ & {} \\ 
 
{} & {} & {} & {} \\  \hline

{} & {} & {} & {} \\

$\{(0,0,1,0), (0,0,-1,0)\}$ & $k$ & $\begin{pmatrix}
1 & 0 & 0  \\
0 & 1 & 0  \\
0 & 0 & 1 \\
\end{pmatrix}$,    
$\begin{pmatrix}
0 \\
0  \\
-2z \\
\end{pmatrix}$ & $\textbf{A}$ ($\mathcal{G}_1$)  \\

{} & {} & {} & {} \\

{} & $d^{-1}c$ & $\begin{pmatrix}
1 & 0 & 0  \\
0 & 1 & 0  \\
0 & 0 & 1 \\
\end{pmatrix}$,    
$\begin{pmatrix}
4x \\
0  \\
0 \\
\end{pmatrix}$ & {} \\

{} & {} & {} & {} \\

{} & $e^{-1}c$ & $\begin{pmatrix}
1 & 0 & 0  \\
0 & 1 & 0  \\
0 & 0 & 1 \\
\end{pmatrix}$,    
$\begin{pmatrix}
0 \\
2y  \\
-2z \\
\end{pmatrix}$ & {} \\ 

{} & {} & {} & {} \\  \hline

{} & {} & {} & {} \\ 

$\{(0,0,0,1), (0,0,0,-1)\}$ & $g$ & $\begin{pmatrix}
1 & 0 & 0  \\
0 & -1 & 0  \\
0 & 0 & -1 \\
\end{pmatrix}$,    
$\begin{pmatrix}
-2x \\
0  \\
0 \\
\end{pmatrix}$ & $\textbf{F}$ ($\mathcal{G}_6$)  \\

{} & {} & {} & {} \\ 

{} & $k^{-1}l$ & $\begin{pmatrix}
1 & 0 & 0  \\
0 & 1 & 0  \\
0 & 0 & 1 \\
\end{pmatrix}$,    
$\begin{pmatrix}
0 \\
0  \\
4z \\
\end{pmatrix}$ & {} \\ 

{} & {} & {} & {} \\ 

{} & $i^{-1}l$ & $\begin{pmatrix}
-1 & 0 & 0  \\
0 & 1 & 0  \\
0 & 0 & -1 \\
\end{pmatrix}$,    
$\begin{pmatrix}
0 \\
-2y  \\
2z \\
\end{pmatrix}$ & {} \\

{} & {} & {} & {} \\  \hline

{} & {} & {} & {} \\ 

$\{(\pm 1/2, \pm 1/2, \pm 1/2, \pm 1/2)\}$ & $e^{-1}g$ & $\begin{pmatrix}
1 & 0 & 0  \\
0 & 1 & 0  \\
0 & 0 & 1 \\
\end{pmatrix}$,    
$\begin{pmatrix}
0 \\
-4y  \\
0 \\
\end{pmatrix}$ & $\textbf{B}$ ($\mathcal{G}_2$)  \\

{} & {} & {} & {} \\ 

{} & $a^{-1}k^{-1}ak$ & $\begin{pmatrix}
1 & 0 & 0  \\
0 & 1 & 0  \\
0 & 0 & 1 \\
\end{pmatrix}$,    
$\begin{pmatrix}
8x \\
0  \\
0 \\
\end{pmatrix}$ & {} \\

{} & {} & {} & {} \\

{} & $a^{-1}k^{-1}j^{-1}b^{-1}fc$ & 
$\begin{pmatrix}
-1 & 0 & 0  \\
0 & -1 & 0  \\
0 & 0 & 1 \\
\end{pmatrix}$,    
$\begin{pmatrix}
-2x \\
-2y \\
-4z \\
\end{pmatrix}$ & {} \\

{} & {} & {} & {} \\  \hline

\end{tabular} \\ \\ 
 
        \end{minipage}}
\end{table}

We point out that if we order the ideal vertices according to our table above then we see that the Euclidean structures on the associated
cusps are given by \textbf{BAAFB}. However, if we go to manifolds no.3 in the Ratcliffe-Tschantz census we see that under the column headed \textbf{LT} the Euclidean structure
on the cusps is given by \textbf{AABBF}. This is because they chose to write out the Euclidean structures in alphabetical order, as opposed to any ordering
on the ideal vertices.

The Euclidean structure on each cusp is found by an appeal to the classification theorem of compact connected orientable flat 3-dimensional Riemannian manifolds
(see \cite{wolf} Thm.3.5.5, p.117). One simply needs to identify the linear holonomy group, and one can do this
from the computations of the associated affine transformations. The matrices of each generator (when viewed as an affine transformation on a horospherical
neighbourhood) generate the linear holonomy group. For example, consider the first entry in the above table, we see that two of the generators have
the identity matrix as their $O(3)$ component, when viewed as an affine transformation, and one of the generators has 
$\biggl(\begin{smallmatrix}
1&0&0\\ 0&-1&0\\ 0&0&-1
\end{smallmatrix} \biggr)$ 
as its $O(3)$ component. Therefore the linear holonomy group is generated by 
$\biggl(\begin{smallmatrix}
1&0&0\\ 0&-1&0\\ 0&0&-1
\end{smallmatrix} \biggr)$, which has order 2. This implies the linear holonomy group associated to this cusp is $\Z_2$. Using the classification theorem
(see \cite{wolf} Thm.3.5.5, p.117) we see that it corresponds to a Euclidean structure of type \textbf{B} (or \textbf{$\mathcal{G}_2$} in Wolf's notation).

Observe that if we look at the cusp corresponding to the ideal vertex class $\{(0,1,0,0), (0,-1,0,0)\}$ we see that the Euclidean structure is given
by $\textbf{A}$, which is a 3-torus. Viewing the 3-torus as an $S^1$-fibre bundle over the 2-torus we see that there is a natural 4-manifold that it bounds. Namely, the
associated disk bundle, which we denote by $\widetilde{\textbf{A}}$. If we then take a horoball neighbourhood about the cusp corresponding to 
$\{(0,1,0,0), (0,-1,0,0)\}$ and chop it off, we obtain a boundary 3-manifold given by $\textbf{A}$. We can then glue in the 4-manifold $\widetilde{\textbf{A}}$ 
by identifying boundary components, thus killing the cusp corresponding to $\{(0,1,0,0), (0,-1,0,0)\}$. At this point there are a few
remarks that we need to make. The first one is that since we are gluing in a disc bundle over a flat surface the Euler characteristic does not change. 
The second remark is that the way we glue in the solid 3-torus $\widetilde{\textbf{A}}$ is not unique, it depends on a choice. Namely, in choosing
how the 3-torus fibres over the 2-torus we made a choice of an $S^1$-fibre, you can think of this as fixing a meridian. When we carry out the gluing
process we are filling in this meridian with a copy of $D^2$. Finally, the 4-manifold we chose to glue in was the solid torus and this choice is not
unique. What we mean by this is that since there are other 4-manifolds bounded by the 3-torus we could have just as well chosen one of them to glue in.
However as our interest is in constructing codimension two link complements we will always focus on the disk bundle case.

From the classification theorem of closed flat 3-manifolds (see \cite{wolf} Thm.3.5.5, p.117) one knows that 
the only orientable closed flat 3-manifolds that are $S^1$-fibre bundles are given by types $\textbf{A}$ and $\textbf{B}$ (or $\mathcal{G}_1$ and $\mathcal{G}_2$ using
Wolf's notation). Furthermore, all the non-orientable ones $\textbf{G}$, $\textbf{H}$, $\textbf{I}$, $\textbf{J}$ ($\mathcal{B}_1$, $\mathcal{B}_2$, $\mathcal{B}_3$, $\mathcal{B}_4$ in Wolf's notation)
are $S^1$ fibre bundles. There are many hyperbolic 4-manifolds in the Ratcliffe-Tschantz census whose cusp structure has type $\textbf{A}$, $\textbf{B}$, 
$\textbf{G}$, $\textbf{H}$, $\textbf{I}$, $\textbf{J}$.
For such manifolds we can obtain a closed 4-manifold by gluing in associated disk bundles as described above, such a gluing procedure does not change
the Euler characteristic. Let $E \rightarrow S$ denote a closed flat 3-manifold written as an $S^1$ fibre bundle over a compact surface $S$ (note that
from the classification theorem we know that $S$ must be either a 2-torus or a Klein bottle). The associated disk bundle can then be written as
\[\pi : \widetilde{E} \rightarrow S\]
The 4-manifold $\widetilde{E}$ has a nice handle decomposition
coming from the handle decomposition of the surface $S$. $S$ is a compact surface, hence can be given a handle decomposition consisting of one 0-handle, 
$n$ 1-handles ($n = 2 - \chi(S)$), and one 2-handle (see \cite{gompf} p.131 last paragraph). The $k$-handles of $\widetilde{E}$ are given by
preimages under $\pi$ of $k$-handles of $B$ (see \cite{gompf} p.131 last paragraph). This means that $\widetilde{E}$ has a handle decomposition
consisting of one 0-handle, two 1-handles and one 2-handle, since $\chi(S) = 0$. 

On the level of a Kirby diagram the gluing in of the disk bundle $\widetilde{E}$ is done by adding one 2-handle, two 3-handles and one 4-handle (see \cite{akbulut} chap.3, p.35). If we are given a hyperbolic 4-manifold with cusps all of which are $S^1$-fibre bundles then we can glue in a disk bundle to each cusp to produce a closed
4-manifold. In our paper \cite{sarat} we showed how to construct a Kirby diagram for any one of the Ratcliffe-Tschantz manifolds. Provided we can understand
where to add the 2-handle corresponding to the gluing of the disk bundle we can then produce a Kirby diagram for the ``filled in'' manifold.
We can then try and apply certain handle slides/cancellations to try and reduce the Kirby diagram of this ``filled in'' manifold. The hope is that we can reduce it to the Kirby diagram of a familiar
closed 4-manifold that we can explicitly identify. Provided we are successful in this reduction process, we would have then found an explicit four dimensional
hyperbolic link complement.

Our general method of constructing four dimensional hyperbolic link complements whose diffeomorphism type can be identified can thus be broken down
in to three general steps. The first step involves finding finite volume non-compact hyperbolic 4-manifolds with cusp structure given by
one of the closed flat 3-manifolds that is an $S^1$-fibre bundle over a flat surface. The second step involves being able to construct
a Kirby diagram for the filling of the hyperbolic 4-manifold. Our approach in this step is to use a construction of a Kirby diagram for the hyperbolic
4-manifold, then explicitly work out how to add 2-handles to the diagram to obtain a Kirby diagram for the filling. Finally the third step involves
carrying out various handle slides/cancellations. We are going to give an explicit example of a filling using a hyperbolic 4-manifold obtained from the Ratcliffe-Tschantz
census soon. The handle slides/cancellations we use are based on three very simply moves from the theory of Kirby calculus. In the next section we will take
the time to explain exactly what these moves are and how they will be used.

The main manifold we are going to be concerned with in terms of the filling process is numbered 1011 in the census. This is a non-orientable five cusped hyperbolic
4-manifold with each cusp structure given by the non-orientable $S^1$ fibre bundle with flat structure given by type $\textbf{G}$ (or $\mathcal{B}_1$ in Wolf's notation). It has orientable double cover the flat manifold given by type $\textbf{A}$ (or $\mathcal{G}_1$ in Wolf's notation), the 3-torus. The structure of each cusp can be found using methods
exactly analogous to what we did for manifold no. 3. The following table summarises this information.

\begin{table}[H]
  \resizebox{0.9\textwidth}{10cm}{\begin{minipage}{\textwidth}

\hskip-2.0cm\begin{tabular}{|l | l | l | l |}
\hline

{} & {} & {} & {} \\

Ideal vertex class & Generators for stabiliser subgroup & Affine transformations & Euclidean structure   \\ 

{} & {} & {} & {} \\ \hline

{} & {} & {} & {} \\

$\{(1,0,0,0), (-1,0,0,0)\}$ & $c$ & $\begin{pmatrix}
1 & 0 & 0  \\
0 & 1 & 0  \\
0 & 0 & 1 \\
\end{pmatrix}$,    
$\begin{pmatrix}
0 \\
-2y  \\
0 \\
\end{pmatrix}$ & \textbf{G} ($\mathcal{B}_1$) \\

{} & {} & {} & {} \\

{}         & $a^{-1}b$ & $\begin{pmatrix}
1 & 0 & 0  \\
0 & 1 & 0  \\
0 & 0 & 1 \\
\end{pmatrix}$,    
$\begin{pmatrix}
4x \\
0  \\
0 \\
\end{pmatrix}$ & {} \\

{} & {} & {} & {} \\

{} & $a^{-1}g$ &   $\begin{pmatrix}
1 & 0 & 0  \\
0 & -1 & 0  \\
0 & 0 & 1 \\
\end{pmatrix}$,    
$\begin{pmatrix}
2x \\
0  \\
-2z \\
\end{pmatrix}$ & {} \\ 

{} & {} & {} & {} \\ \hline 

{} & {} & {} & {} \\ 

$\{(0,1,0,0), (0,-1,0,0)\}$ & $a$ & $\begin{pmatrix}
1 & 0 & 0  \\
0 & 1 & 0  \\
0 & 0 & 1 \\
\end{pmatrix}$,    
$\begin{pmatrix}
-2x \\
0  \\
0 \\
\end{pmatrix}$ & \textbf{G} ($\mathcal{B}_1$) \\ 

{} & {} & {} & {} \\

{} & $e^{-1}f$ & $\begin{pmatrix}
1 & 0 & 0  \\
0 & 1 & 0  \\
0 & 0 & 1 \\
\end{pmatrix}$,    
$\begin{pmatrix}
0 \\
4y  \\
0 \\
\end{pmatrix}$ & {} \\

{} & {} & {} & {} \\

{} & $e^{-1}i$ & $\begin{pmatrix}
-1 & 0 & 0  \\
0 & 1 & 0  \\
0 & 0 & 1 \\
\end{pmatrix}$,    
$\begin{pmatrix}
0 \\
2y  \\
-2z \\
\end{pmatrix}$ & {} \\ 
 
{} & {} & {} & {} \\  \hline

{} & {} & {} & {} \\

$\{(0,0,1,0), (0,0,-1,0)\}$ & $k$ & $\begin{pmatrix}
1 & 0 & 0  \\
0 & 1 & 0  \\
0 & 0 & 1 \\
\end{pmatrix}$,    
$\begin{pmatrix}
0 \\
0  \\
-2z \\
\end{pmatrix}$ & \textbf{G} ($\mathcal{B}_1$) \\

{} & {} & {} & {} \\

{} & $c^{-1}d$ & $\begin{pmatrix}
1 & 0 & 0  \\
0 & 1 & 0  \\
0 & 0 & 1 \\
\end{pmatrix}$,    
$\begin{pmatrix}
4x \\
0  \\
0 \\
\end{pmatrix}$ & {} \\

{} & {} & {} & {} \\

{} & $c^{-1}e$ & $\begin{pmatrix}
1 & 0 & 0  \\
0 & 1 & 0  \\
0 & 0 & -1 \\
\end{pmatrix}$,    
$\begin{pmatrix}
2x \\
-2y  \\
0 \\
\end{pmatrix}$ & {} \\ 

{} & {} & {} & {} \\  \hline

{} & {} & {} & {} \\ 

$\{(0,0,0,1), (0,0,0,-1)\}$ & $j$ & $\begin{pmatrix}
1 & 0 & 0  \\
0 & 1 & 0  \\
0 & 0 & 1 \\
\end{pmatrix}$,    
$\begin{pmatrix}
0 \\
-2y  \\
0 \\
\end{pmatrix}$ & \textbf{G} ($\mathcal{B}_1$) \\

{} & {} & {} & {} \\ 

{} & $g^{-1}h^{-1}$ & $\begin{pmatrix}
1 & 0 & 0  \\
0 & 1 & 0  \\
0 & 0 & 1 \\
\end{pmatrix}$,    
$\begin{pmatrix}
4x \\
0  \\
0 \\
\end{pmatrix}$ & {} \\ 

{} & {} & {} & {} \\ 

{} & $g^{-1}k$ & $\begin{pmatrix}
1 & 0 & 0  \\
0 & -1 & 0  \\
0 & 0 & 1 \\
\end{pmatrix}$,    
$\begin{pmatrix}
2x \\
0  \\
-2z \\
\end{pmatrix}$ & {} \\

{} & {} & {} & {} \\  \hline

{} & {} & {} & {} \\ 

$\{(\pm 1/2, \pm 1/2, \pm 1/2, \pm 1/2)\}$ & $e^{-1}g$ & $\begin{pmatrix}
1 & 0 & 0  \\
0 & 1 & 0  \\
0 & 0 & 1 \\
\end{pmatrix}$,    
$\begin{pmatrix}
0 \\
-4y  \\
0 \\
\end{pmatrix}$ & \textbf{G} ($\mathcal{B}_1$) \\

{} & {} & {} & {} \\ 

{} & $a^{-1}k^{-1}ak$ & $\begin{pmatrix}
1 & 0 & 0  \\
0 & 1 & 0  \\
0 & 0 & 1 \\
\end{pmatrix}$,    
$\begin{pmatrix}
8x \\
0  \\
0 \\
\end{pmatrix}$ & {} \\

{} & {} & {} & {} \\

{} & $a^{-1}k^{-1}j^{-1}fc$ & 
$\begin{pmatrix}
1 & 0 & 0  \\
0 & -1 & 0  \\
0 & 0 & 1 \\
\end{pmatrix}$,    
$\begin{pmatrix}
-4x \\
2y \\
-4z\\
\end{pmatrix}$ & {} \\

{} & {} & {} & {} \\  \hline

\end{tabular} \\ \\ 
 
        \end{minipage}}
\end{table}

\section{Elementary Moves}\label{elementary}

In this section we are going to go through the three main reduction moves we will be using to simplify the Kirby diagram of a filling of
one of the Ratcliffe-Tschantz manifolds. We call these three main reductions moves the ``three elementary moves''.
In \cite{sarat} we showed how to construct the Kirby diagram of manifold no. 3 in the Ratcliffe-Tschantz census. The reader might want to refer to that
paper for the following discussion.

We showed that the Kirby diagram of manifold no. 3 can be viewed via the four following diagrams:

\centerline {\graphicspath{ {manifold_3/2-handle_cycles/} }\includegraphics[width=11cm, height=11cm]{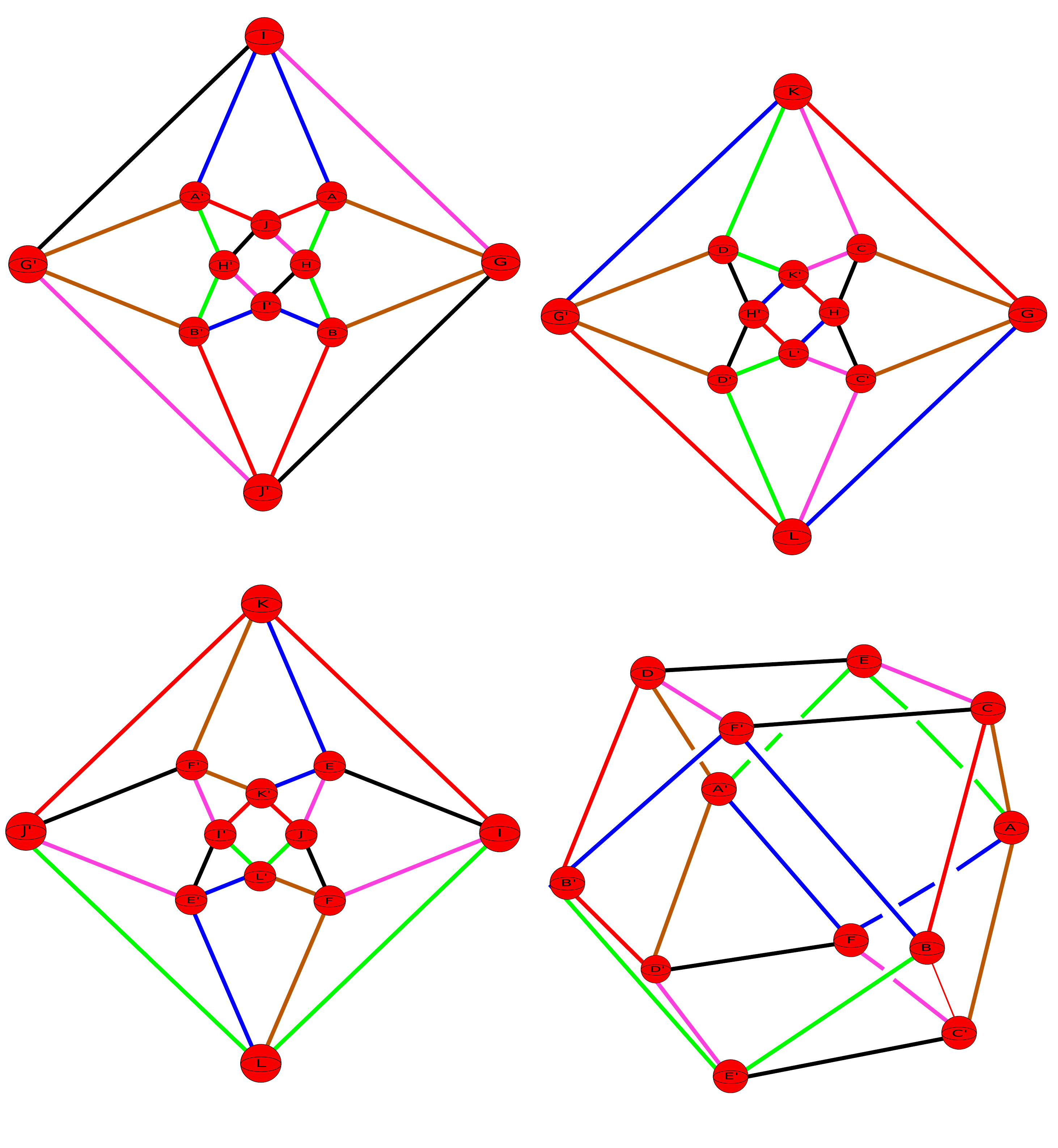}} 

The top two diagrams represent that part of the Kirby diagram that lies in the $x-y$ and $x-z$ planes, reading left to right. The bottom left diagram
represents that part of the diagram contained in the $y-z$ plane. Finally, the bottom right diagram represents the six 2-handles that did not all lie
in a single 2-plane. The attaching circles of the 2-handles are colour coded, each colour representing the attaching circle of one 2-handle. The reader
can refer to \cite{sarat} for an explanation on how we obtain the 2-handle structure. 
 
Consider the boundary component given by the code $\textbf{A}$ (or $\mathcal{G}_1$ in Wolf's notation), this is the 3-torus and it bounds the solid
4-manifold $S^1 \times S^1 \times D^2$. When we glue in the solid 3-torus we need to add one 4-handle, two 3-handles 
(as the Euler characteristic of $T^2 = 0$) and one 2-handle. 
If we go back to the table outlining the generators of each parabolic subgroup corresponding
to each boundary component, we see that the first boundary component labelled $\textbf{A}$ has a generator given by the translation $a$. This means that algebraically
the translation $a$ represents an $S^1$-fibre of the boundary component. Therefore when we glue in a solid torus to this boundary component we can do so along the
$S^1$-fibre corresponding to the translation $a$. Algebraically this means we are killing the transformation 
$a$, hence the 2-handle we are attaching must be a straight line segment running between $A-A'$ once. As the 1-handle $A-A'$ lies in the
$x-y$ plane we draw this 2-handle as a dashed line segment lying in the $x-y$ plane, it must also be added to the diagram showing the six 2-handles that do not
all lie in a single 2-plane.

\centerline {\graphicspath{ {manifold_3/filling_boundary/} }\includegraphics[width=10cm, height=10cm]{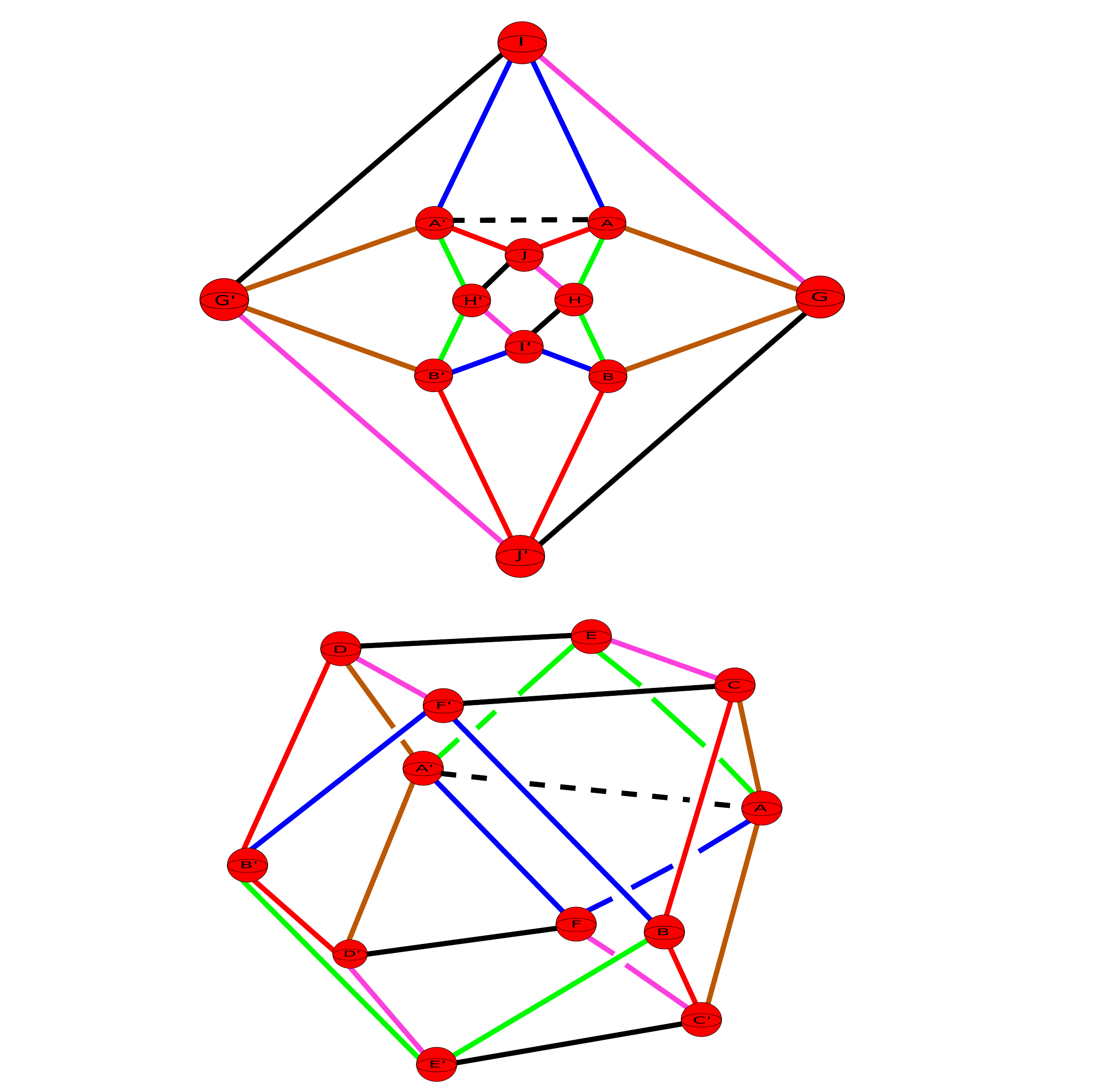}}

There are two subtle points we need to address with this gluing procedure. First of all, we did not explain how the normal bundle to the added 2-handle
looks in the Kirby diagram. In order to understand this one must carry out a similar analysis as was done when trying to understand how
the normal bundles of the 2-handles in the Kirby diagram look like (see Sect.5 in \cite{sarat}). The point is that these added 2-handles
will, most of the time, lie in a single plane, hence the trivialisation of their normal bundle is easy to understand. When the added 2-handle does
not lie in a single plane we will find that it has a planar framing (see Sect.5 in \cite{sarat}), in other words a parallel curve to the 2-handle behaves as if the added 
2-handle was lying in a plane. The second point to address has to do with how exactly we know where to put the added 2-handle in our Kirby diagram. In the above
diagram the added 2-handle lies in the  $x-y$ plane, just before the diagram we said that we can draw this 2-handle as a straight line running between $A,A'$ in the $x-y$ plane because both $A$ and $A'$ lie in the $x-y$ plane. The question is, why is this the right place for the added 2-handle? The basic idea of why this
is the right place has to do with how the fundamental domain of this boundary component looks. It consists of two cubes, one coming from
the ideal vertex $(0,1,0,0)$, and another coming from $(0,-1,0,0)$. The cube centred at $(0,1,0,0)$ has two of its sides being $A$ and $A'$ and the
$S^1$-fibre that we are filling along corresponds to a straight line joining $A$ and $A'$.
When we formed the Kirby diagram we did so by taking the dual polyhedron to
the 24-cell $P$, when we do this the cubes will look like octahedrons in the Kirby diagram, as the dual of a cube
is an octahedron. The point is that we can then go to our Kirby diagram find the associated dual octahedron, and then identify the corresponding $S^1$-fibre
as a straight line running between some 1-handles. If we do this for the above mentioned boundary component we find that the added 2-handle does indeed lie
where we have drawn it in the above diagram. This explanation may seem convoluted, but for now we insist the reader to not pay too much attention to it as the
primary aim of this section is to explain the ``elementary moves'' we will be using to reduce a Kirby diagram, and we don't want the reader to get bogged down
with some minor details. In the next section we will deal with an explicit example, and give full details of exactly how the dual octahedra look like in the Kirby diagram, which in turn will tell us exactly how the added 2-handles look like, and how they should be added.

Coming back to the above diagram we can see that the added 2-handle only passes over $A-A'$ once, in other words the attaching sphere of the added 2-handle transversely intersects the belt sphere of the 1-handle $A-A'$ once.
This means this 1-handle and 2-handle pair form a cancelling pair and can be erased from the diagram. Any other 2-handles that pass over the 1-handle $A-A'$ must first be
slid over the added 2-handle, and then we can erase the pair of handles from the diagram. Observe that because all the added 2-handles, coming from boundary fillings, are 
unknotted and have parallel curves that do not twist around the 2-handle in any way (i.e. they are planar framed), whenever we slide 2-handles over them nothing non-trivial will happen, making the sliding process very straightforward.
This standard handle cancellation move is the \textbf{first elementary move} we will be using to try and reduce our Kirby diagram. 
Let us show how the diagram in the $x-y$ plane changes when we carry out such a cancellation.

\centerline {\graphicspath{ {manifold_3/elementary_moves/} }\includegraphics[width=6cm, height=6cm]{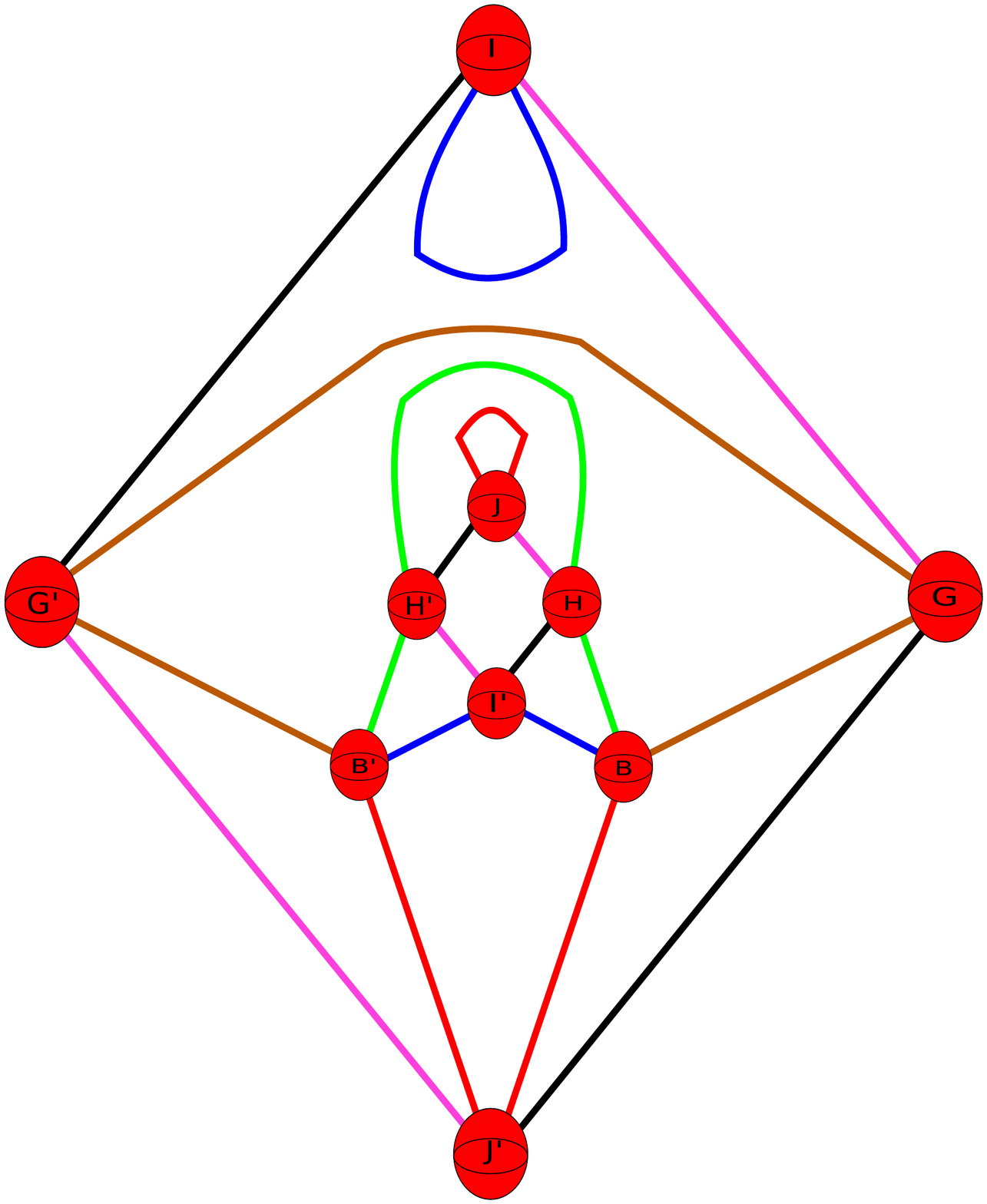}}

After this cancellation has taken place we can see that a few 2-handles have slid into new positions. In particular notice how the blue 2-handle has a component
that now loops back into $I$, similarly the red 2-handle has a component that loops back into $J$. 

The \textbf{second elementary move} we will be making use of is to push such 2-handles through the 1-handle piece they loop back into. For example if
we consider the component of the blue 2-handle that loops back into $I$, we see that we can push it through $I$ to come out as a component that loops back
into $I'$.

\centerline {\graphicspath{ {manifold_3/elementary_moves/} }\includegraphics[width=6cm, height=6cm]{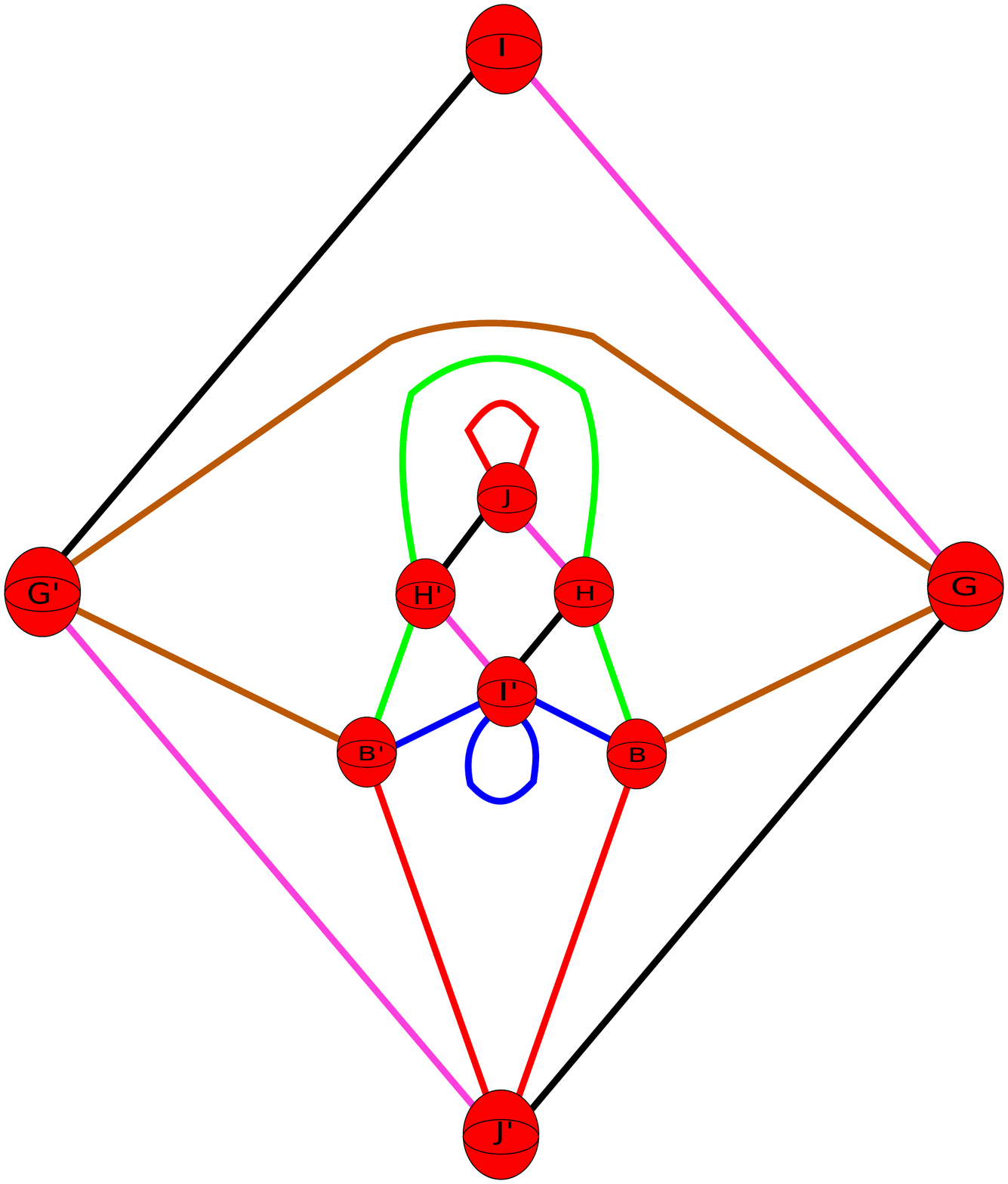}}

We can then slide the blue 2-handle off $I'$, giving a blue 2-handle that runs between $B-B'$ once.

\centerline {\graphicspath{ {manifold_3/elementary_moves/} }\includegraphics[width=6.5cm, height=5.8cm]{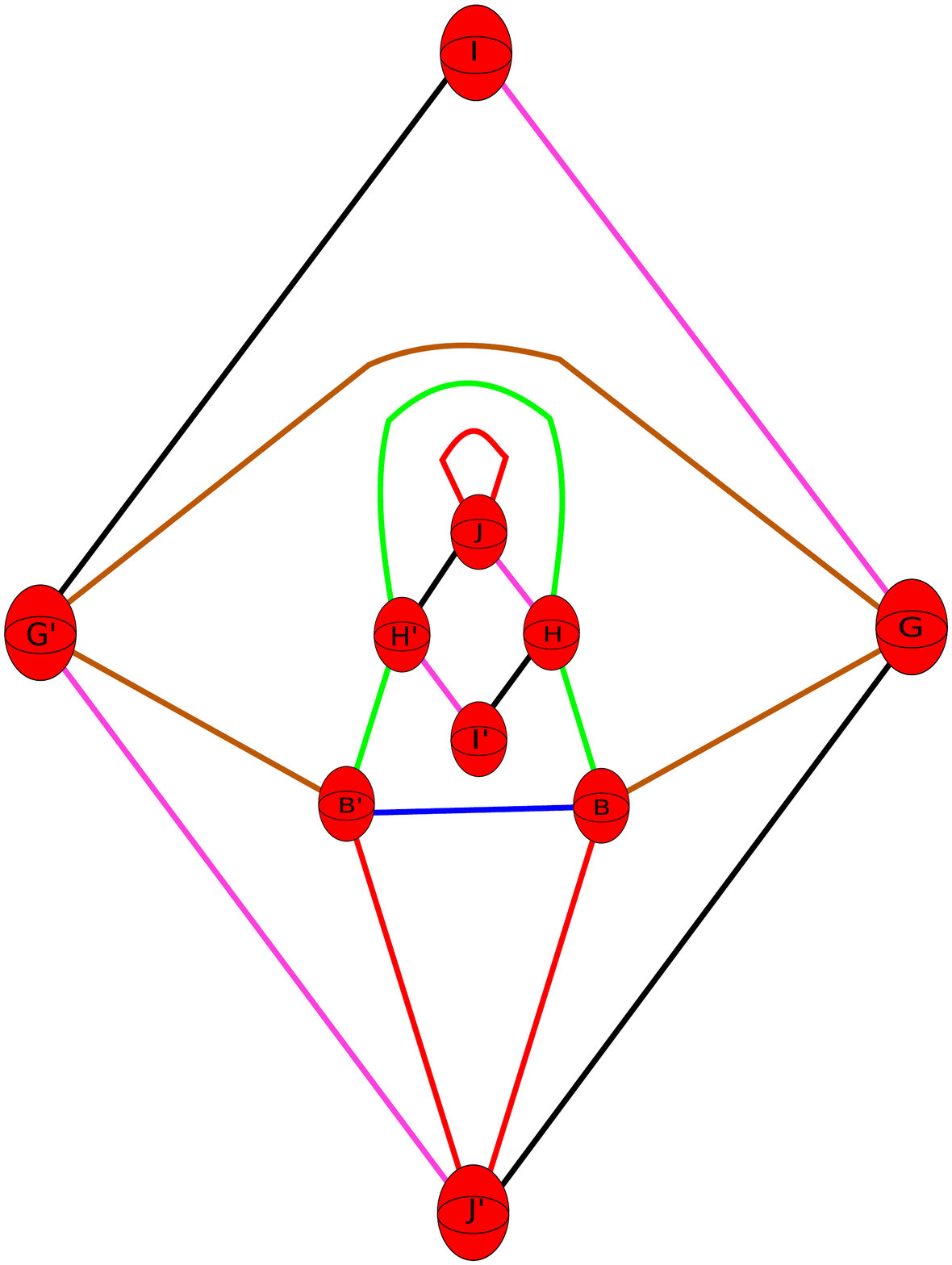}}

We can also carry out the same move for the component of the red 2-handle that loops back into $J$:

\centerline {\graphicspath{ {manifold_3/elementary_moves/} }\includegraphics[width=6cm, height=6cm]{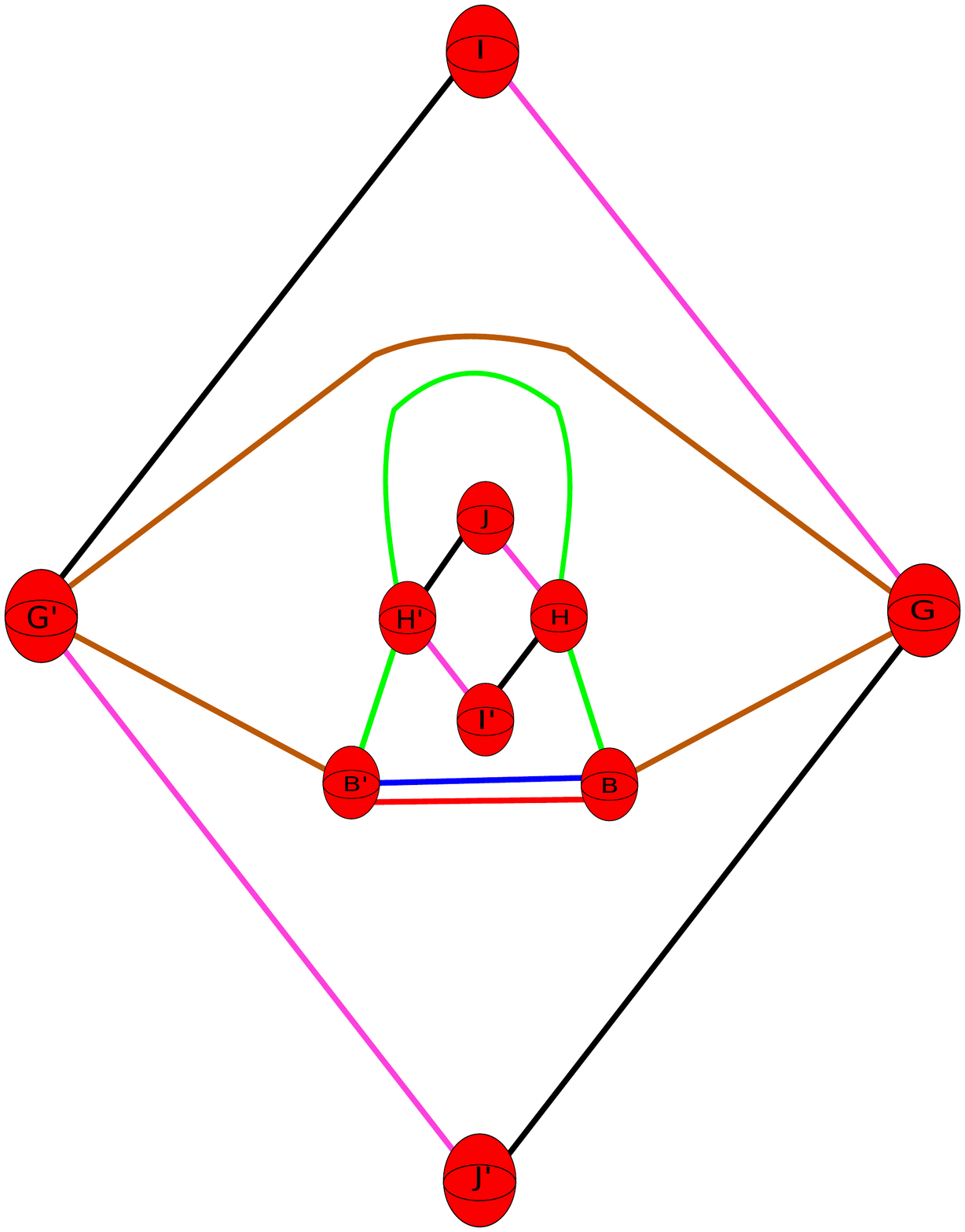}}

In the above move it is very important that we have that the 2-handle we are pushing through
being unknotted. This is because some of our 1-handles are being identified via orientation preserving diffeomorphisms, hence pushing 2-handles through will cause
bits of knots to get mirrored, however if the 2-handles are all unknotted then one does not have to worry about such minor
technicalities.

We now have a blue 2-handle and a red 2-handle passing over the 1-handle $B-B'$ once. The \textbf{third elementary move} we will be making use of can be
described as follows. We can cancel the 1-handle $B-B'$ using the blue 2-handle that now runs over it once, in so doing the red 2-handle
will form an unknotted circle:

\centerline {\graphicspath{ {manifold_3/elementary_moves/} }\includegraphics[width=6cm, height=6cm]{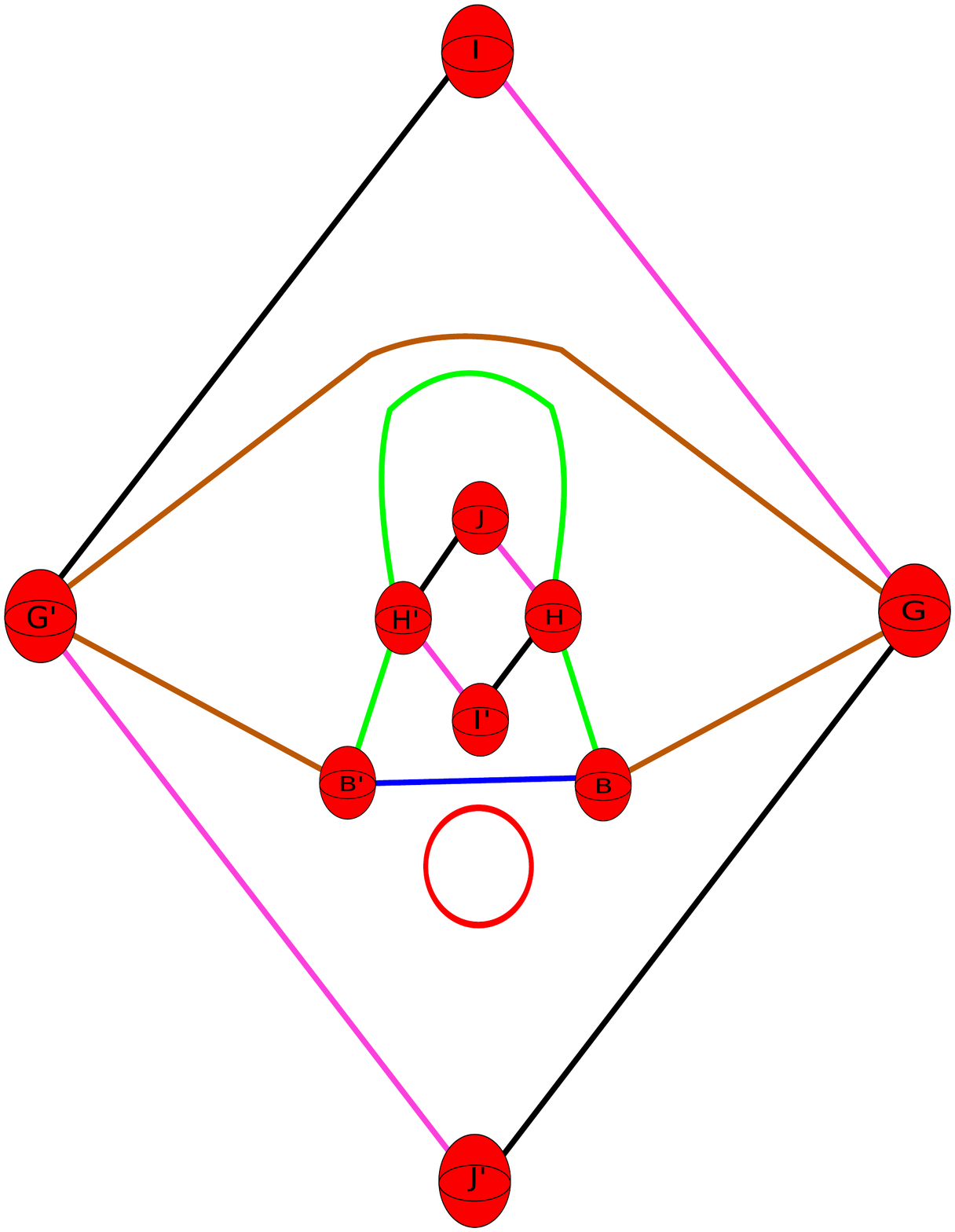}}

Note that this unknotted circle now has a well-defined notion of a framing number. You can see that because the 2-handles in the $x-y$ plane all
have parallel curves that do not twist around them in any way the unknotted circle must have framing zero.
This unknotted circle then cancels a 3-handle and can be deleted from the diagram, a proof of this fact can be found in \cite{gompf} prop.5.1.9, p.148.

The three moves described above will be heavily used in various situations in the sections to come, this is why we took the time to explain each one
carefully. On top of this we will carry out various handle slides, as all our 2-handle are unknotted and have a planar framing the handle slides
we carry out will always be straightforward. If at any point we exploit the use of a non-trivial move we will take the time to carefully explain
how one proceeds.

Filling in the boundary components of any one of the Ratcliffe-Tschantz manifolds produces a compact 4-manifold (possibly non-orientable) for which we know
how to build a Kirby diagram for. Using the elementary moves outlined above one can try and reduce the Kirby diagram of this smooth 4-manifold, the hope is that
after sufficiently many reductions the end Kirby diagram is that of a compact smooth 4-manifold that we can identify. This gives a way of trying to identify
which compact smooth 4-manifolds have smooth complements that are given by the Ratcliffe-Tschantz manifolds, and in turn allows one to obtain explicit
examples of smooth hyperbolic link complements.
In general this proves to be a very difficult task, the reason being is
that we have at least twenty four 2-handles to deal with and we view many of them in certain 2-planes. The added 2-handles (coming from filling boundary components)
will in general pass through some of these 2-planes giving rise to intersection points that have to be carefully tracked when carrying out various elementary moves.
In the next section we will give an explicit example of a situation in which we can identify the filled in 4-manifold up to diffeomorphism.

\section{An explicit example}

In this section we will show that the double cover of the Ratcliffe-Tschantz manifold numbered 1011 is a smooth complement in the standard smooth 4-sphere.
In other words if we perform a boundary filling on the orientable double cover of manifold 1011 we get a closed
smooth 4-manifold that is diffeomorphic to $S^4$.

For the remainder of this section we are going to denote the manifold numbered 1011 in the Ratcliffe-Tschantz census by $M$. The structure of the Kirby diagram
for $M$ was computed in \cite{sarat}.
We remind the reader of the structure of the Kirby diagram for $M$:

The following picture shows those 2-handles lying in the $x-y$ plane, with the table following explaining the colour coding.

\centerline {\graphicspath{ {2-handle_cycles/} }\includegraphics[width=10cm, height=10cm]{2cycle_x-y}}

\begin{tabular}{|l | l | }
\hline
colour & equivalence class  \\ \hline

green & $\xymatrix{
A \cap H \ar[r]^a & A'\cap G'\ar[r]^{g^{-1}} & B\cap G \ar[r]^{b} & B'\cap H' \ar[r]^{h^{-1}} & A\cap H }$ \\ \hline

red & $\xymatrix{
A \cap J \ar[r]^a & A'\cap J\ar[r]^j & B'\cap J' \ar[r]^{b^{-1}} & B\cap J' \ar[r]^{j^{-1}} & A\cap J }$  \\ \hline  

brown & $\xymatrix{
A \cap G \ar[r]^a & A'\cap H' \ar[r]^{h^{-1}} & B\cap H \ar[r]^{b} & B'\cap G' \ar[r]^{g^{-1}} & A\cap G }$  \\ \hline  

blue & $\xymatrix{
A \cap I \ar[r]^a & A'\cap I\ar[r]^i & B'\cap i' \ar[r]^{b^{-1}} & B\cap I' \ar[r]^{i^{-1}} & A\cap I }$ \\ \hline

pink & $\xymatrix{
G \cap I \ar[r]^g & G'\cap J' \ar[r]^{j^{-1}} & G'\cap J \ar[r]^{g^{-1}} & G\cap I' \ar[r]^{i^{-1}} & G\cap I }$  \\   \hline

black & $\xymatrix{
H \cap J \ar[r]^h & H'\cap I' \ar[r]^{i^{-1}} & H'\cap I \ar[r]^{h^{-1}} & H\cap J' \ar[r]^{j^{-1}} & H\cap J }$  \\   \hline

\end{tabular} \\ \\

The 2-handles that lie in the $x-z$ plane can be seen in the following picture.

\centerline {\graphicspath{ {2-handle_cycles/} }\includegraphics[width=10cm, height=10cm]{2cycle_x-z}}

\begin{tabular}{|l | l | }
\hline
colour & equivalence class  \\ \hline

green & $\xymatrix{
D \cap K \ar[r]^d & D'\cap L \ar[r]^{l} & D'\cap L' \ar[r]^{d^{-1}} & D\cap K' \ar[r]^{k^{-1}} & D\cap K }$  \\   \hline

red & $\xymatrix{
G \cap K \ar[r]^g & G'\cap L' \ar[r]^{l^{-1}} & H'\cap L \ar[r]^{h^{-1}} & H\cap K' \ar[r]^{k^{-1}} & G\cap K }$  \\   \hline

brown & $\xymatrix{
C \cap G \ar[r]^c & C'\cap G\ar[r]^g & D\cap G' \ar[r]^{d} & D'\cap G' \ar[r]^{g^{-1}} & C\cap G }$  \\ \hline  

blue & $\xymatrix{
C \cap H \ar[r]^c & C'\cap H \ar[r]^h & D\cap H' \ar[r]^{d} & D'\cap H' \ar[r]^{h^{-1}} & C\cap H }$  \\ \hline  

pink & $\xymatrix{
C \cap K \ar[r]^c & C'\cap L \ar[r]^{l} & C'\cap L' \ar[r]^{c^{-1}} & C\cap K' \ar[r]^{k^{-1}} & C\cap K }$  \\   \hline

black & $\xymatrix{
G \cap L \ar[r]^g & G'\cap K' \ar[r]^{k^{-1}} & H'\cap K \ar[r]^{h^{-1}} & H\cap L' \ar[r]^{l^{-1}} & G\cap L }$  \\   \hline

\end{tabular} \\ \\

The 2-handles that lie in the $y-z$ plane can be seen in the following picture.

\centerline {\graphicspath{ {2-handle_cycles/} }\includegraphics[width=10cm, height=10cm]{2cycle_y-z}}

\begin{tabular}{|l | l | }
\hline
colour & equivalence class   \\ \hline

green & $\xymatrix{
I \cap L \ar[r]^i & I'\cap L \ar[r]^{l} & J'\cap L' \ar[r]^{j^{-1}} & J\cap L' \ar[r]^{l^{-1}} & I\cap L }$  \\   \hline

red & $\xymatrix{
I \cap K \ar[r]^i & I'\cap K \ar[r]^{k} & J'\cap K' \ar[r]^{j^{-1}} & J\cap K' \ar[r]^{k^{-1}} & I\cap K }$  \\   \hline

brown & $\xymatrix{
F \cap L \ar[r]^f & F'\cap K' \ar[r]^{k^{-1}} & F'\cap K \ar[r]^{f^{-1}} & F\cap L' \ar[r]^{l^{-1}} & F\cap L }$  \\   \hline

blue & $\xymatrix{
E \cap K \ar[r]^e & E'\cap L' \ar[r]^{l^{-1}} & E'\cap L \ar[r]^{e^{-1}} & E\cap K' \ar[r]^{k^{-1}} & E\cap K }$  \\   \hline

pink & $\xymatrix{
E \cap J \ar[r]^e & E'\cap I' \ar[r]^{i^{-1}} & F\cap I \ar[r]^{f} & F'\cap J' \ar[r]^{j^{-1}} & E\cap J }$  \\   \hline

black & $\xymatrix{
E \cap I \ar[r]^e & E'\cap J' \ar[r]^{j^{-1}} & F\cap J \ar[r]^{f} & F'\cap I' \ar[r]^{i^{-1}} & E\cap I }$  \\   \hline

\end{tabular} \\ \\

Finally, the 2-handles that do not lie in any one of the above three planes can be seen in the following picture.

\centerline {\graphicspath{ {2-handle_cycles/} }\includegraphics[width=10cm, height=10cm]{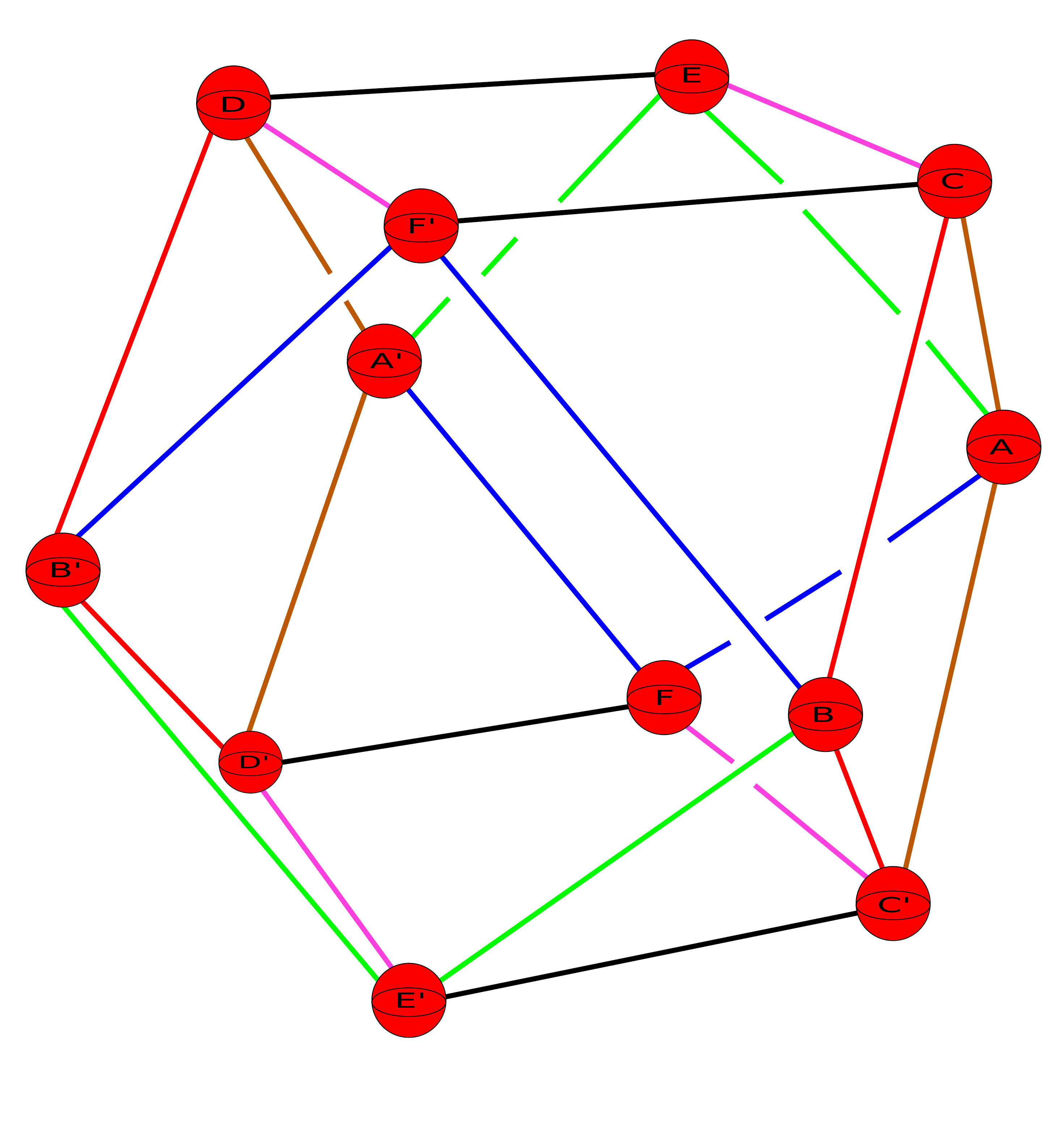}}

\begin{tabular}{|l | l | }
\hline
colour & equivalence class   \\ \hline

green & $\xymatrix{
A \cap E \ar[r]^a & A'\cap E\ar[r]^e & B\cap E' \ar[r]^{b} & B'\cap E' \ar[r]^{e^{-1}} & A\cap E }$  \\ \hline  

red & $\xymatrix{
B \cap C \ar[r]^b & B'\cap D\ar[r]^d & B'\cap D' \ar[r]^{b^{-1}} & B\cap C' \ar[r]^{c^{-1}} & B\cap C }$  \\ \hline  

brown & $\xymatrix{
A \cap C \ar[r]^a & A'\cap D\ar[r]^d & A'\cap D' \ar[r]^{a^{-1}} & A\cap C' \ar[r]^{c^{-1}} & A\cap C }$  \\ \hline

blue &  $\xymatrix{
A \cap F \ar[r]^a & A'\cap F\ar[r]^f & B\cap F' \ar[r]^{b} & B'\cap F' \ar[r]^{f^{-1}} & A\cap F }$  \\ \hline  

pink & $\xymatrix{
C \cap E \ar[r]^c & C'\cap F \ar[r]^f & D\cap F' \ar[r]^{d} & D'\cap E' \ar[r]^{e^{-1}} & C\cap E }$  \\ \hline  

black & $\xymatrix{
C \cap F' \ar[r]^c & C'\cap E' \ar[r]^{e^{-1}} & D\cap E \ar[r]^{d} & D'\cap F \ar[r]^{f} & C\cap F' }$  \\   \hline

\end{tabular} \\ \\

Recall the manifold $M$ is non-orientable and has five cusps each of which has the type \textbf{G} (or in Wolf's notation type $\mathcal{B}_1$). 
We denote the 4-manifold that bounds this 3-manifold by $\widetilde{\textbf{G}}$, remember this is the associated disc bundle to \textbf{G}.
We need to choose a translation 
in each parabolic subgroup that is going to correspond to an $S^1$-fibre which we will fill in. We have already shown the parabolic information corresponding
to each cusp in the table at the end of section \ref{parabolics}. From that table the reader can see that we can take
the first generator associated to each stabiliser subgroup as the translation corresponding to an $S^1$-fibre. That is, we take the transformations
$c$, $a$, $k$, $j$ and $e^{-1}g$. 
We are going to fill in each of these five $S^1$-fibres by gluing in the manifold $\widetilde{\textbf{G}}$, on the level of the above Kirby diagrams
we need to show where the added 2-handles go. 

In order to understand where these added 2-handles will go in our Kirby diagram let us go back and see how we obtained these parabolic
translations. The idea was to take each ideal vertex then take a horospherical neighbourhood about each of these vertices, and then by applying
various isometries we could work out the parabolic translations. For example when we took the ideal vertex $\{(1,0,0,0)\}$ we found that
its equivalence class consisted of two points $\{(1,0,0,0), (-1,0,0,0)\}$, and hence a fundamental domain for the parabolic subgroup corresponding to
this class were two cubes, one centred at the ideal vertex $(1,0,0,0)$ and one centred at $(-1,0,0,0)$. The cube around the ideal vertex
$(1,0,0,0)$ took the form:

\centerline {\graphicspath{ {parabolic_dual/} }\includegraphics[width=5cm, height=4cm]{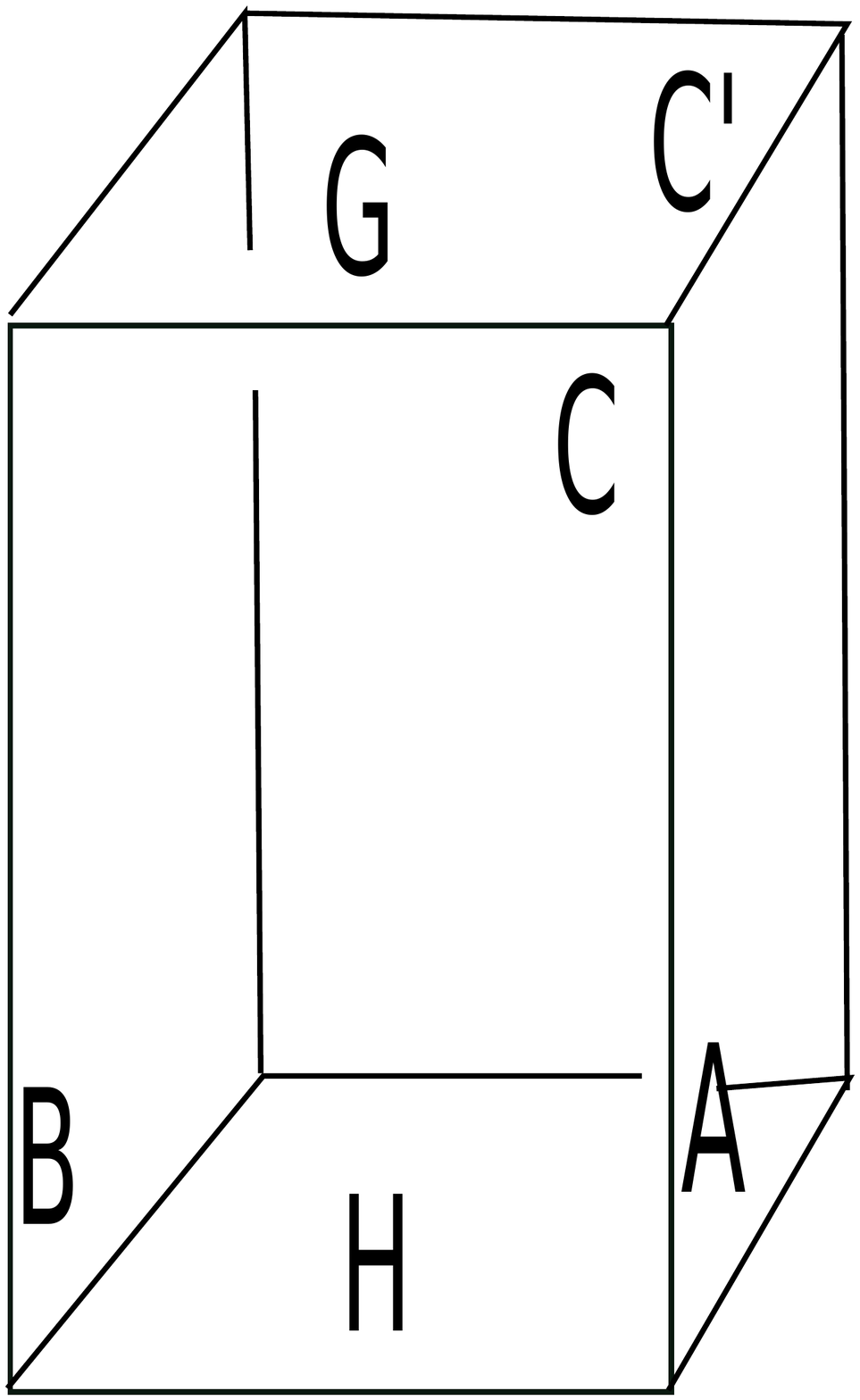}}

We then found that the isometry $c$ was a parabolic translation in the parabolic subgroup corresponding to the class $\{(1,0,0,0), (-1,0,0,0)\}$.
When we formed the handle decomposition of the 24-cell what we were really doing was taking a dual cell structure. This means that if we take the dual
of the above box we will be getting that part of the fundamental domain in the handle decomposition of the 24-cell. As the dual of a cube is an octahedron,
taking the dual of the above cube gives:

\centerline {\graphicspath{ {parabolic_dual/} }\includegraphics[width=6cm, height=6cm]{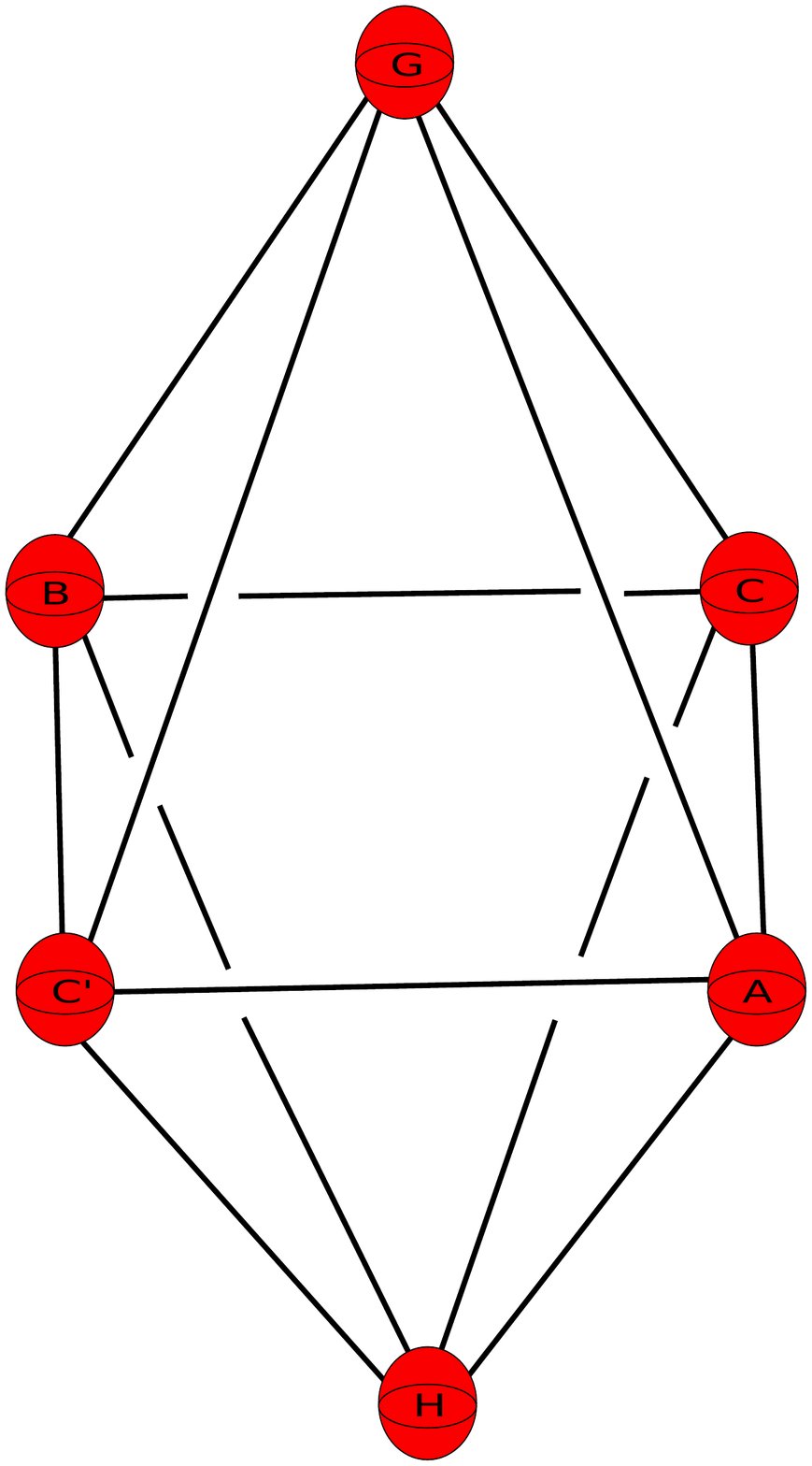}}

Filling in the $S^1$-fibre of the associated boundary component involves attaching a 2-handle from $C$ to $C'$. In our Kirby
diagram this involves drawing an attaching circle between $C$ and $C'$ in our octahedron. Note that the $S^1$-fibre is identified in the fundamental domain
by a straight line joining side $C$ to side $C'$ and running inside the rectangular box making up the fundamental domain. Therefore, in our Kirby diagram
the added 2-handle, representing a filling of the $S^1$-fibre corresponding to the translation $c$, will be a straight line running from
$C$ to $C'$ and contained in the dual octahedron.

\centerline {\graphicspath{ {parabolic_dual/} }\includegraphics[width=6cm, height=6cm]{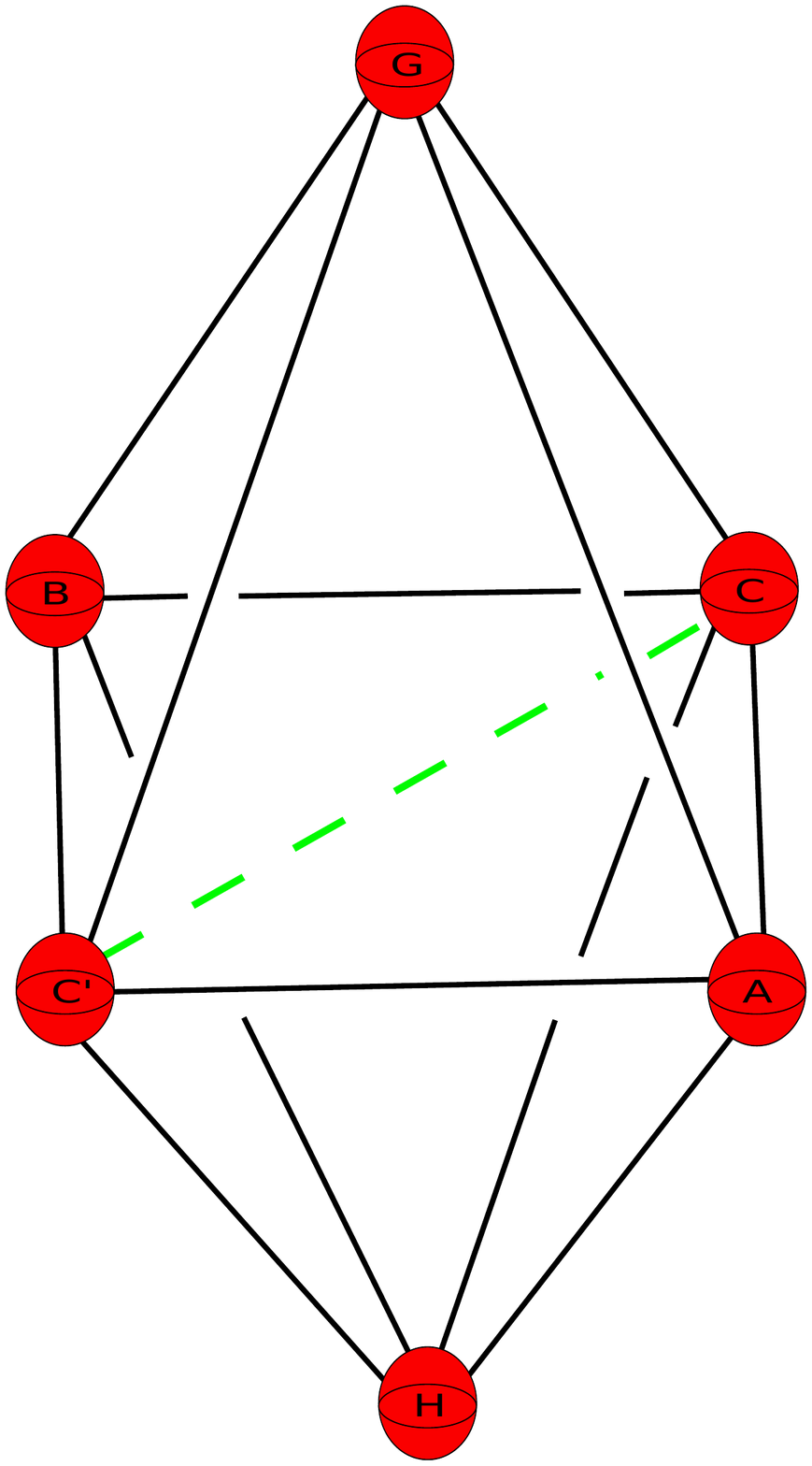}}

We can also see how this looks in part of our ambient handle decomposition diagram. Recall, the coordinates of $C$, $C'$, $A$, $B$, $G$ and $H$ are:

\begin{tabular}{|l|l|l||l|l|l|}
$A$ &  $S_{(+1,+1,0,0)}$ &  $(\frac{1}{\sqrt{2}} ,\frac{1}{\sqrt{2}}, 0)$ &  $B$ & $S_{(+1,-1,0,0)}$ &  $(\frac{1}{\sqrt{2}} ,\frac{-1}{\sqrt{2}}, 0)$ \\ 
 $C$ &  $S_{(+1,0,+1,0)}$ &   $(\frac{1}{\sqrt{2}}, 0, \frac{1}{\sqrt{2}})$ &  $C'$ & $S_{(+1,0,-1,0)}$  &  $(\frac{1}{\sqrt{2}}, 0, \frac{-1}{\sqrt{2}})$ \\ 
 $G$ &  $S_{(+1,0,0,+1)}$ &  $(1 + \sqrt{2}, 0, 0)$ &  $H$ &   $S_{(+1,0,0,-1)}$ &  $(-1 + \sqrt{2}, 0, 0)$ \\ 
 \end{tabular}

The following shows a picture of this octahedron (in purple) in part of the ambient handle decomposition diagram.

\centerline{\graphicspath{ {parabolic_dual/} }\includegraphics[width=9cm, height=6cm]{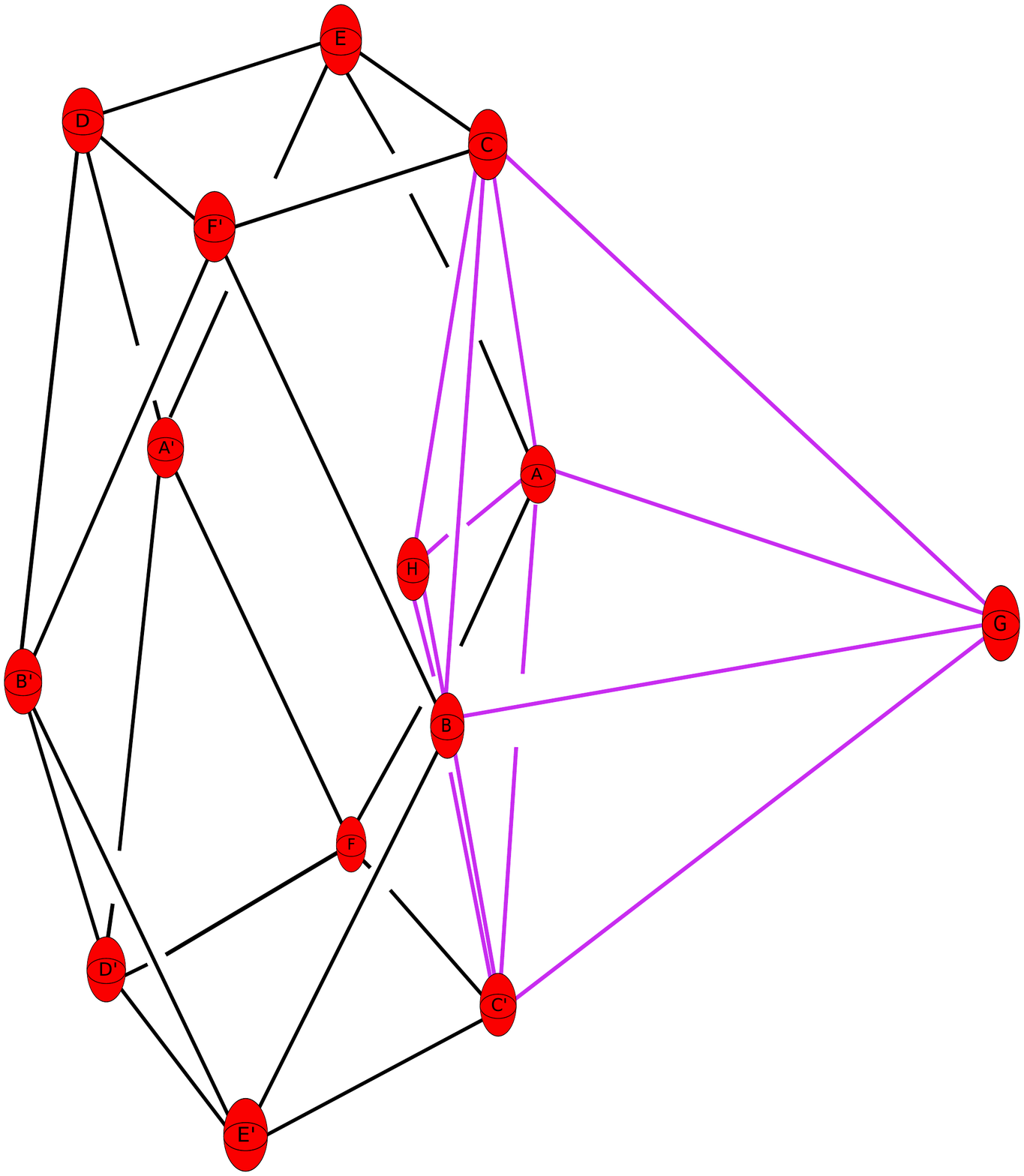}}

The filling of $C-C'$ can be seen in the following picture, with the attaching circle of the 2-handle we are using to do the filling being shown
in green. You can clearly see that it lies in the $x-z$ plane.

\centerline {\graphicspath{ {parabolic_dual/} }\includegraphics[width=9cm, height=6cm]{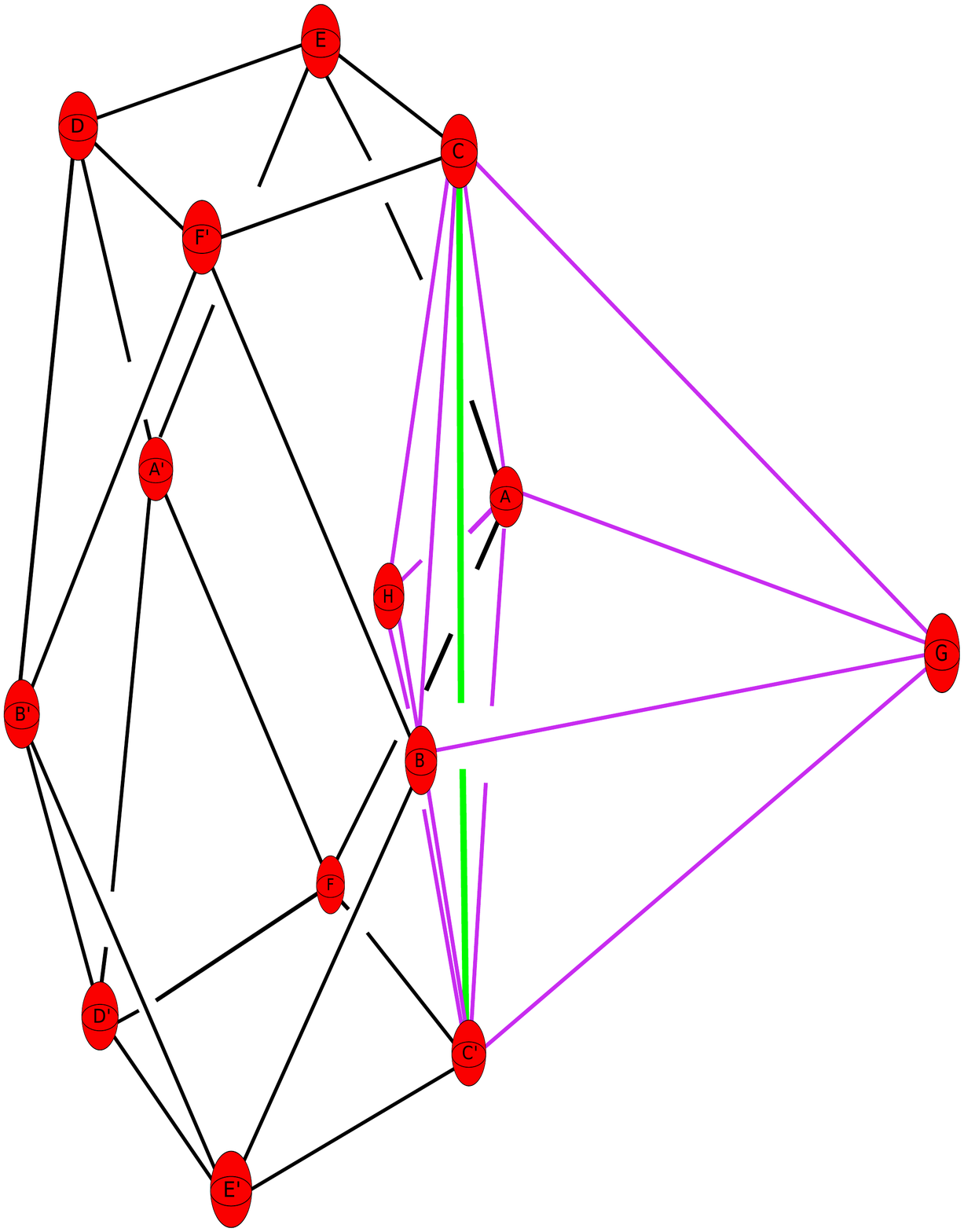}}

We can do the same for all the other fillings we are going to carry out. The reader should keep in mind that each attaching circle of an added 2-handle,
corresponding to a filling, comes from a straight line joining two sides of the rectangular box making up the fundamental domain of the boundary component
we are filling.

If we consider the ideal vertex $\{(0,1,0,0)\}$ we saw that its equivalence class consisted of the two vertices $\{(0,1,0,0), (0,-1,0,0)\}$. The horosphere
about the vertex $\{(0,1,0,0)\}$ looks like:

\centerline {\graphicspath{ {parabolic_dual/} }\includegraphics[width=5cm, height=4cm]{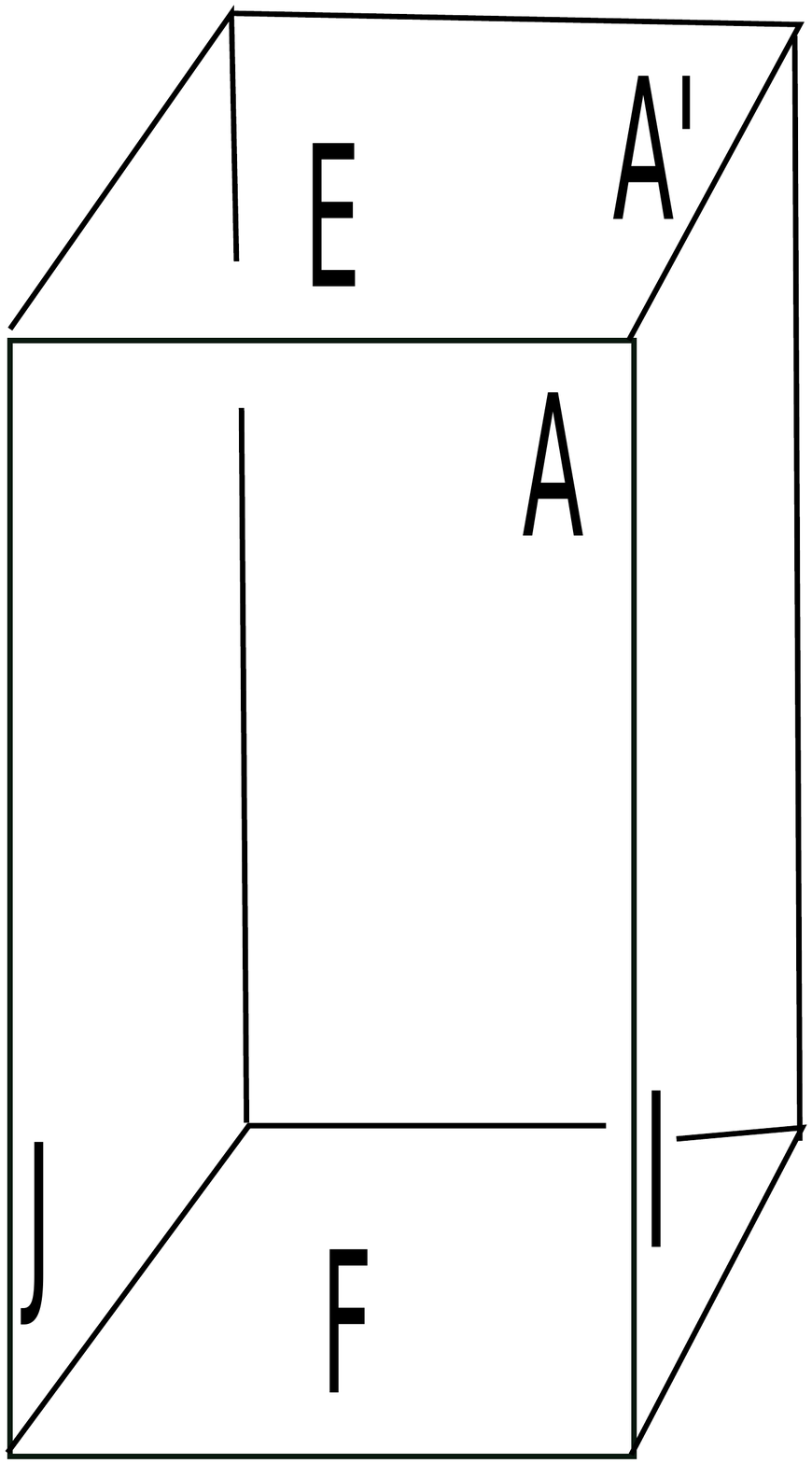}}

Taking the dual gives the following octahedron:

\centerline {\graphicspath{ {parabolic_dual/} }\includegraphics[width=6cm, height=6cm]{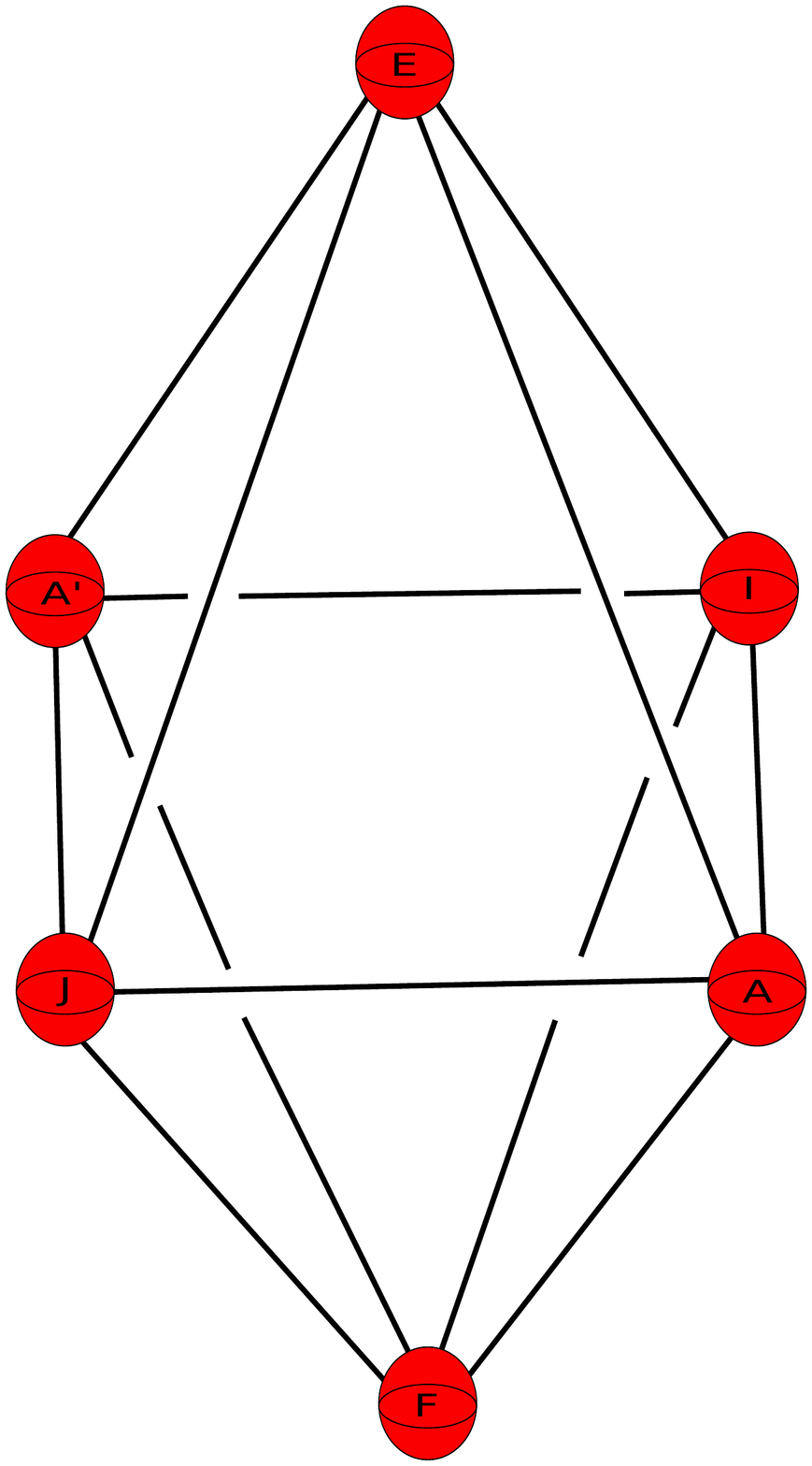}}

Filling in the $S^1$ fibre of the associated boundary component involves attaching a 2-handle from $A$ to $A'$ in the following way.

\centerline{\graphicspath{ {parabolic_dual/} }\includegraphics[width=6cm, height=6cm]{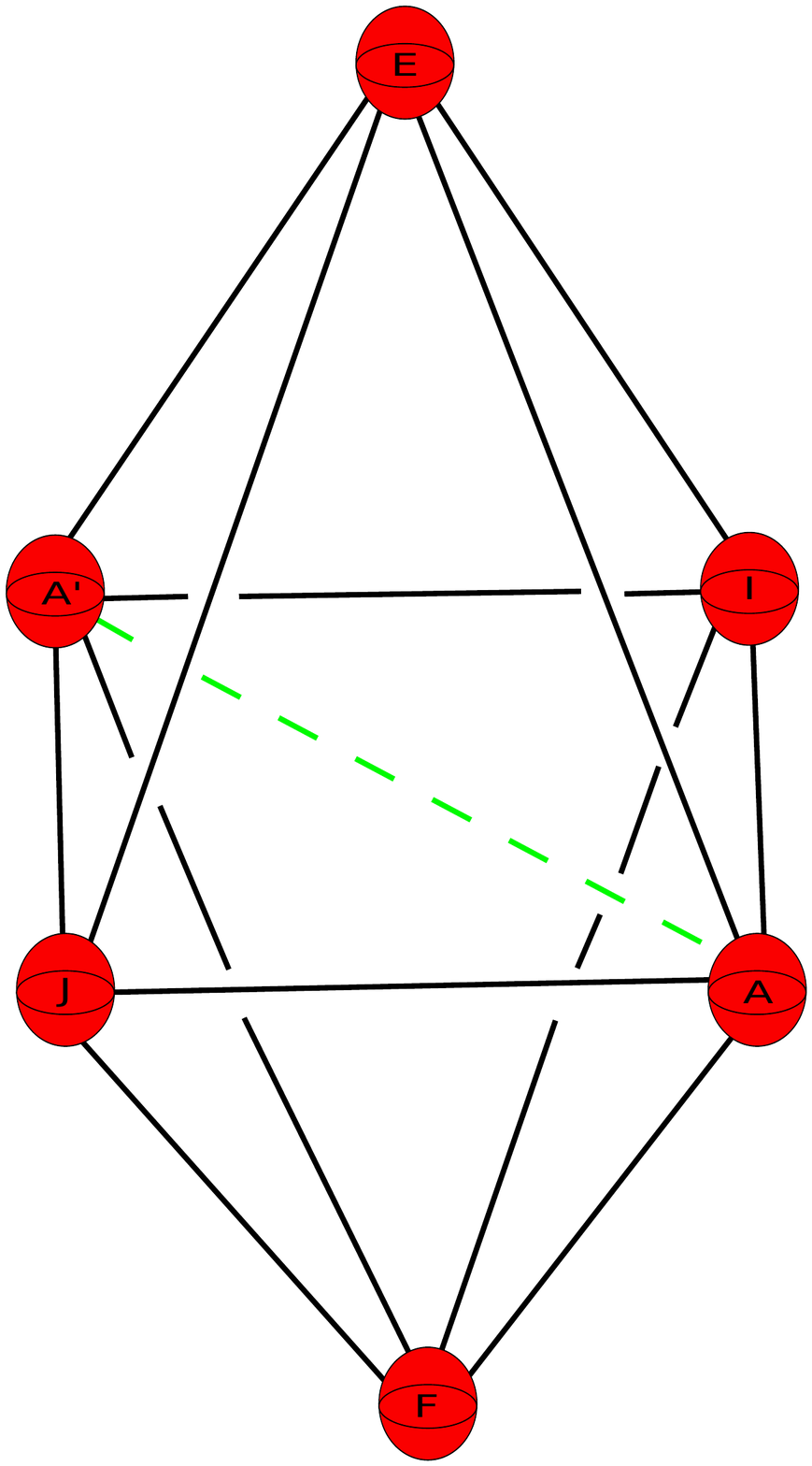}}

A picture of this dual octahedron in part of the three dimensional handle decomposition diagram is shown in the following picture, with the octahedron given by the purple edges.

\centerline {\graphicspath{ {parabolic_dual/} }\includegraphics[width=6cm, height=6cm]{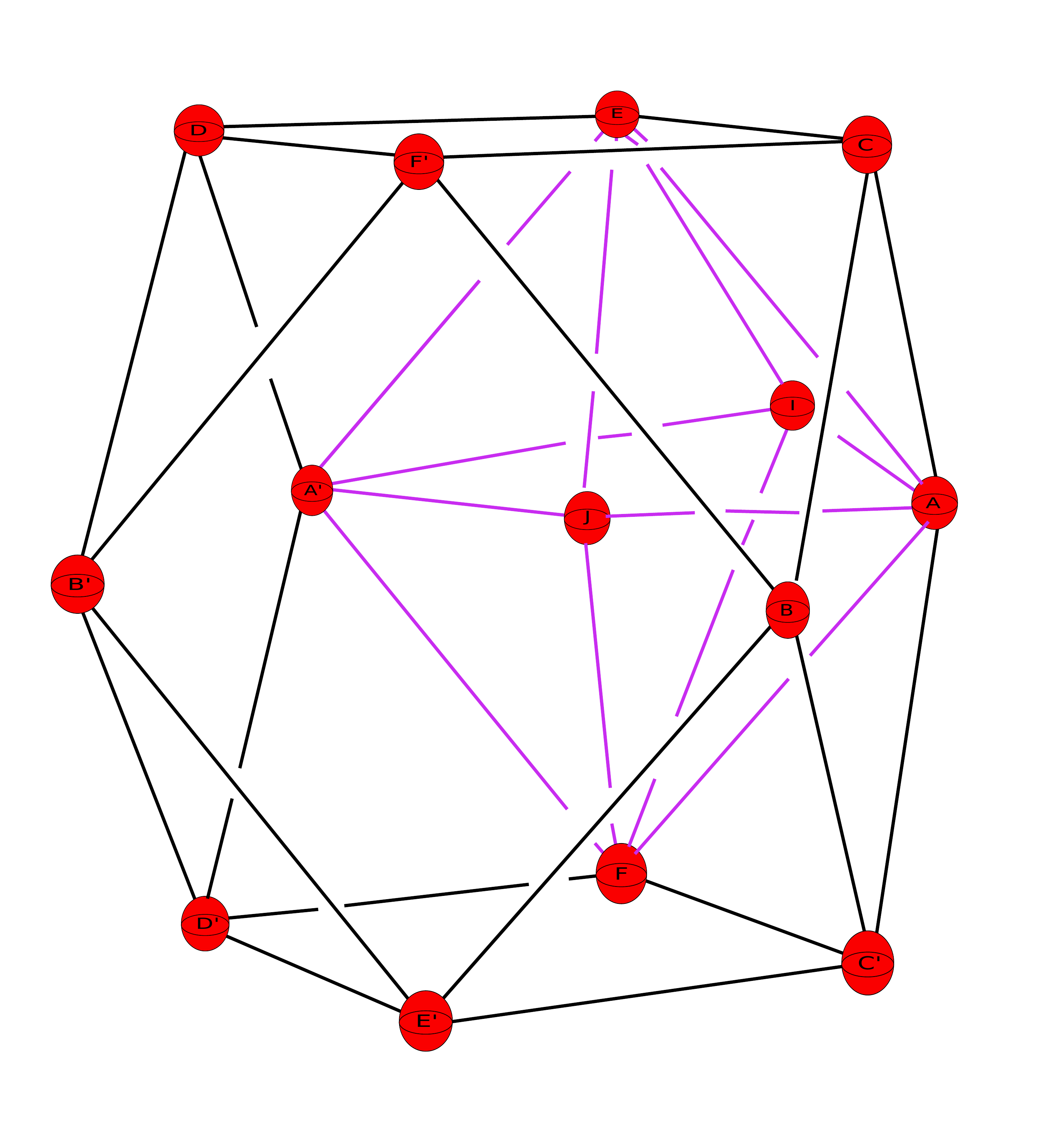}}

The filling of $A-A'$ can be seen in the following picture, with the added 2-handle in green. One can clearly see that it lies in the $x-y$ plane, and does
not interfere with the other 2-handles.

\centerline {\graphicspath{ {parabolic_dual/} }\includegraphics[width=6cm, height=6cm]{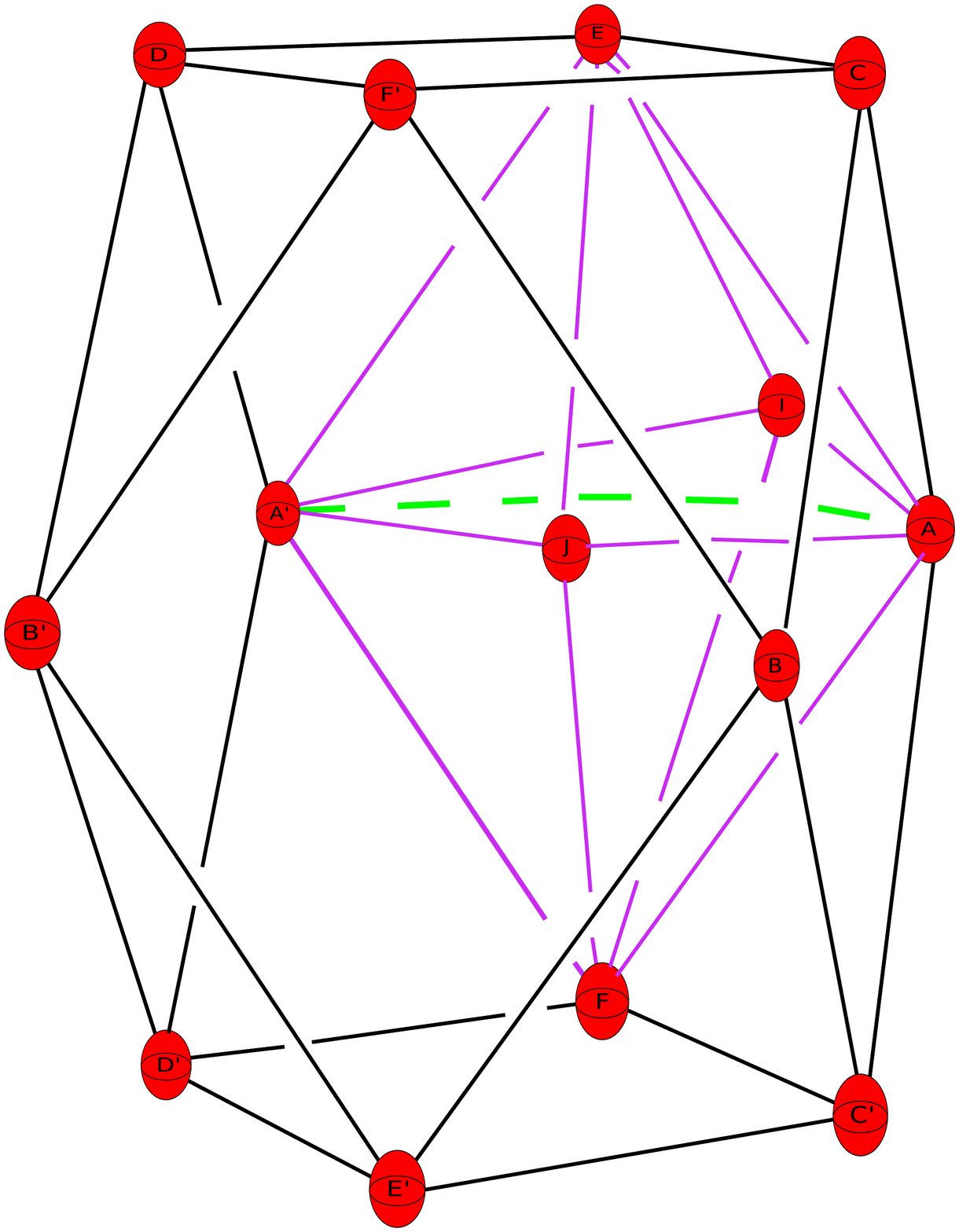}}

For the ideal vertex $(0,0,1,0)$ we found that the equivalence class was $\{(0,0,1,0), (0,0,-1,0)\}$, the horosphere about the vertex $(0,0,1,0)$ looks like:

\centerline{\graphicspath{ {parabolic_dual/} }\includegraphics[width=5cm, height=4cm]{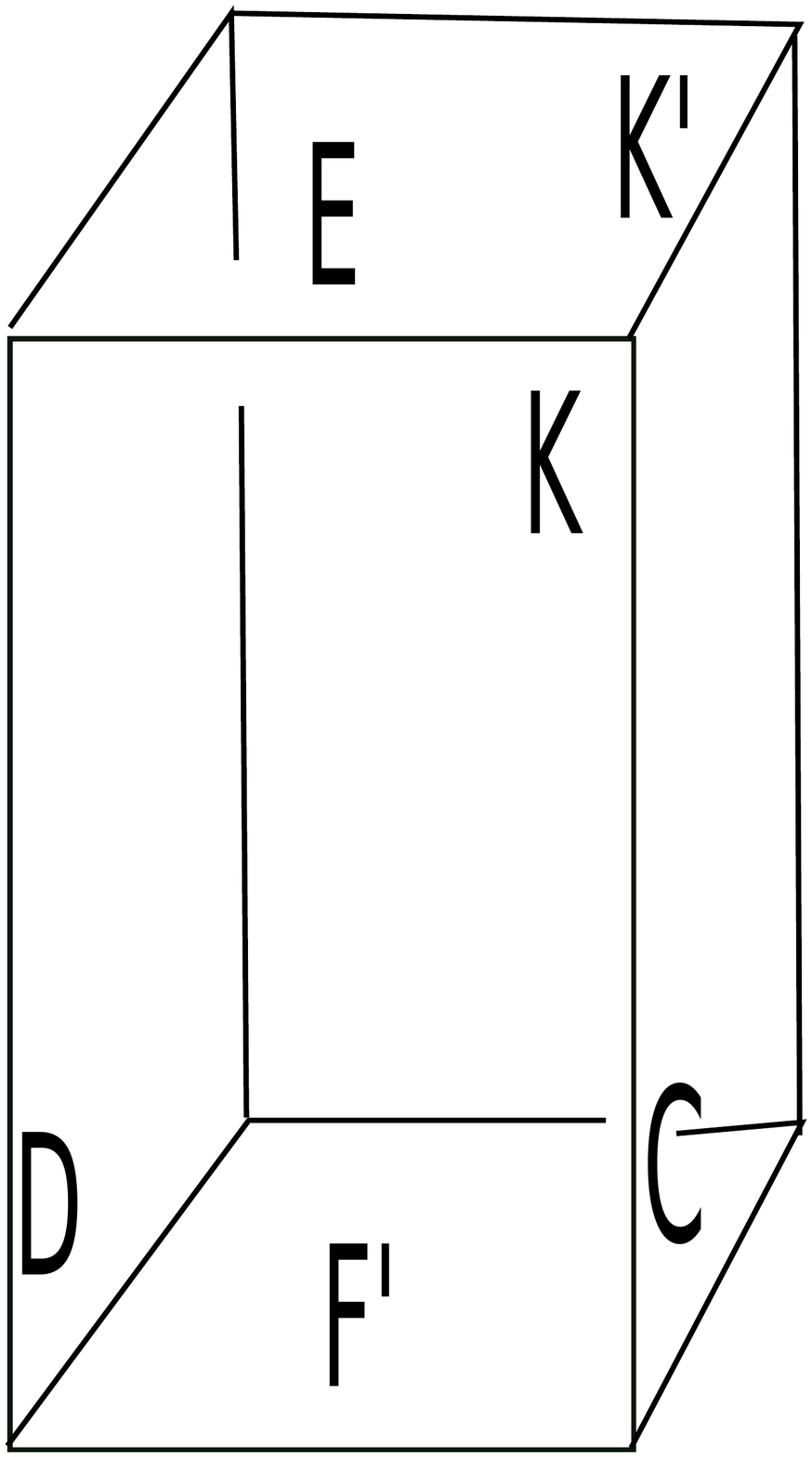}}

The dual octahedron then takes the form:

\centerline {\graphicspath{ {parabolic_dual/} }\includegraphics[width=6cm, height=6cm]{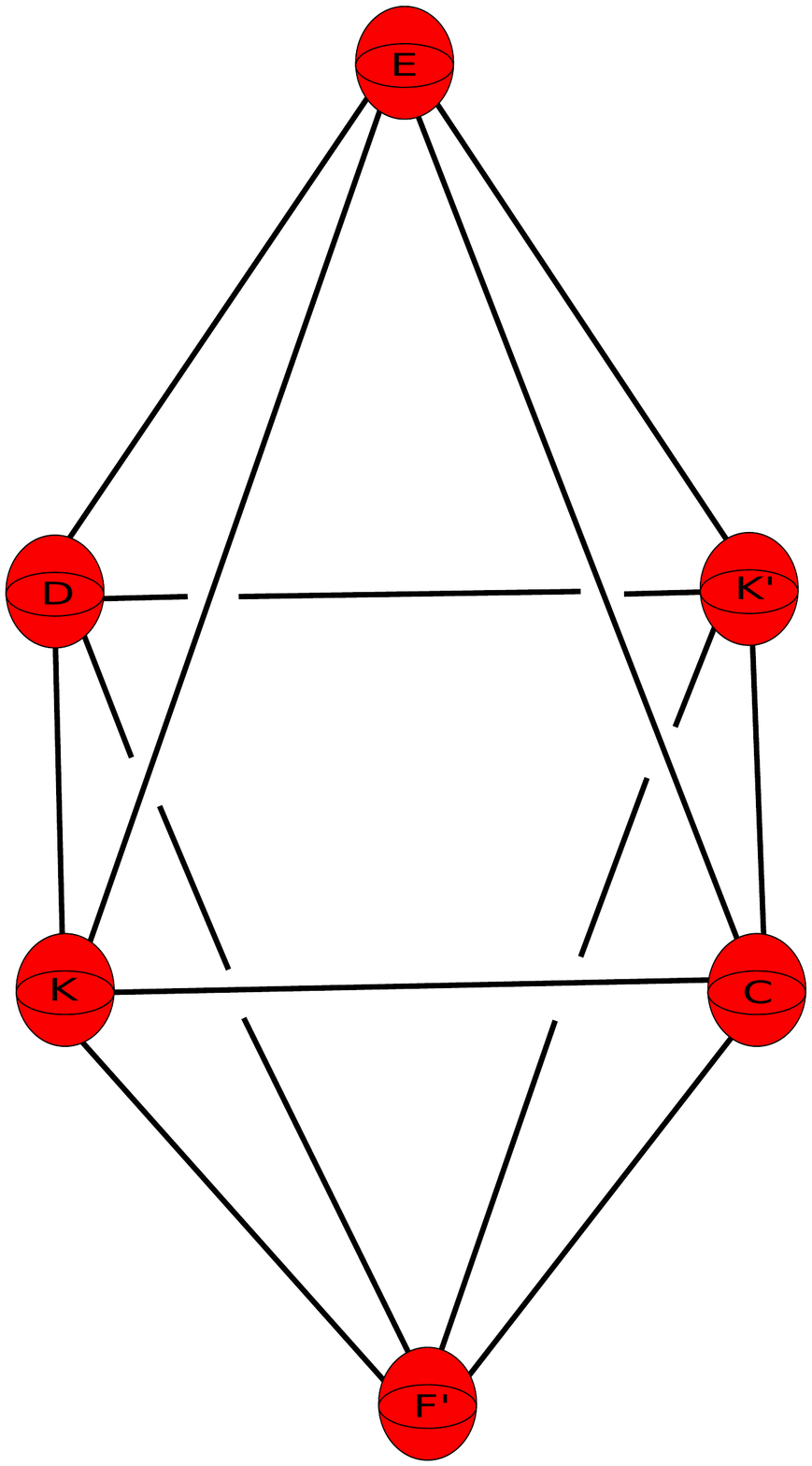}}

Filling in the $S^1$ fibre of the associated boundary component involves attaching a 2-handle from $K$ to $K'$, and in this case we take the following
attaching circle.

\centerline {\graphicspath{ {parabolic_dual/} }\includegraphics[width=6cm, height=6cm]{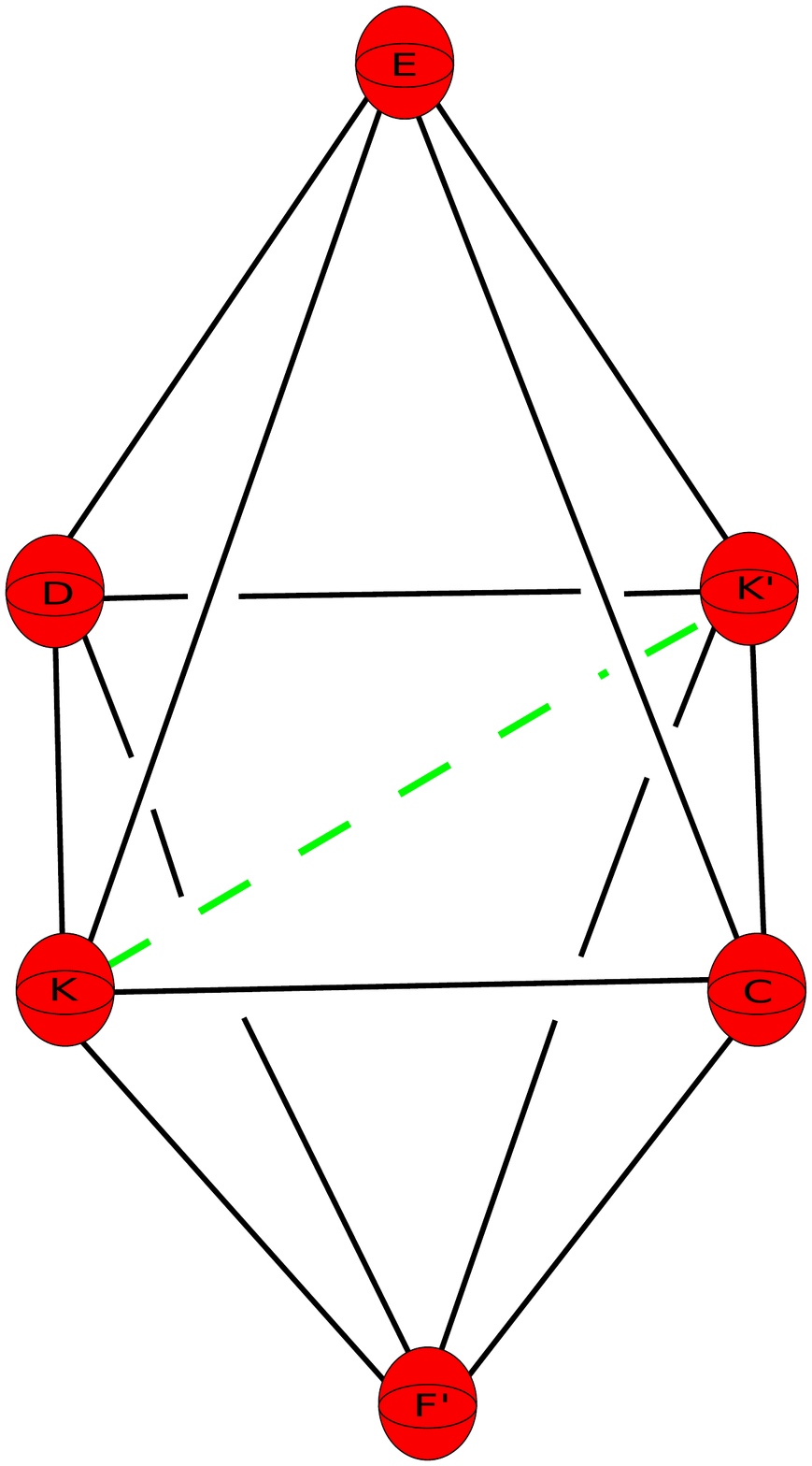}}

A picture of this dual octahedron in part of the handle decomposition is shown in the following picture, with the octahedron given by the purple edges

\centerline {\graphicspath{ {parabolic_dual/} }\includegraphics[width=9cm, height=7cm]{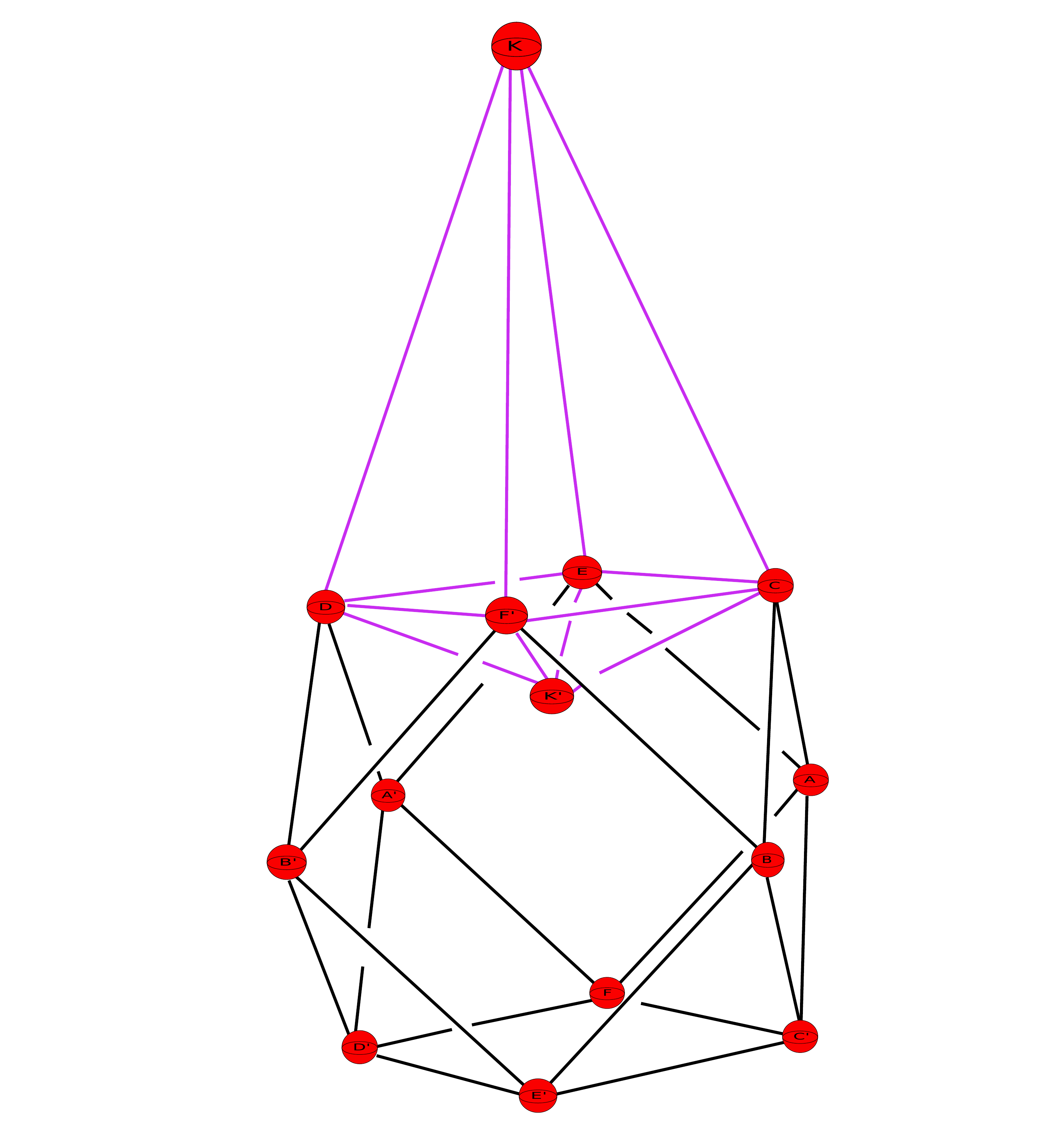}}

The attaching circle of the 2-handle we are attaching is shown in the following picture, it is clear that this 2-handle lies in the $x-z$ and $y-z$ planes, and
does not interfere with any of the other 2-handles.

\centerline{\graphicspath{ {parabolic_dual/} }\includegraphics[width=9cm, height=8cm]{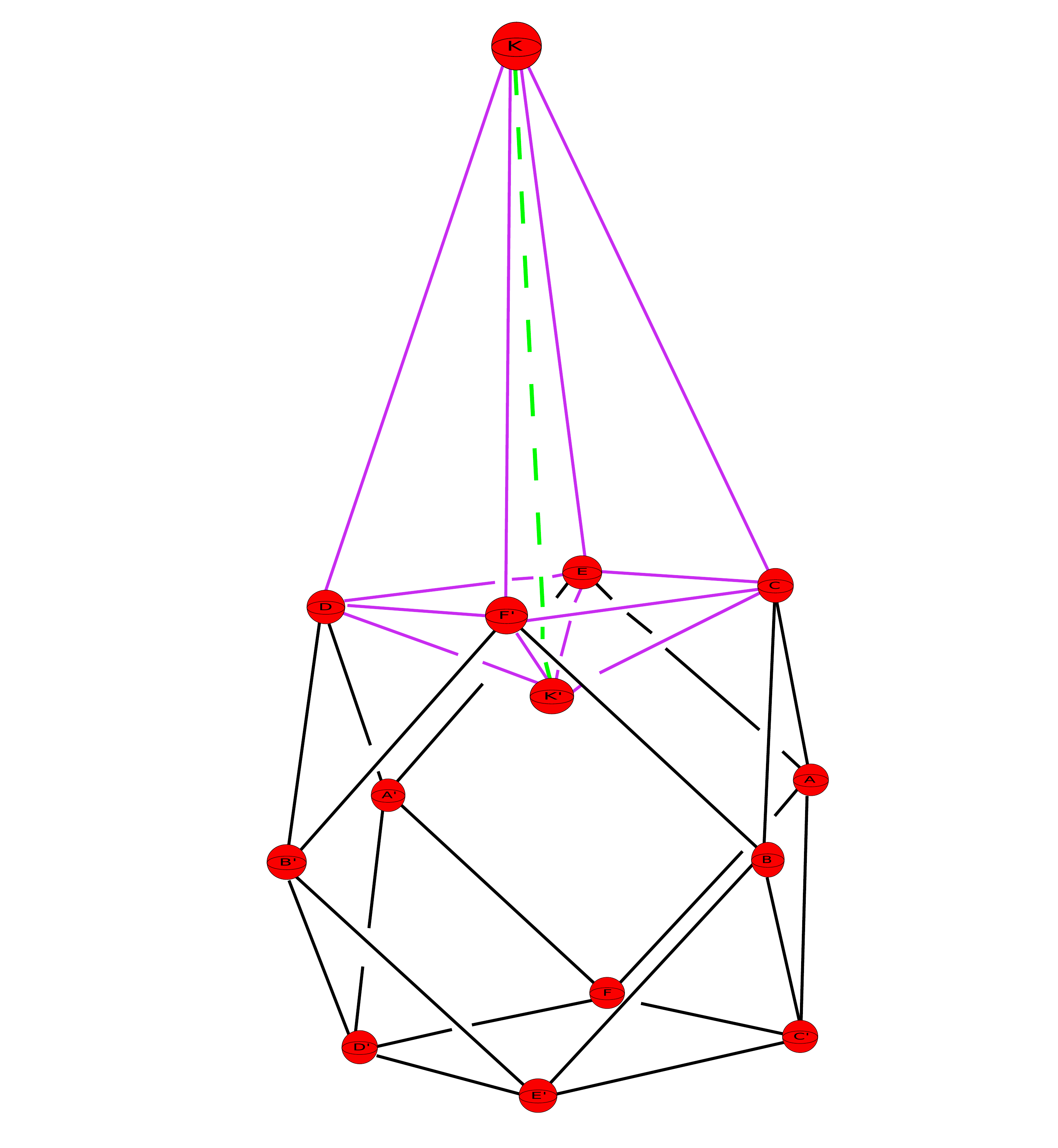}}

For the vertex class $\{(0,0,0,1), (0,0,0,-1)\}$ we found that the isometry $j$ was a parabolic translation. In this case we need to look at the
horosphere about the ideal vertex $(0,0,0,-1)$.

\centerline{\graphicspath{ {parabolic_dual/} }\includegraphics[width=5cm, height=3.5cm]{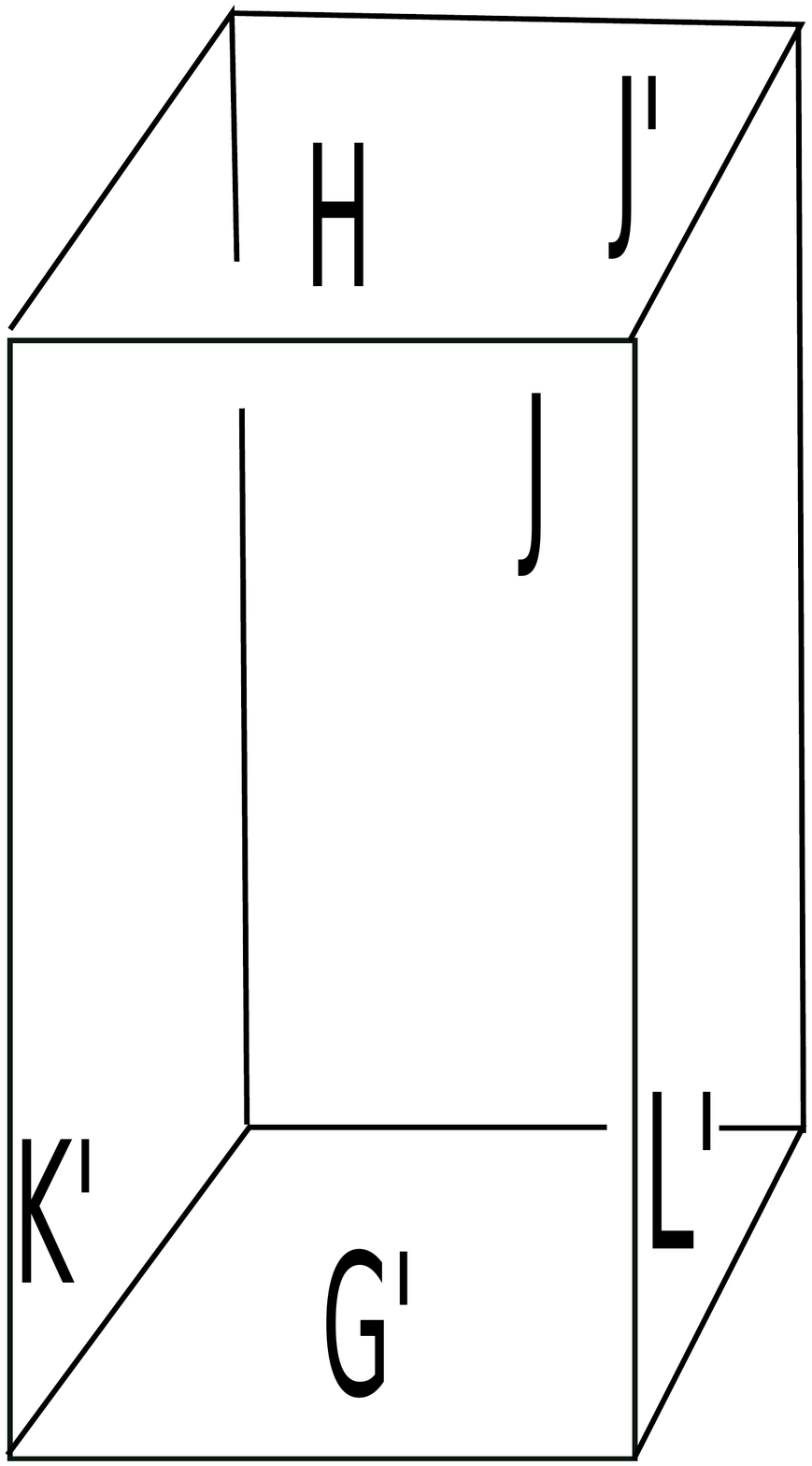}}

The dual octahedron is then given by

\centerline{\graphicspath{ {parabolic_dual/} }\includegraphics[width=6cm, height=6cm]{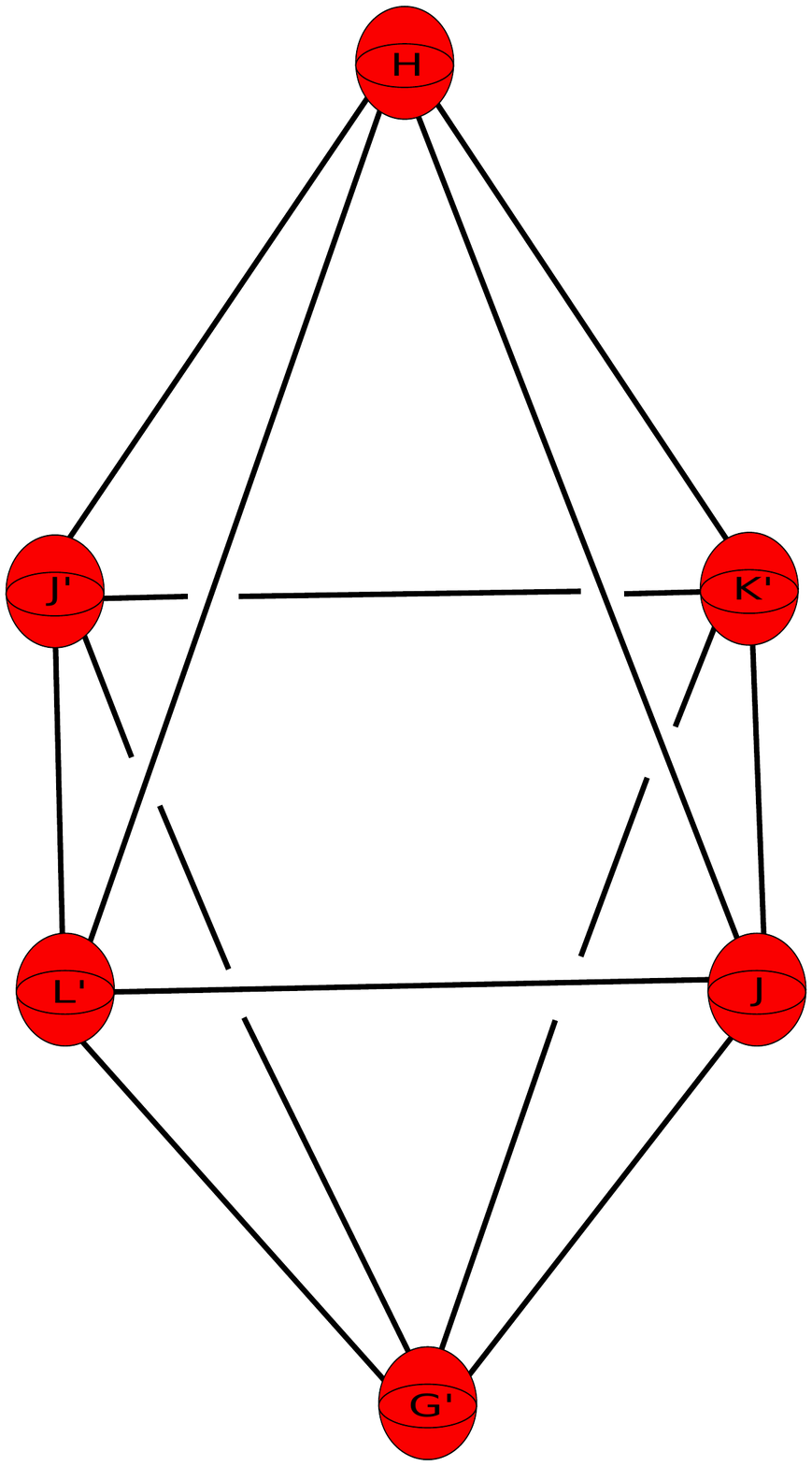}}

and the added 2-handle can be seen in the following picture.

\centerline{\graphicspath{ {parabolic_dual/} }\includegraphics[width=6cm, height=6cm]{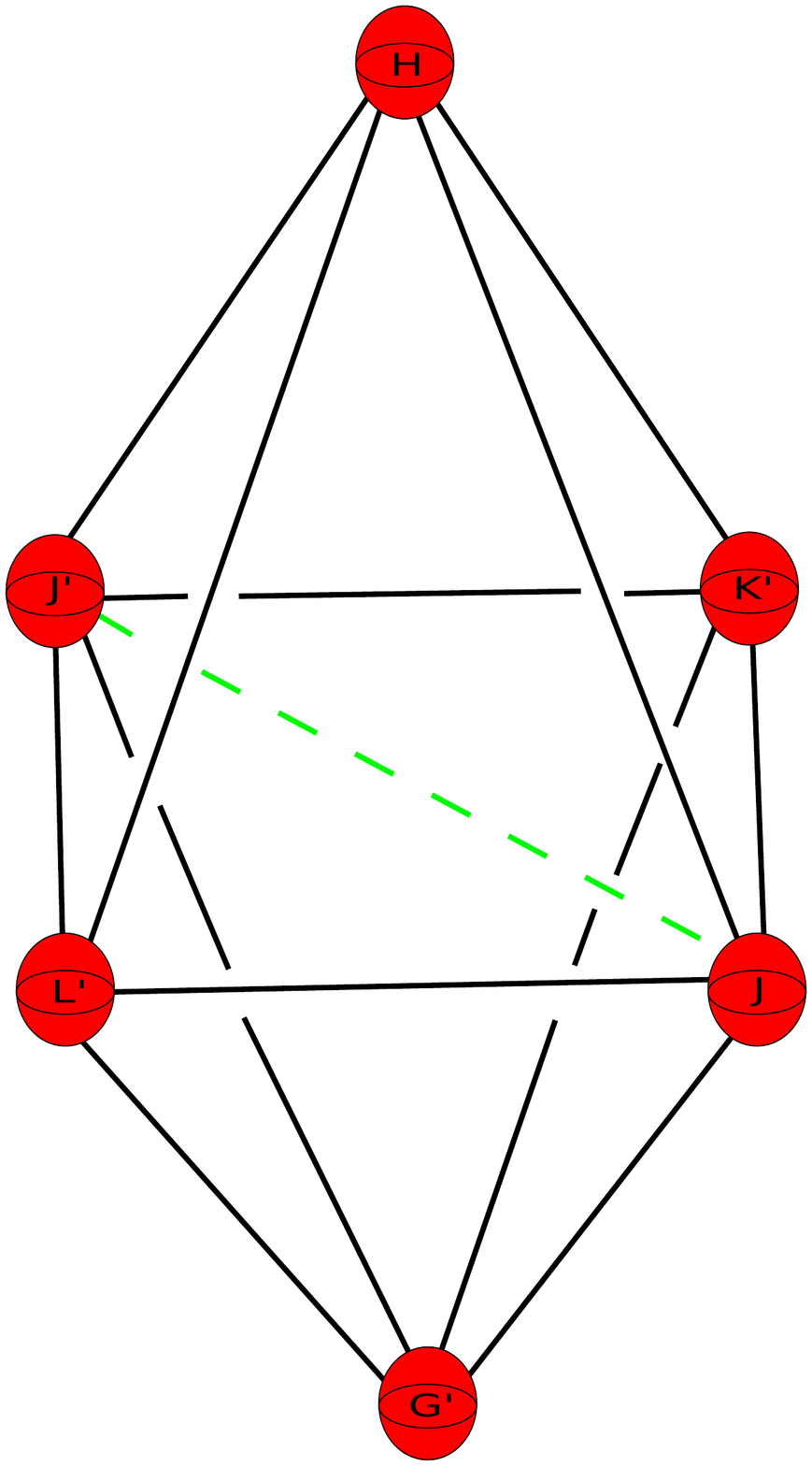}}

A picture of this dual octahedron in part of the handle decomposition is shown in the following picture, with the octahedron given by the purple edges.

\centerline{\graphicspath{ {parabolic_dual/} }\includegraphics[width=7cm, height=6cm]{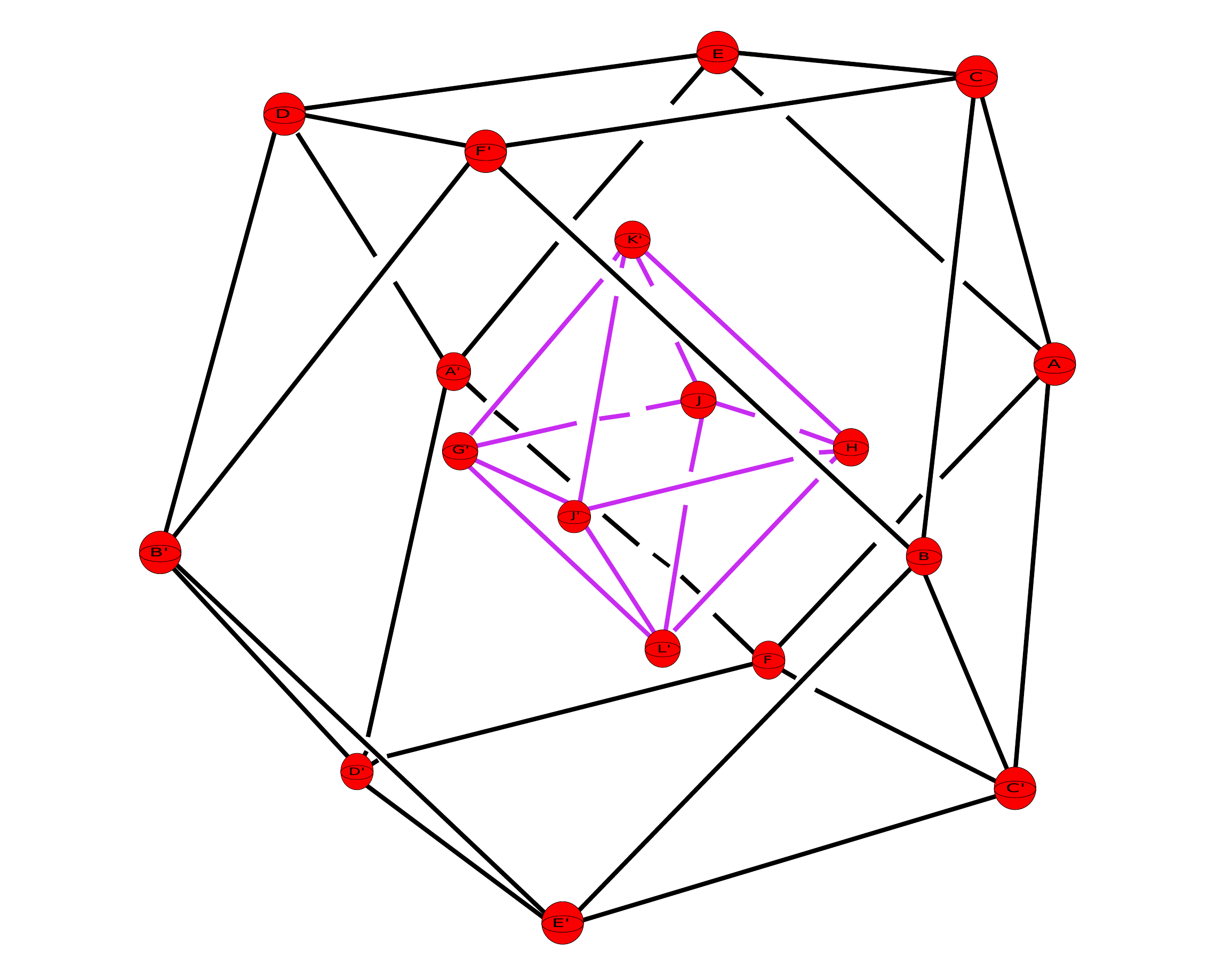}}

The attaching circle of the 2-handle we are attaching is shown in the following picture, it is clear that this 2-handle lies in the $x-y$ and $y-z$ planes, and
does not interfere with any of the other 2-handles.

\centerline{\graphicspath{ {parabolic_dual/} }\includegraphics[width=8cm, height=7cm]{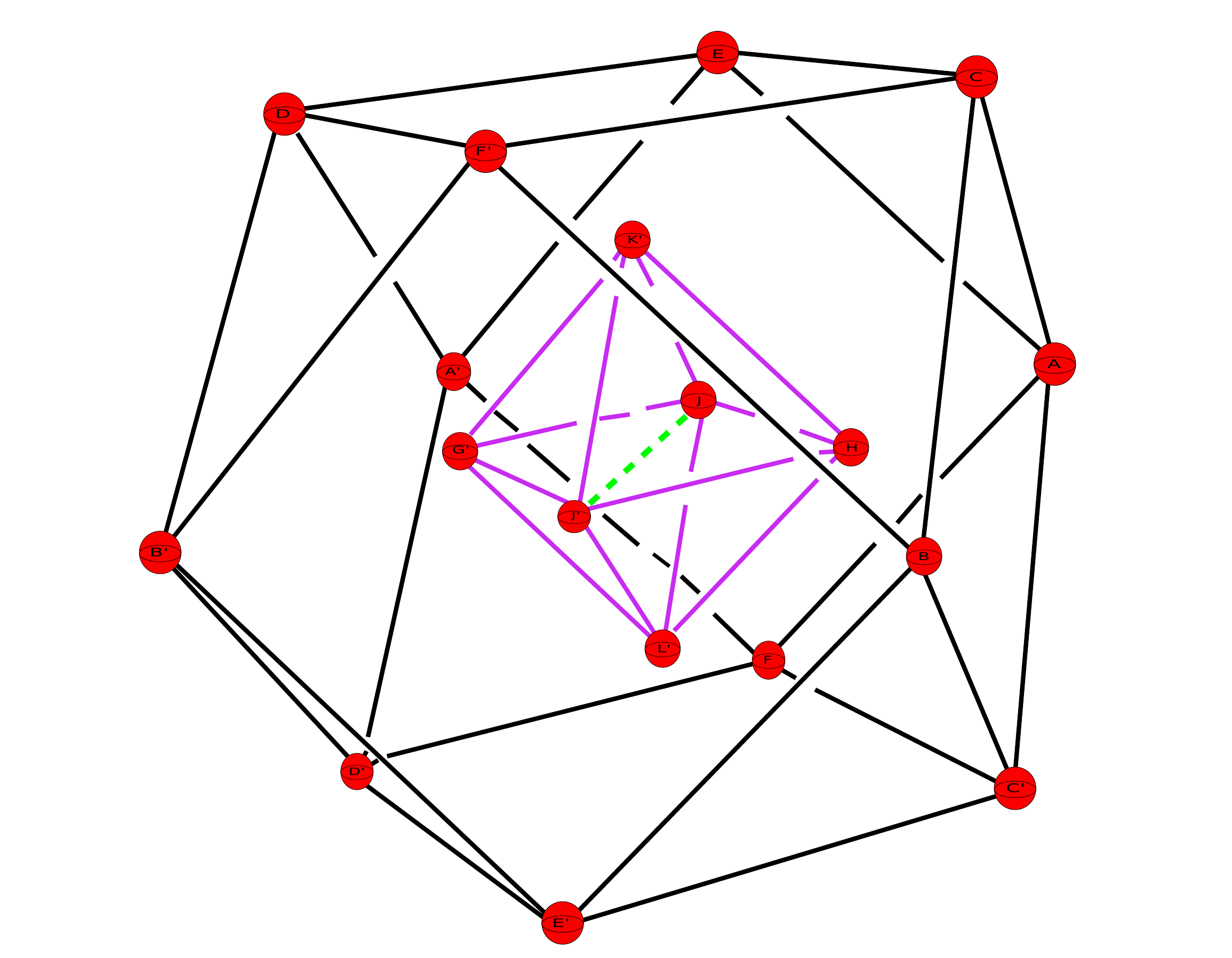}}

The final vertex class to consider consists of the ideal vertices $\{(\pm 1/2, \pm 1/2, \pm 1/2, \pm 1/2)\}$. In this case we found that
the isometry $e^{-1}g$ corresponded to a parabolic translation. Therefore, when we fill the fibre corresponding to this isometry we need to add a 2-handle component
from $E$ to $G$ and one from $E'$ to $G'$. In this case we need to consider two horospheres, the one about the vertex $(1/2, 1/2, 1/2, 1/2)$ and the one
about the vertex $(-1/2, -1/2, -1/2, -1/2)$. Both these horosphere are shown in the following picture, the one on the left corresponding to
$(1/2, 1/2, 1/2, 1/2)$.

\centerline{\graphicspath{ {parabolic_dual/} }\includegraphics[width=5cm, height=3.5cm]{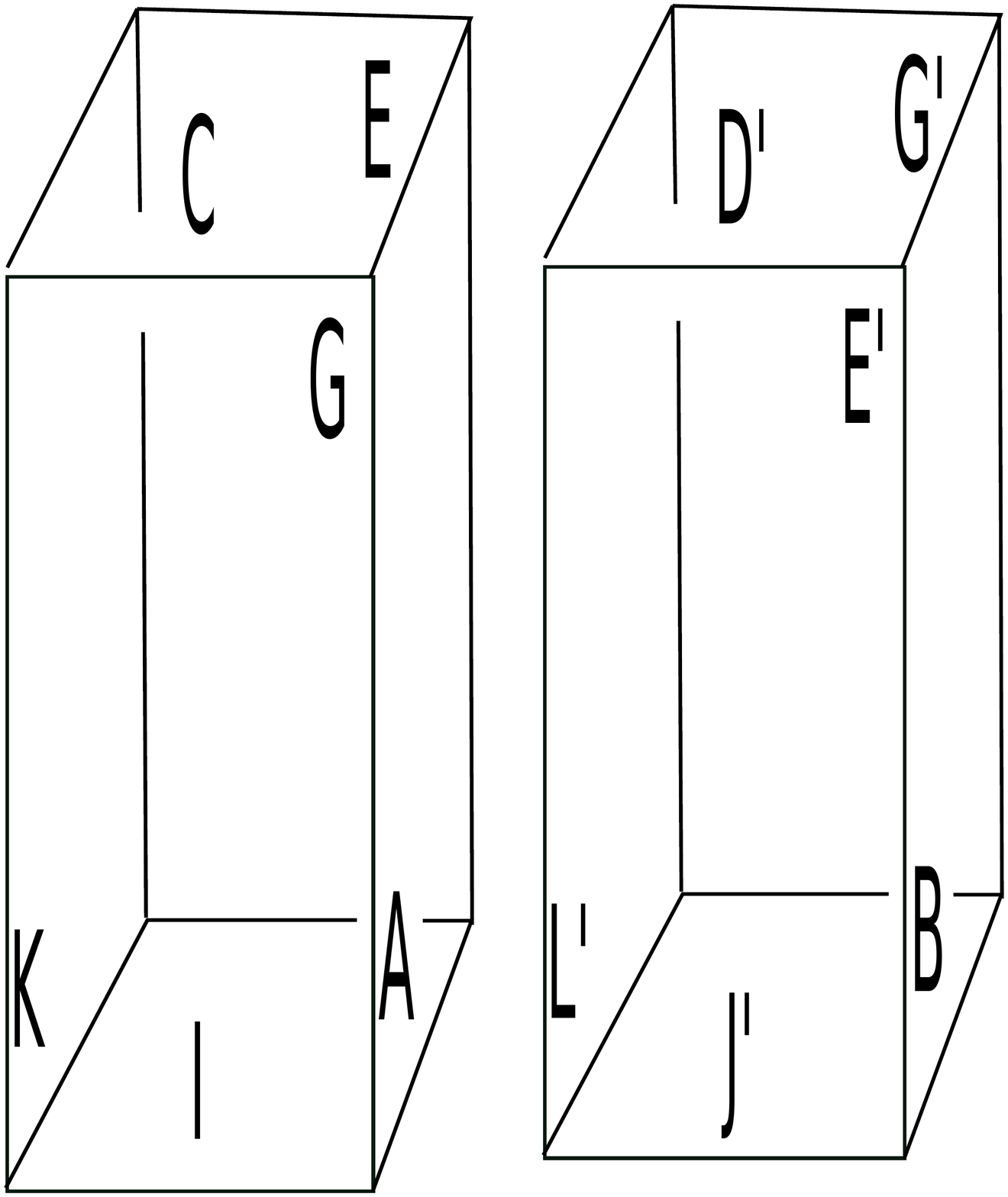}}

The dual octahedrons are then given by:

\centerline{\graphicspath{ {parabolic_dual/} }\includegraphics[width=6cm, height=5cm]{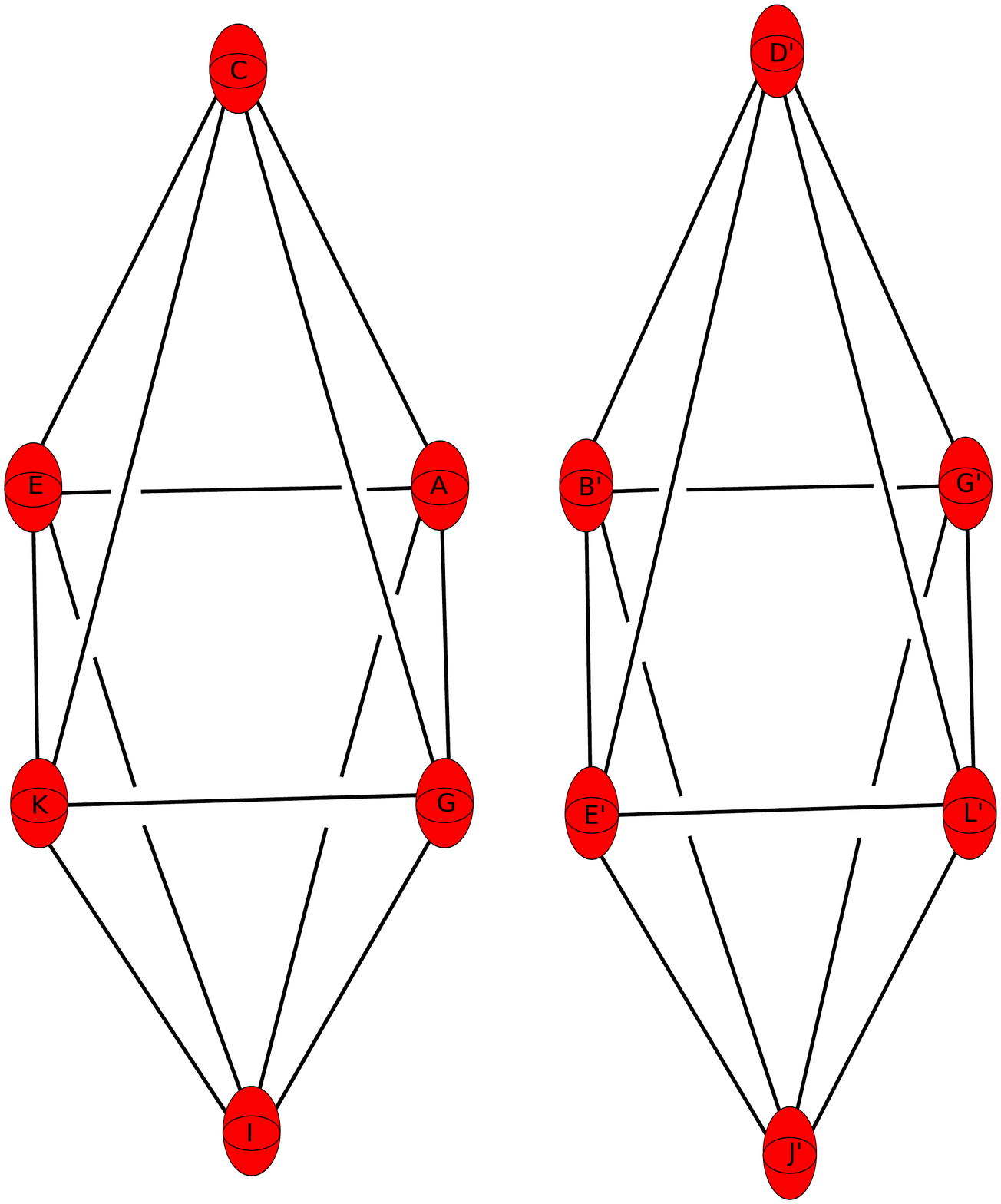}}

The following shows a three dimensional picture of the first octahedron above.

\centerline{\graphicspath{ {parabolic_dual/} }\includegraphics[width=7cm, height=10cm]{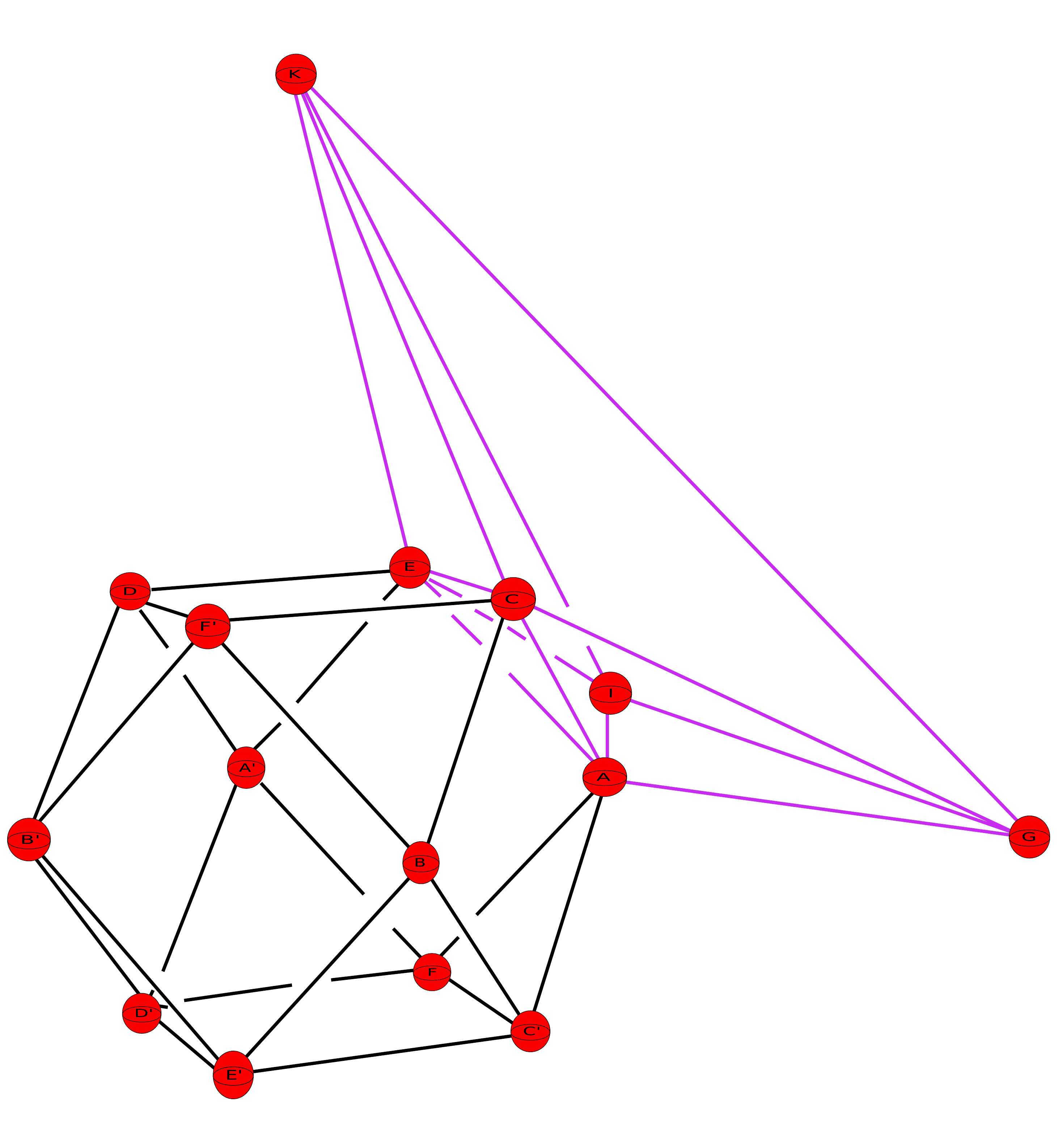}}

The 2-handle that we are going to add has one component running from $E$ to $G$, and can be seen in the following diagram. The reader should note that
we have drawn the added component as a curved arc when it really should be a straight line running from $E$ to $G$. The reason for drawing it as a curved arc
is simply because it is easier to view in the diagram, the reader should really picture this as a straight line.

\centerline{\graphicspath{ {parabolic_dual/} }\includegraphics[width=8cm, height=10cm]{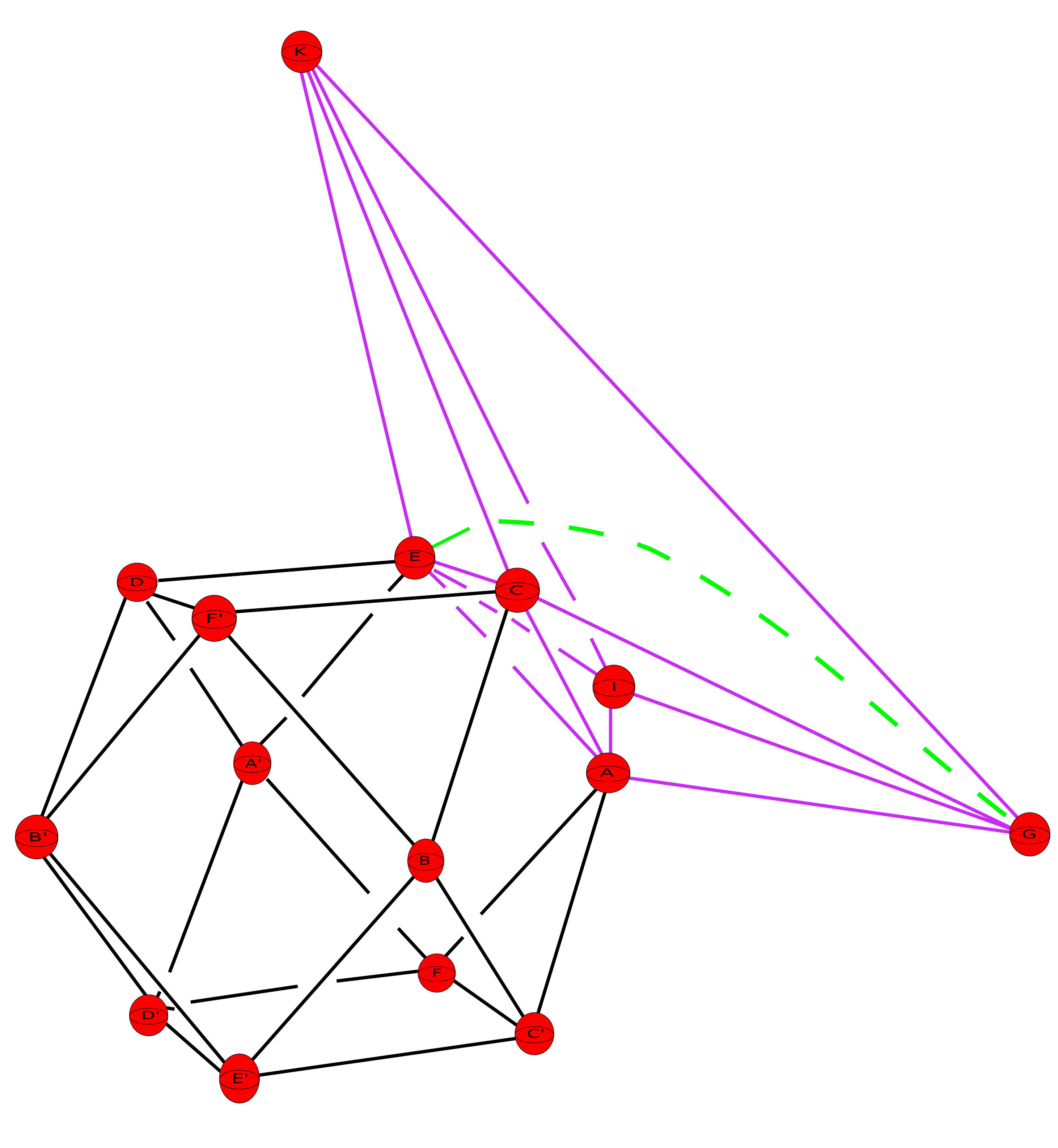}}

Observe that the added component falls outside of the diagrams showing the $x-y$, $x-z$, $y-z$ planes and the diagrams showing the six 2-handles that do not
all lie in a single 2-plane. Due to this, any handle cancellation/slides we do which take place in the $x-y$, $x-z$ and $y-z$ planes or the diagram
corresponding to the six 2-handles that do not all lie in a single 2-plane, will not interfere with this added 2-handle component running from $E$ to $G$. 

The following shows a three dimensional picture of the second octahedron above

\centerline{\graphicspath{ {parabolic_dual/} }\includegraphics[width=6cm, height=7cm]{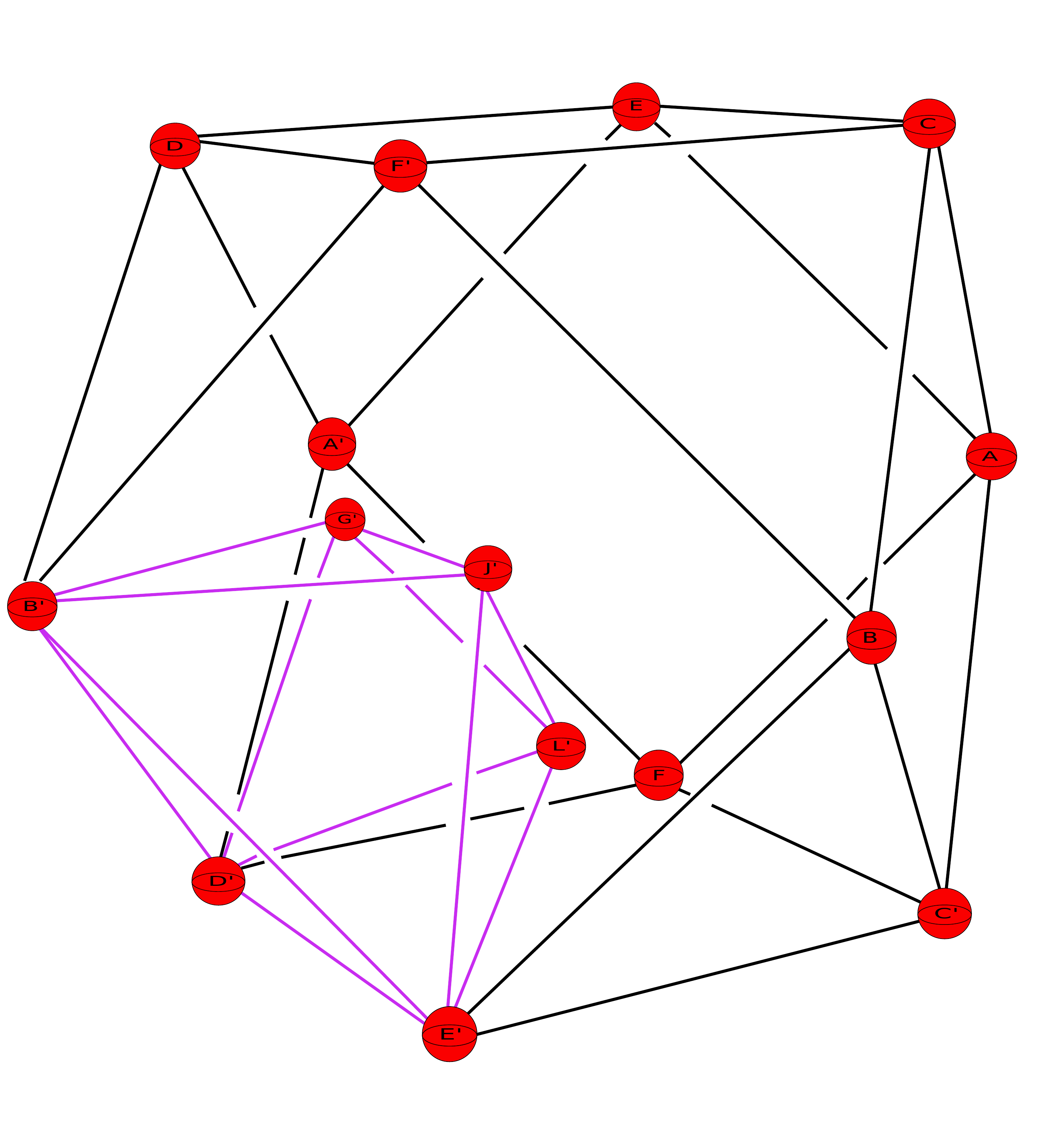}}

The 2-handle component that we are adding runs from $E'$ to $G'$ and can be seen in the following diagram:

\centerline{\graphicspath{ {parabolic_dual/} }\includegraphics[width=6cm, height=7cm]{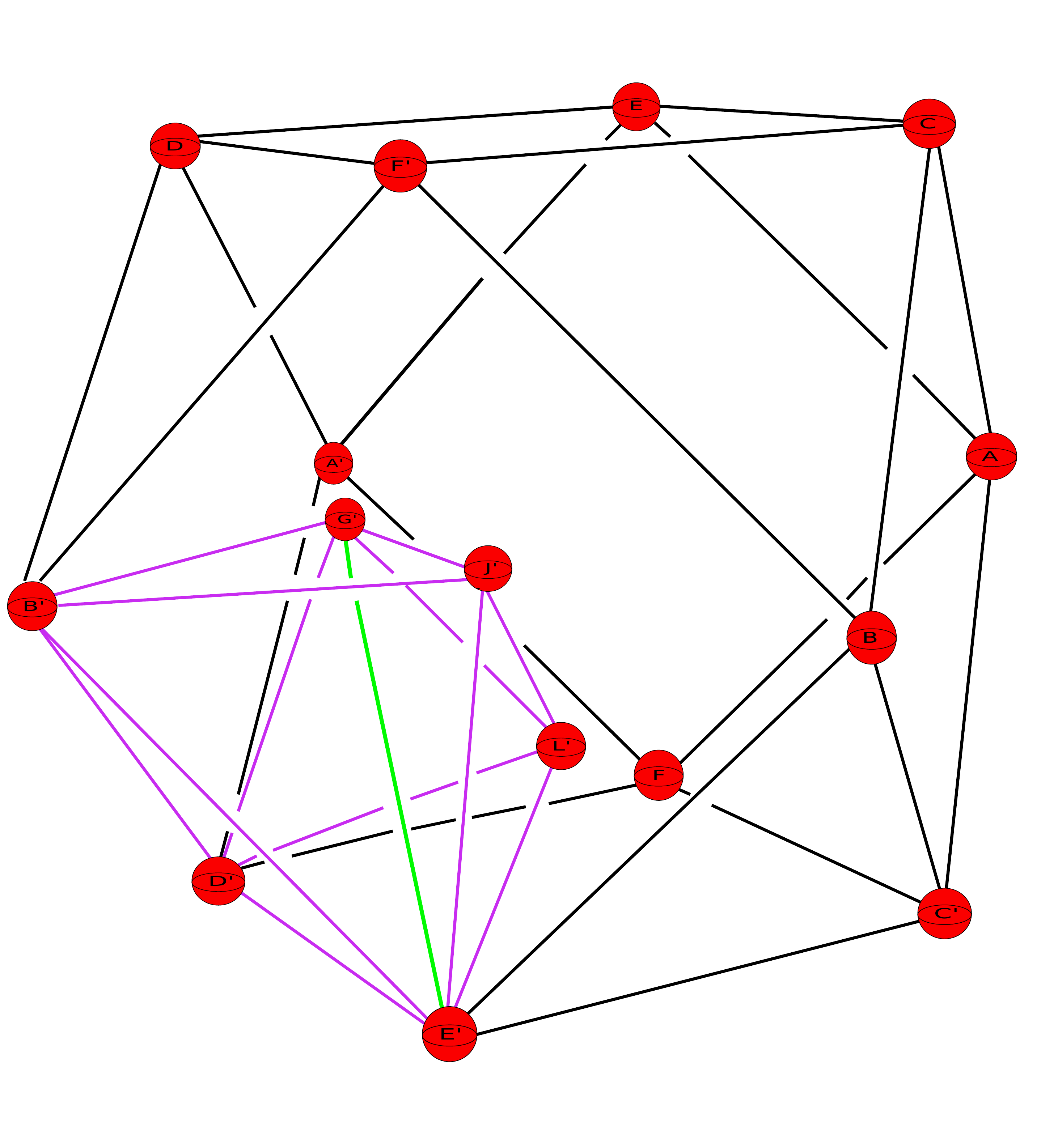}}

In this case the added 2-handle runs inside the diagram corresponding to the six 2-handles that do not all lie in a single 2-plane. However, it does not
lie in any of the planes $x-y$, $x-z$ or $y-z$. Furthermore, all the handle cancellations/slides we do to begin with will not interfere with this 2-handle.

The following picture shows the Kirby diagrams with the added attaching circle of the added $2$-handles that pass over $A-A'$, $J-J'$, $K-K'$ and $C-C'$, they
are the black dashed lines. We note that these 2-handles pass over the 1-handles in question once. Also, the added 2-handle running between $C-C'$ in the $x-z$ plane
is drawn as a curved arc, as opposed to a straight line, purely for ease of viewing.

\centerline{\graphicspath{ {filling_cusps/} }\includegraphics[width=11cm, height=10cm]{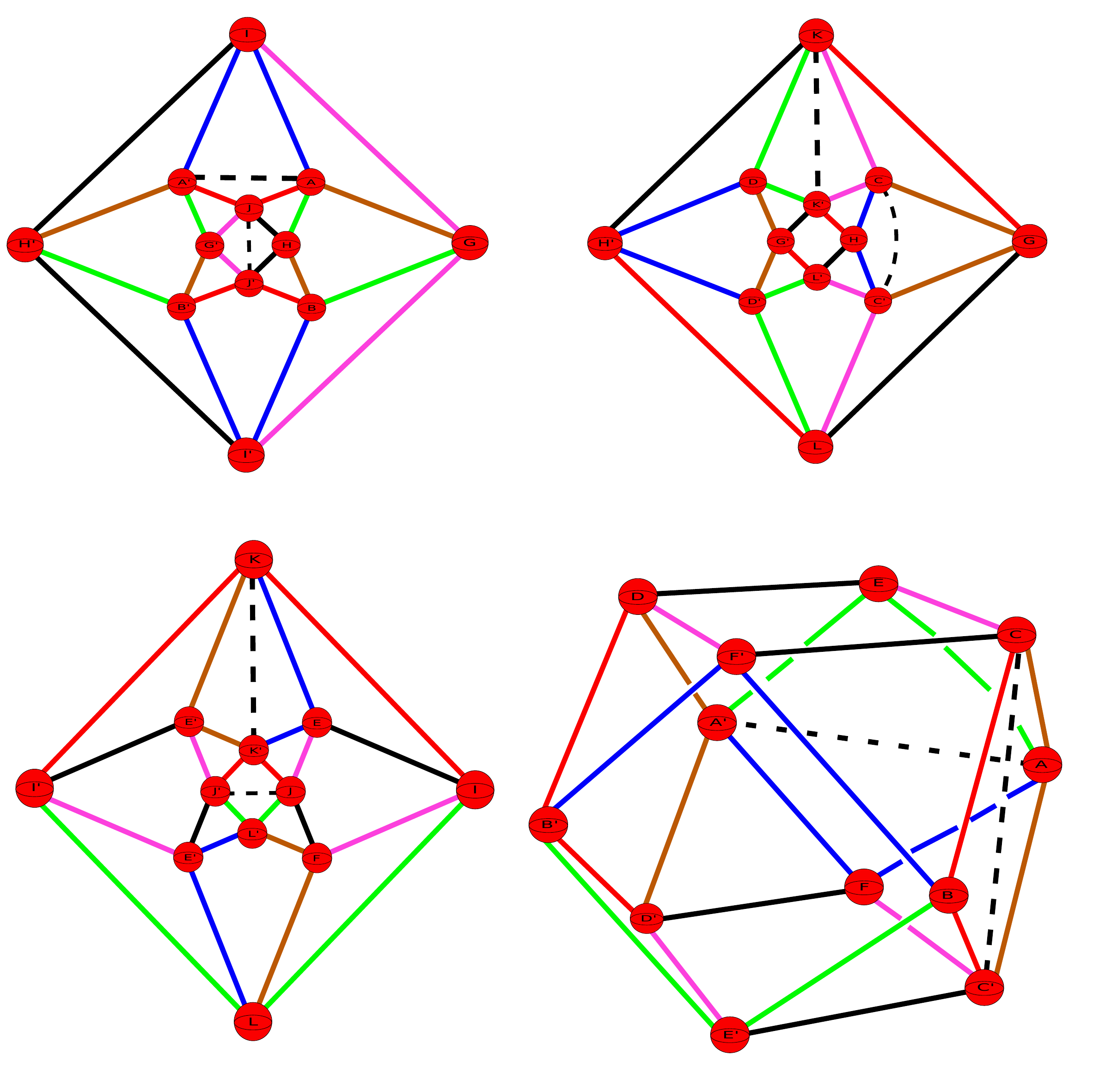}}

The added 2-handle that runs from $E$ to $G$ and $E'$ to $G'$ are shown in the following picture.

\centerline{\graphicspath{ {filling_cusps/} }\includegraphics[width=8cm, height=9cm]{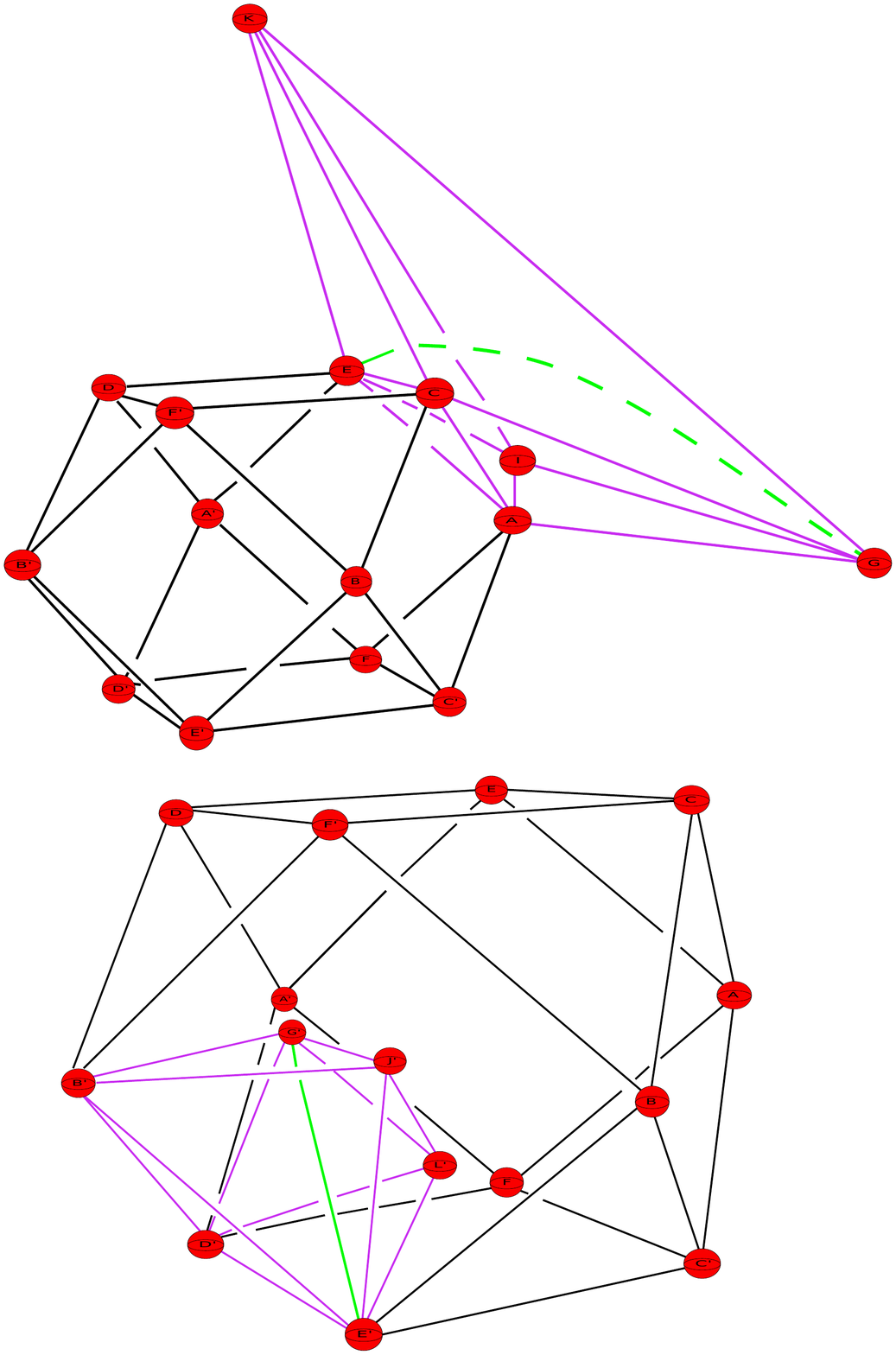}}

The above shows the Kirby diagram associated to the smooth closed 4-manifold that is obtained from filling in the boundary components of $M$. 
In sect.5 of \cite{sarat} we explained how the 2-handles associated to the Kirby diagram of $M$ had what we called a planar framing.
This had the feature that whenever we took a parallel curve to such a 2-handle it would never cross over the 2-handle. We want to briefly show that the same
feature is possessed by the added 2-handles making up the closed boundary filling of $M$.

We need to start by outlining exactly how the attaching regions of the 1-handles are being identified. We gave explicit computations of the identifying diffeomorphisms
in sect.5 of \cite{sarat}. The following table shows what the identifying diffeomorphism for the attaching regions of the 1-handles are.

\begin{tabular}{|l|l|}
\hline
{} & {} \\
1-handle & Identifying diffeomorphism \\
{} & {} \\
\hline
{} & {} \\
$A,A'$ \& $B, B'$ & $(x,y,z) \mapsto (-x,y,z)$ \\
{} & {} \\
\hline
{} & {} \\
$C,C'$ \& $D, D'$ & $(x,y,z) \mapsto (x,y,-z)$ \\
{} & {} \\
\hline
{} & {} \\
$E,E'$ \& $F, F'$ & $(x,y,z) \mapsto \bigg(-\frac{x}{x^2+y^2+z^2},-\frac{y}{x^2+y^2+z^2},-\frac{z}{x^2+y^2+z^2}\bigg)$ \\
{} & {} \\
\hline
{} & {} \\
$G,G'$ \& $H, H'$ & $(x,y,z) \mapsto \bigg(-\frac{x}{x^2+y^2+z^2},-\frac{y}{x^2+y^2+z^2},-\frac{z}{x^2+y^2+z^2}\bigg)$ \\
{} & {} \\
\hline
{} & {} \\
$I,I'$ \& $J, J'$ & $(x,y,z) \mapsto (x,-y,z)$ \\
{} & {} \\
\hline
{} & {} \\
$K,K'$ \& $L, L'$ & $(x,y,z) \mapsto \bigg(\frac{x}{x^2+y^2+z^2},\frac{y}{x^2+y^2+z^2},\frac{z}{x^2+y^2+z^2}\bigg)$ \\
{} & {} \\
\hline

\end{tabular} \\

The first four added 2-handles all lie in at least one single 2-plane, it is easy to check that a parallel curve to such a 2-handle cannot cross over the 2-handle
in any way, and lies in the same plane.
We will show this for
the added 2-handle that runs between $A, A'$. This added 2-handle lies in the $x-y$ plane, we take a parallel curve just above it that is going into
$A$, and that also lies in the $x-y$ plane. It is shown as the orange curve going into $A$ in the following diagram.

\centerline {\graphicspath{ {framing/} } \includegraphics[width=7cm, height=7cm]{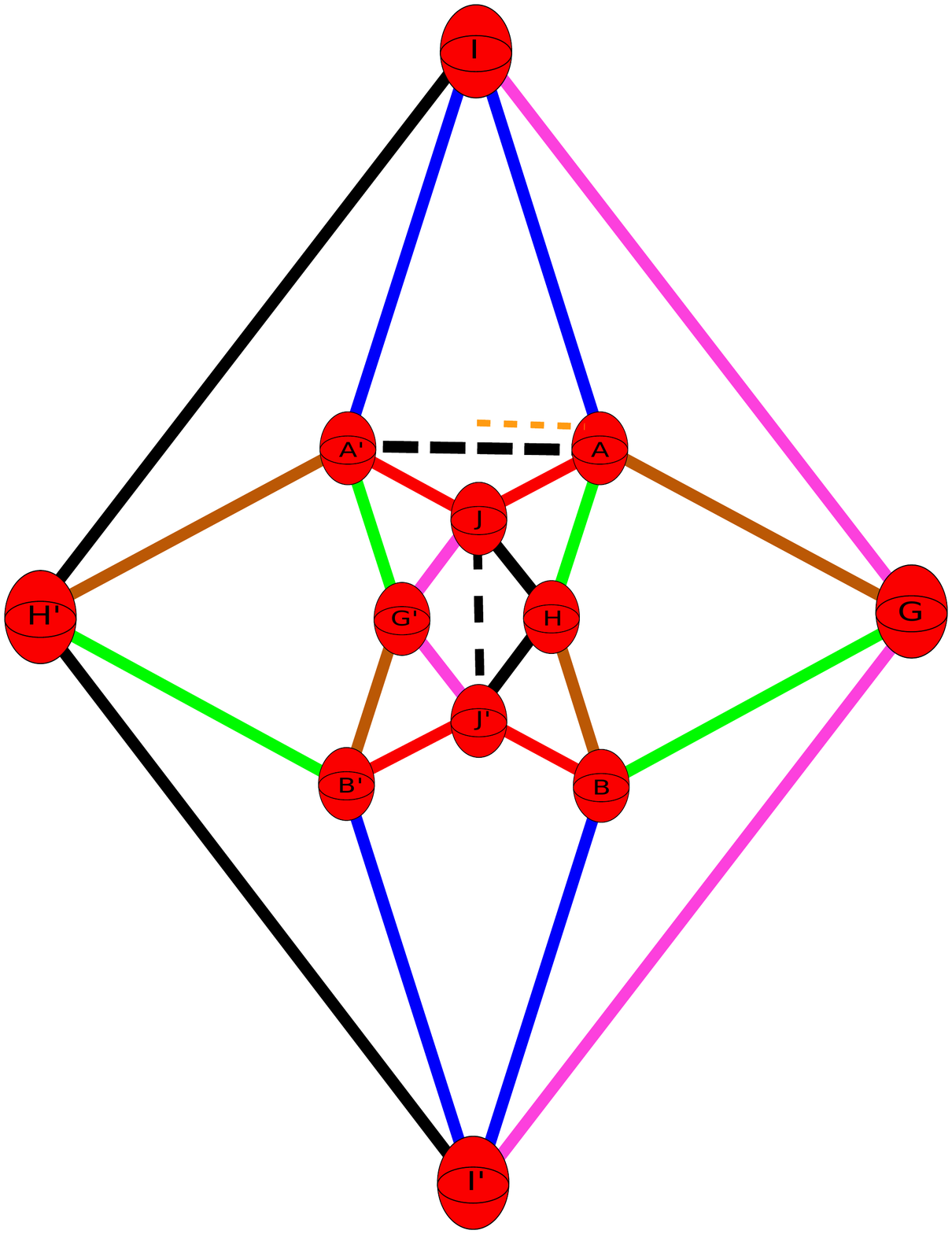}}

The attaching map for $A, A'$ is given by the reflection $(x,y,z) \mapsto (-x,y,z)$, it is therefore clear that when the parallel curve goes into
$A$ it comes out of $A'$ also above the added 2-handle. The following diagram shows what the parallel curve looks like.

\centerline {\graphicspath{ {framing/} } \includegraphics[width=7cm, height=7cm]{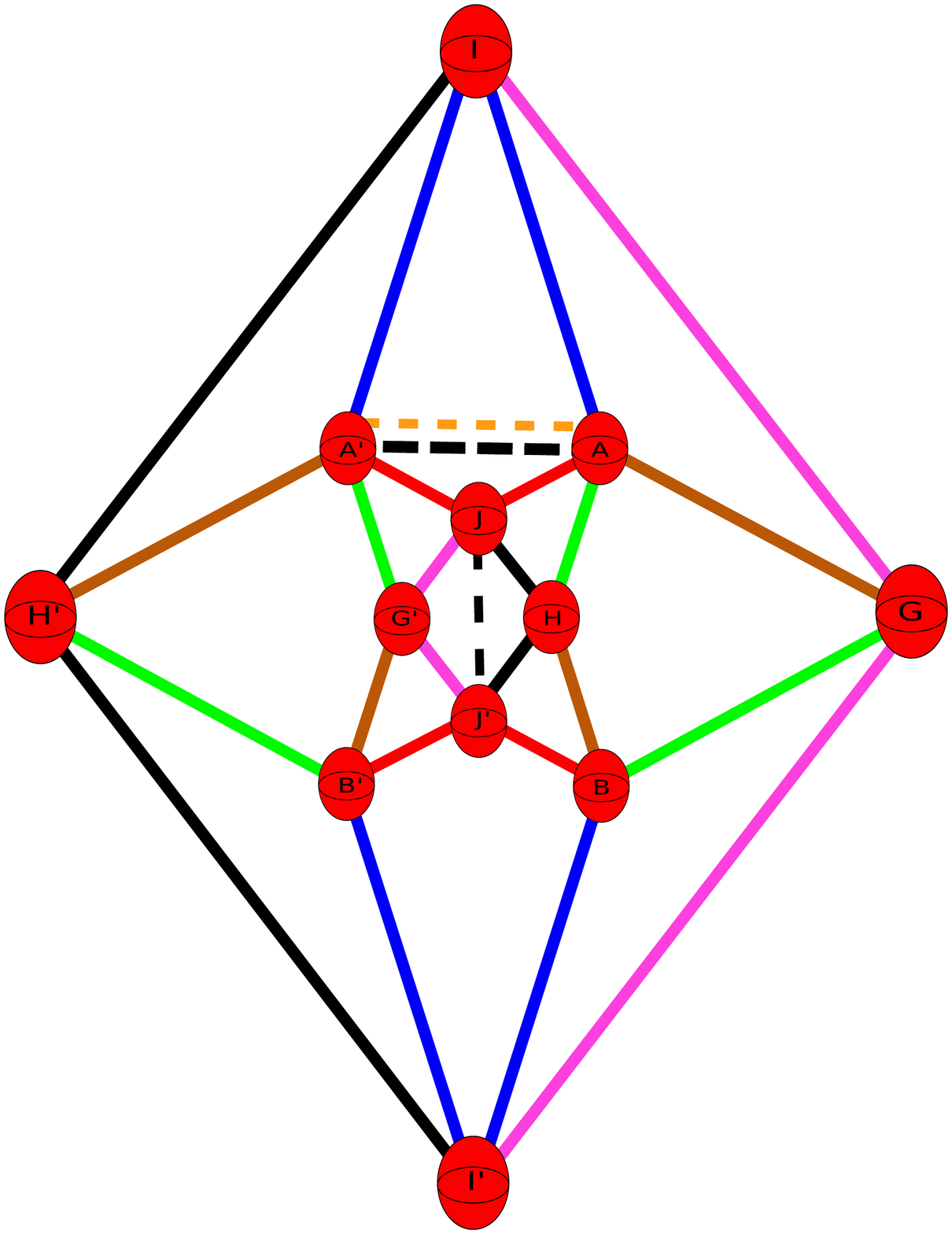}}

We see that a parallel curve to the added 2-handle running between $A,A'$ does not at any point cross over the added 2-handle. A similar argument for the added 2-handles running between $C,C'$, $J,J'$, $K,K'$ shows that
these added 2-handles have parallel curves behaving in a similar way. 

Let us show that the added 2-handle running from $E$ to $G$ and from $E'$ to $G'$ has parallel curve exhibiting a similar behaviour.
We start by taking a parallel curve component to the piece going into $E$, we choose the parallel curve component to lie just above this added 2-handle piece.
The following picture shows the added 2-handle (corresponding to the boundary filling) in green, and the parallel curve component in orange. Remember
that we are drawing the added 2-handle component, in green, as a curved arc purely for ease of viewing. The reader should really think of this as a straight
line running from $E$ to $G$. Due to this we have had to draw the parallel curve as a curved arc, again the reader should think of this as a straight line
running into $E$ and that is parallel to what should be a straight line running between $E$ and $G$.

\centerline {\graphicspath{ {framing/} } \includegraphics[width=8cm, height=7cm]{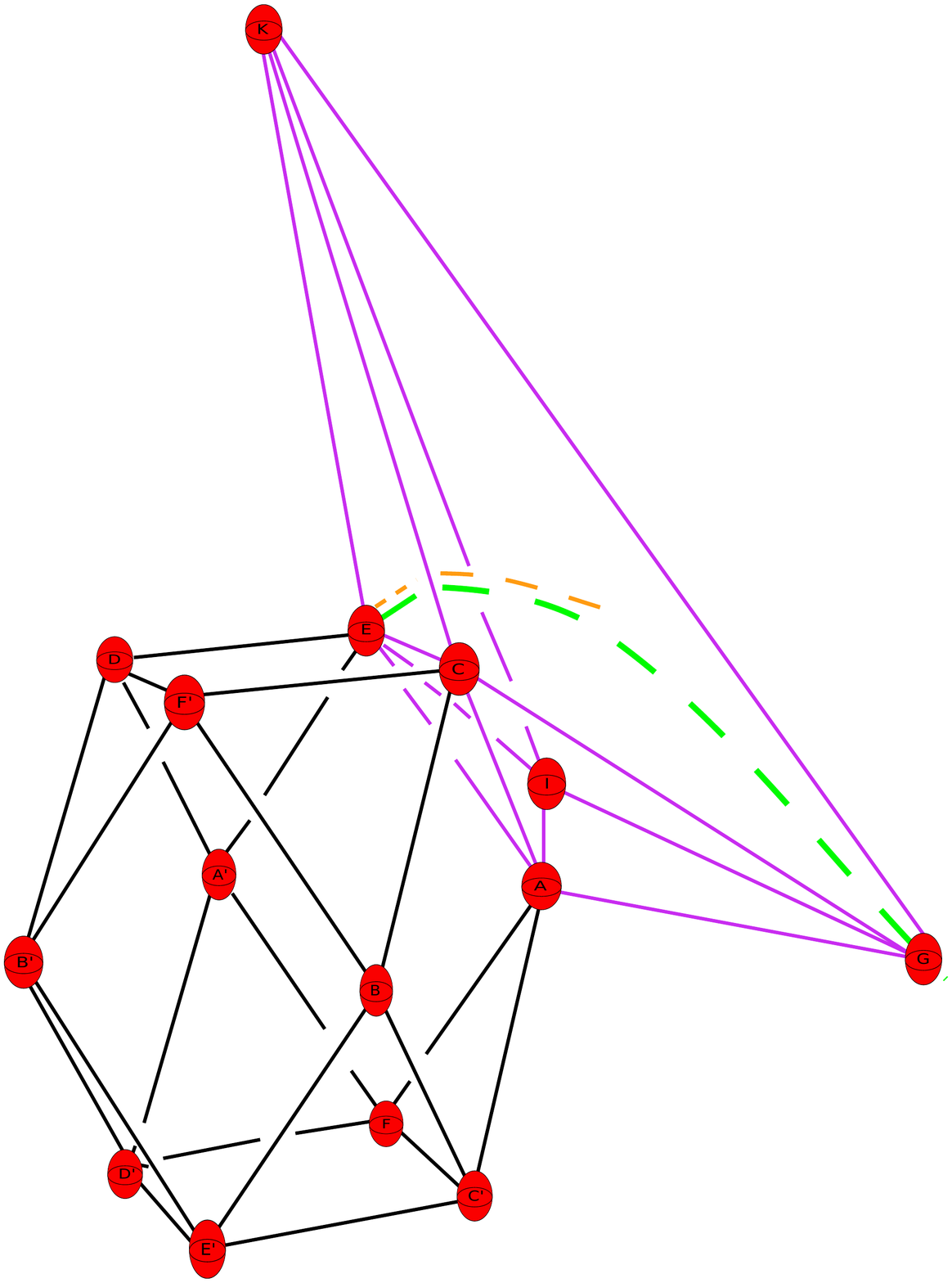}}

The attaching map for $E, E'$ was given by the composition of the inversion in $S^2$ followed by the antipodal map. 
A rough computation shows that the orange parallel curve comes out of $E'$ slightly above the added green 2-handle (we suggest the reader to look
at sect.5 in \cite{sarat} for an example of this type of computation).

\centerline {\graphicspath{ {framing/} } \includegraphics[width=7cm, height=5cm]{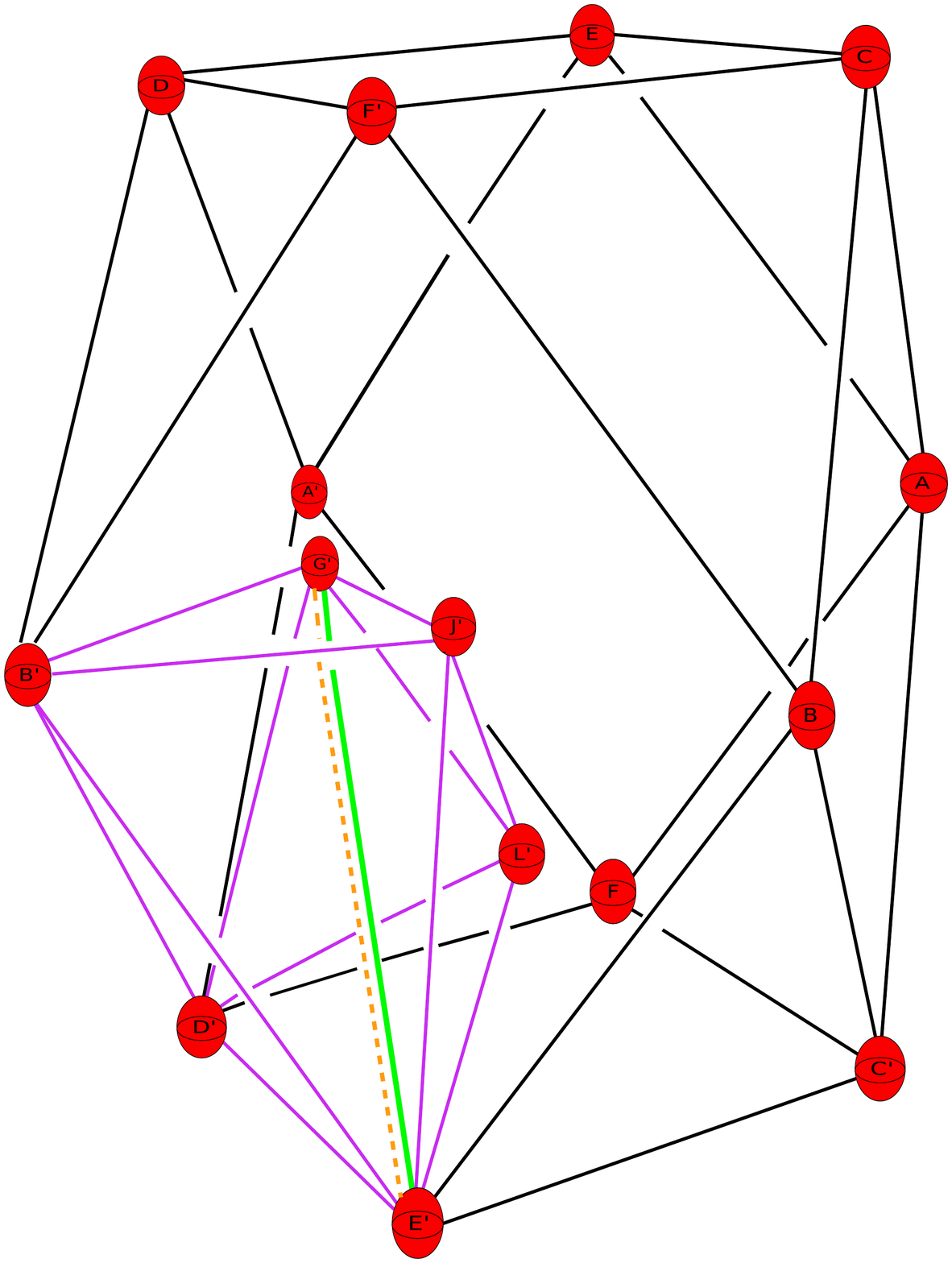}}

The attaching map for $G, G'$ is the same as the one for $E, E'$, and a little rough analysis shows that when the parallel curve comes 
out of $G'$ it does so just above the added green 2-handle.

\centerline {\graphicspath{ {framing/} } \includegraphics[width=10cm, height=8cm]{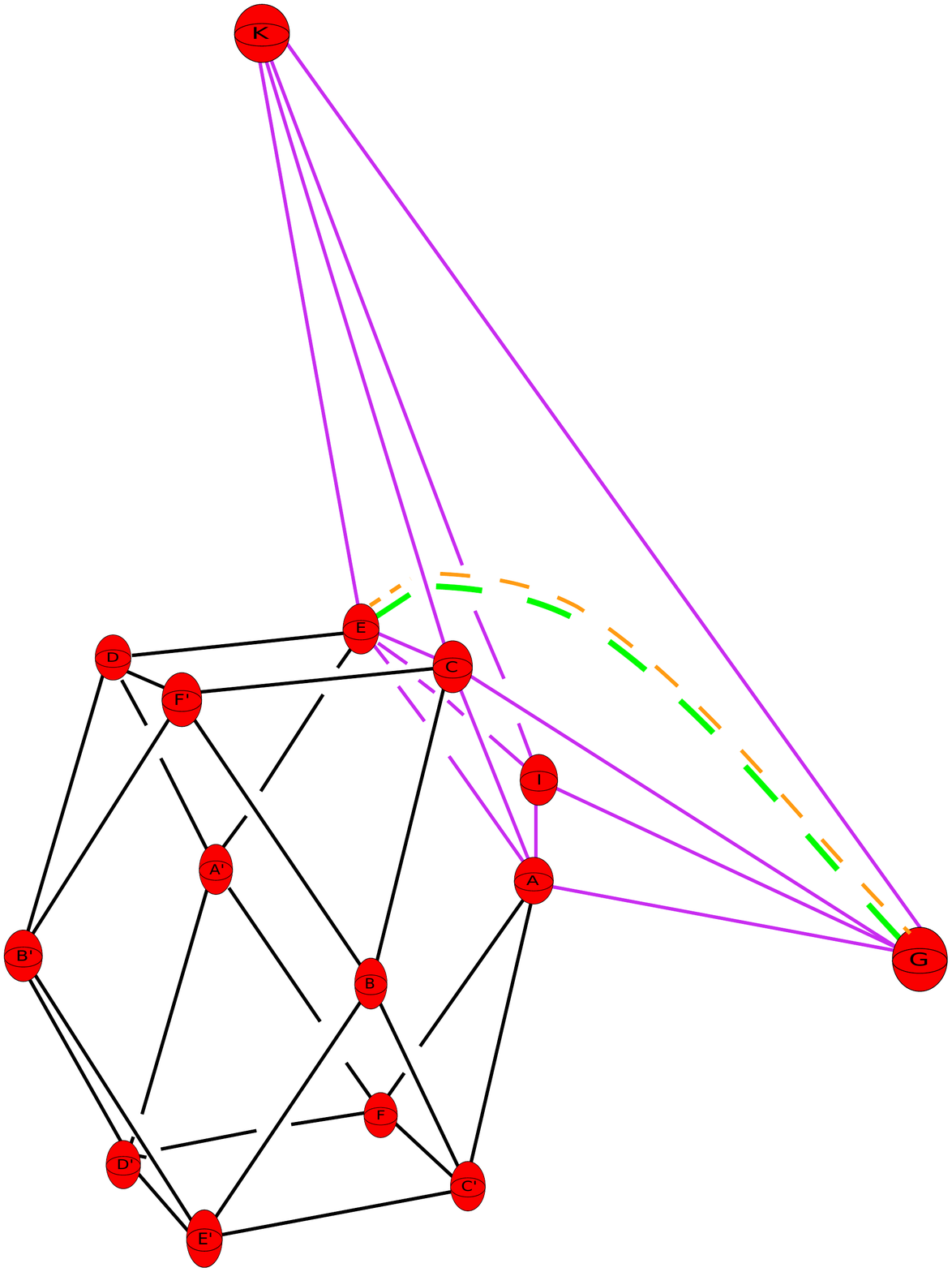}}

We thus see that the parallel curve never crosses over the added 2-handle, and in this regard has a planar framing. The importance of this observation is that
the closed filled in manifold, obtained from filling in the boundary of $M$, has a Kirby diagram with every 2-handle having a planar framing.
In turn this allows us to carry out handle slides in a very straight forward manner. In particular, for those 2-handles that reside in the diagram
corresponding to the six 2-handles that do not all lie in a single plane the handle slides can be carried out as if these 2-handles were all lying
in a single 2-plane.

There is one technicality to do with the \textbf{second elementary move} we described in section \ref{elementary}. This was the handle slide that involved pushing
a component of 2-handle through one attaching sphere of a 1-handle so that it came out of the second attaching sphere. 
When we described this handle slide, we
showed how it took place in the case that the 2-handles were all residing in a single 2-plane. In this situation it was straightforward to understand how
the 2-handle component behaved as it went through one attaching sphere and then out the other. We also want to carry out such handle slides in the
diagram corresponding to the six 2-handles that do not all lie in a single 2-plane. In this case one has to be a bit careful, the reason being that
if the attaching map identifying the attaching spheres is some ``wild'' diffeomorphism then it may well be the case that when we push a component of 2-handle
through one attaching sphere it could come out the other in a very non-trivial way. The point is that for all the Kirby diagrams we consider this will
never be a problem as the attaching maps between attaching spheres are very straightforward. Let us explain this in a bit more detail.

Suppose we have a 1-handle $S, S'$ and we have a component of a 2-handle that goes in to $S$, as shown in the following diagram.

\centerline {\graphicspath{ {framing/} } \includegraphics[width=4cm, height=2cm]{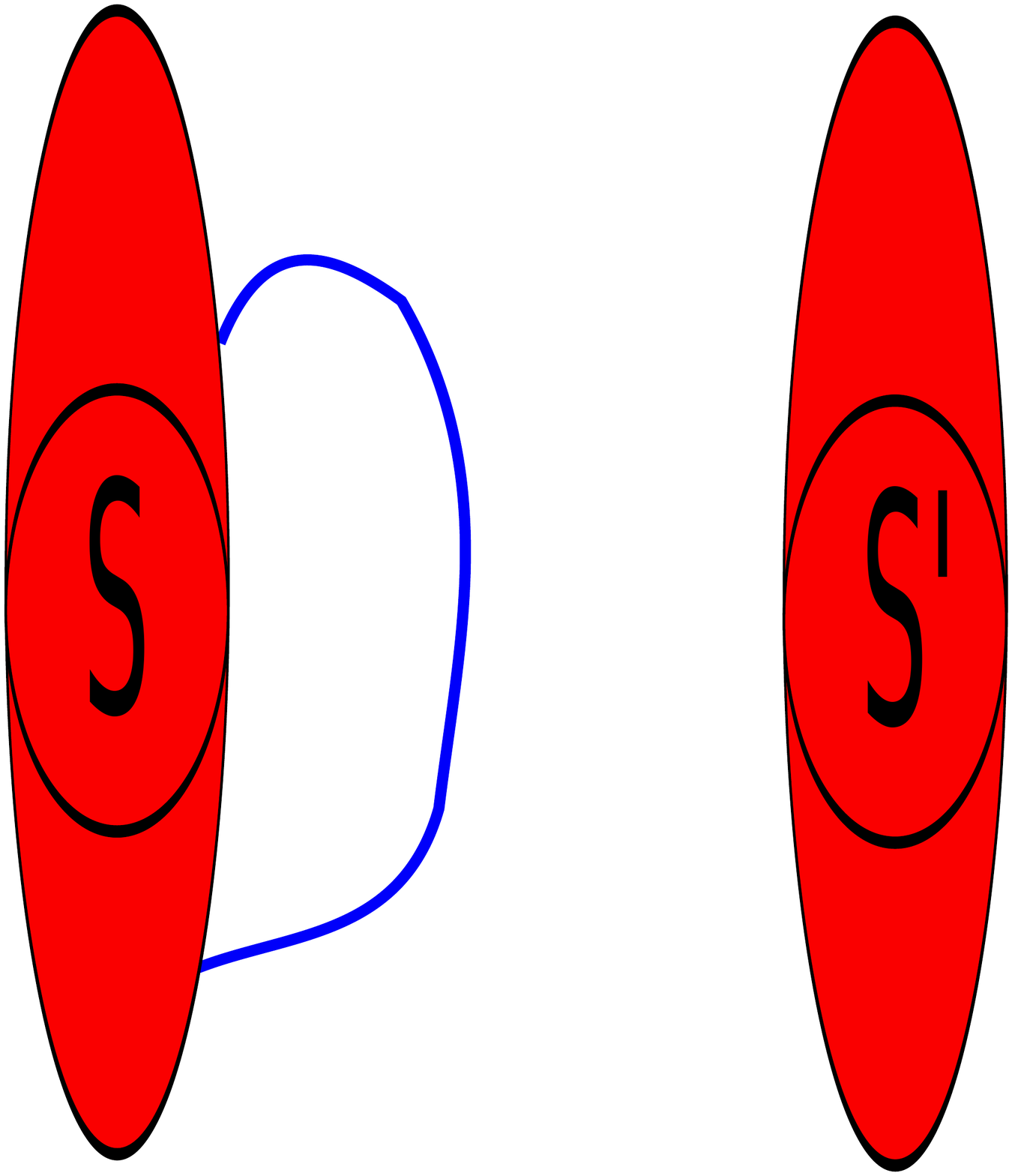}}

We can then push the blue 2-handle component through $S$ so that it comes out of $S'$. What we are really doing when we carry out such a move is we are
isotoping the blue curve onto the attaching sphere $S$, producing a curve on $S$. For example, the following diagram shows one such isotopy in which
the blue 2-handle component has been moved to give a curve on the attaching sphere $S$.

\centerline {\graphicspath{ {framing/} } \includegraphics[width=4cm, height=2cm]{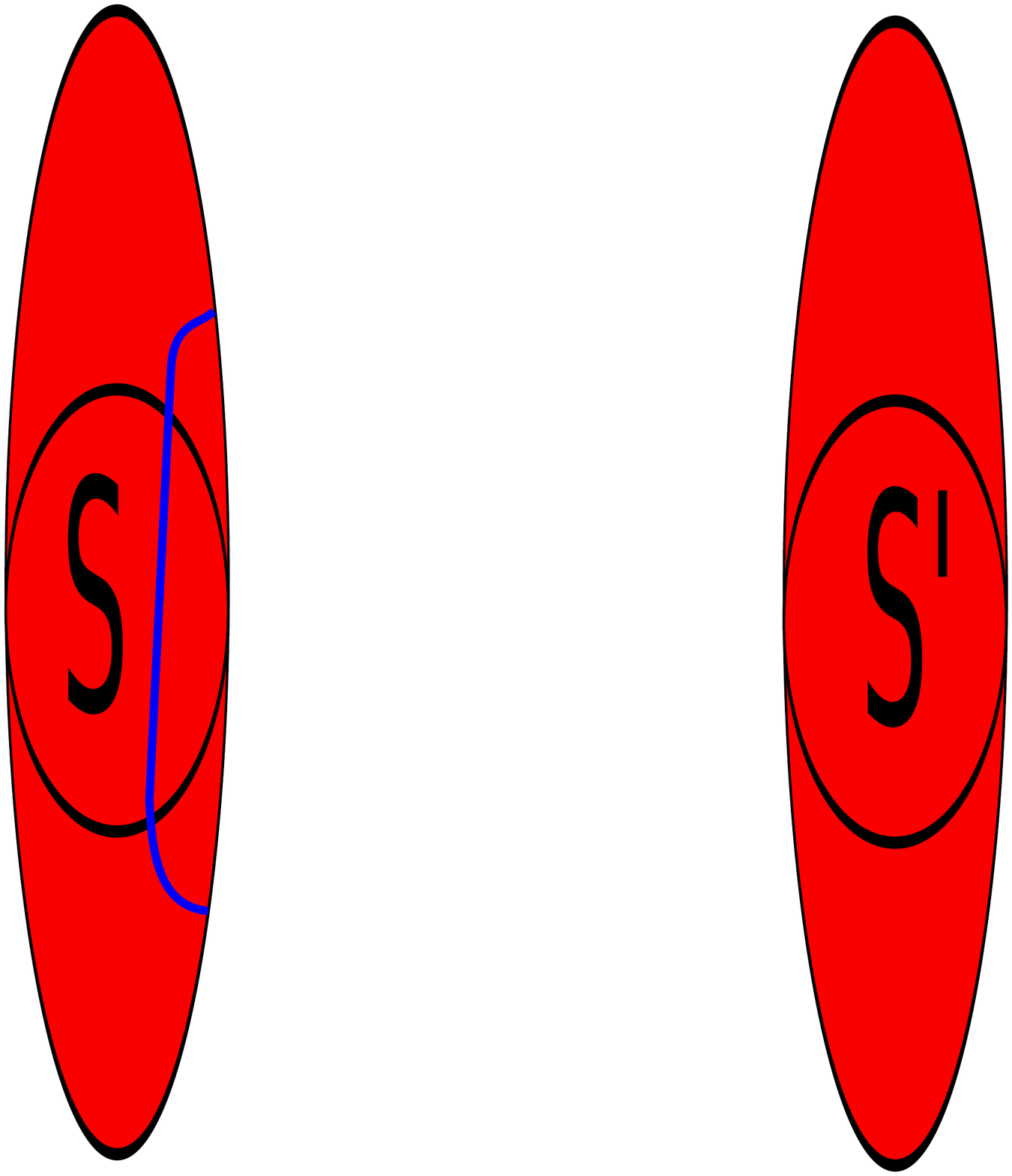}}

Now, we have an attaching map from $S$ to $S'$, which is some diffeomorphism identifying the spheres. Therefore, the blue curve on $S$ will map to some
image curve on $S'$, once we understand how this image curve on $S'$ looks we can then understand what the blue 2-handle component looks like
after we have pushed it through $S$. It is at this point that one must be very careful, for the attaching map identifying $S$ to $S'$ could be
some ``wild'' diffeomorphism, which in turn may map the blue curve on $S$ to some curve on $S'$ that winds around $S'$ several times. If there are
other bits of 2-handle coming out of $S'$, then this image curve could wrap around some of these other bits of 2-handle, and this is precisely
why it is important to know what the attaching map between the attaching spheres is.
In our case all attaching maps are reflections or compositions of reflections with inversion in $S^2$, therefore it is a straightforward process
to work out what a 2-handle component looks like when we push it through a 1-handle.

We give an explicit example of this sort of computation. Consider the diagram consisting of the six 2-handles that did not all lie in a single plane together
with the added 2-handle running from $A$ to $A'$, corresponding to a filling of a boundary component.

\centerline {\graphicspath{ {framing/} } \includegraphics[width=6cm, height=6cm]{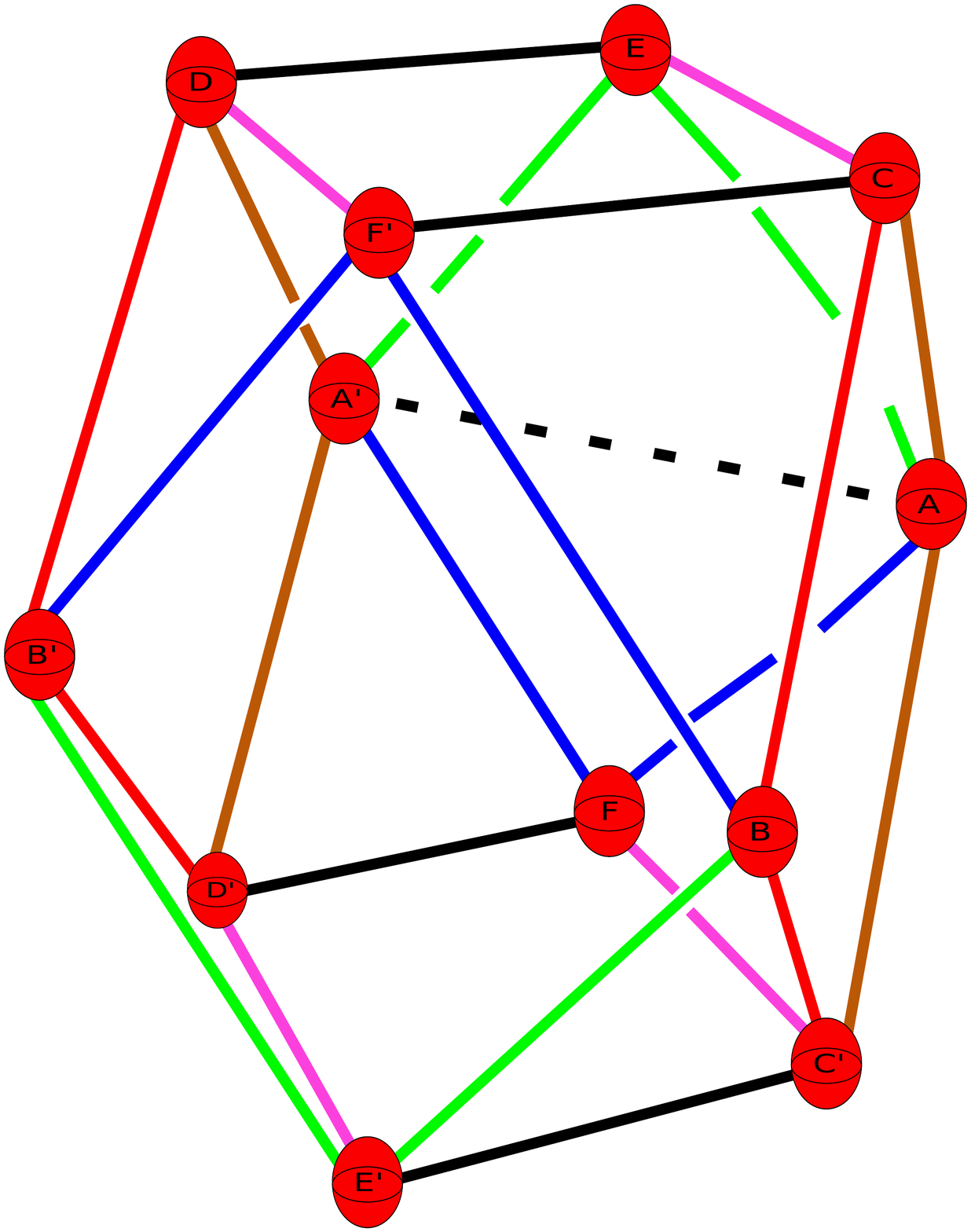}}

When we use the added 2-handle running between $A$ and $A'$ (black dashed line) to cancel $A,A'$ the green 2-handle component moves in to the following
position.

\centerline {\graphicspath{ {framing/} } \includegraphics[width=6cm, height=5cm]{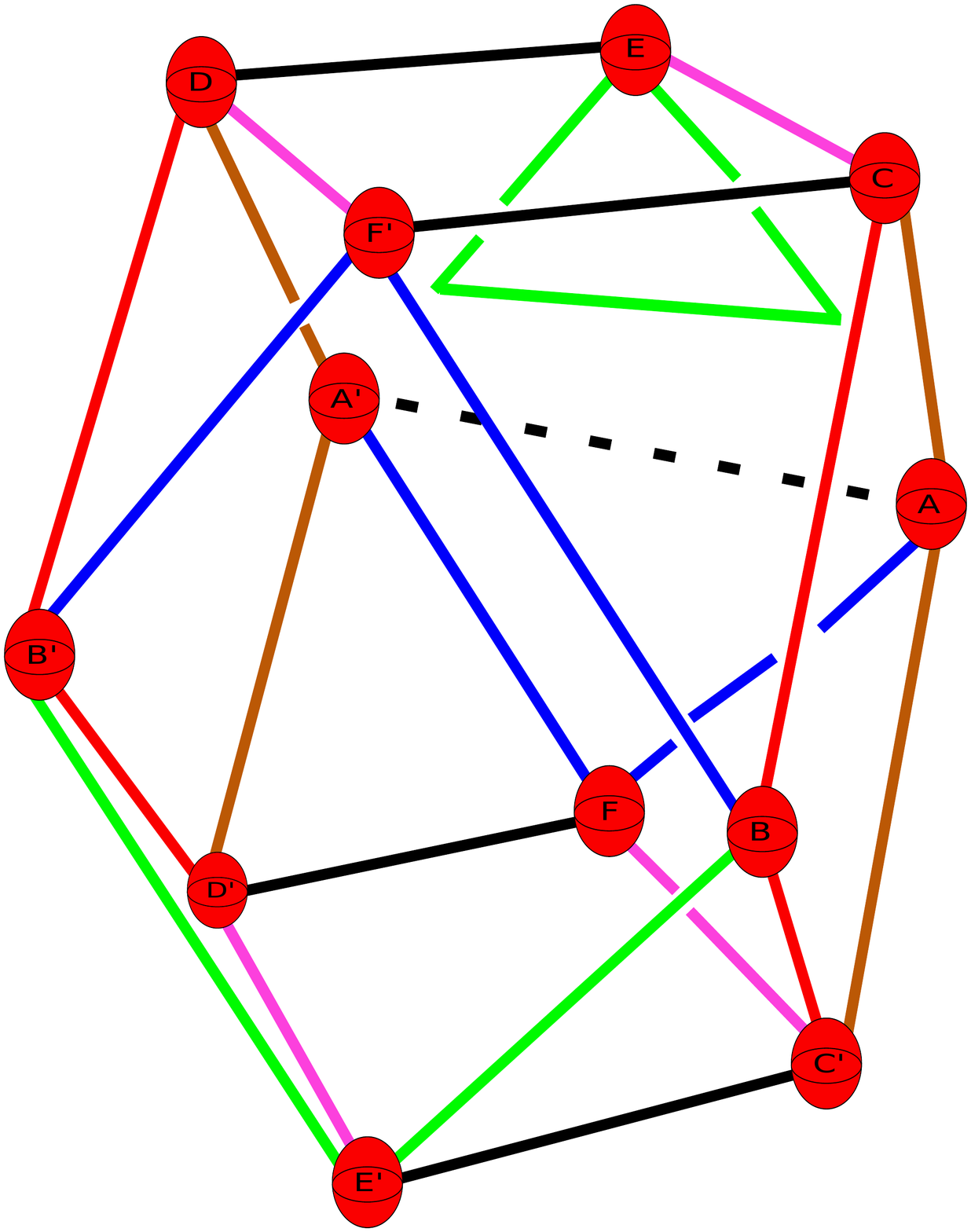}}

We can then slide this green 2-handle component through $E$ so that it comes out of $E'$. To start with this involves isotoping the green curve looping
back into $E$ to the following curve on $E$.

\centerline {\graphicspath{ {framing/} } \includegraphics[width=6cm, height=6cm]{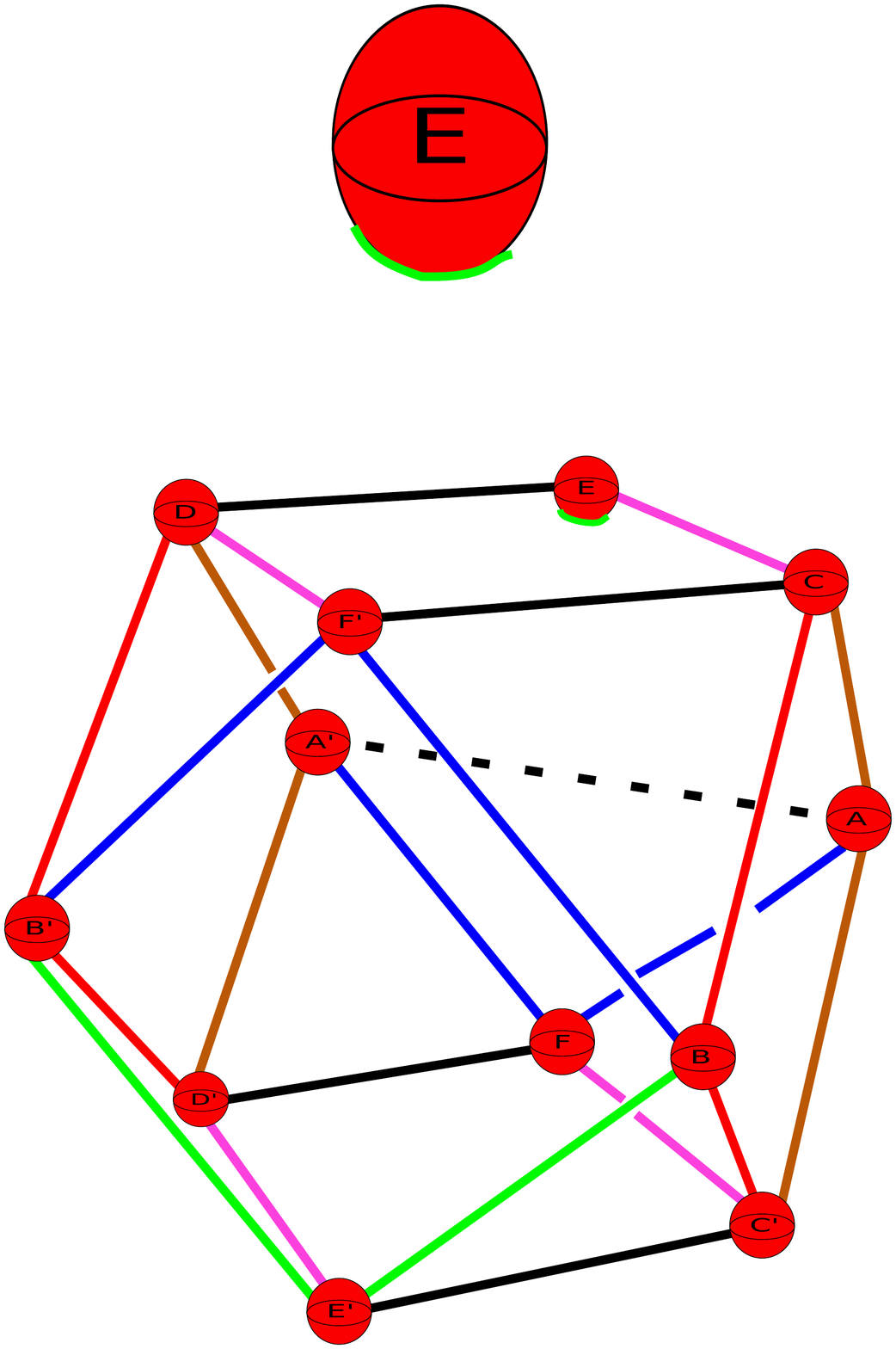}}

We have shown a close up of the curve on $E$ and how it looks in the whole diagram. The curve sits in a plane parallel to the x-z plane

If we want to push it through so that it comes out of $E'$ we need to understand what the attaching map for $E, E'$ is. Observe that there are other
2-handles that are hitting $E'$, we have a pink and black 2-handle component in the above diagram, but $E'$ also sits in the y-z plane, and there
are blue, black and pink 2-handle components in that plane meeting $E'$ (see diagram showing 2-handles in the y-z plane). Therefore it is possible
that when we send the green curve on $E$, using the attaching map of $E, E'$, we could end up with an image curve on $E'$ that winds around some of these
other 2-handle components meeting $E'$. 

In this situation we know exactly how $E$ is identified to $E'$, we showed that the identifying map was given by
\[ (x,y,z) \mapsto \bigg(-\frac{x}{x^2+y^2+z^2},-\frac{y}{x^2+y^2+z^2},-\frac{z}{x^2+y^2+z^2}\bigg) \]
which is the composition of the inversion in $S^2$ followed by the antipodal map. Recall that the attaching sphere corresponding to $E$ is a sphere 
(of some small arbitrary radius) centred at the point $(0, \frac{1}{\sqrt{2}}, \frac{1}{\sqrt{2}})$. 
The green curve on $E$ lies inside $S^2$:

\centerline {\graphicspath{ {framing/} } \includegraphics[width=3cm, height=3cm]{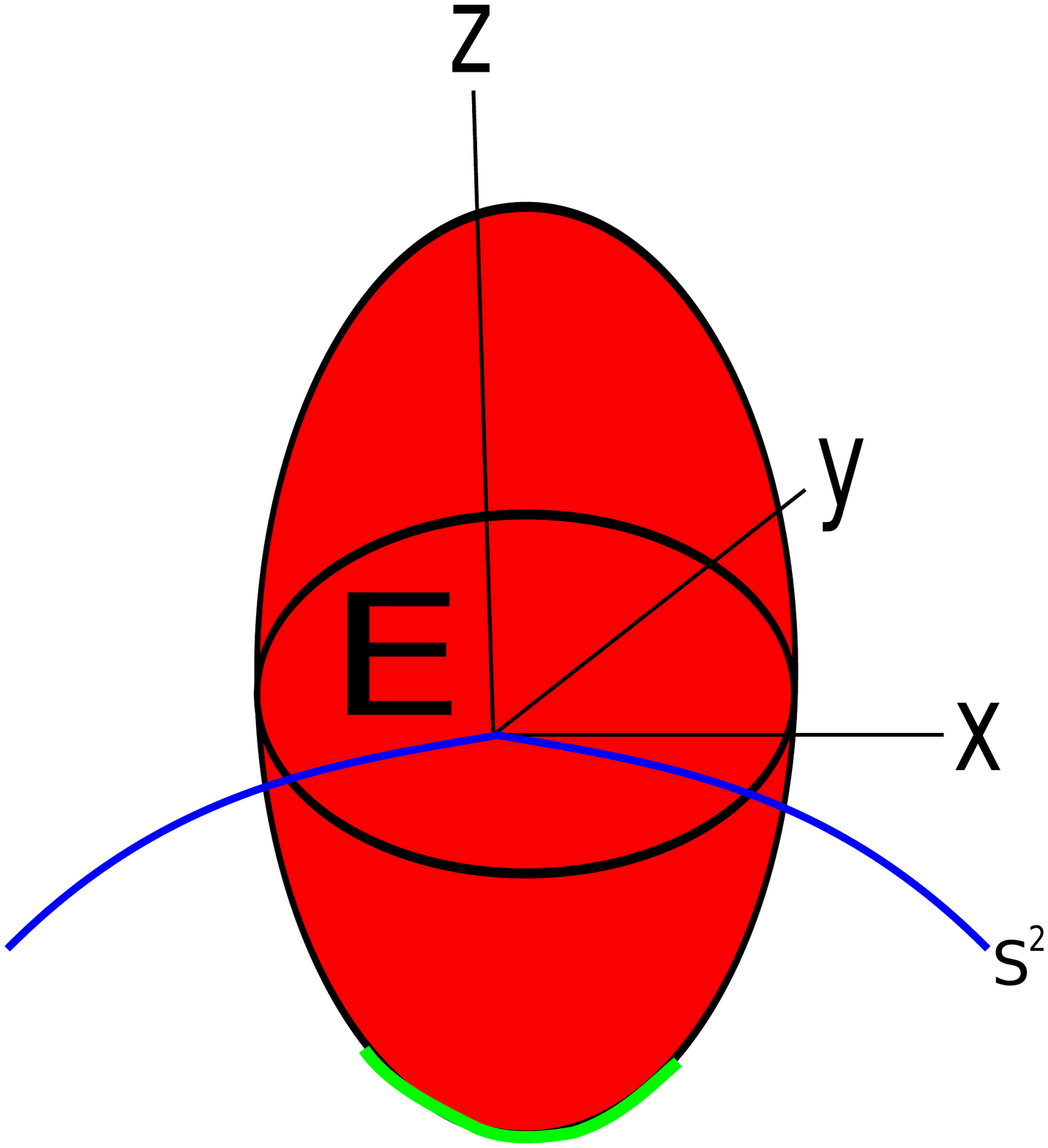}}

When we apply the
inversion map it will map to a curve outside of $S^2$. The following diagram shows what it looks like when we apply the inversion map.

\centerline {\graphicspath{ {framing/} } \includegraphics[width=4cm, height=4cm]{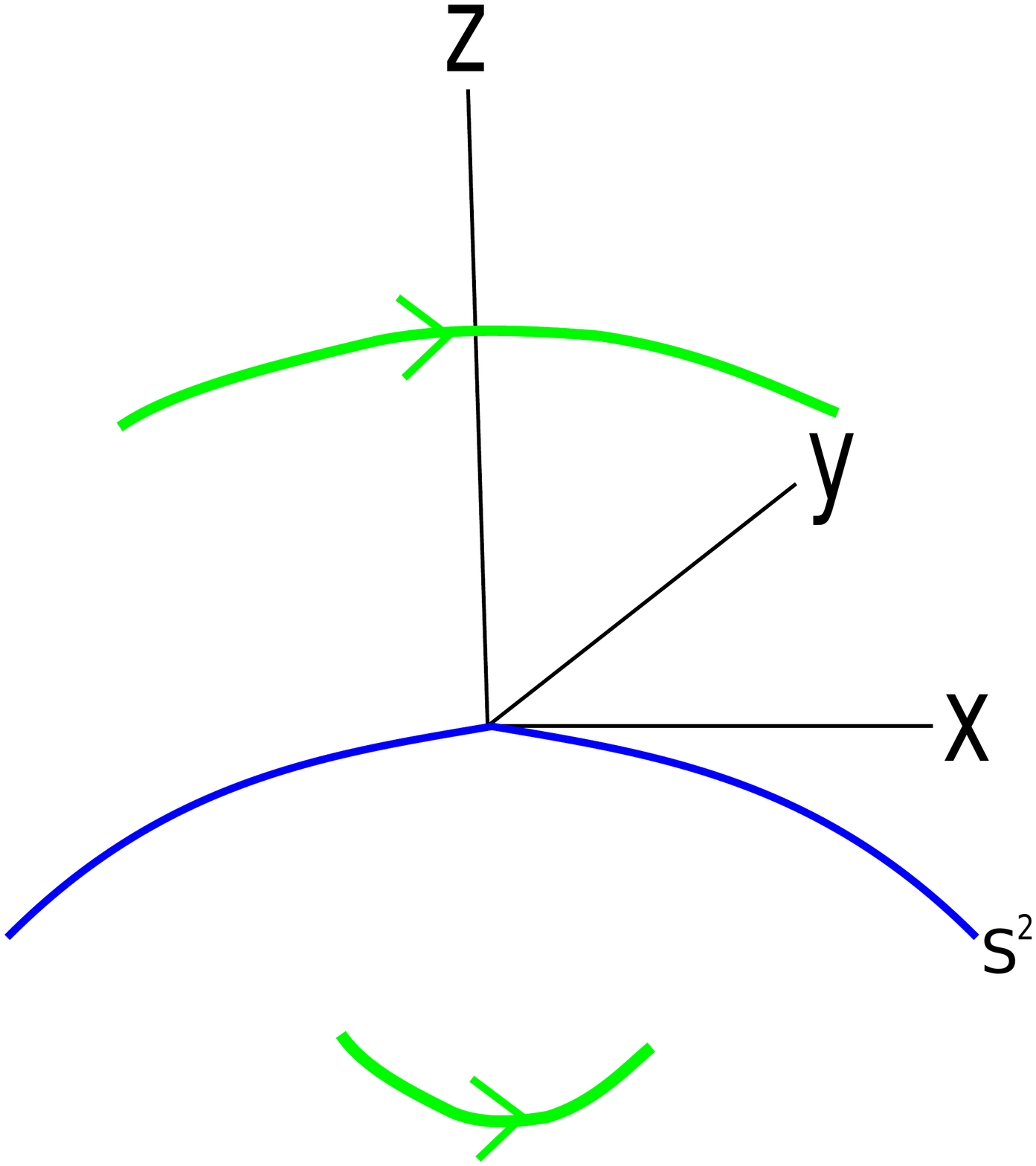}}

We then apply the antipodal map to get the following curve on $E'$.

\centerline {\graphicspath{ {framing/} } \includegraphics[width=2cm, height=2cm]{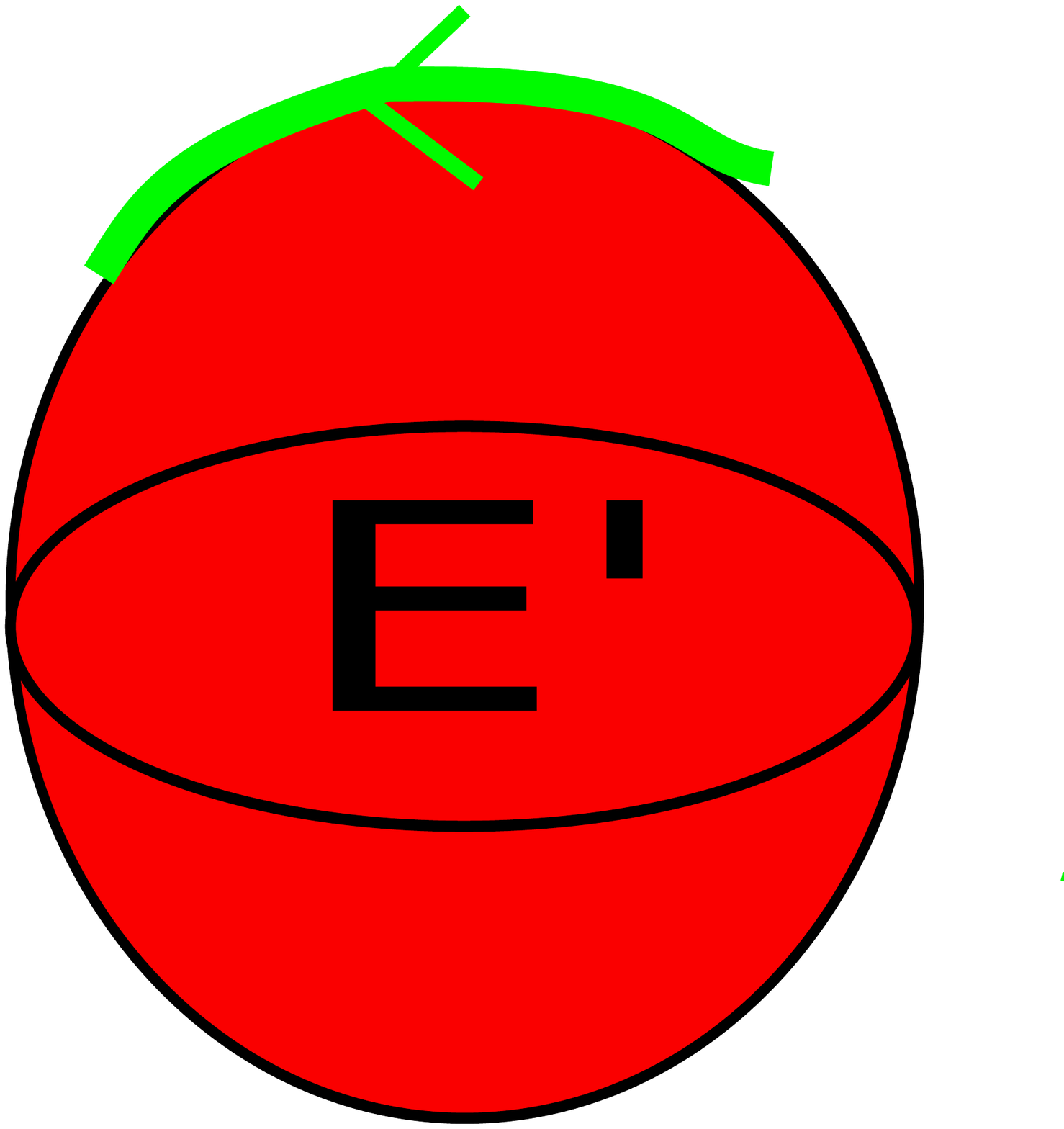}}

In the whole diagram we then see that the image of the green curve, which was residing on $E$, under the attaching map from $E$ to $E'$ looks like:

\centerline {\graphicspath{ {framing/} } \includegraphics[width=6cm, height=6cm]{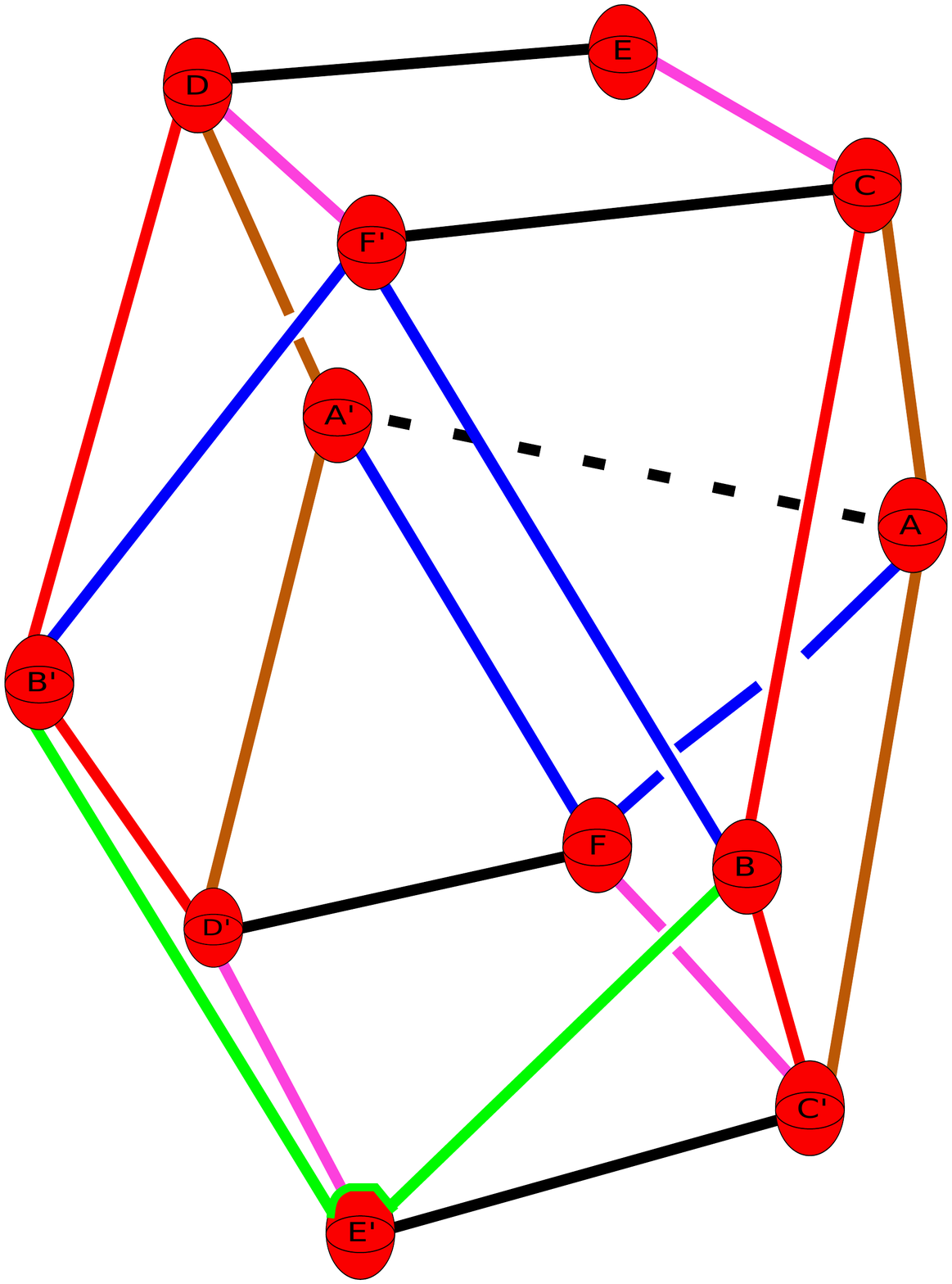}}

Thus, what we are able to conclude is that when we push the original green curve through $E$ it comes out of $E'$ in straightforward manner, and does not
``wind around'' any of the other 2-handle components meeting $E'$. We can then proceed to push the green 2-handle component off $E'$ to obtain a 
green 2-handle component passing over $B, B'$ once and that does not interfere with any of the other 2-handles

\centerline {\graphicspath{ {framing/} } \includegraphics[width=6cm, height=6cm]{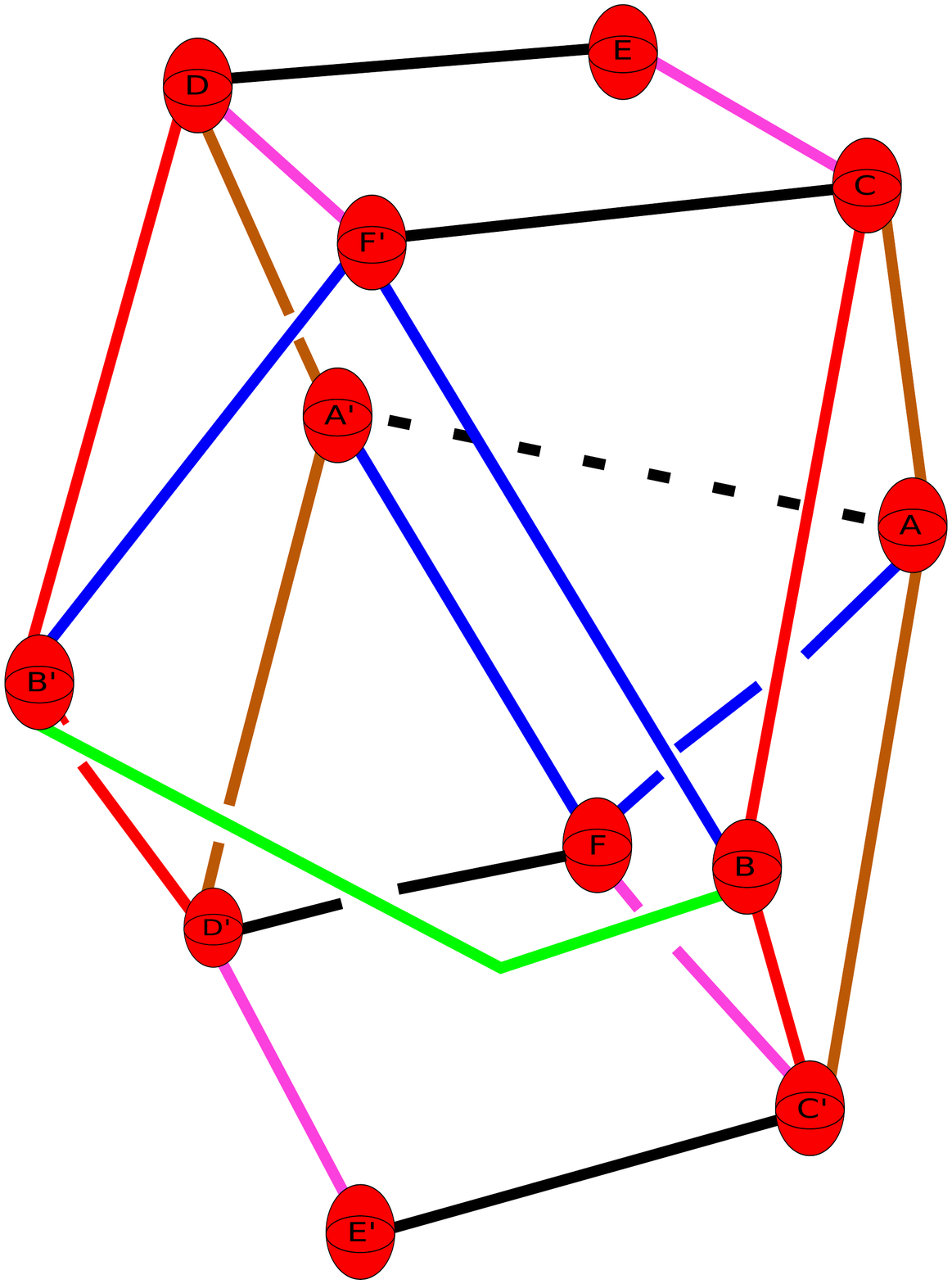}}

In general we will carry out several moves where we push a 2-handle component through the attaching sphere of a 1-handle. As all our attaching maps
are reflections or reflections composed with inversion in $S^2$, we will find that we will never be in a situation where the curve on one attaching
sphere is pushed to a curve in the image attaching sphere in such a way that it winds around other bits of 2-handle.

Before we proceed to doing some handle cancellations we remark that some of these added 2-handles will intersect some of the other planes. For example
consider the added 2-handle that runs between $C-C'$, as the co-ordinate of $C = (\frac{1}{\sqrt{2}}, 0, \frac{1}{\sqrt{2}})$ and 
$C' = (\frac{1}{\sqrt{2}}, 0, \frac{-1}{\sqrt{2}})$, we can see that in the bottom right picture above, this 2-handle must intersect the $x-y$ plane.
We have drawn this intersection point  as a black dot in the top left picture below. Similarly the added 2-handle that runs over $A-A'$ intersects
the $y-z$ plane in a point, and this has been drawn as a black dot in the bottom left picture below. Finally, the 2-handle passing over
$J-J'$ intersects the $x-z$ plane in a point, and this can be seen as a black dot in the top right picture below.

\centerline{\graphicspath{ {filling_cusps/} }\includegraphics[width=12cm, height=13cm]{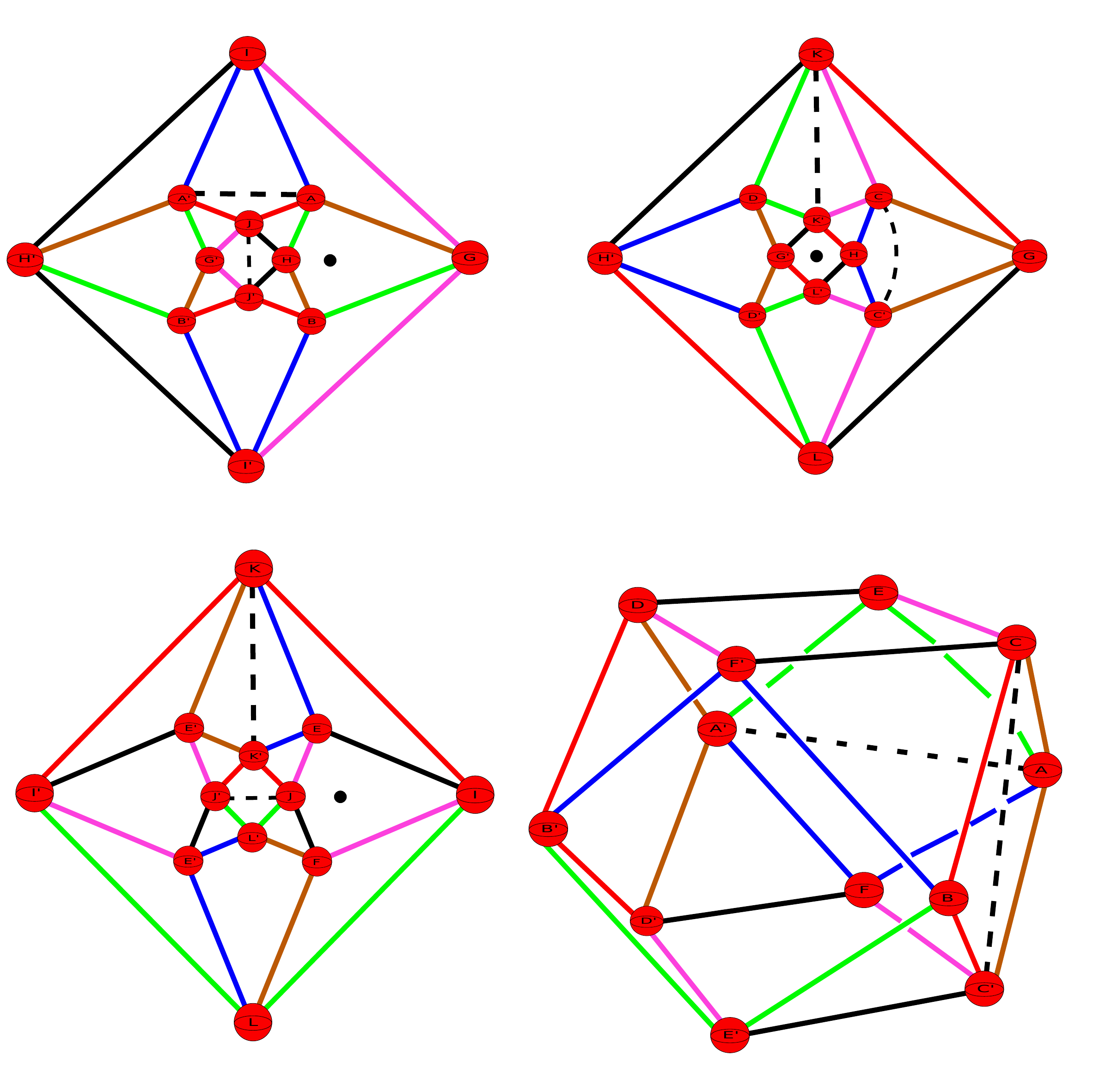}}

When we start cancelling handles and doing handle slides we have to be careful that we do not cut through any of the other 2-handles for
then we would be doing an ``illegal'' move, and that if one of our 2-handles goes around another, then a handle slide on one may cause the other to change position
and we have to keep track of this. It turns out that none of this really creates any problems as all our moves are primarily done in the planes in question, and
it can be clearly seen that they do not tangle around any of the added 2-handles that intersect these planes.

The above discussion may seem confusing so we go through a few handle cancellations in detail to make our point.

The added 2-handles corresponding to boundary fillings pass over a 1-handle once and hence form a handle cancellation pair, therefore we may cancel them from our diagram. We recall that
if we have a handle cancelling pair that does not meet any other 2-handles in a Kirby diagram then we can simply delete them from the diagram. However, if
there are some 2-handles that pass over the 1-handle that we are cancelling, we must push those 2-handles through the 2-handle that is cancelling the
1-handle.

We start by cancelling $A-A'$ (the reader should compare this with the first elementary move we showed in the previous section) using the added 2-handle that passes over it once. This only affects the 2-handles in the $x-y$ plane and the six extra
2-handles that do not lie in any one plane. The following pictures show the Kirby diagrams of these handles once we have done the cancellation.

\centerline{\graphicspath{ {handle_cancellation/cancelling_A-A'/} }\includegraphics[width=10cm, height=10cm]{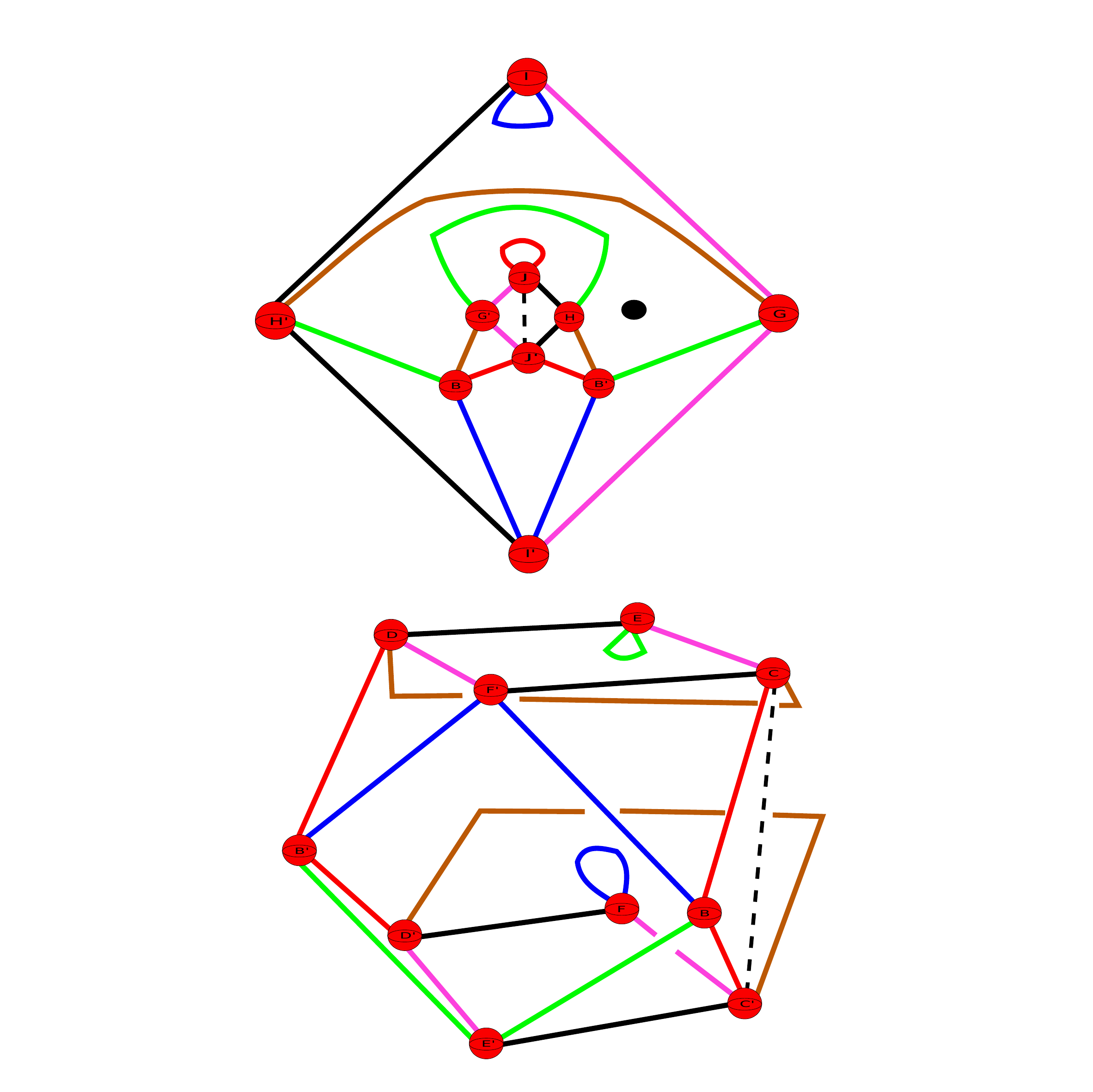}}

Recall that this 2-handle that we just used to cancel $A-A'$ intersected the $y-z$ plane, hence when we use it to cancel $A-A'$ we will introduce new intersection
points in the $y-z$ plane, and it is important to keep track of these new intersection points. 

We start by giving an analysis of what happens when we cancel $A-A'$ from the $x-y$ plane. First of all, the cancellation from the $x-y$ plane introduces some
new 2-handles in the $x-y$ plane. We have the blue 2-handle that starts at $I$ and loops back into $I$. This creates an intersection point in the $y-z$ plane
close to $I$, you can see it as the black dot in the following diagram:

\centerline{\graphicspath{ {handle_cancellation/cancelling_A-A'/} }\includegraphics[width=6cm, height=5cm]{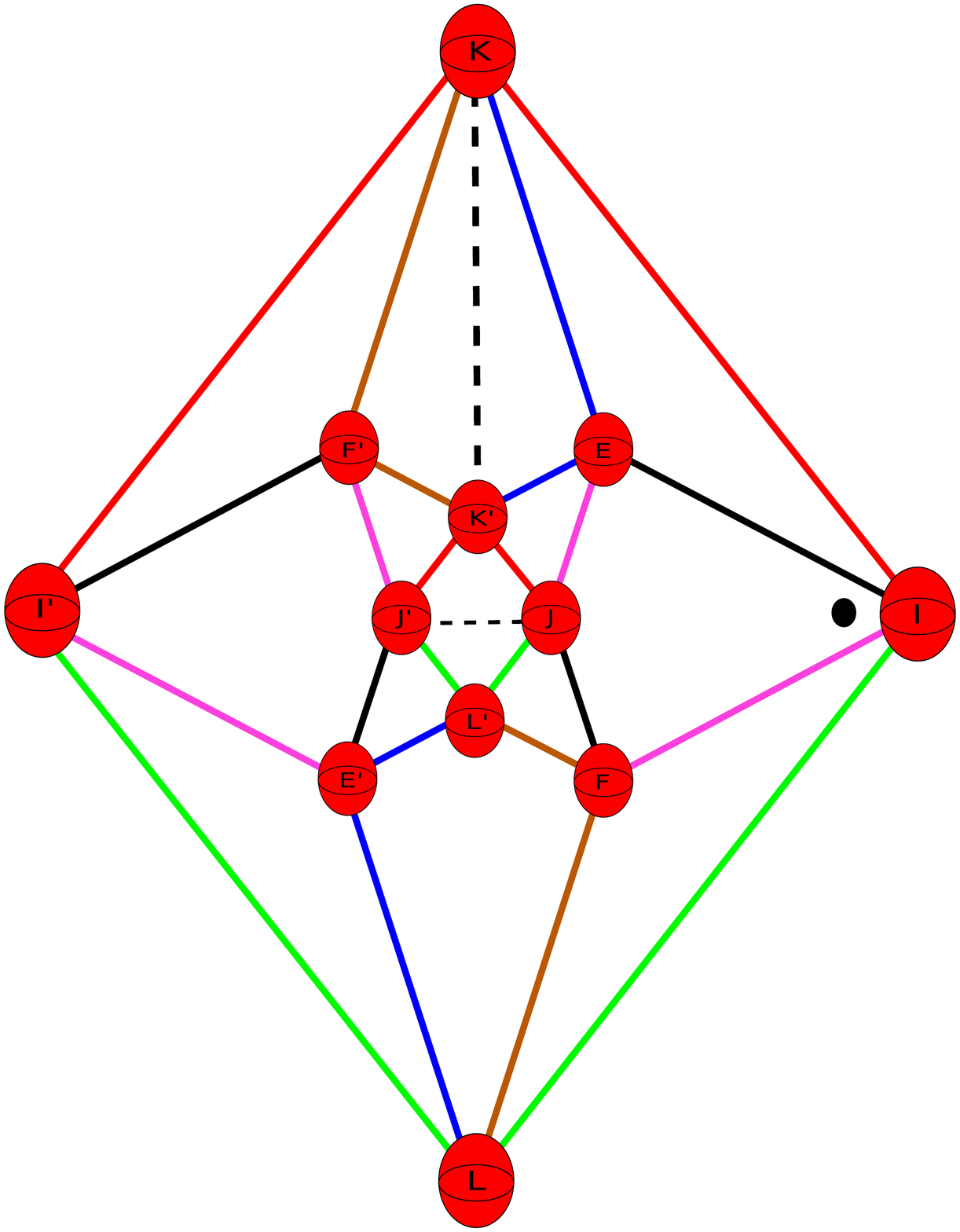}}

We also have the new brown 2-handle in the $x-y$ plane that goes between $G$ and $H'$, this creates a second intersection point with the $y-z$ plane:

\centerline{\graphicspath{ {handle_cancellation/cancelling_A-A'/} }\includegraphics[width=6cm, height=6cm]{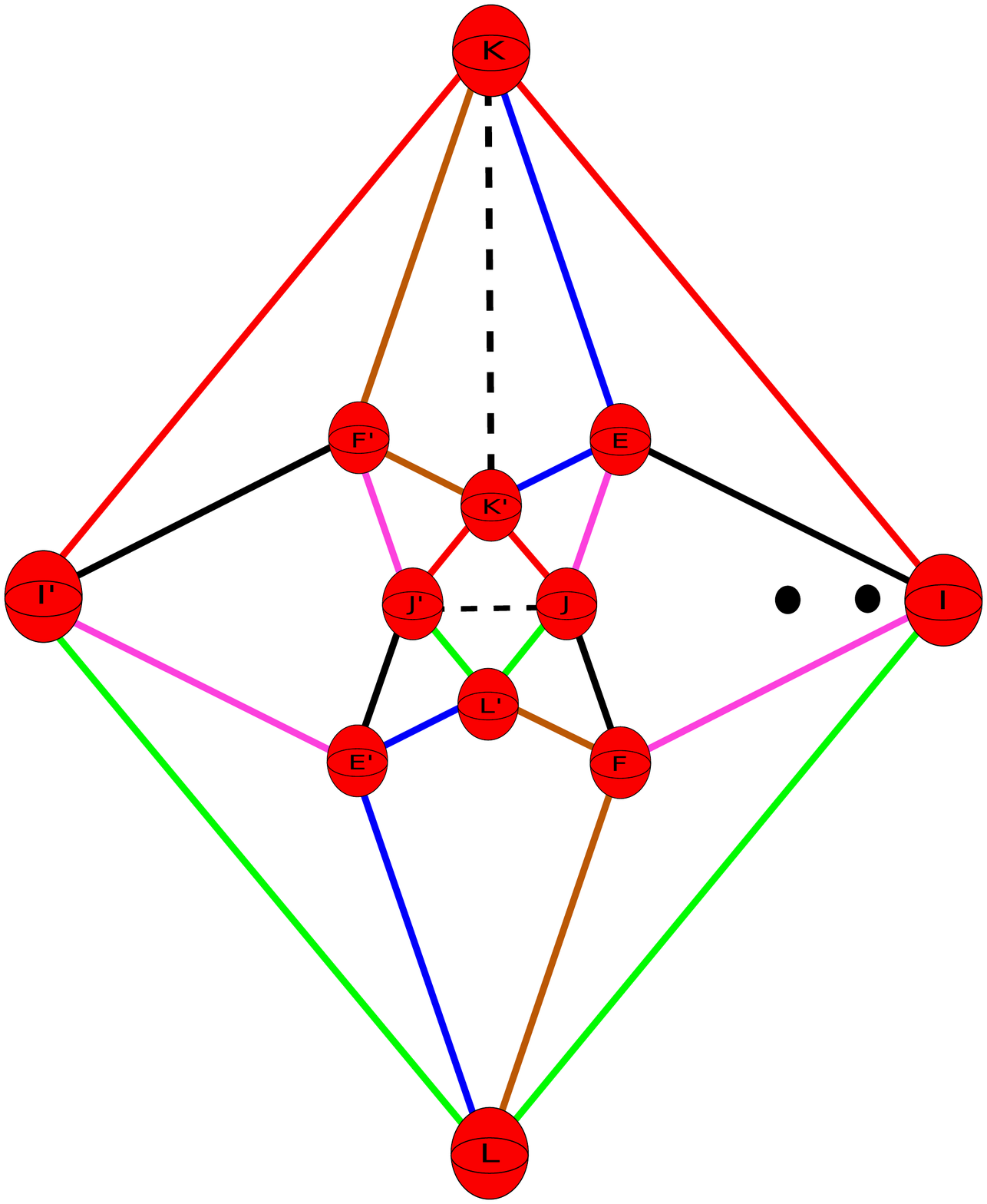}}

Remaining in the $x-y$ plane we also have the new green 2-handle that passes between $G'$ and $H$, this introduces a third intersection point in the $y-z$ plane.

\centerline{\graphicspath{ {handle_cancellation/cancelling_A-A'/} }\includegraphics[width=6cm, height=6cm]{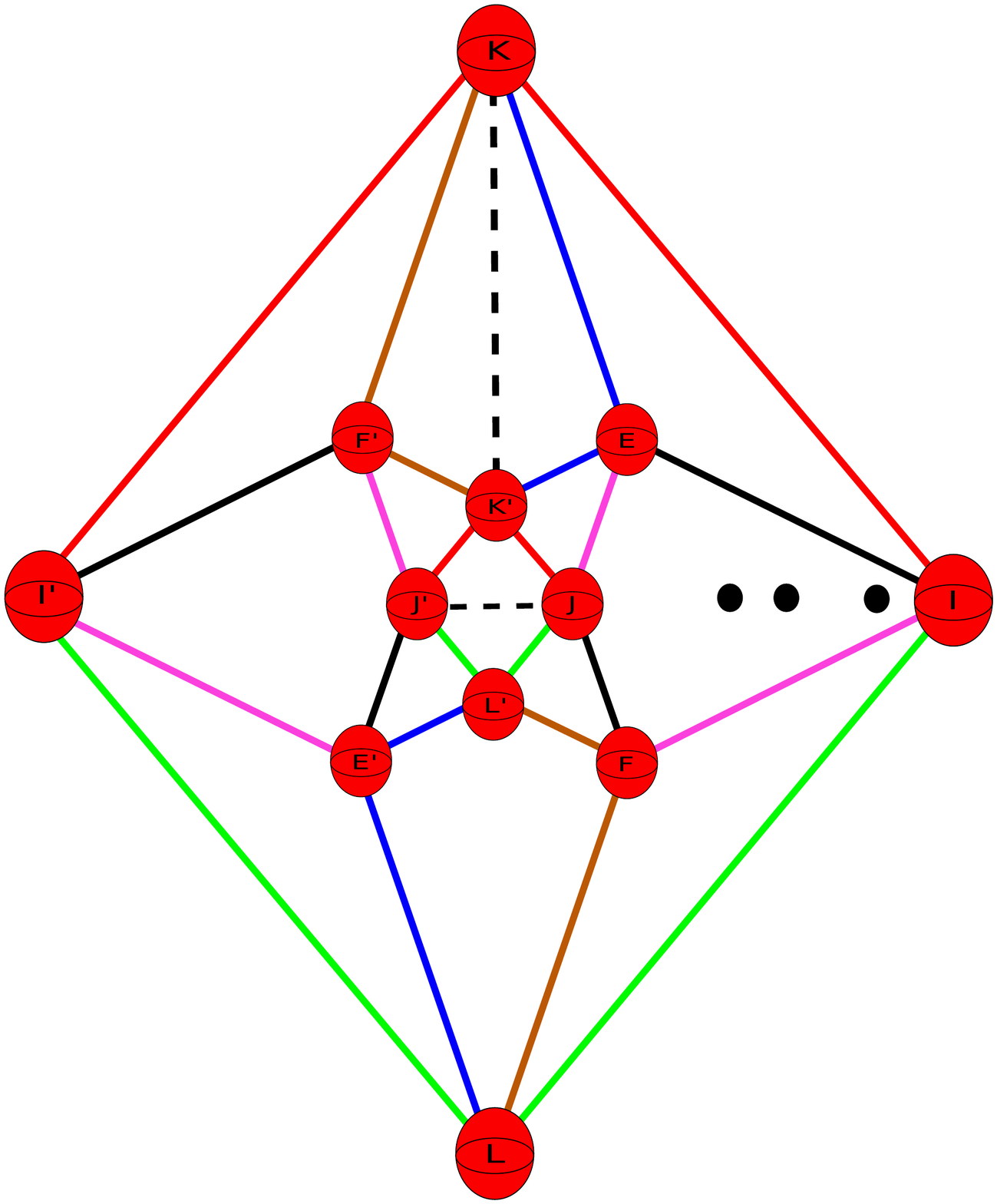}}

Finally, we have the red 2-handle that starts at $J$ and loops back into it, this gives a fourth intersection point.

\centerline{\graphicspath{ {handle_cancellation/cancelling_A-A'/} }\includegraphics[width=7cm, height=6.5cm]{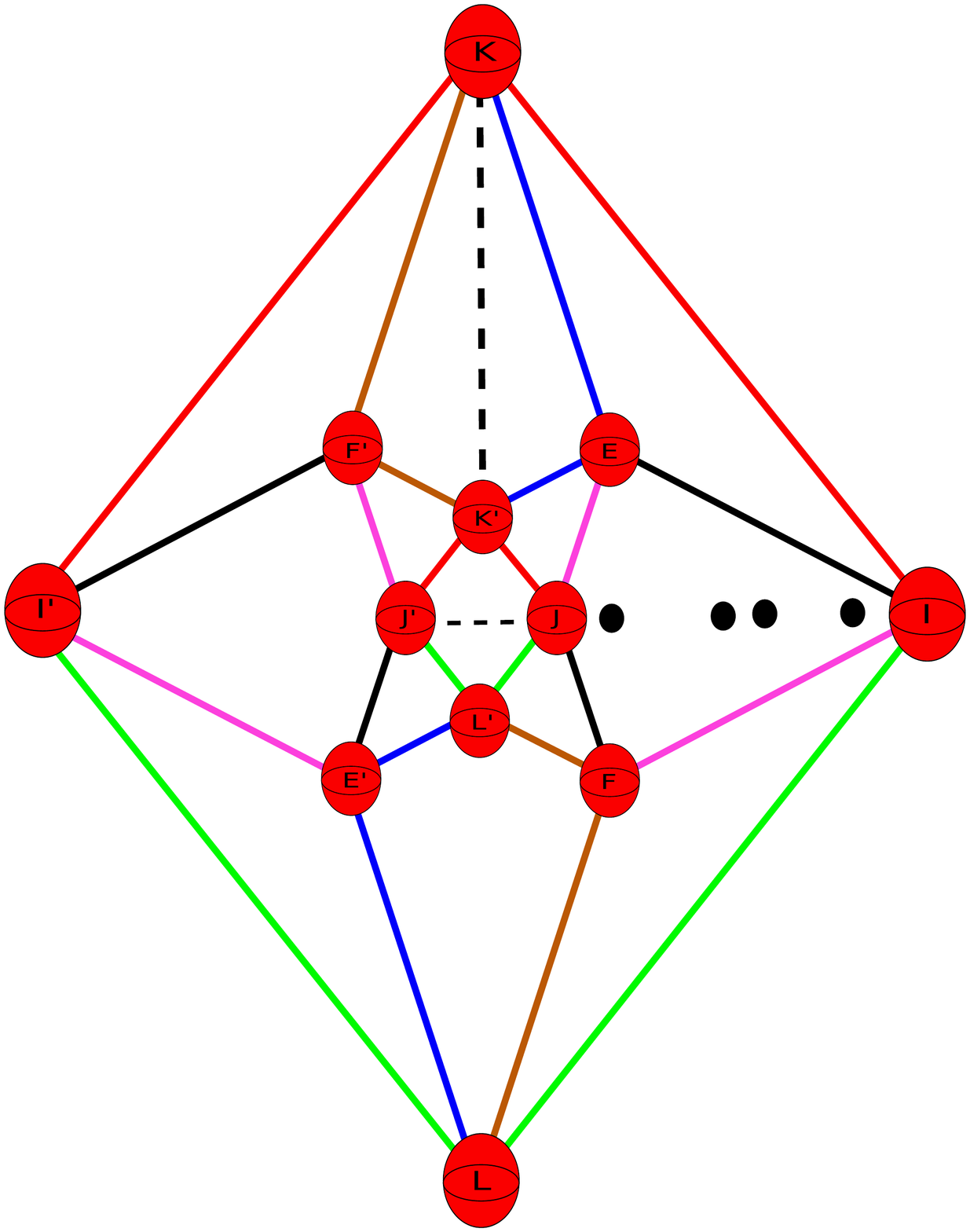}}

In total we obtain four new intersection points from cancelling $A-A'$ from the $x-y$ plane. However we are not done yet, there will also be some
intersection points arising from the cancellation of $A-A'$ from the diagram containing the six 2-handles that do not all lie in any 2-plane.

The analysis of this situation is exactly analogous to what we did above. We start with the green 2-handle that starts at $E$ and loops back into it. This gives
an intersection point in the $y-z$ plane near $E$ as shown in the following diagram.

\centerline{\graphicspath{ {handle_cancellation/cancelling_A-A'/} }\includegraphics[width=7cm, height=7cm]{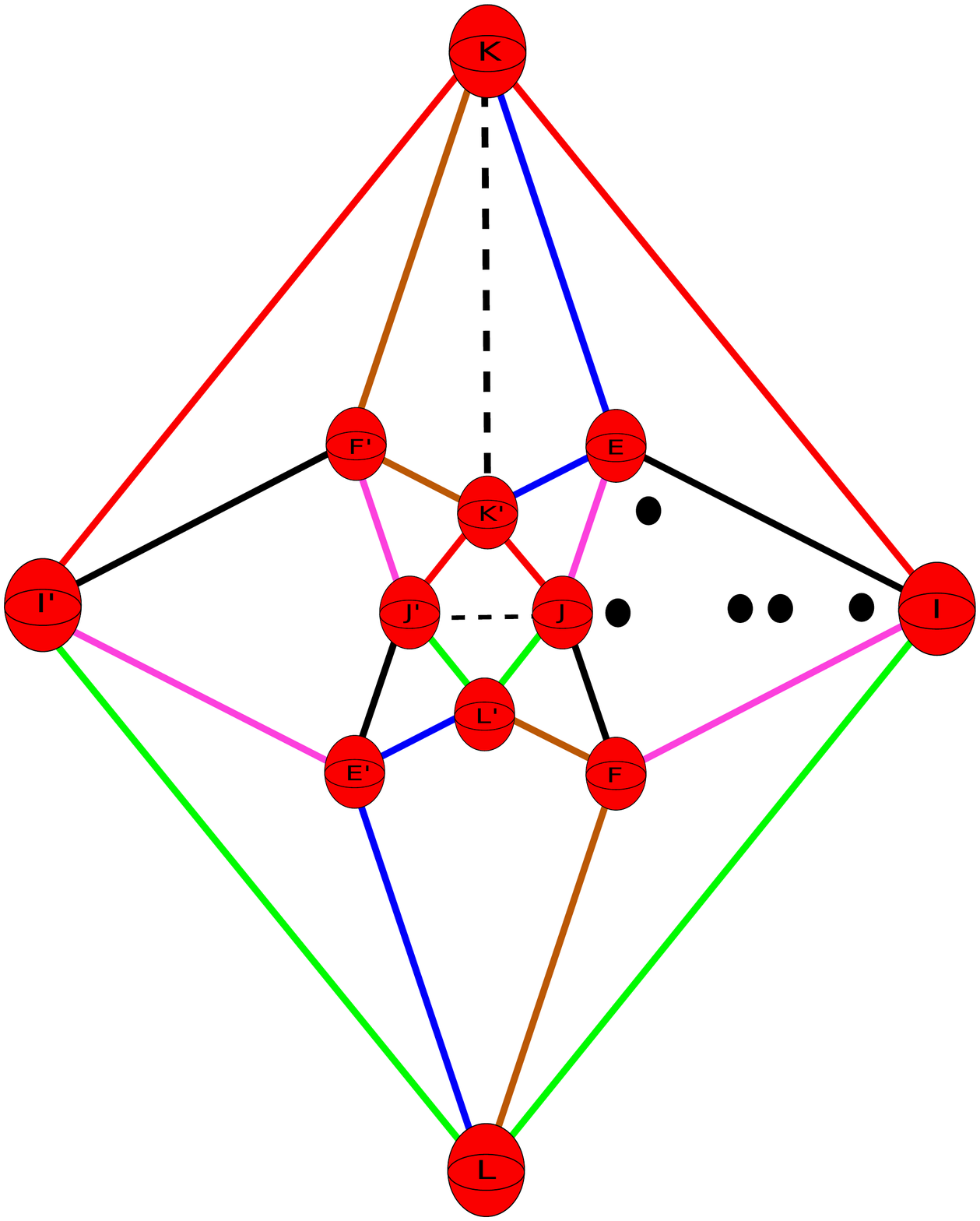}}

The brown 2-handle that runs between $C$ and $D$ gives an intersection point that looks like:

\centerline{\graphicspath{ {handle_cancellation/cancelling_A-A'/} }\includegraphics[width=7cm, height=7cm]{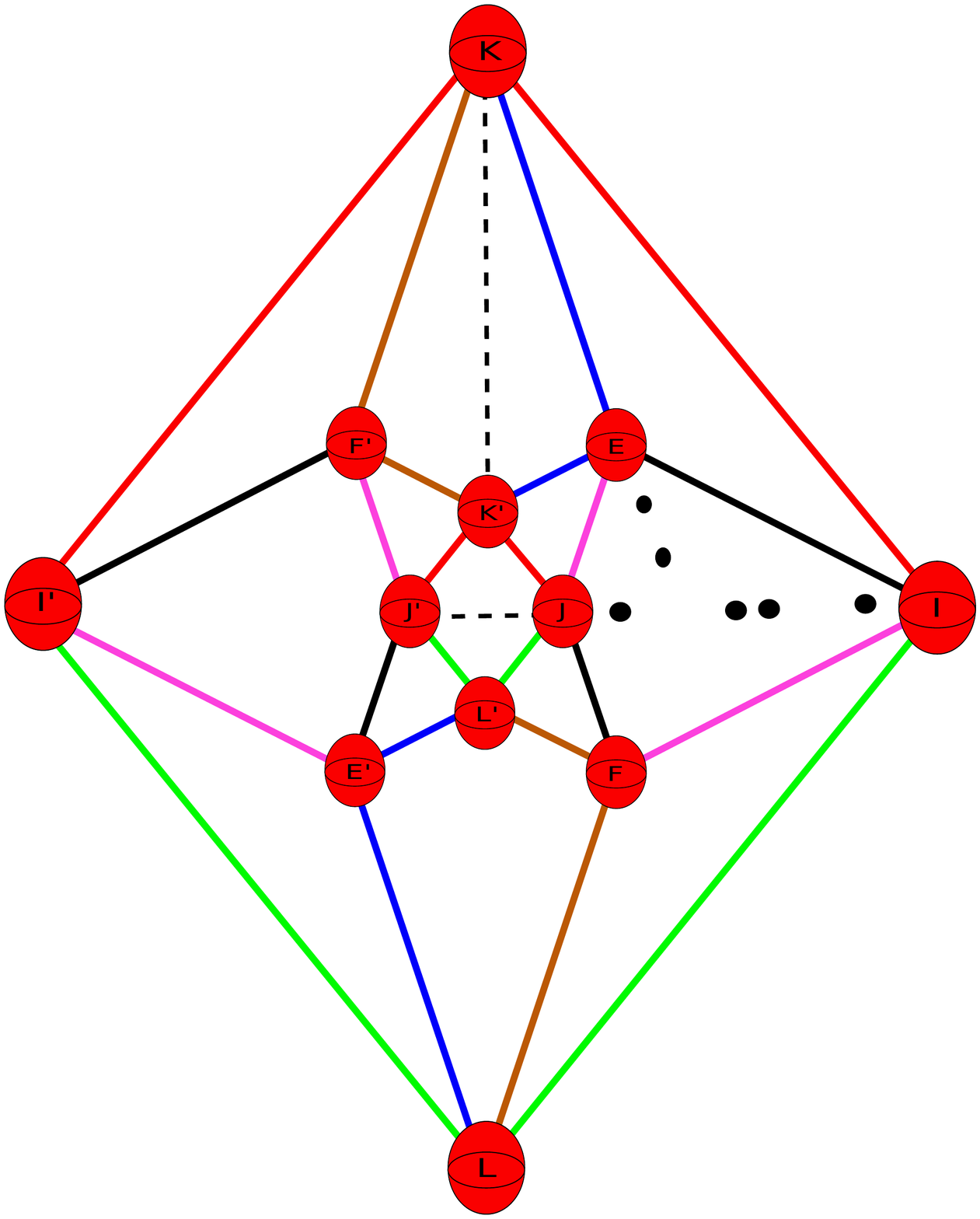}}

There is also the brown 2-handle that runs between $C'$ and $D'$, this gives a point of intersection that looks like:

\centerline{\graphicspath{ {handle_cancellation/cancelling_A-A'/} }\includegraphics[width=7cm, height=7cm]{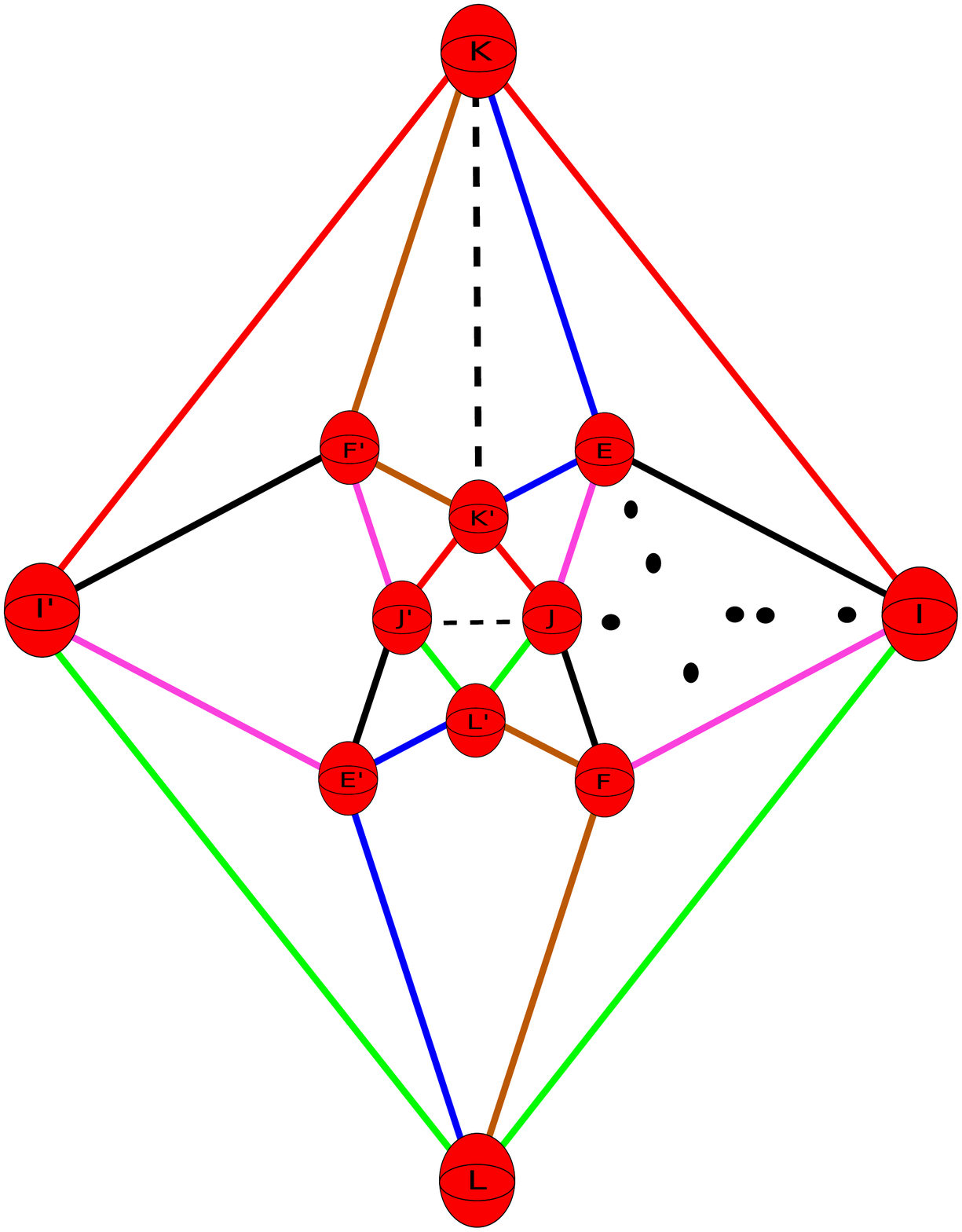}}

Finally, we have the 2-handle that starts at $F$ and loops back into it, it gives an intersection point near $F$:

\centerline{\graphicspath{ {handle_cancellation/cancelling_A-A'/} }\includegraphics[width=7cm, height=7cm]{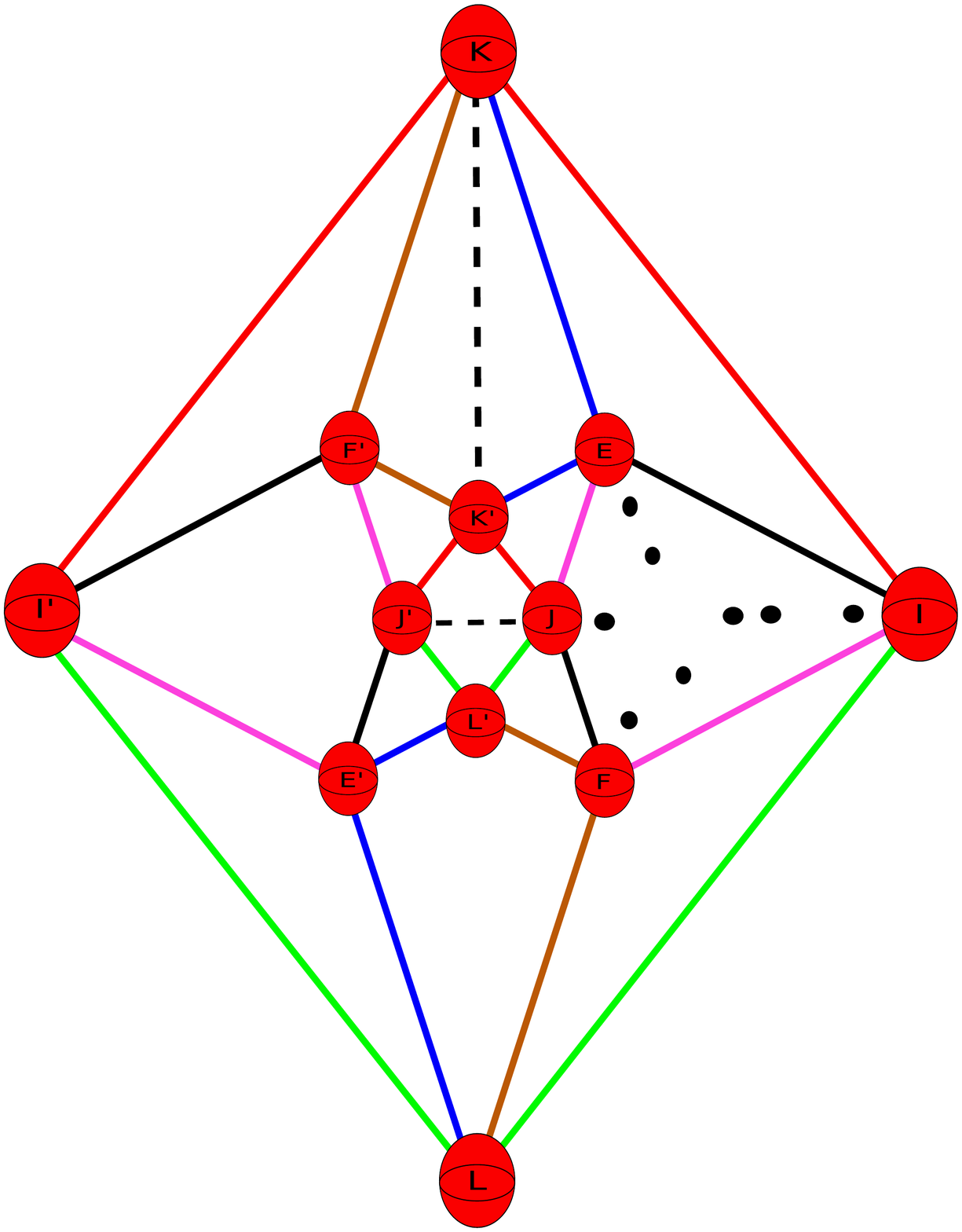}}

In total we get 8 new intersection points, to make things easier we want to keep track of how each intersection came about. Namely, we want to keep track
of which 2-handle in which diagram gave us a particular intersection point. To do this we will introduce a simple coding system that will allow us to keep track
of where these points of intersection are coming from. Each intersection point will be labelled with either four letters or two letters. In the case of
four letters, the first two will tell us that the intersection point is coming from the $x-y$, $x-z$ or $y-z$ planes, and the next two letters will tell us the
1-handles that the 2-handle causing the intersection point is passing between, the actual 2-handle will be clear from context. In the case of a two letter code, we 
are to interpret that the 2-handle causing the intersection point is coming from the diagram representing the six 2-handles that do not all lie in a single 2-plane.
The two letters denote the two 1-handles that this 2-handle passes over, again the actual 2-handle in question will be clear from the context.

At times it will be difficult for us to show the labels of all intersection points in the whole diagram simply because there may be too many intersection points
and not enough space. In these instances we will always show two diagrams, a whole diagram with unlabelled points of intersection and a close up diagram
of the labelled intersection points. It will be easy to tell which point corresponds to which point in the two diagrams.

It is time to show an explicit example of what exactly we mean by this code, and how exactly it looks. The following picture shows the $y-z$ plane
with a black dot near $I$ labelled \textbf{XY\underline{\space}II}, the first two letters being XY tell us that the 2-handle giving this point of intersection
is coming from the $x-y$ plane, the second two letters are $II$ and this tells us that the 2-handle in the $x-y$ plane is starting at $I$ and looping back into
$I$.

\centerline{\graphicspath{ {handle_cancellation/cancelling_A-A'/} }\includegraphics[width=6cm, height=5.5cm]{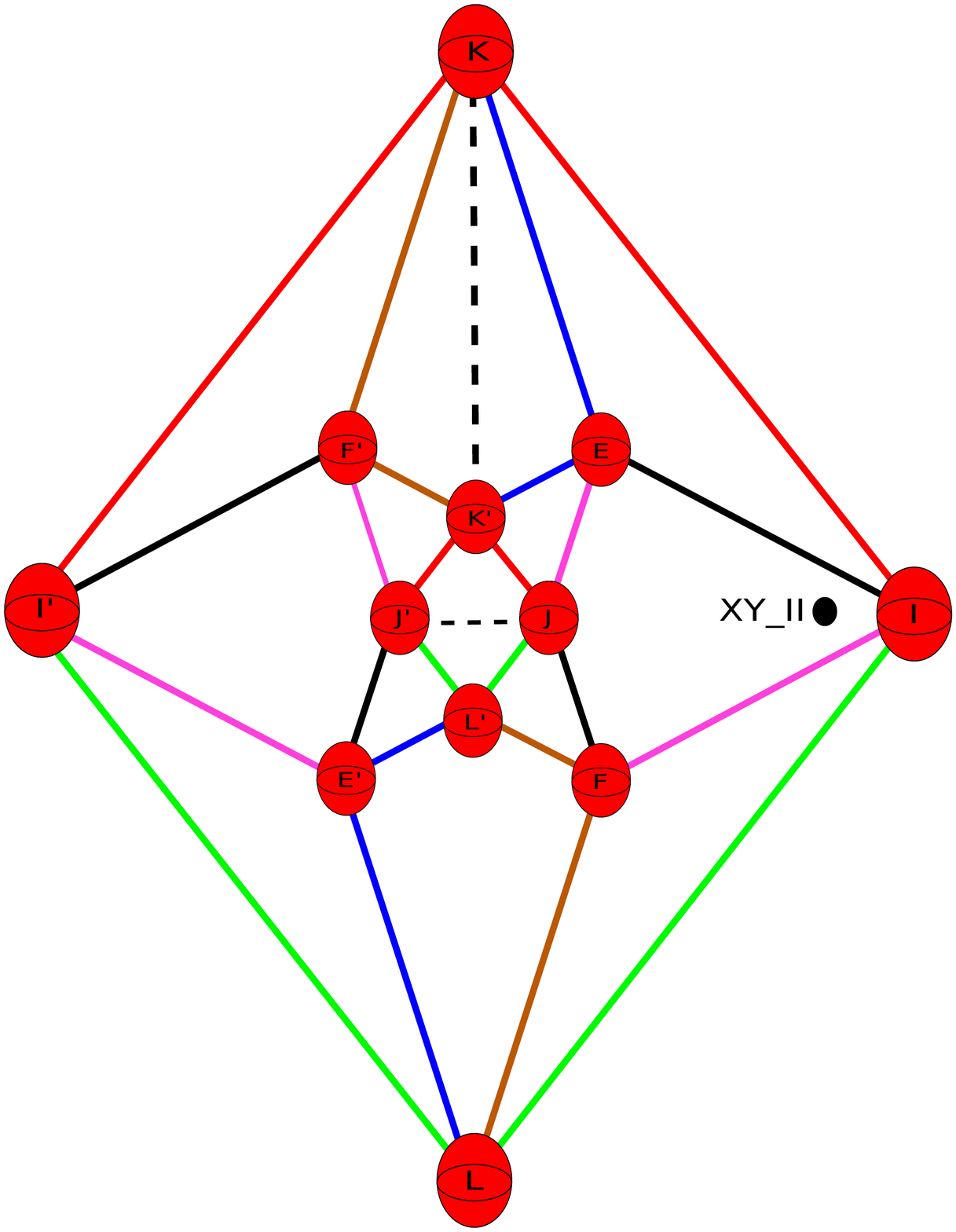}}

This is precisely the point of intersection we talked about above which arose through cancelling $A-A'$ in the $x-y$ plane, and we showed a picture
of it previously. There were also many other intersection points (remember there were 8 in total!), as another example the following picture shows a labelled 
intersection point coming from cancelling $A-A'$ in the diagram corresponding to the six 2-handles that do not all lie in a single 2-plane. The labelling
consists of just two letters, which immediately tells us that the 2-handle contributing this point of intersection is coming from the diagram corresponding
to the six 2-handles that do not all lie in a single 2-plane. The two letters are CD, which tell us that the 2-handle in question is running between 
the 1-handle $C$ and the 1-handle $D$.

\centerline{\graphicspath{ {handle_cancellation/cancelling_A-A'/} }\includegraphics[width=6cm, height=5.5cm]{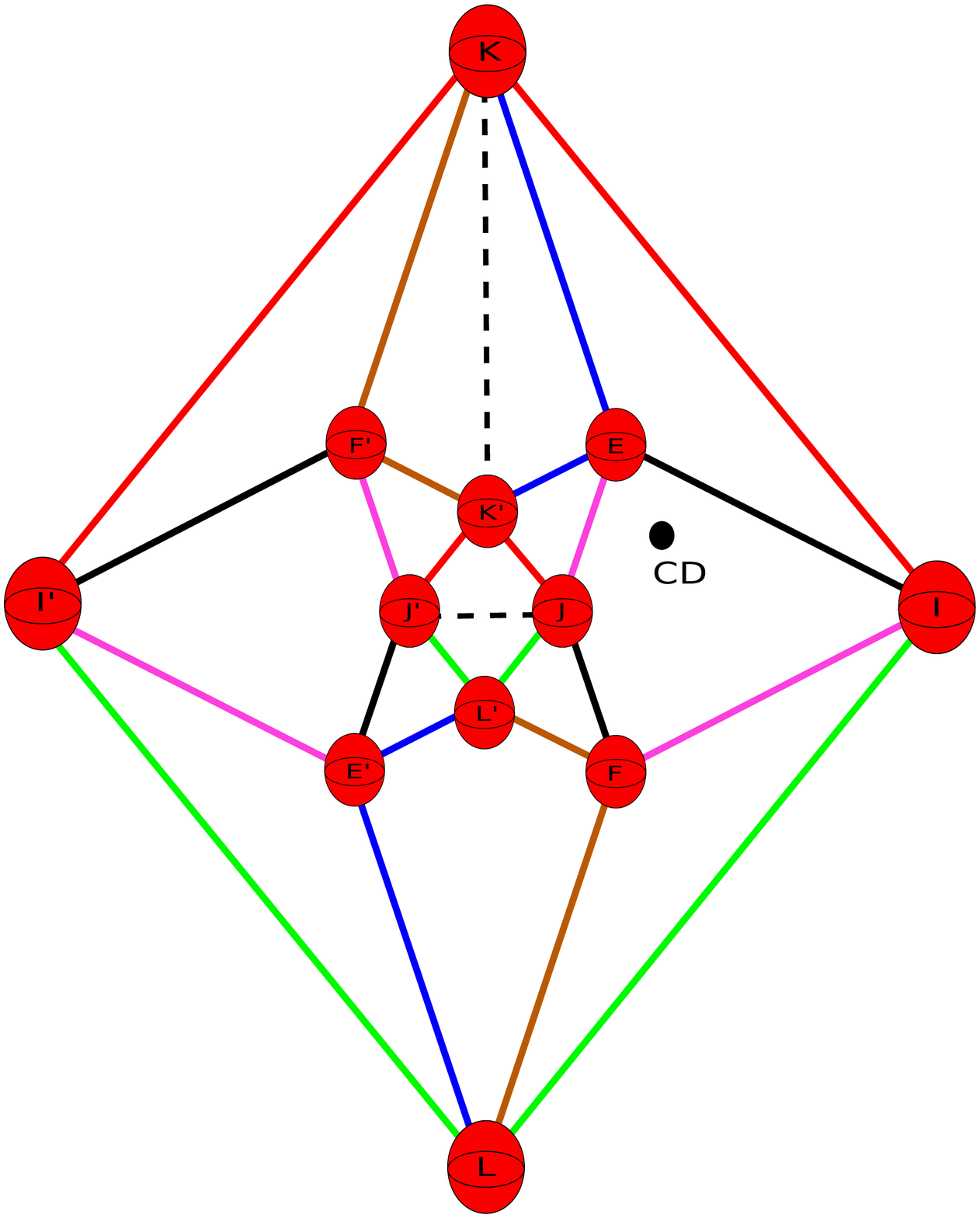}}

When we cancelled $A-A'$ we saw that there were a total of 8 new points of intersection. The following shows a close up of all these points of intersection
with their corresponding labels.

\centerline{\graphicspath{ {handle_cancellation/cancelling_A-A'/} }\includegraphics[width=7.5cm, height=6cm]{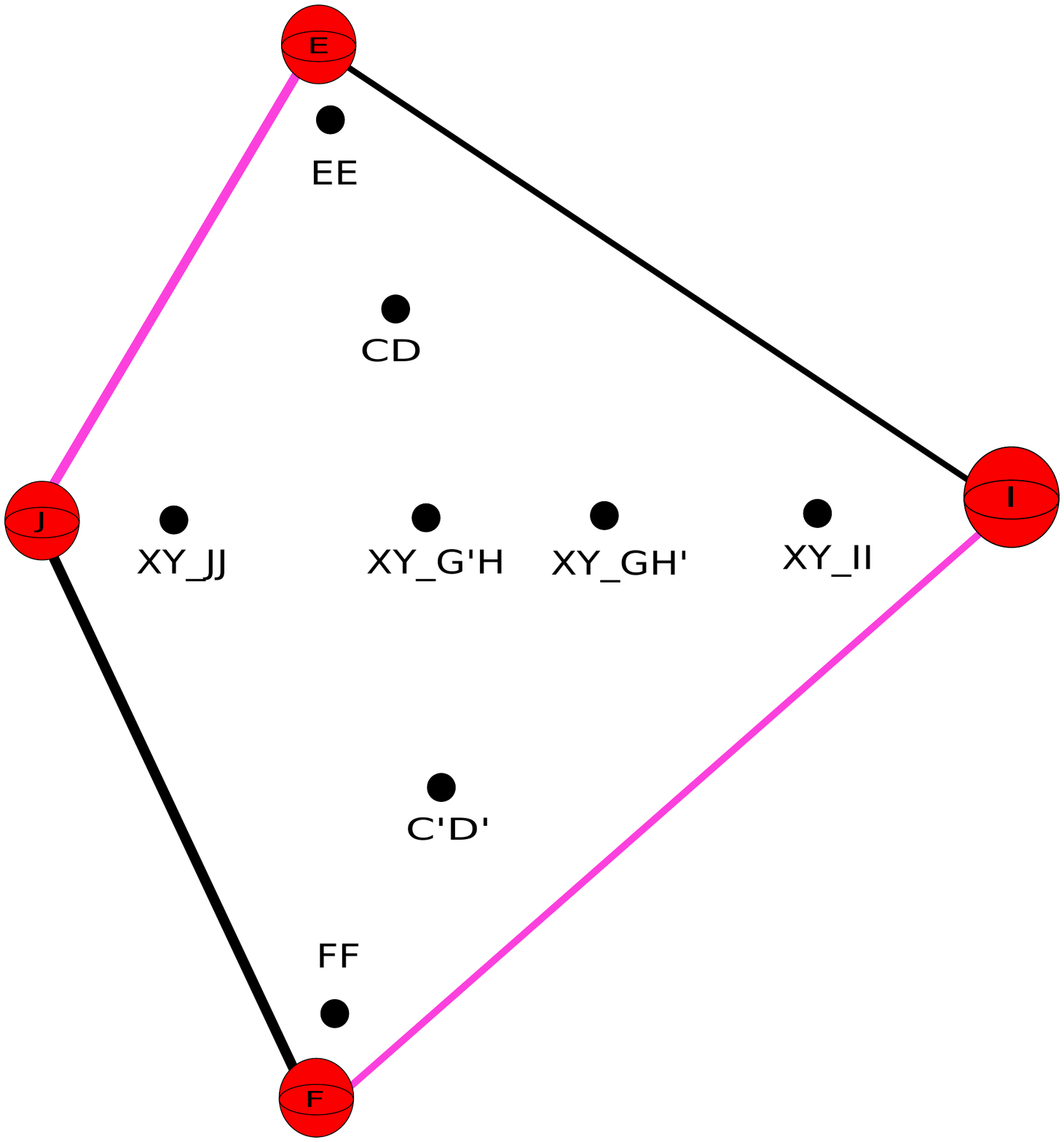}}

So far we have done one handle cancellation, and in some detail explained how we get extra points of intersection. The following picture shows
what we have done.

\centerline{\graphicspath{ {handle_cancellation/cancelling_A-A'/} }\includegraphics[width=9cm, height=10cm]{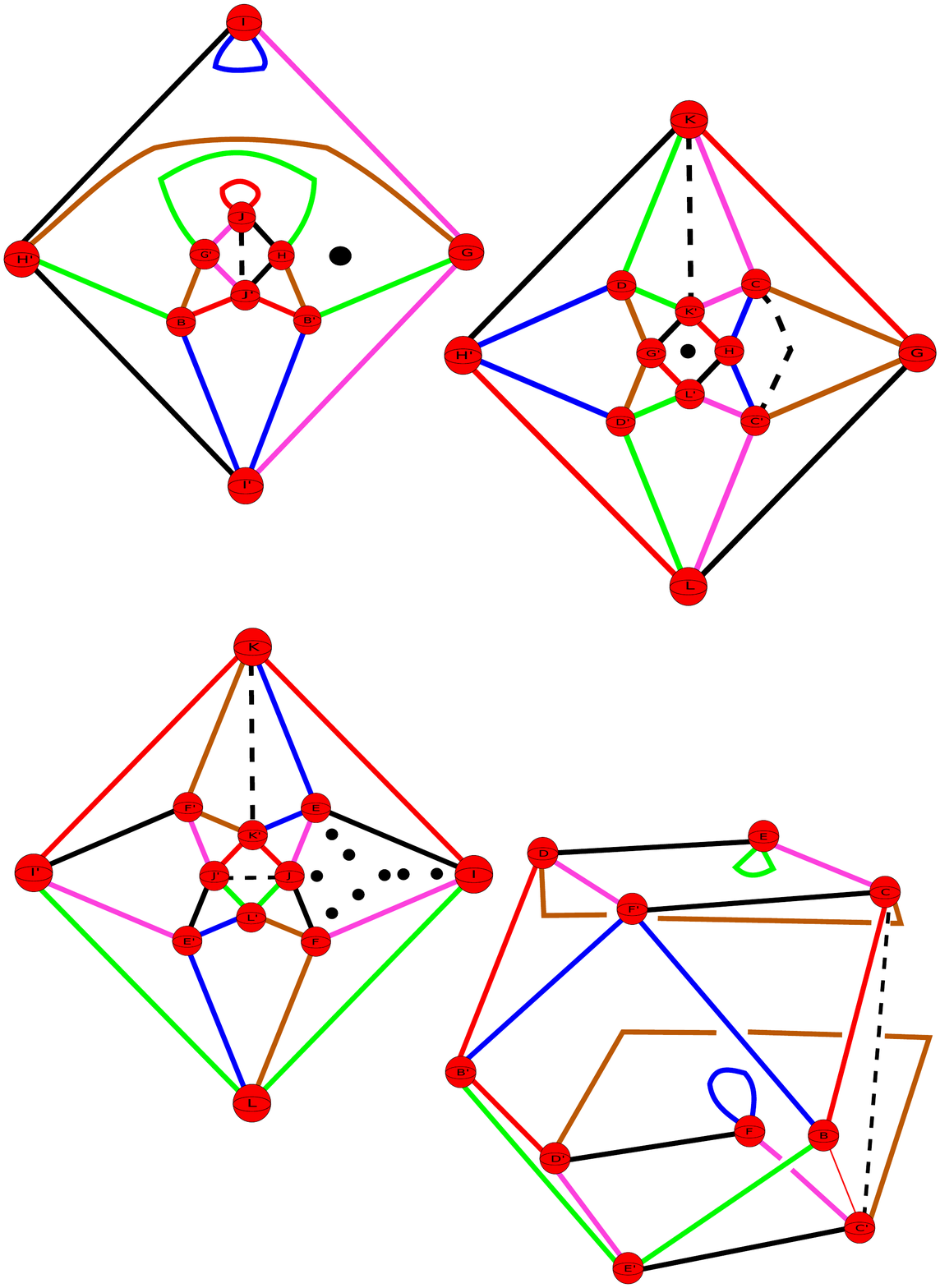}}

In the top left diagram we have a 2-handle, corresponding to (the one in blue)

$\xymatrix{
A \cap I \ar[r]^a & A'\cap I\ar[r]^i & B'\cap i' \ar[r]^{b^{-1}} & B\cap I' \ar[r]^{i^{-1}} & A\cap I }$

which has a component that loops back into $I$. In the above picture this is the blue arc at the top meeting $I$ twice. 
We can then slide it over the 1-handle $I-I'$ to give a 2-handle that meets $B-B'$ once.

There is also a second handle slide we can do in the $x-y$ plane. Namely, we have the 2-handle in red that starts at $J$ and loops back into it. We can then slide it
over the 1-handle $J-J'$ to give a 2-handle in red that also meets $B-B'$ once.

\centerline{\graphicspath{ {handle_cancellation/cancelling_A-A'/} }\includegraphics[width=6cm, height=5.5cm]{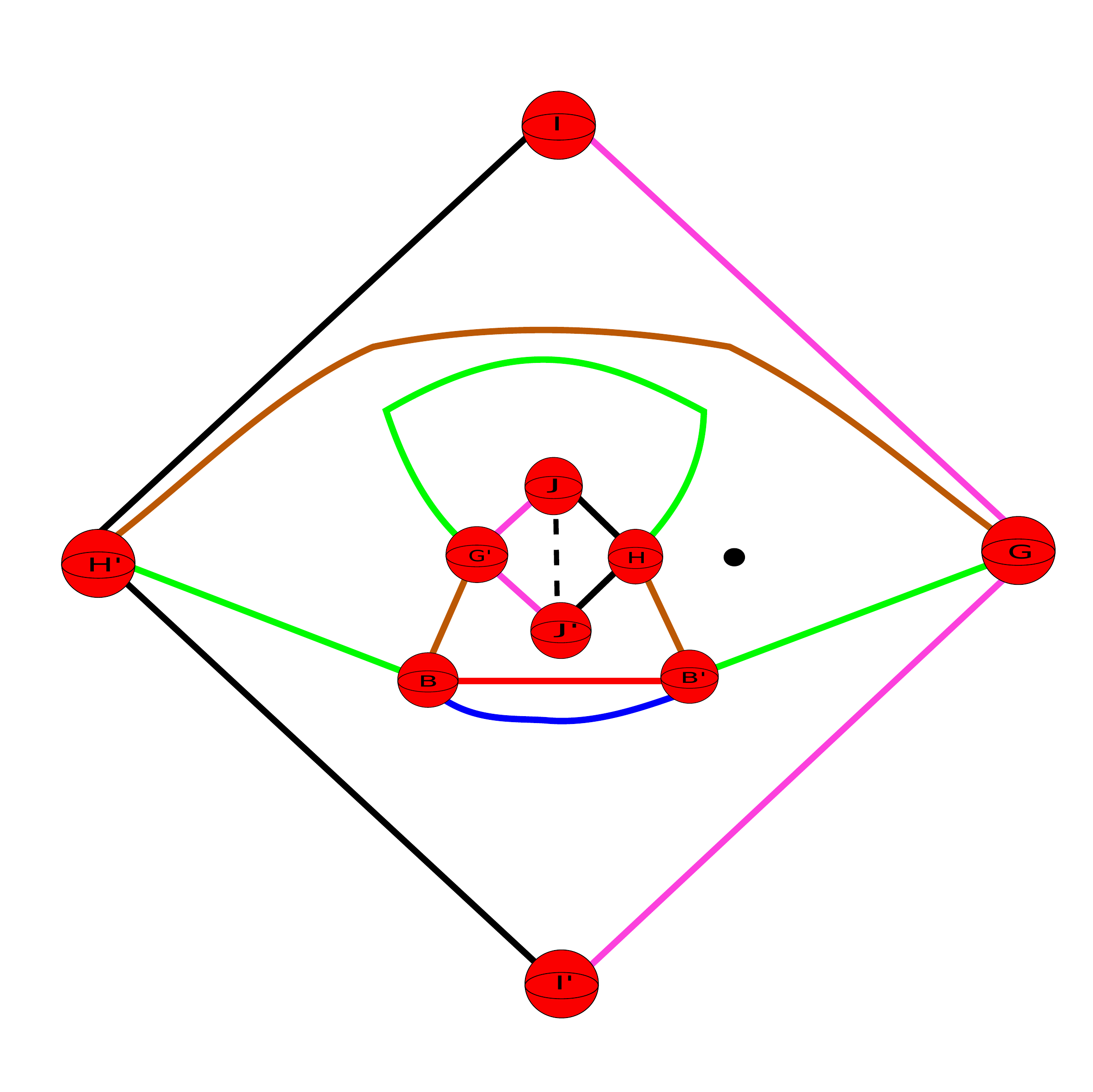}}

In carrying out these two handle slides the intersection points labelled \textbf{XY\underline{\space}II} and \textbf{XY\underline{\space}JJ} in the $y-z$ plane
will disappear. However, in the $x-y$ plane the 2-handle (in blue) that meets $B-B'$ once will give a point of intersection
in the $y-z$ plane in the region bounded by the 1-handles $I'$, $F'$, $J'$ and $E'$ and the 2-handles that run between them, and the red 2-handle
that meets $B-B'$ will also give an intersection point in the $y-z$ plane in the region bounded by the 1-handles $I'$, $F'$, $J'$, $E'$ and the 2-handles that run between them.

The picture below shows these new intersection points along with the other intersection points that we have encountered so far. 
The two diagrams at the bottom of the picture show close ups of the regions that these intersection points lie in and the labels of the intersection
points.

\centerline{\graphicspath{ {handle_cancellation/cancelling_A-A'/} }\includegraphics[width=8cm, height=9cm]{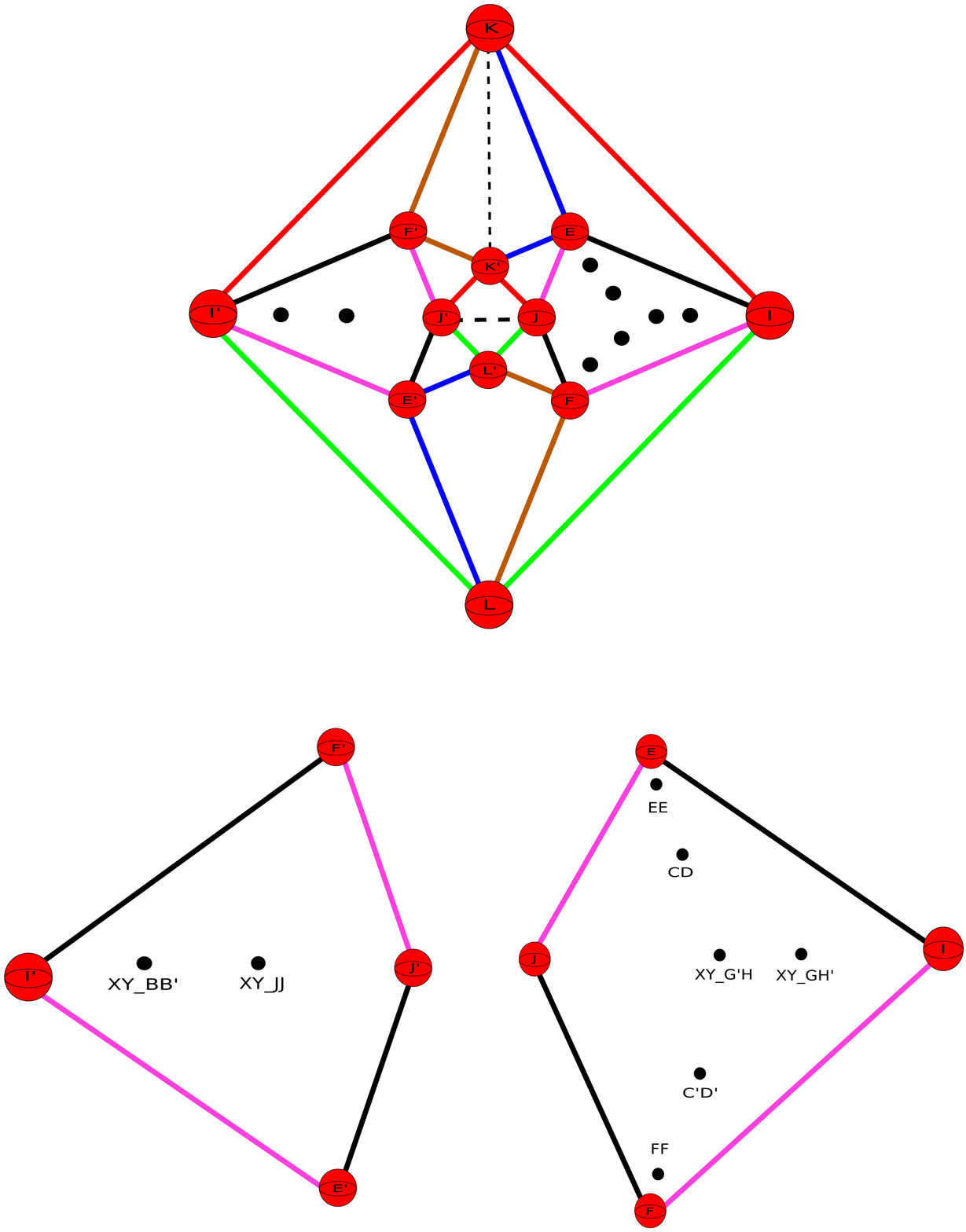}}

We can also carry out two handle slides in the diagram consisting of the 2-handles that do not all lie in a single 2-plane. We can slide the green
2-handle that starts and ends at $E$ through $E$ and then off $E'$ to obtain a 2-handle that meets $B-B'$ once. Similarly we can slide the blue 2-handle
that starts and ends at $F$ through $F$ and then off $F'$ to give another 2-handle that meets $B-B'$ once.

\centerline{\graphicspath{ {handle_cancellation/cancelling_A-A'/} }\includegraphics[width=6cm, height=5.5cm]{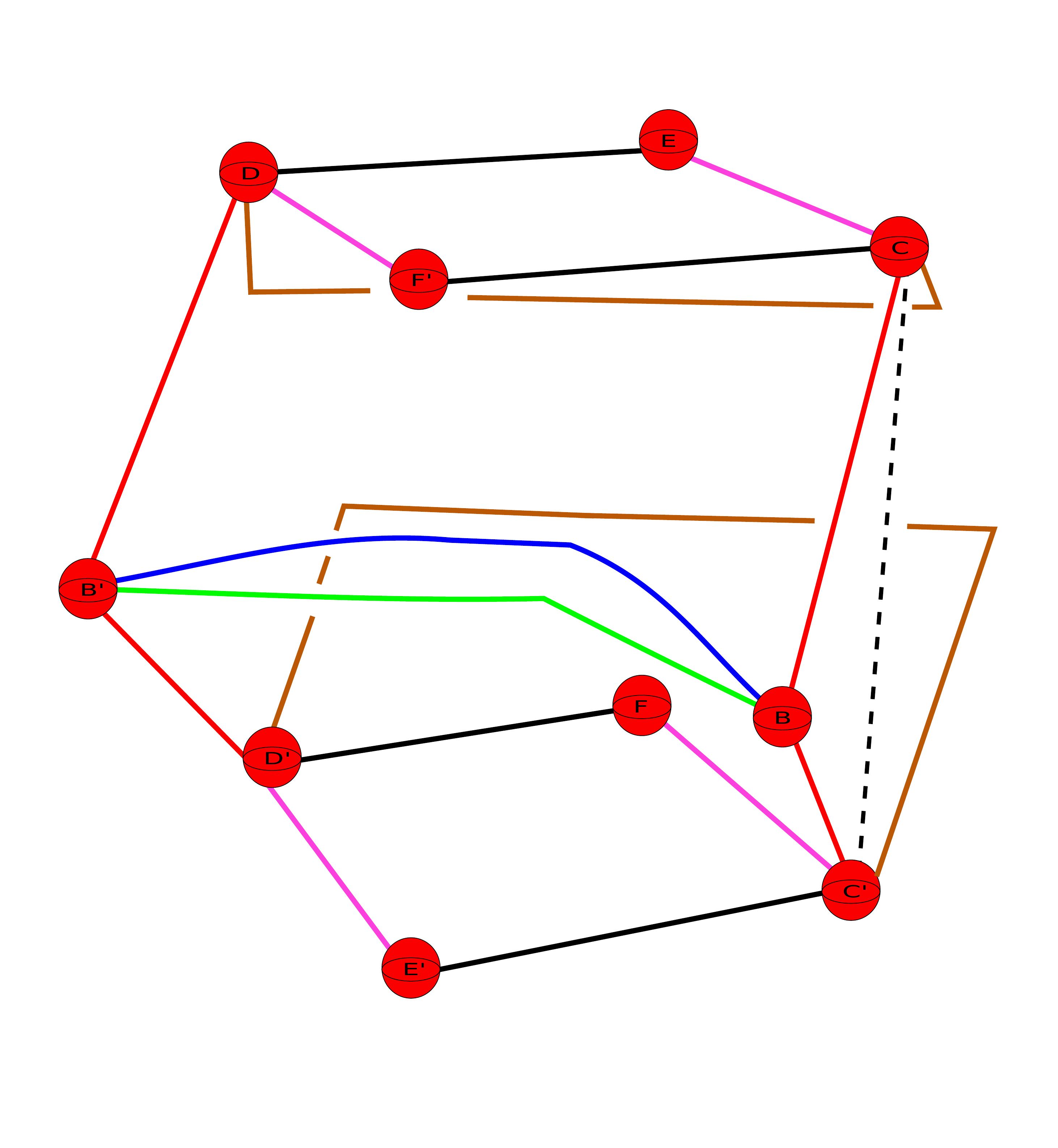}}

These handle slides that we have carried out have not interfered with the added 2-handles running from $E$ to $G$ and $E'$ to $G'$. In the case of
the added 2-handle that runs between $E$ and $G$ this is easy to see as it lies outside of the diagrams for which these handle cancellations and slides are
being done. In the case of the 2-handle that runs between $E'$ and $G'$ it is also easy to see, using some three dimensional insight, that none
of the cancellations and slides we have done so far affect the handle. The following picture shows how this 2-handle sits after we have carried out the above
cancellations and slides, it is drawn in yellow. The picture is supposed to give the reader some insight in to why it is the case that 
carrying out the above cancellations and slides does not affect the added 2-handle that runs between $E'$ and $G'$.

\centerline{\graphicspath{ {handle_cancellation/cancelling_A-A'/} }\includegraphics[width=6cm, height=5cm]{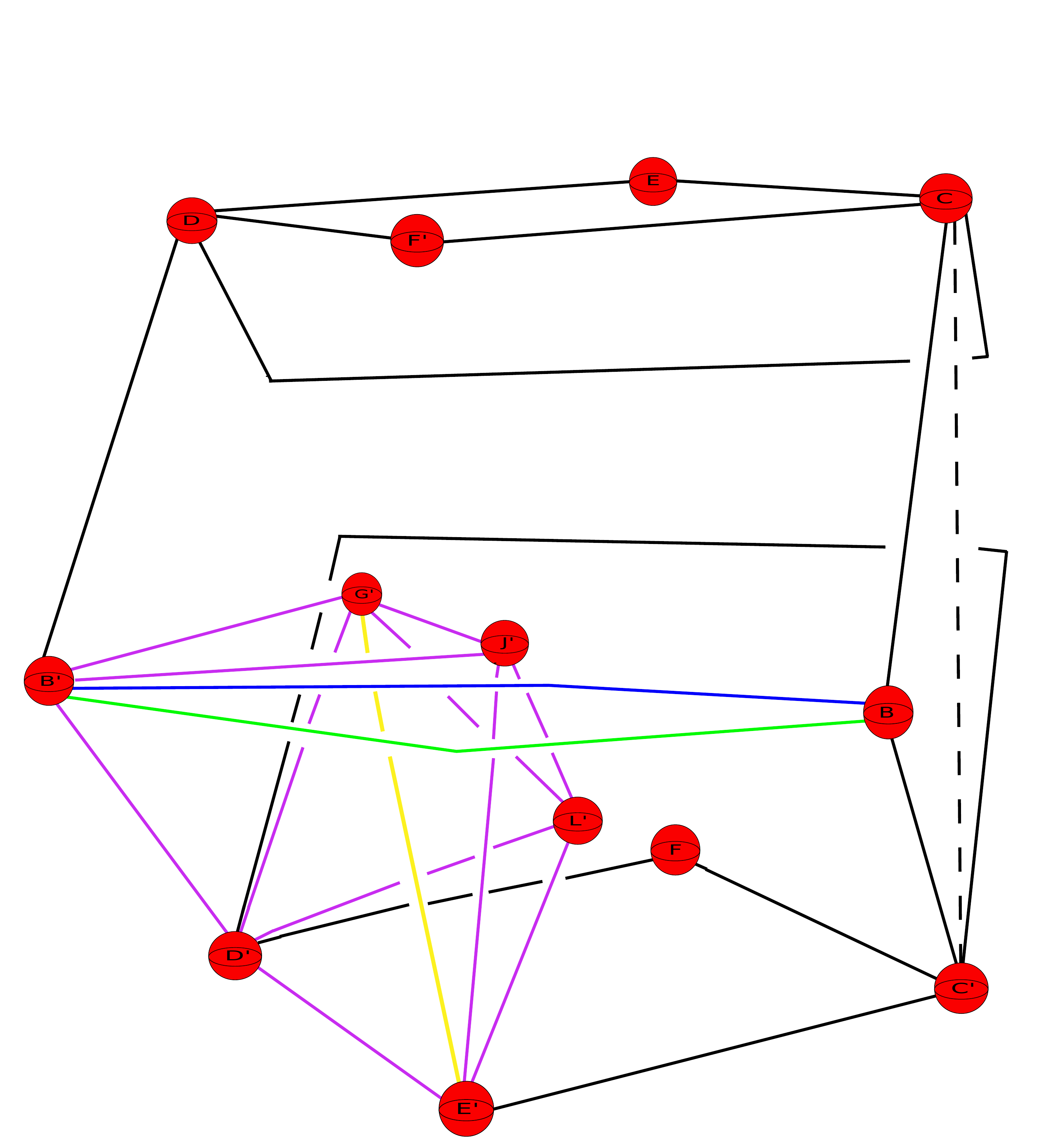}}

The two handle slides just carried out above will cause the intersection points in the $y-z$ plane labelled  \textbf{EE} and  \textbf{FF} to disappear, but the two new 2-handles
between $B-B'$, that arose through these handle slides, will give two new points of intersection in the $y-z$ plane in the region bounded by
$I'$, $F'$, $J'$, $E'$ and the 2-handles that run between them.

\centerline{\graphicspath{ {handle_cancellation/cancelling_A-A'/} }\includegraphics[width=8cm, height=9cm]{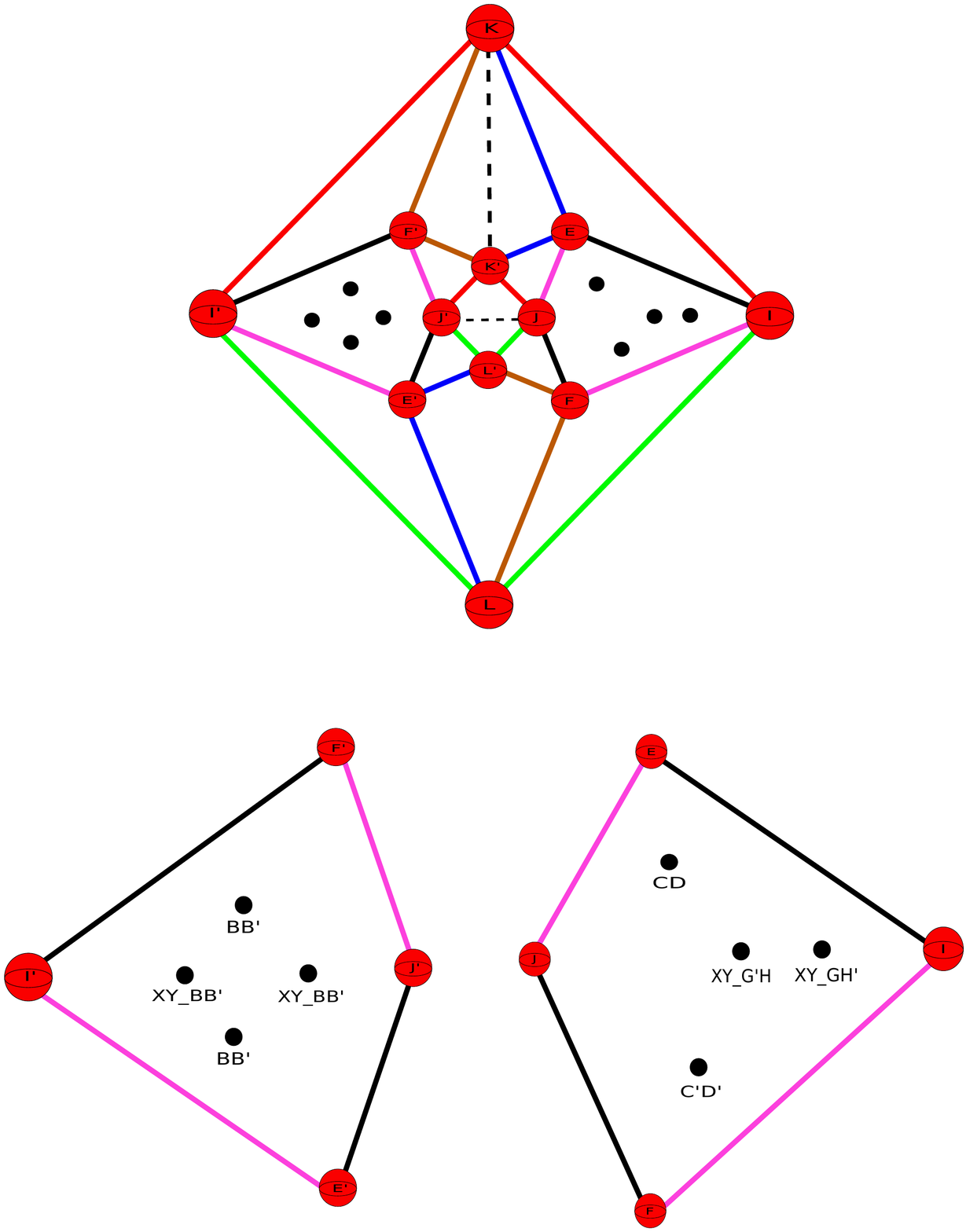}}

The following picture shows the end result of cancelling $A-A'$, along with the various handle slides we under took after this cancellation

\centerline{\graphicspath{ {handle_cancellation/cancelling_A-A'/} }\includegraphics[width=8cm, height=10cm]{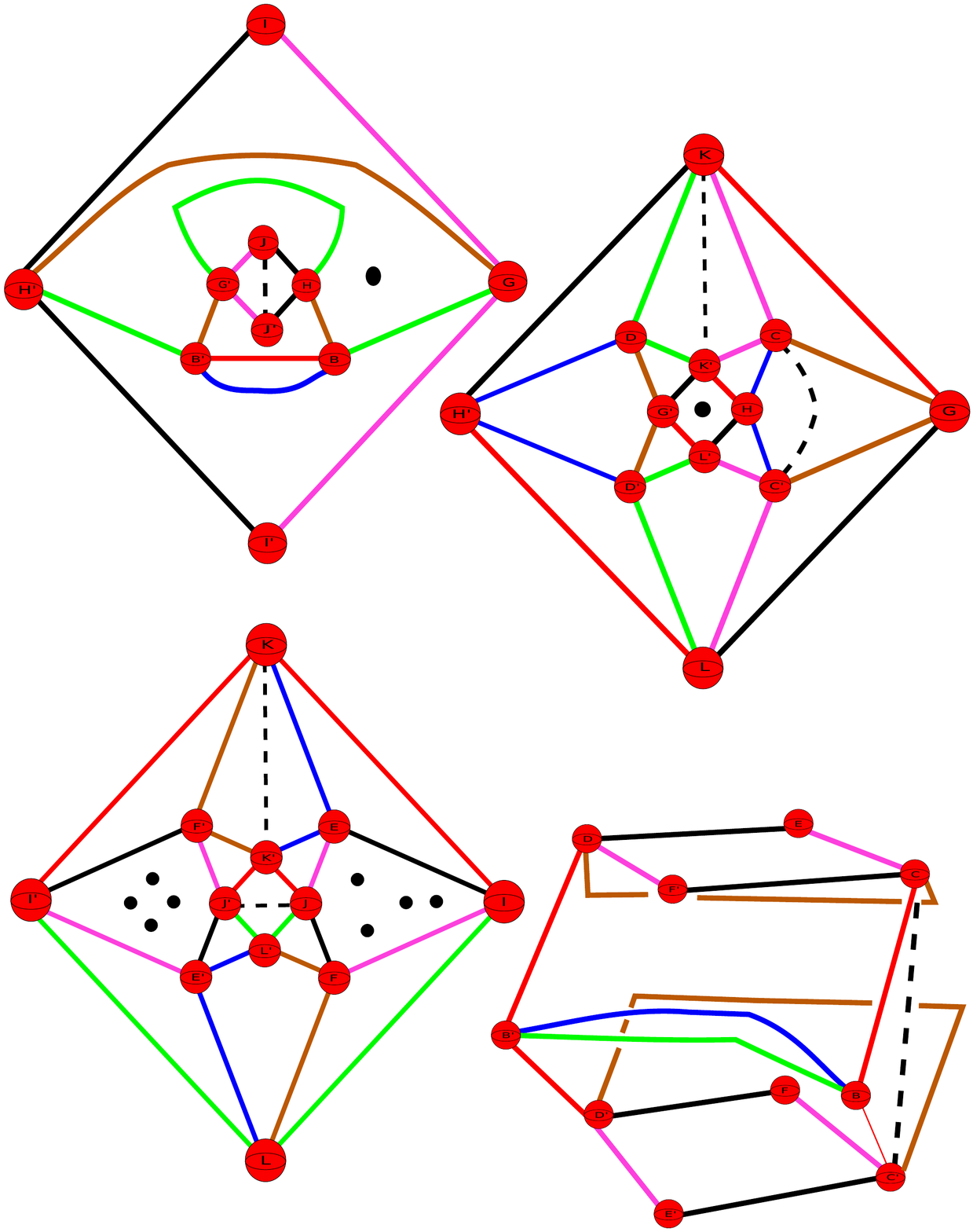}}

and the picture below shows how the added 2-handles running between $E$ and $G$, and $E'$ and $G'$ sit in our Kirby diagram so far.

\centerline{\graphicspath{ {handle_cancellation/cancelling_A-A'/} }\includegraphics[width=09cm, height=10cm]{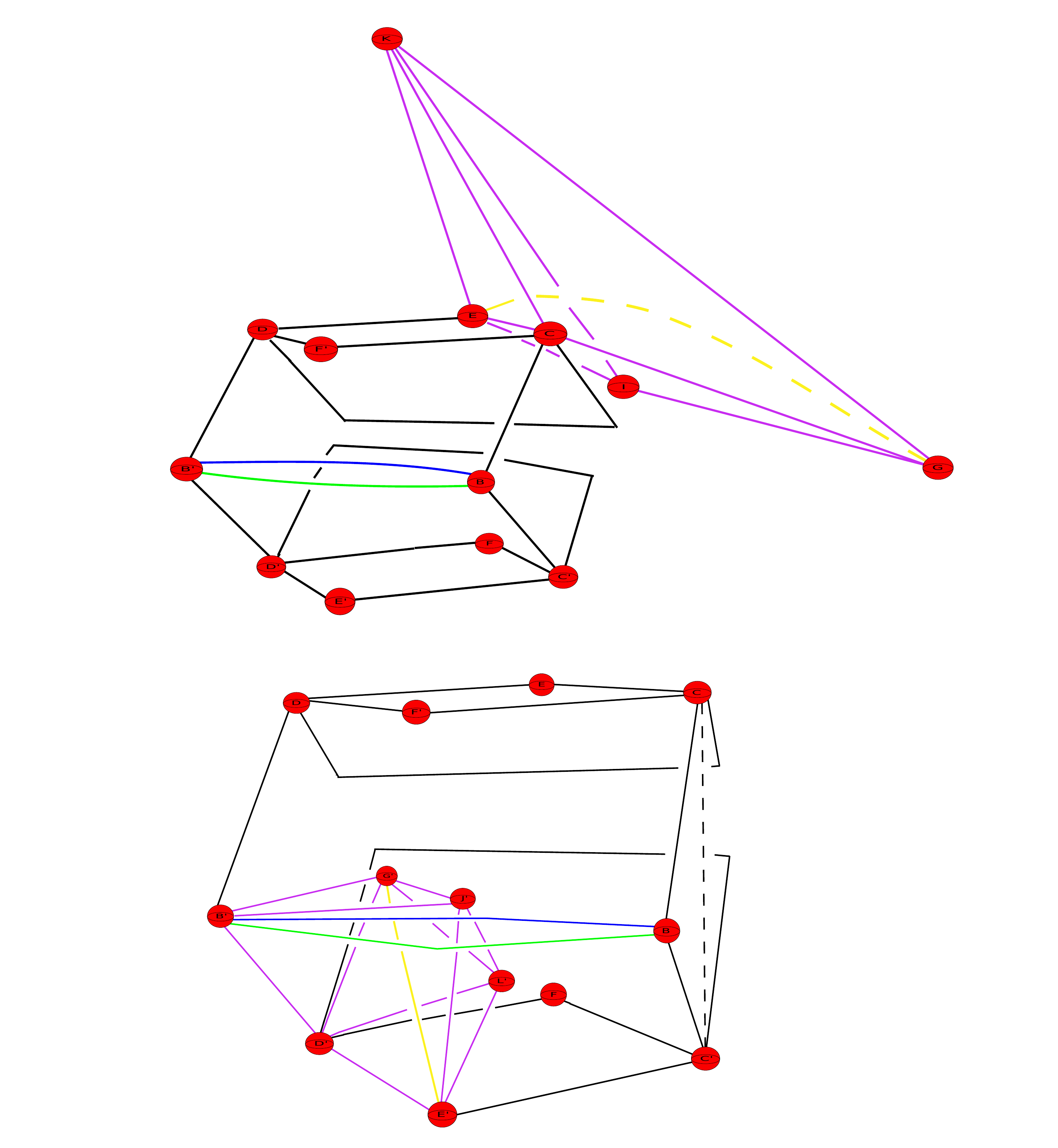}}

We now cancel $K-K'$, in this case only the 2-handles in the $x-z$ and $y-z$ planes change. 

The following picture shows how the Kirby diagrams change after we have carried out this cancellation.

\centerline{\graphicspath{ {handle_cancellation/cancelling_K-K'/} }\includegraphics[width=8cm, height=8cm]{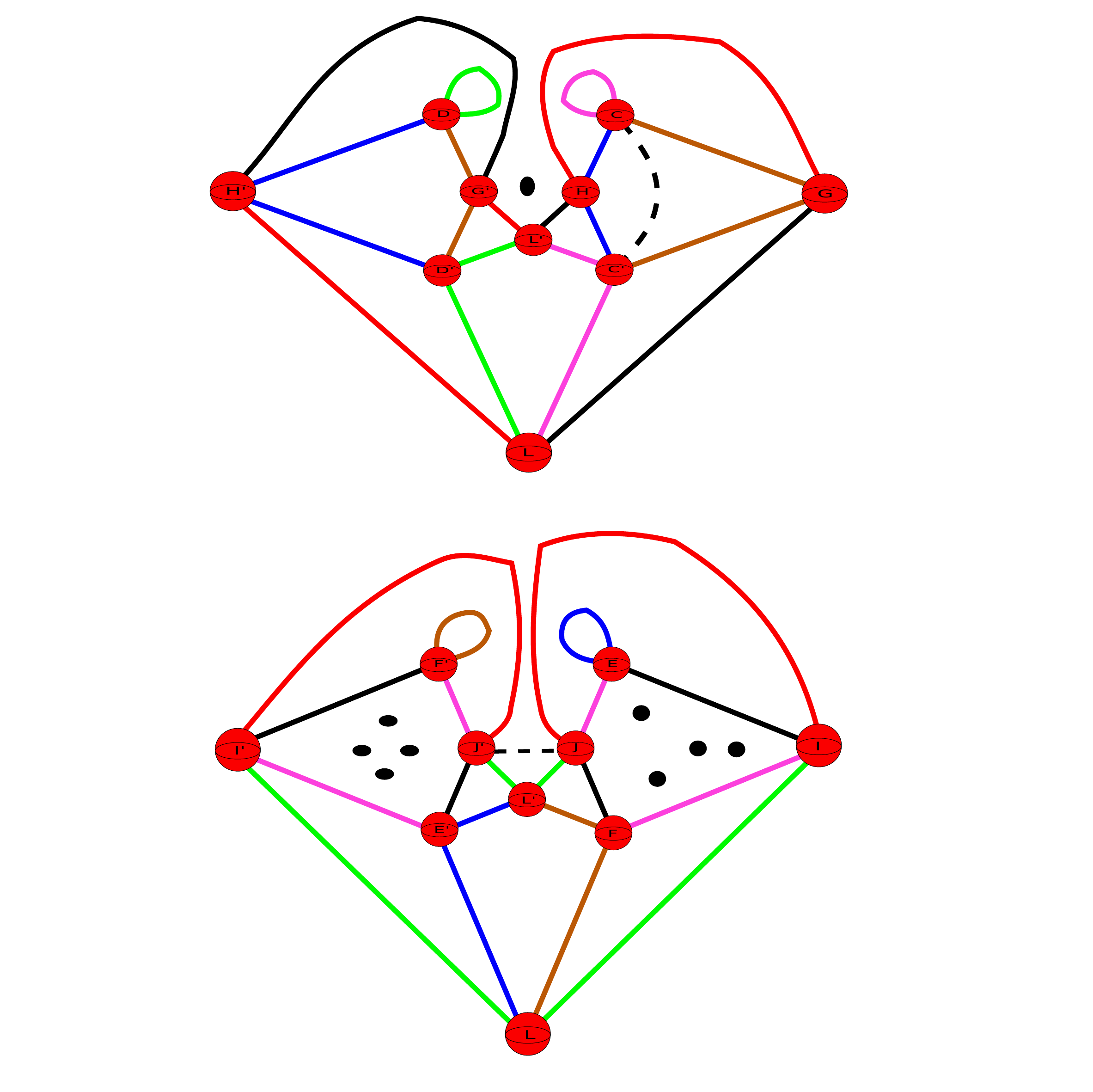}}

At this point one can ask if cancelling $K-K'$ has introduced any new intersection points in any of the other diagrams ? The only candidate is the
$x-y$ plane. In this case it is easy to see that no intersection points are created. Let us give a brief explanation of why this is the case.
Start with the $x-z$ plane, this meets the $x-y$ plane along a straight line that passes through the one handles $H'$, $G'$, $H$ and $G$.
In the following picture you can see this line as the horizontal yellow line through $H'$, $G'$, $H$ and $G$, where the top diagram is that of
the $x-y$ plane, and the bottom diagram is that of the $x-z$.

\centerline{\graphicspath{ {handle_cancellation/cancelling_K-K'/} }\includegraphics[width=8cm, height=9cm]{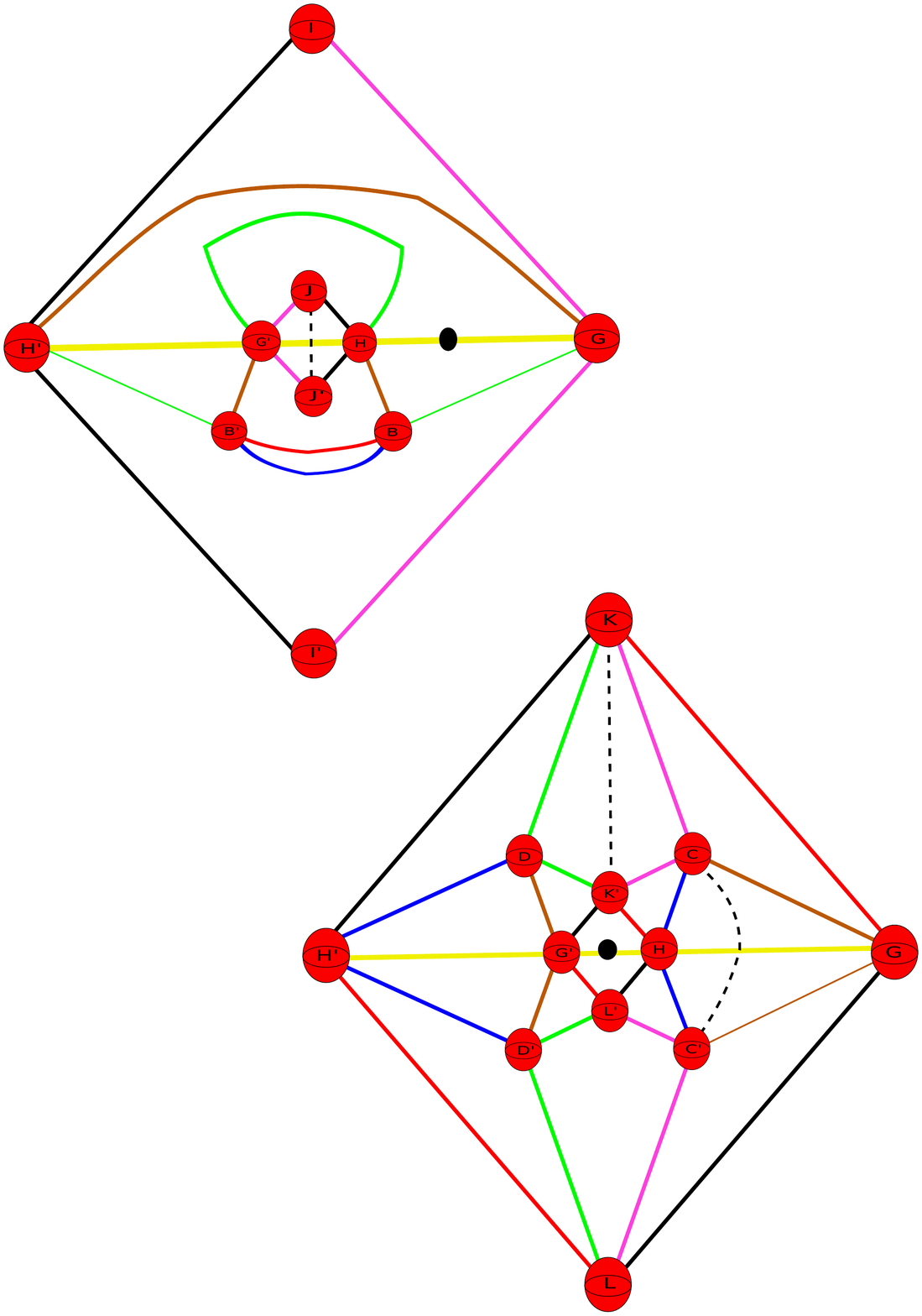}}

When we cancel $K-K'$ from the $x-z$ plane it is easy to see that none of the cancellations interfere with the yellow line, meaning that in cancelling
$K-K'$ from the $x-z$ plane none of the 2-handles that emerge from this cancellation cross the yellow line. This means none of these 2-handles
contribute to any intersection points in the $x-y$ plane.

In the case of the $y-z$ plane we have that it intersects the $x-y$ plane along a vertical line through $I'$, $J'$, $J$ and $I$, as shown in the
following picture.

\centerline{\graphicspath{ {handle_cancellation/cancelling_K-K'/} }\includegraphics[width=8cm, height=8cm]{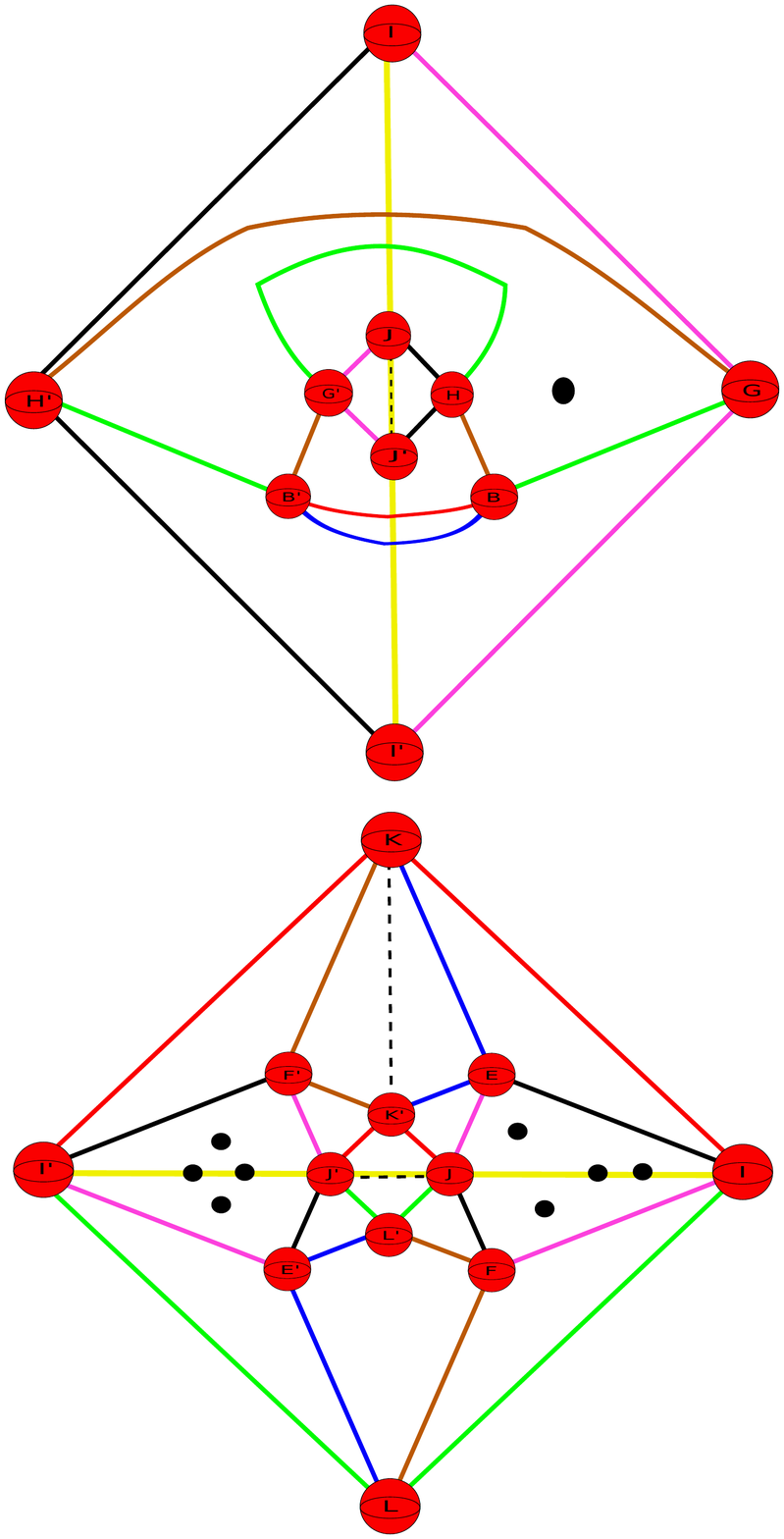}}

In this case you can also clearly see that when we cancel $K-K'$ from the $y-z$ plane none of the new 2-handles that emerge from this cancellation
cross the horizontal yellow line in the $y-z$ plane. This means none of these 2-handles intersect the $x-y$ plane, and hence we do not get any points of
intersection from cancelling $K-K'$ from the $y-z$ plane.

Coming back to the pictures of the $x-z$ and $y-z$ planes, after cancelling $K-K'$ we see that we can do some handle slides. First of all, in the
$x-z$ plane we have a green 2-handle that starts at $D$ and loops back into it, and a pink 2-handle that starts at $C$ and loops back into it. We can slide
the green 2-handle through $D$ to $D'$ and then off $D'$ to give a 2-handle meeting $L-L'$ once. We can do the same
with the pink 2-handle to get another 2-handle meeting $L-L'$ once.

In the case of the $y-z$ plane we have the red 2-handle that starts at $F'$ and loops back into it, and we have the blue 2-handle that starts at $E$ and
loops back into it. We can then slide these (just as we did above) to give two 2-handles that meet $L-L'$ once.

\centerline{\graphicspath{ {handle_cancellation/cancelling_K-K'/} }\includegraphics[width=9cm, height=8cm]{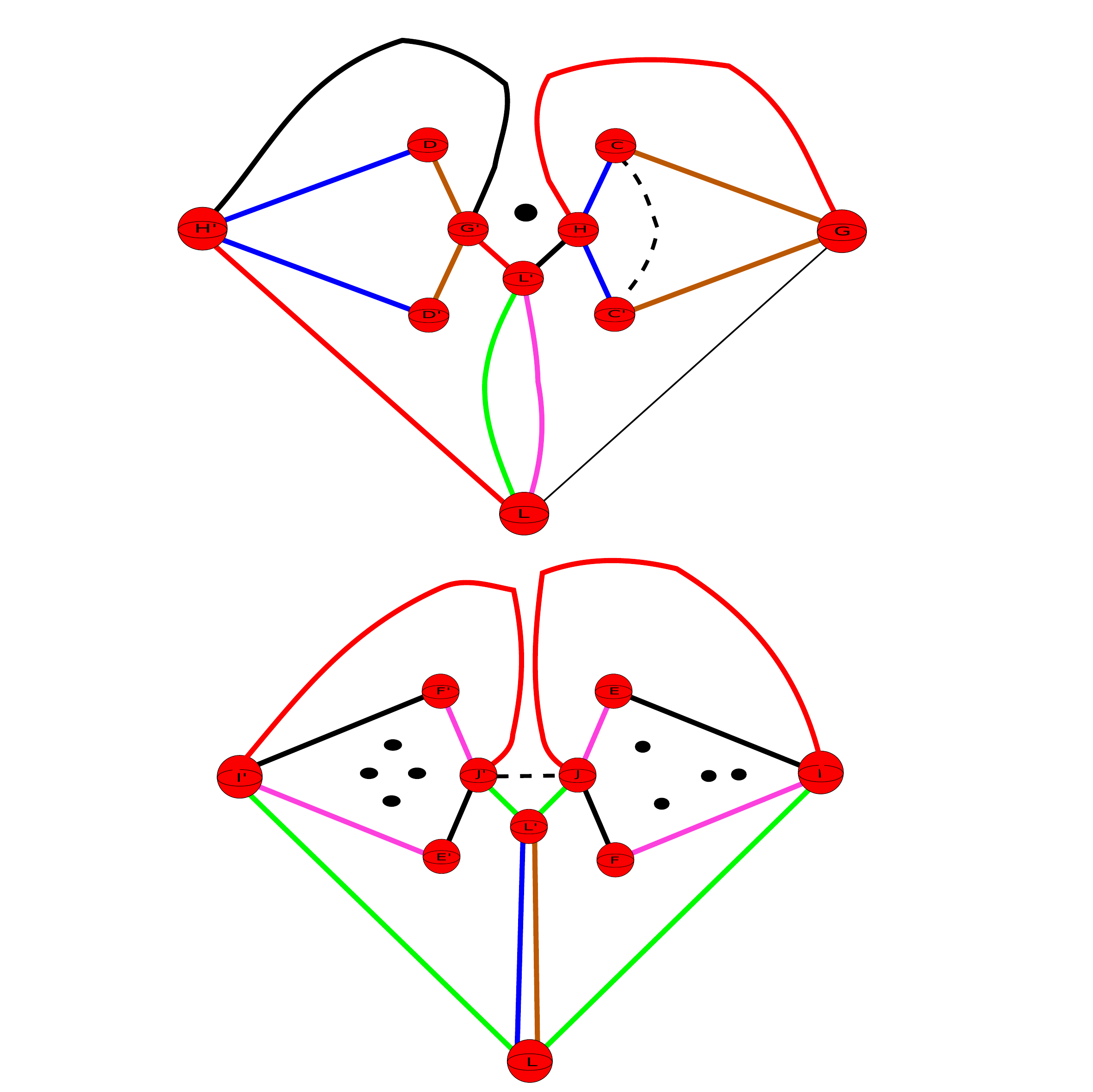}}

Observe that in carrying out these handle slides we have not added any new points of intersection to the $x-y$ plane. For example the two handle slides
we did in the $x-z$ plane gave us two 2-handles that met $L-L'$ once. These two 2-handles do not cross the intersection line of the $x-y$ plane with
the $x-z$ plane, hence these new 2-handles cannot intersect the $x-y$ plane. A similar analysis for the handle slides in the $y-z$ plane shows that
they do not add any points of intersection as well.

It is clear that in cancelling the 1-handle $K-K'$ we have introduced no new intersection points in the other diagrams. For example when we cancelled
$K-K'$ in the $x-z$ plane, we introduced a 2-handle that loops back in to $C$, it is the pink 2-handle in the first cancellation diagram. This 2-handle does not introduce any
intersection points in the $x-y$ plane because $C$ lies above the $x-y$ plane. Similarly we also introduced a 2-handle that goes from $H$ to $G$, it is 
drawn as the red 2-handle in the picture above. This 2-handle also does not create any points of intersection with the $x-y$ plane. A similar analysis shows that
when we first cancel $K-K'$ none of the new 2-handles formed create any new intersection points.  

We then proceeded to doing some handle slides, and we can ask whether these handle slides create any new intersection points. It turns out that they do not, and one
can easily see this by simply thinking of what happens in a 3-dimensional picture. Let us give an explicit analysis of this, we will work with the $y-z$ plane and
the 2-handle that loops back in to $E$, it is drawn in blue in the first picture showing what happens when we cancel $K-K'$ (see the picture before the three above pictures).
We push this 2-handle through $E-E'$ to get a 2-handle meeting $L-L'$, note that $E'$, $L'$ and $L$ all lie below the $x-y$ plane. Hence when we push
this 2-handle through $E-E'$ to give a 2-handle between $L-L'$ we create no intersection points with the $x-y$ plane. Similarly, this new 2-handle does not create
any intersection points with the $x-z$ plane.

So far we have filled in two boundary components of the manifold $M$ corresponding to the fibres given by the isometries $a$ and $k$.
We then carried out various handle cancellations and handle slides. The following picture puts all that we have done so far together.

\centerline{\graphicspath{ {handle_cancellation/cancelling_K-K'/} }\includegraphics[width=9cm, height=10cm]{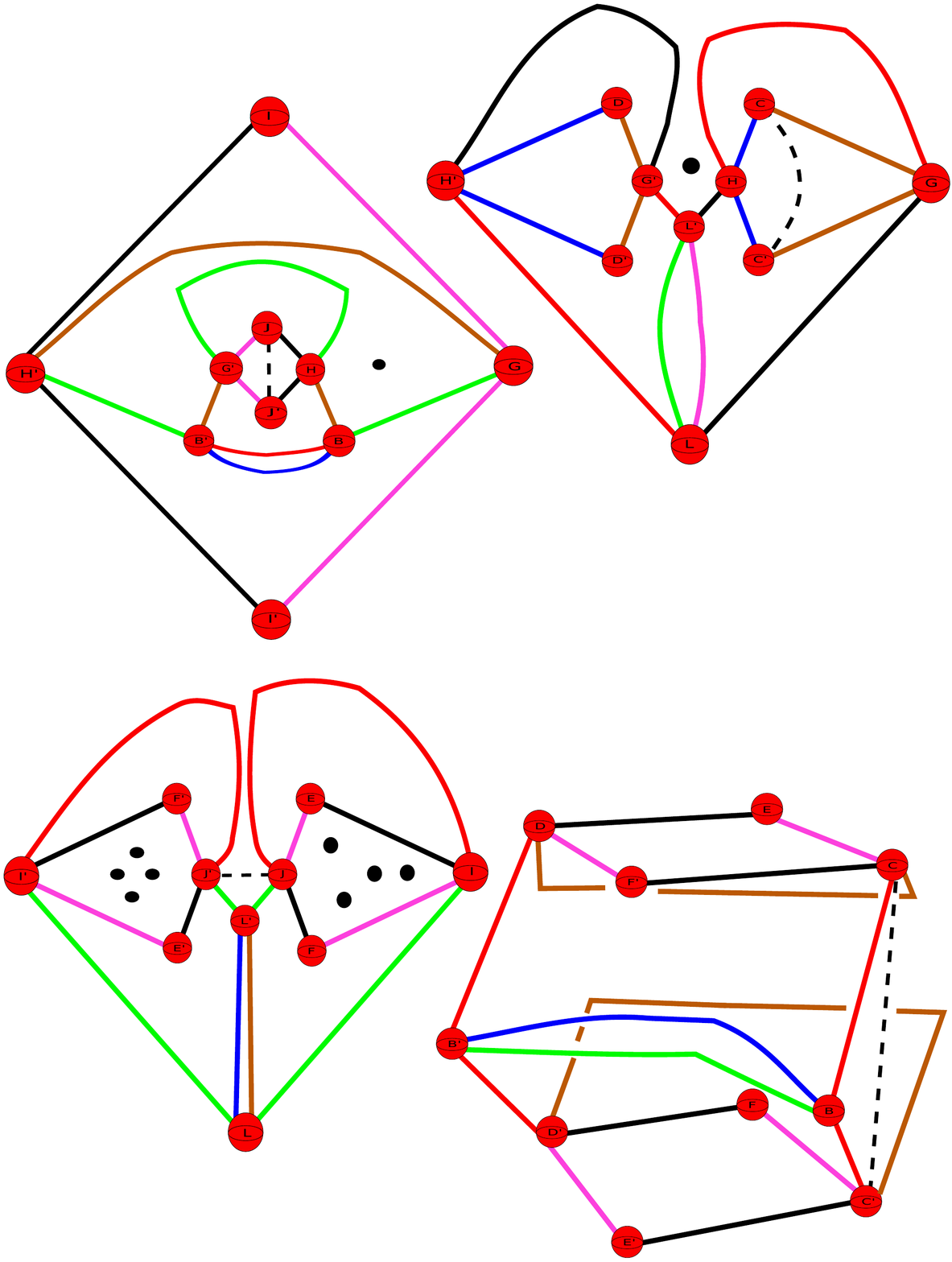}}

The labelling of the intersection points are shown in the following diagram. The top two diagrams show close ups of the regions around the intersection
points in the $x-y$ and $x-z$ planes respectively (viewing from left to right). The bottom two diagrams show close ups of the intersection points in the
$y-z$ plane.

\centerline{\graphicspath{ {handle_cancellation/cancelling_K-K'/} }\includegraphics[width=9cm, height=8cm]{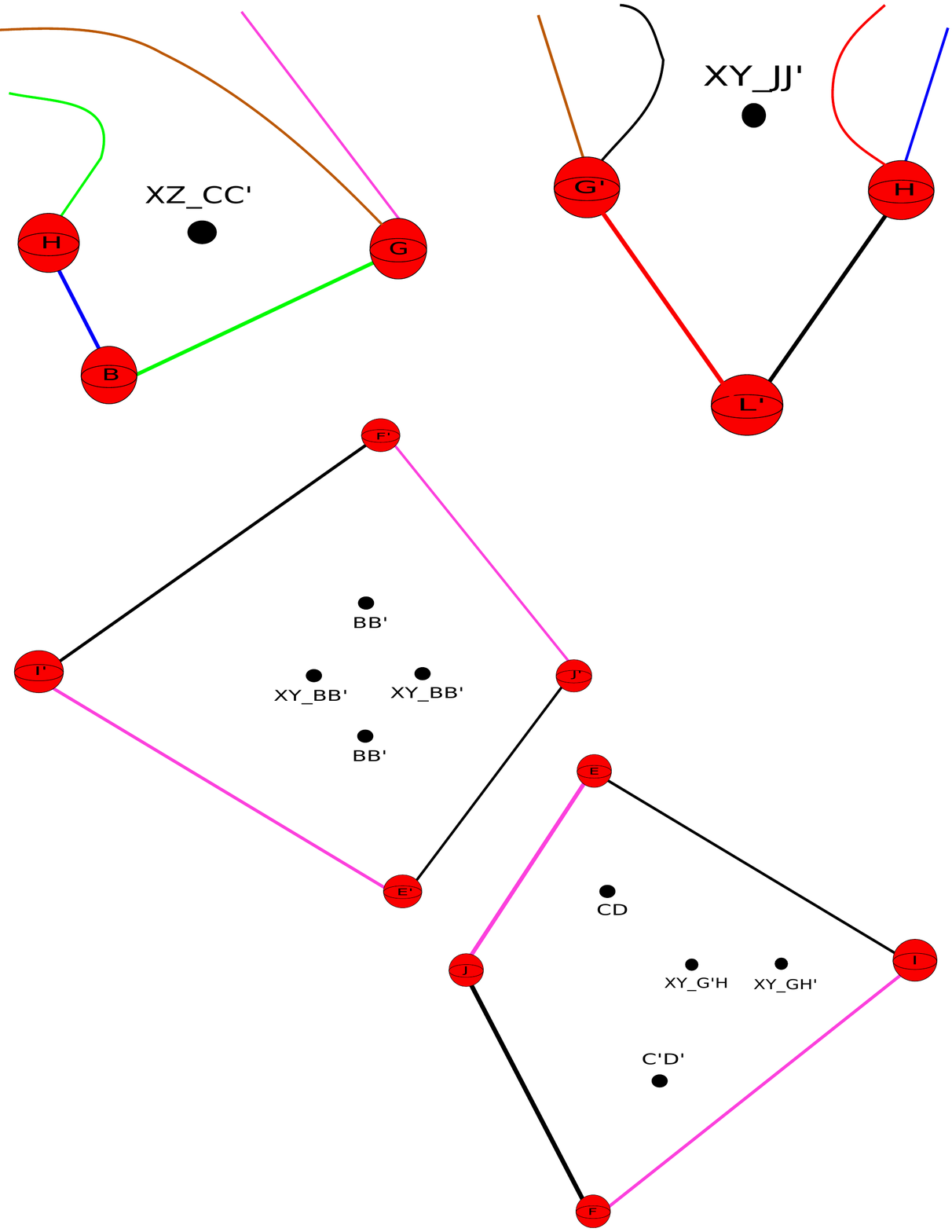}}

The next step is to cancel $J-J'$, this cancellation will only affect the $x-y$ and $y-z$ planes.
The diagrams on the right correspond to sliding the obvious handles in the left picture over 1-handles and then off them to obtain 2-handles
that only meet certain 1-handles once. For example in cancelling $J-J'$ from the $x-y$ plane we obtain a pink 2-handle that starts at $G'$ and loops
back into it (see the top left diagram in the picture below), we can then slide this through $G'$ and off $G$ to obtain a pink 2-handle that meets $I-I'$ once (see the top
right diagram in the picture below).

\centerline{\graphicspath{ {handle_cancellation/cancelling_J-J'/} }\includegraphics[width=11cm, height=10cm]{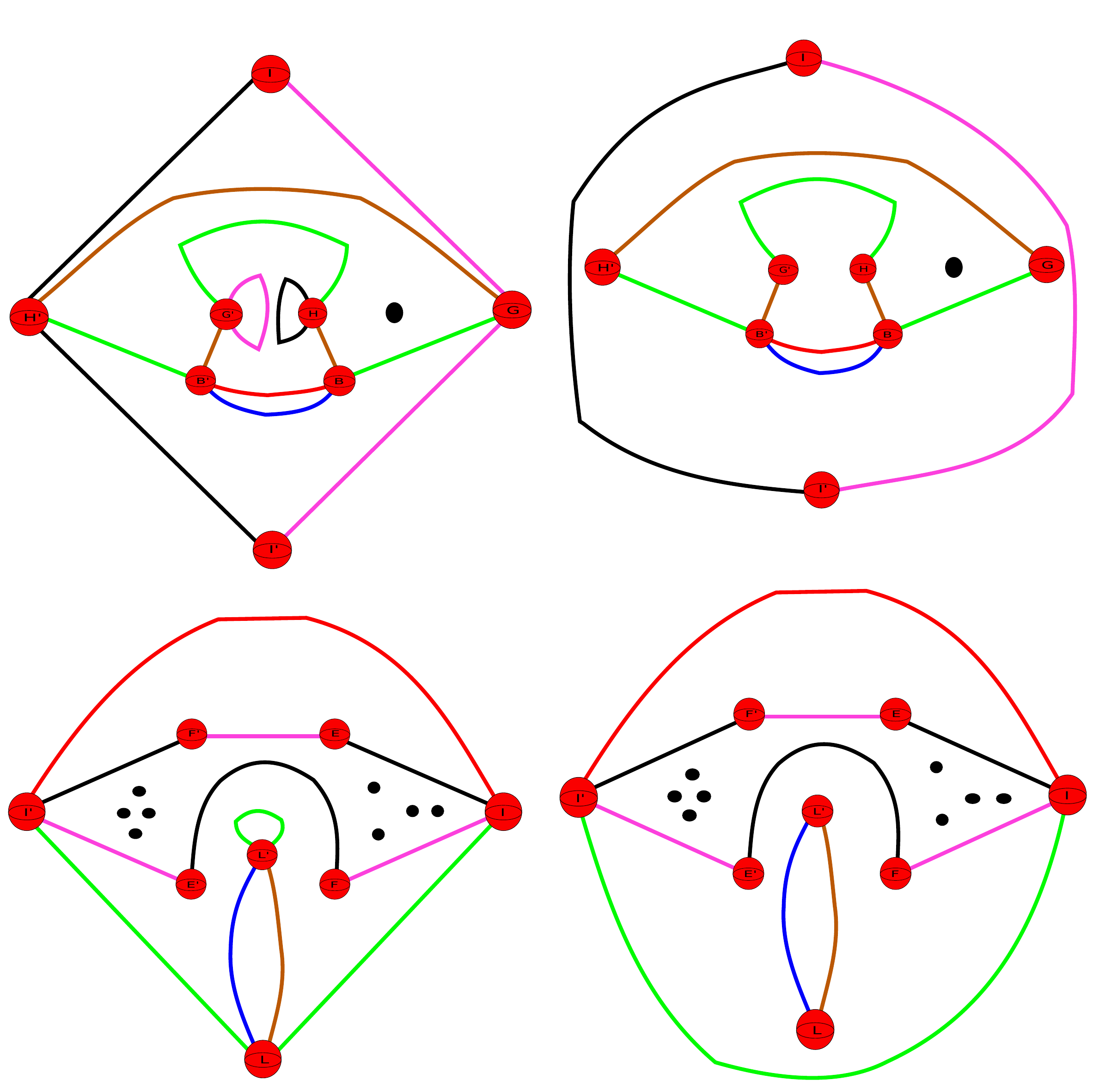}}

When we first cancel $J-J'$ in the $x-y$ plane we introduce some new 2-handles, namely the one in pink that loops back into $G'$ and the one in black
that does the same with $H$. These two 2-handles create two new points of intersection with the $x-z$ plane:

\centerline{\graphicspath{ {handle_cancellation/cancelling_J-J'/} }\includegraphics[width=7cm, height=5cm]{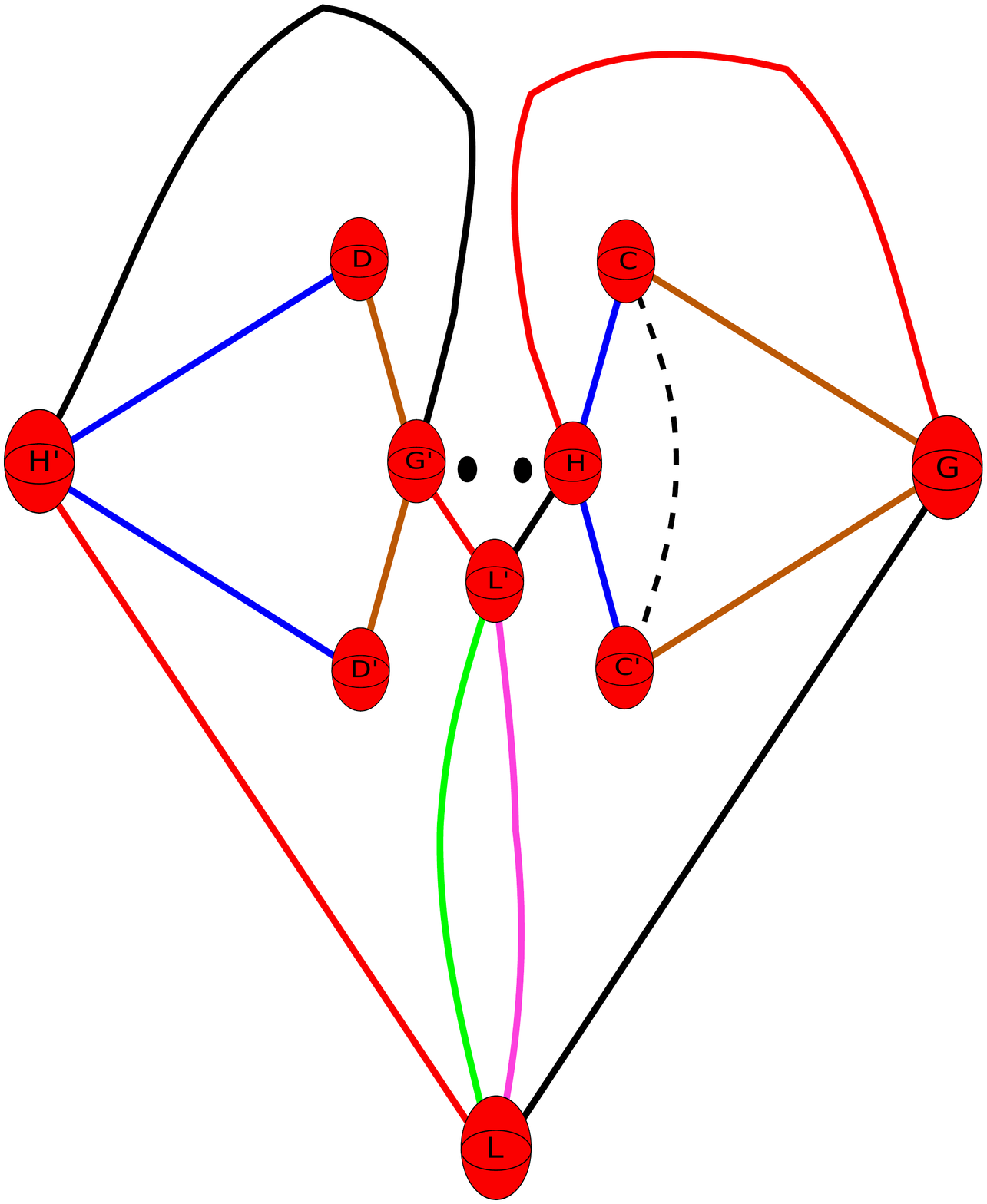}}

When we cancel $J-J'$ in the $y-z$ plane we introduce a 2-handle between $I-I'$, $E-F'$, $E'-F$, and one that loops back into $L'$. These four 2-handles
each create intersection points with $x-z$ plane as well, you can see them as the four vertical black dots:

\centerline{\graphicspath{ {handle_cancellation/cancelling_J-J'/} }\includegraphics[width=7cm, height=5cm]{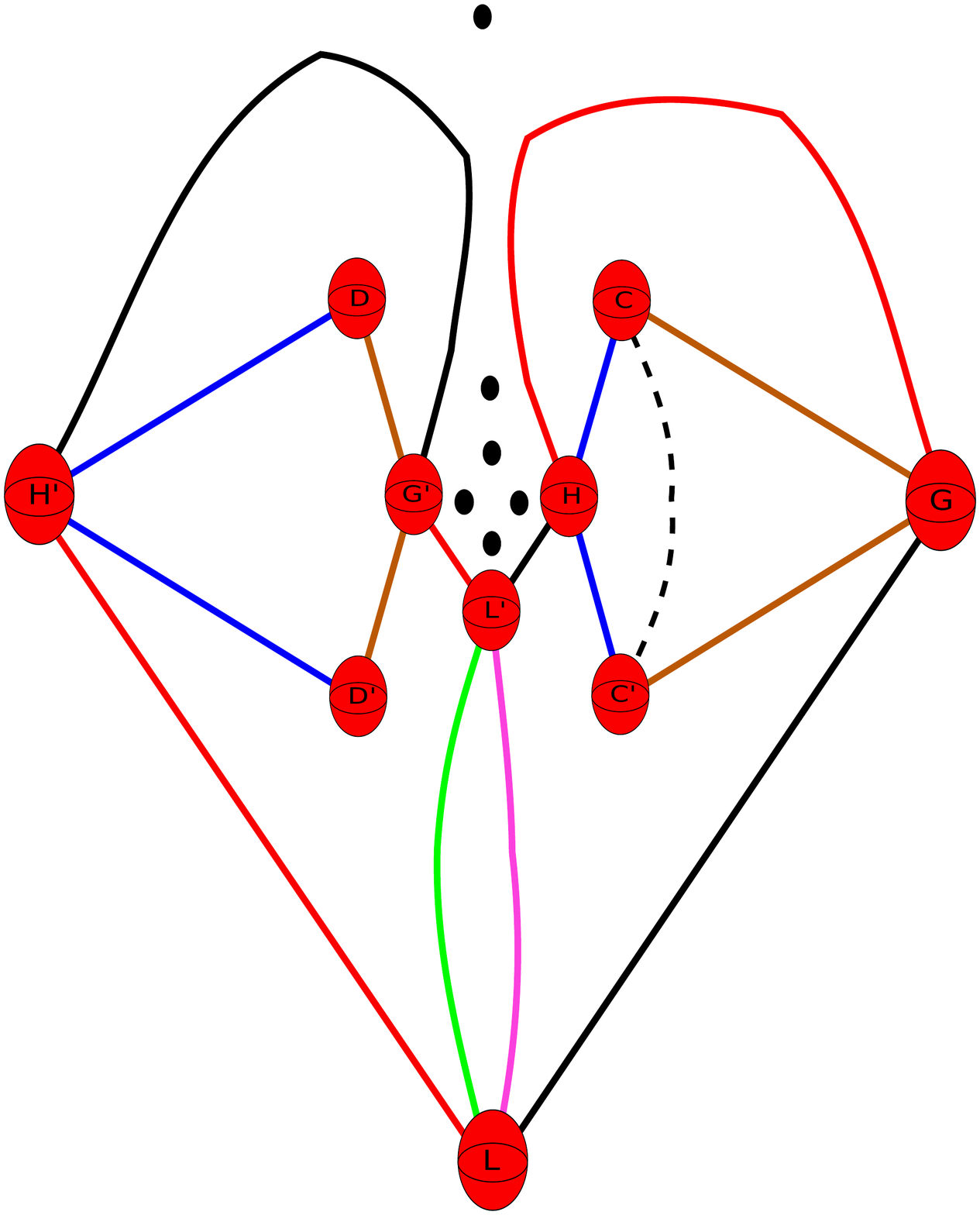}}

The labelling of these 6 new intersection points is given in the following picture.

\centerline{\graphicspath{ {handle_cancellation/cancelling_J-J'/} }\includegraphics[width=5.5cm, height=5.5cm]{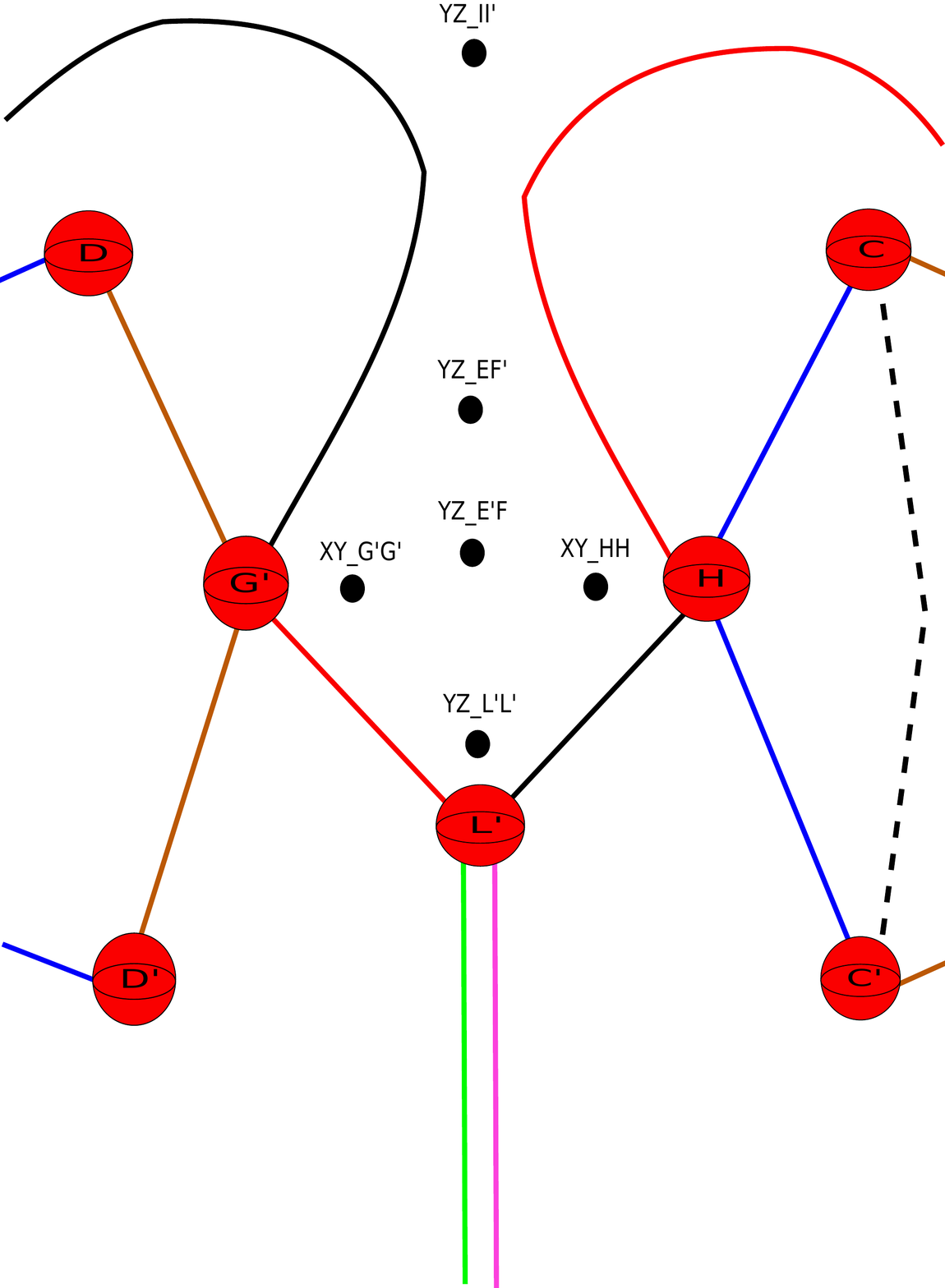}}

After cancelling $J-J'$ we carried out some some handles slides. In the $x-y$ plane we slid the new 2-handle that starts at $G'$ and loops back into it (the one in pink)
through $G'$ and then off $G$ to give a 2-handle that passes over $I-I'$ once. This will cause the intersection point in the $x-z$ plane, labelled 
\textbf{XY\underline{\space}G'G'}, to disappear. However, a new intersection point will arise corresponding to the new 2-handle that passes over $I-I'$ once.
This new intersection point is shown in the following diagram as the black dot to the right of $G$, you can also see that the intersection point labelled
\textbf{XY\underline{\space}G'G'} is no longer in the picture.

\centerline{\graphicspath{ {handle_cancellation/cancelling_J-J'/} }\includegraphics[width=5.5cm, height=5.5cm]{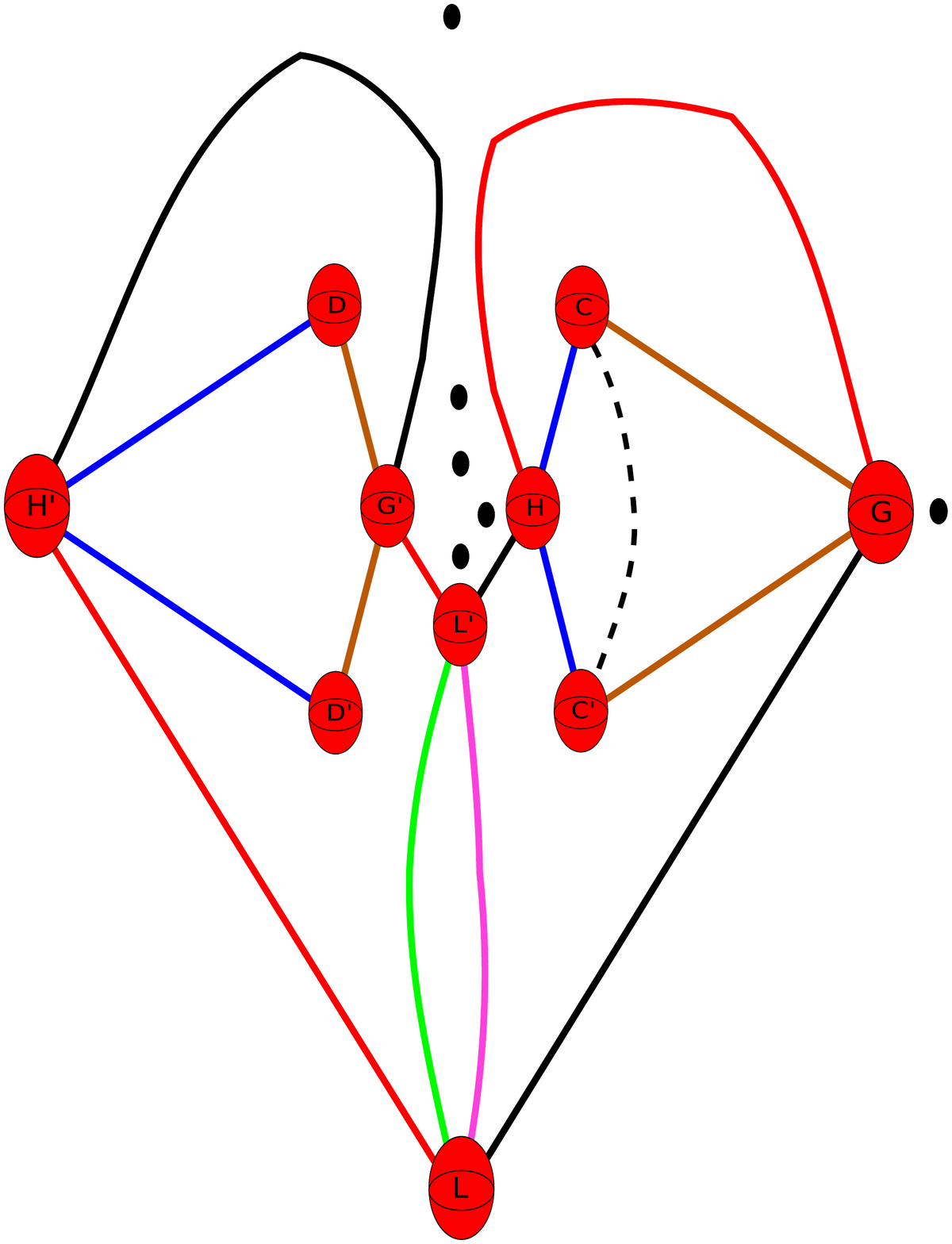}}

The other handle slide we did in the $x-y$ plane was to take the black 2-handle that starts at $H$ and loops back into it and slide it through $H$ and then off
$H'$ to give a new 2-handle that passes over $I-I'$ once. This handle slide causes the intersection point labelled \textbf{XY\underline{\space}HH} to disappear with
a new intersection point appearing to the left of $H'$ (viewed in the $x-z$ plane).

\centerline{\graphicspath{ {handle_cancellation/cancelling_J-J'/} }\includegraphics[width=6cm, height=4cm]{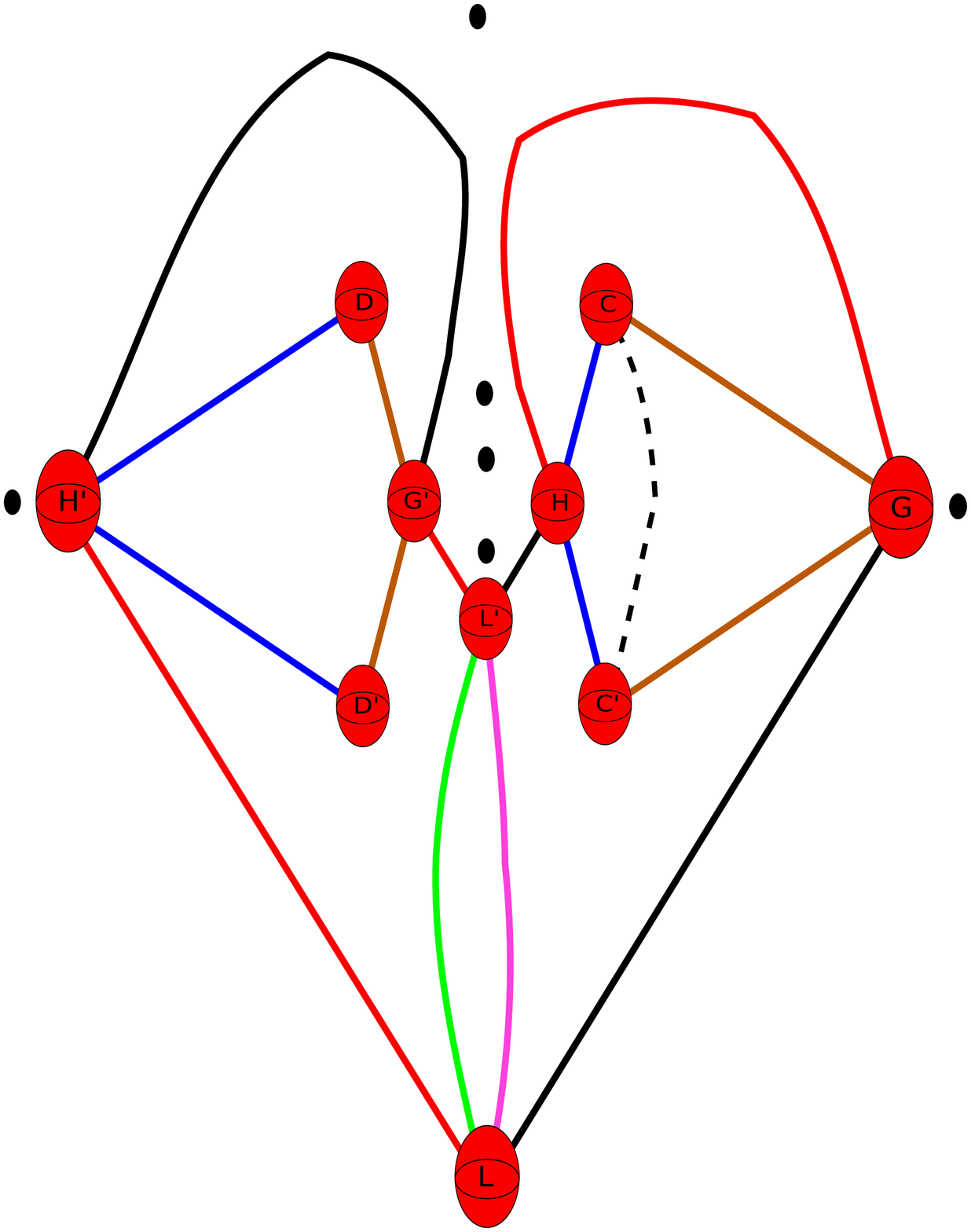}}

We also carried out one handle slide in the $y-z$ plane. This handle slide corresponded to taking the green 2-handle that starts at $L'$ and loops back into it, pushing
it through $L'$ so it comes out at $L$, and then sliding it off to give a new 2-handle that passes over $I-I'$ once. This causes the intersection point
labelled \textbf{YZ\underline{\space}L'L'} to disappear, with a new intersection point just below $L$ to appear. The following picture shows this new intersection
point and all the others with their corresponding label.

\centerline{\graphicspath{ {handle_cancellation/cancelling_J-J'/} }\includegraphics[width=6cm, height=4.5cm]{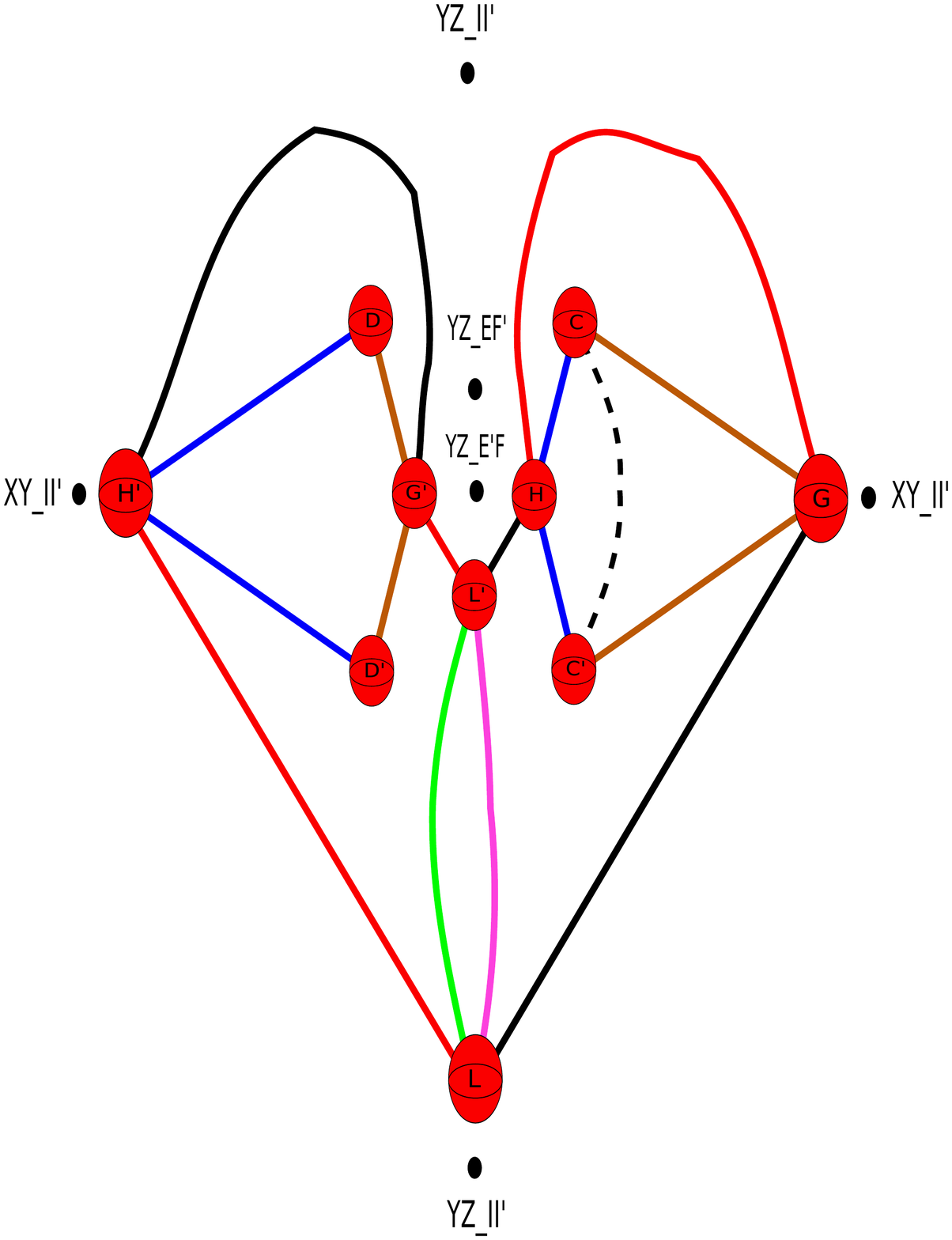}}

If we put everything we have done so far we get the following:

\centerline{\graphicspath{ {handle_cancellation/cancelling_J-J'/} }\includegraphics[width=8cm, height=6.5cm]{all_together_int}}

The pictures showing how the added 2-handle running between $E$ and $G$, and $E'$ and $G'$ is the same as before.

We move on to cancelling $C-C'$, this will only affect the $x-z$ plane and the six 2-handles not lying in any one plane.

\centerline{\graphicspath{ {handle_cancellation/cancelling_C-C'/} }\includegraphics[width=5cm, height=6cm]{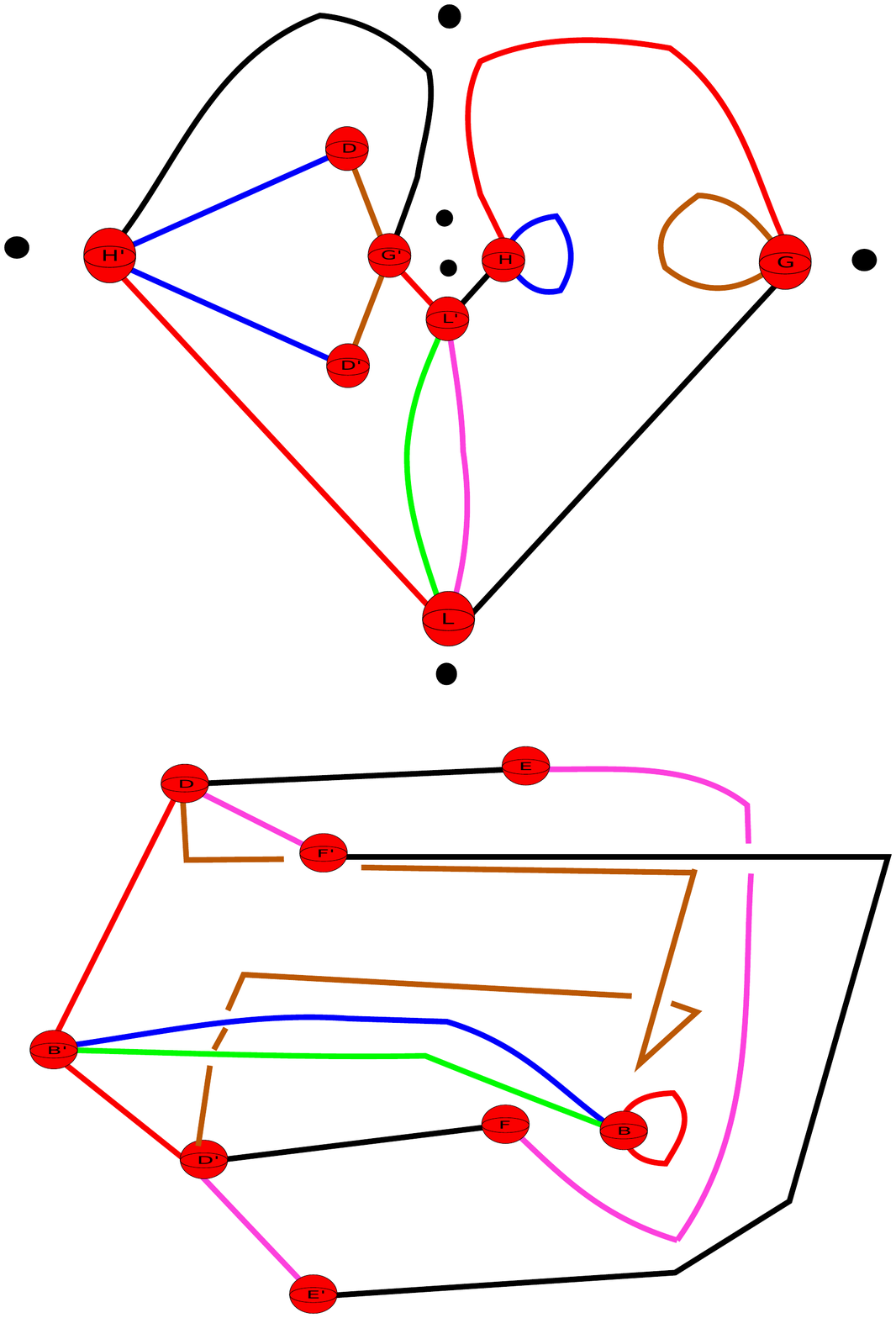}}

Recall that the added 2-handle that passed between $C-C'$ (the one we just used above to cancel $C-C'$) intersected the $x-y$ plane. Therefore
the above handle cancellation will give rise to some points of intersection in the $x-y$ plane. We take the time to explain how these points of intersection
look like.

Start with the $x-z$ plane, when we cancel $C-C'$ we got two new 2-handles one in blue which starts at $H$ and loops back into it, and one in brown that starts
at $G$ and loops back into it. These will give two new intersection points in the $x-y$ plane, as shown below:

\centerline{\graphicspath{ {handle_cancellation/cancelling_C-C'/} }\includegraphics[width=6cm, height=6cm]{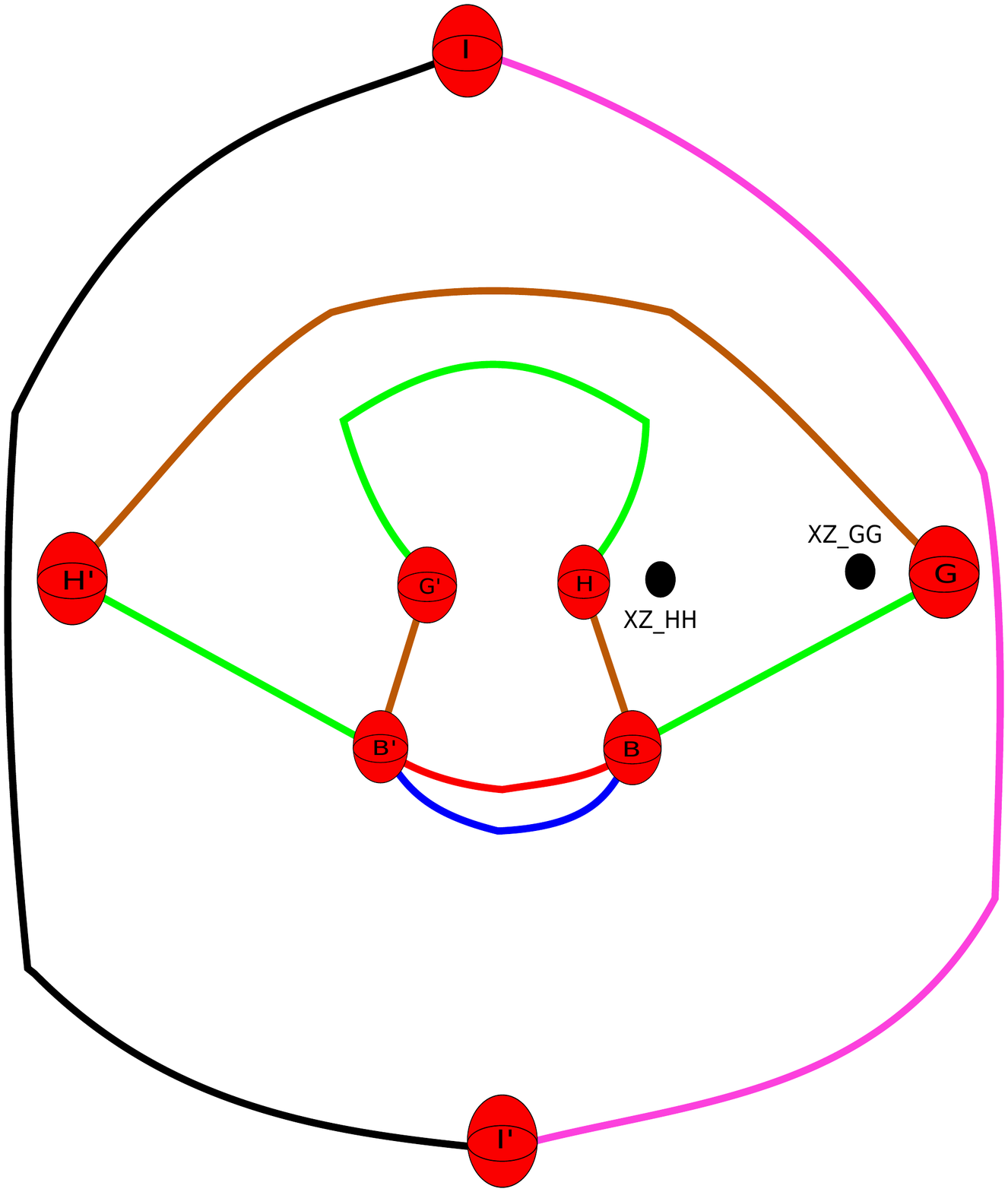}}

Moving to the case of the six 2-handles that did not all lie in a single 2-plane, we see that we get four new points of intersection. First of all,
we have the red 2-handle that starts and ends at $B$. Secondly, we have the black 2-handle that runs from $E'$ to $F'$ and the pink 2-handle that runs from
$F$ to $E$. Finally, we have the brown 2-handle that runs from $D$ to $D'$. All four of these 2-handles give four points of intersection in the $x-y$ plane.
The following picture show these points of intersections in the $x-y$ plane, with the bottom diagram being a close up showing the labelling of these intersection
points.

\centerline{\graphicspath{ {handle_cancellation/cancelling_C-C'/} }\includegraphics[width=7cm, height=8cm]{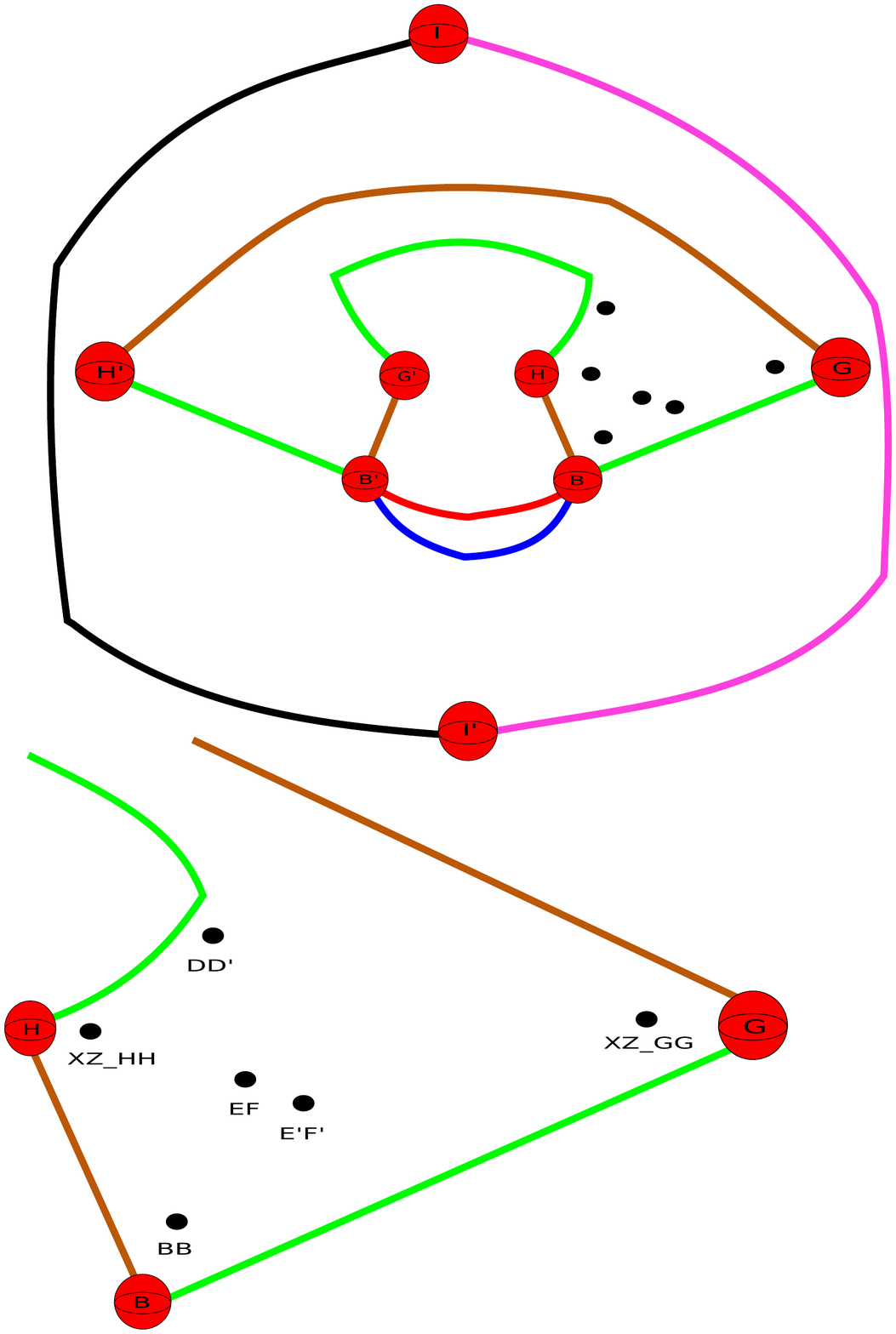}}

We also have to deal with the points of intersection in the $y-z$ plane. We have two points of intersection labelled \textbf{C'D'} and \textbf{CD}, when
we cancel $C-C'$ from the diagram that consisted of the six 2-handles that did not all lie in a single 2-plane the part of the brown 2-handle component
running from $C'$ to $D'$ and the part running from $C$ to $D$ join together to give a brown 2-handle running from $D$ to $D'$. Thus in the $y-z$ plane
we will still see two points of intersection but their labelling will be \textbf{DD'} because this new brown 2-handle running from $D$ to $D'$ intersects
the $y-z$ plane in two distinct points.

\centerline{\graphicspath{ {handle_cancellation/cancelling_C-C'/} }\includegraphics[width=6cm, height=7cm]{y-z_int}}

We can then carry out three handle slides. In the $x-z$ plane we can slide the blue 2-handle that starts and ends at $H$ through $H$ and then off $H'$ to give a blue 2-handle 
between $D-D'$. Similarly, we can slide the red 2-handle that starts and ends at $G$ off $G'$ to give a red 2-handle that meets $D-D'$ once. Finally, in the
diagram corresponding to the six 2-handles that did not all lie in a single 2-plane, we have the red 2-handle that starts and ends at $B$. We can slide it off $B'$
to give a red 2-handle meeting $D-D'$ once.

\centerline{\graphicspath{ {handle_cancellation/cancelling_C-C'/} }\includegraphics[width=9cm, height=9cm]{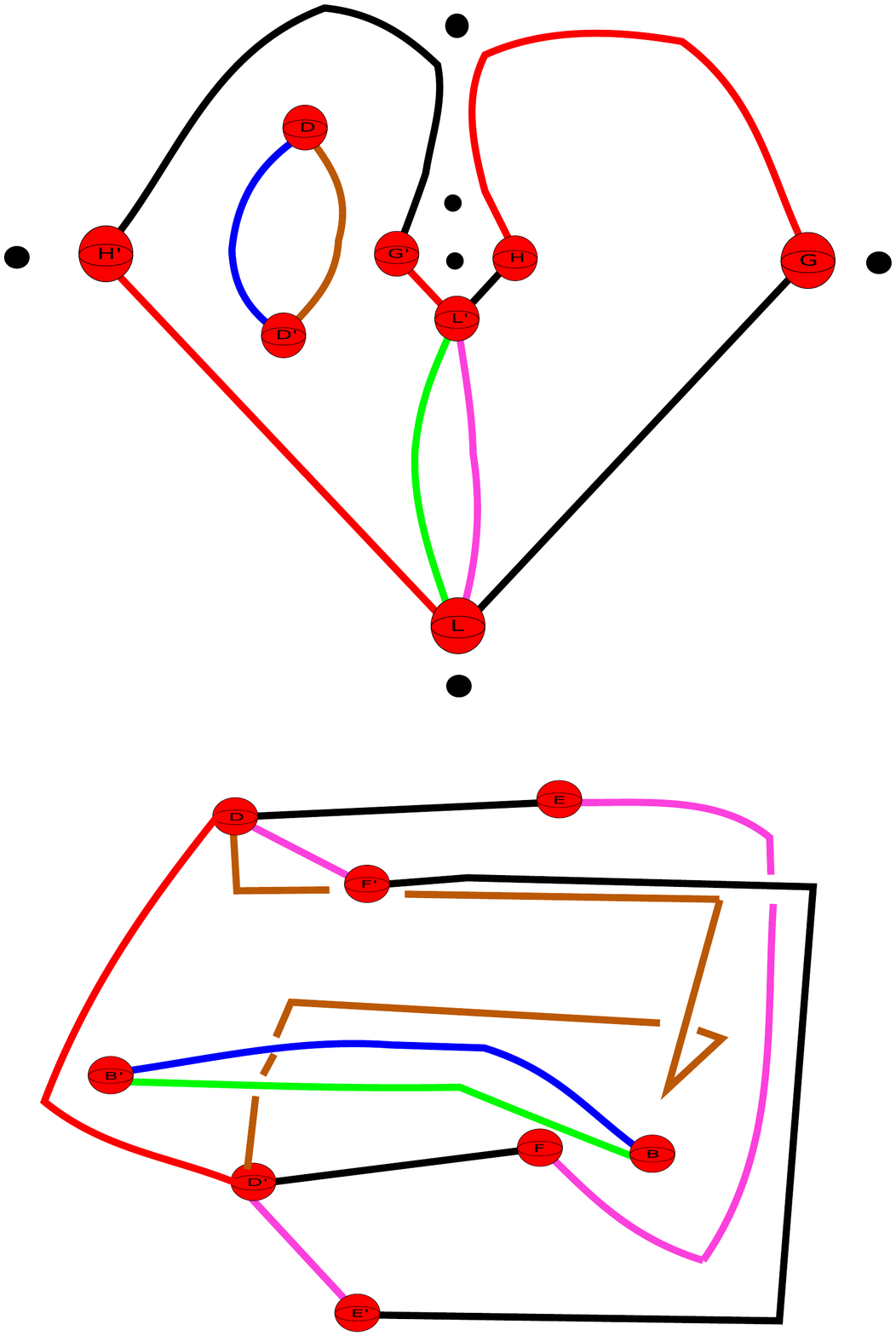}}

The first two handle slides we undertook will cause the points of intersection labelled \textbf{XZ\underline{\space}HH} and \textbf{XZ\underline{\space}GG}
to disappear from the $x-y$ plane, however we will get two new points of intersection coming from the new 2-handles that run between $D-D'$, these are shown
in the picture below with labelling \textbf{XZ\underline{\space}DD'}. The second handle slide we undertook will cause the point of intersection labelled
\textbf{BB} to disappear, but the new 2-handle that we obtained running from $D$ to $D'$ will give a new intersection point labelled  \textbf{DD'}.
These new points of intersection can be seen in the picture below.

\centerline{\graphicspath{ {handle_cancellation/cancelling_C-C'/} }\includegraphics[width=6cm, height=6cm]{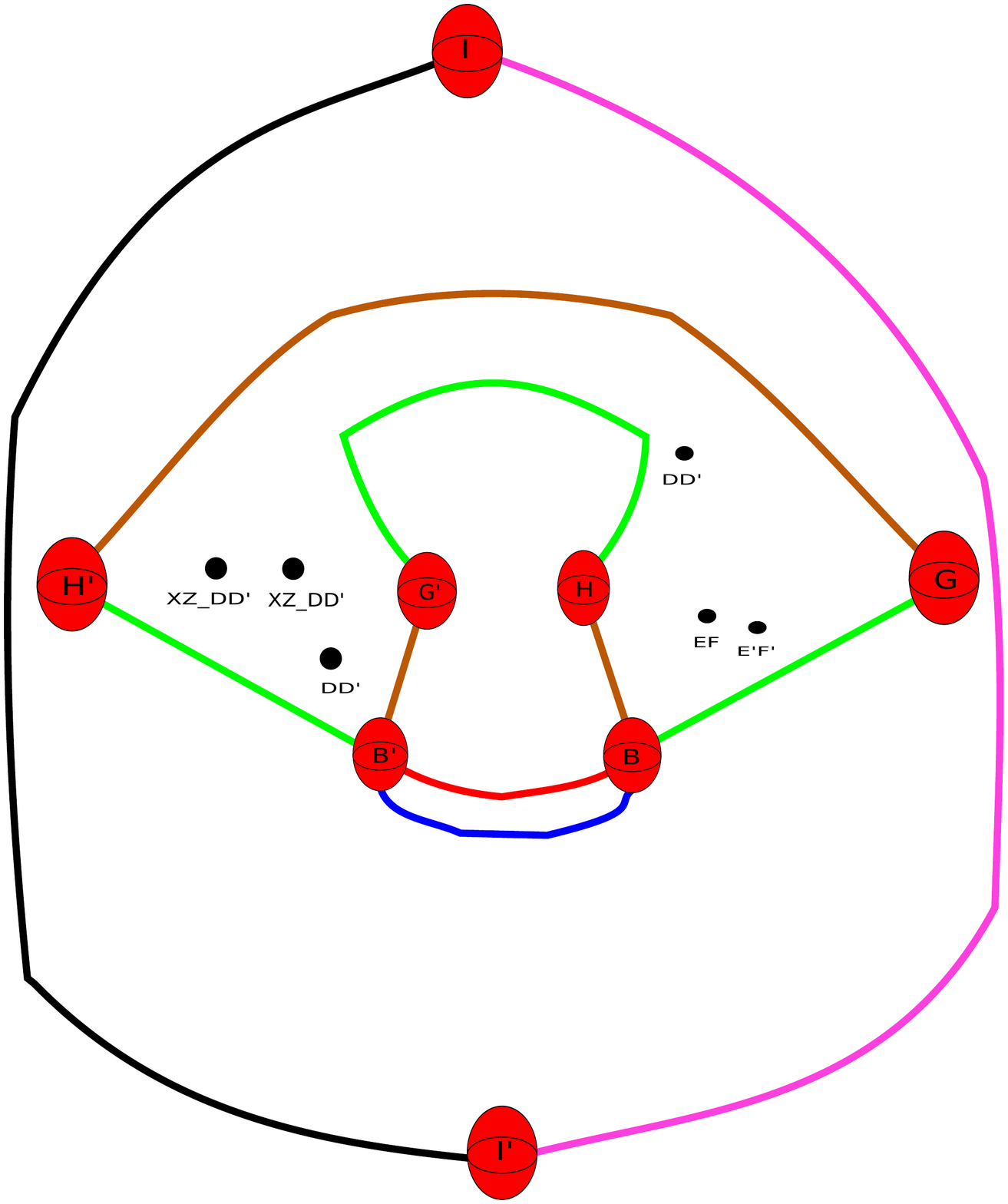}}

These handle slides do not affect the $y-z$ plane.

So far we have we have cancelled the 1-handles $A-A'$, $K-K'$, $J-J'$ and $C-C'$, the result of all these cancellations and various handle slides are shown
in the picture below along with points of intersections that arise when we cancel/slide handles.

\centerline{\graphicspath{ {handle_cancellation/cancelling_C-C'/} }\includegraphics[width=11cm, height=9cm]{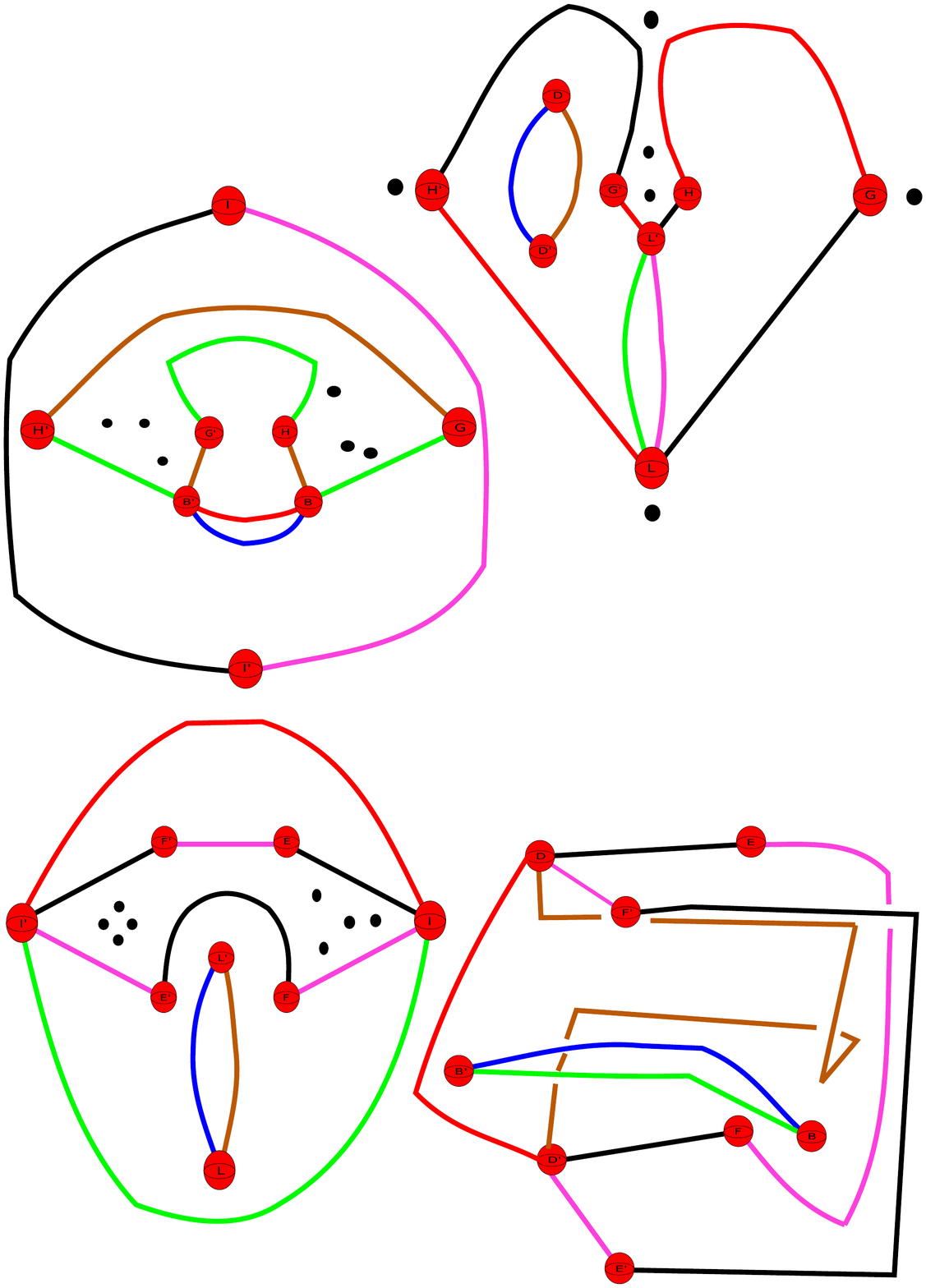}}

It is clear that none of the handle cancellations and slides we have done so far have interfered with the added 2-handles running from $E$ to $G$ and
$E'$ to $G'$. This will be the case for many of the handle cancellations and slides we do in the following, and because of this we will often omit
drawing the 2-handles between $E$ to $G$ and $E'$ to $G'$. However, it is recommended that the reader keep a mental image of these two 2-handles
so as to help convince themselves that none of the handle cancellations/slides we carry out do indeed affect these two 2-handles.

In the above picture you can see several 2-handles that pass over certain 1-handles only once. For example, if we look at the $x-y$ plane we can see that
there is a red 2-handle and a blue 2-handle that passes between the 1-handle $B-B'$ once, hence we may use either of these to cancel $B-B'$.

Cancelling $B-B'$ only affects the handle diagram in the $x-y$ plane and the diagram corresponding to the six $2$-handles not lying in any one plane.

\centerline{\graphicspath{ {handle_cancellation/cancelling_B-B'/} }\includegraphics[width=10cm, height=7.2cm]{x-y_x-y-z}}

The unknotted circles in the above picture denote 2-handles not meeting any 1-handles and have zero framing. Hence they can be used to cancel a 3-handle, and so we may
simply erase them from the diagram (recall elementary move number 3, which we outlined in the previous section).

Recall that the 2-handles running between $B$ and $B'$ gave points of intersection in the $y-z$ plane. To remind the reader the picture below shows all 
the points of intersection in the $y-z$ plane.

\centerline{\graphicspath{ {handle_cancellation/cancelling_B-B'/} }\includegraphics[width=8cm, height=6cm]{y-z_int}}

When we cancel $B-B'$ from the $x-y$ plane using the red 2-handle, the intersection point corresponding to this 2-handle in the $y-z$ plane will disappear. However,
a new one corresponding to the green 2-handle running from $G$ to $H'$ will appear, and one corresponding to the brown 2-handle running from $G'$ to $H$ will appear. 
The blue 2-handle running between $B$ and $B'$ (still staying in the $x-y$ plane)
also corresponds to an intersection point in the $y-z$ plane, when we cancel $B-B'$ this 2-handle forms a closed loop with framing 0. It then
gives a point of intersection in the $y-z$ plane, however since the framing of this 2-handle is 0 it cancels a 3-handle, and so we can simply delete it from
our diagram. This means that we can simply delete the corresponding point of intersection from the $y-z$ plane. Similarly, when we cancelled $B-B'$ from
the diagram that consisted of the 2-handles that did not all lie in a single 2-plane we obtained two 2-handles, one in blue and the other in green which
are loops with zero framing. Hence they each cancel a 3-handle respectively and can be deleted from the diagram. This means that the points of intersection
they give in the $y-z$ plane can be deleted.

\centerline{\graphicspath{ {handle_cancellation/cancelling_B-B'/} }\includegraphics[width=6cm, height=6cm]{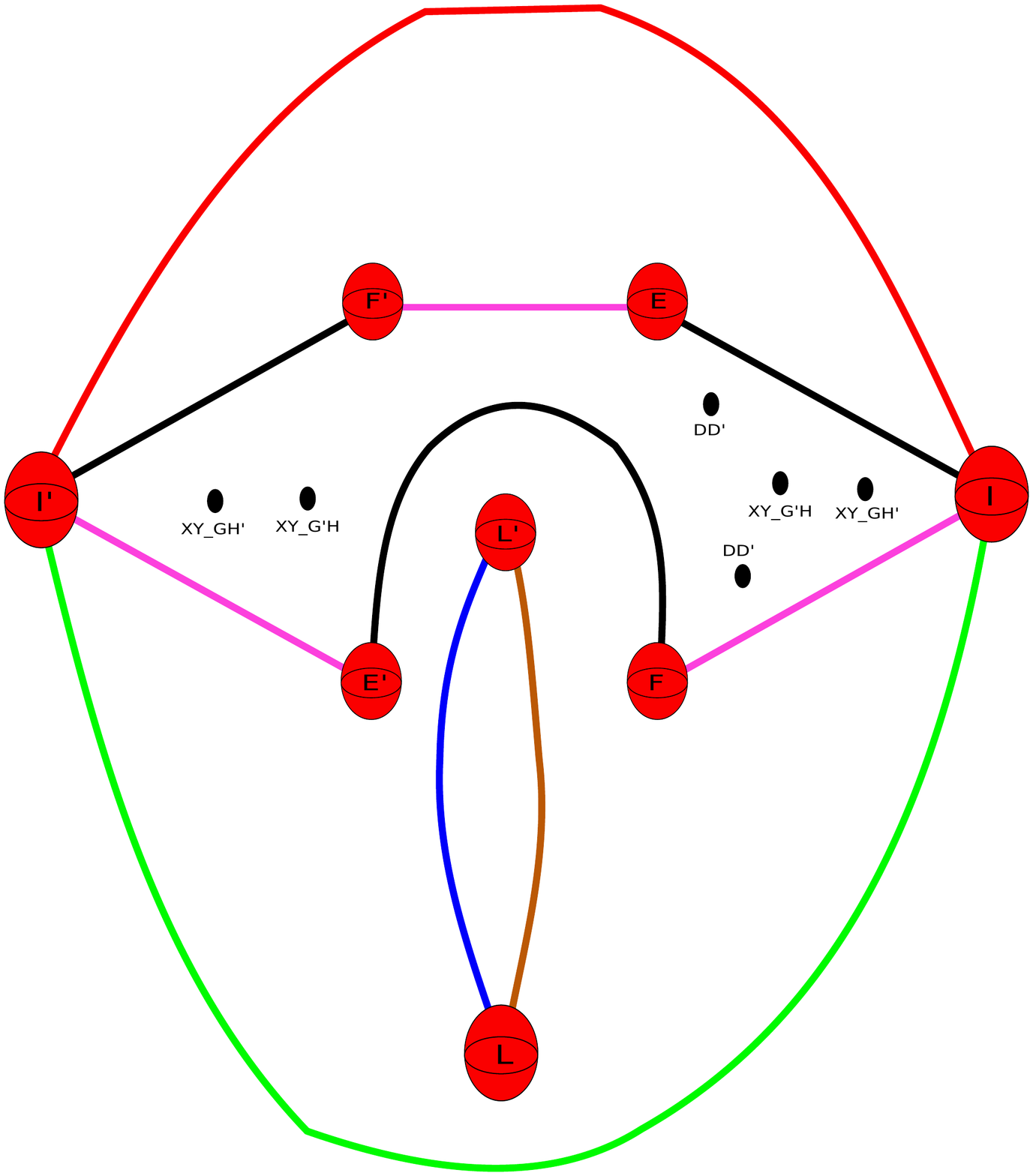}}

Observe that there are two points of intersection labelled \textbf{XY\underline{\space}GH'}, and another two labelled \textbf{XY\underline{\space}G'H}.
It is easy to work out which one corresponds to which 2-handle
in the $x-y$ plane. For example, the point of intersection between $E'$, $I'$ and $F'$ labelled \textbf{XY\underline{\space}GH'} corresponds to the green 2-handle running between $G$ and $H'$ in the $x-y$ plane, and
the point of intersection labelled \textbf{XY\underline{\space}GH'} between $E$, $I$ and $F$ corresponds to the brown 2-handle running between $G$ and $H'$ in the
$x-y$ plane.

Before we show how all diagrams look like after the cancellation of $B-B'$, we want to go back to the diagram that corresponded to the six 2-handles
that did not all lie in a single 2-plane.

\centerline{\graphicspath{ {handle_cancellation/cancelling_B-B'/} }\includegraphics[width=5cm, height=4.5cm]{x-y-z}}

Look at the brown 2-handle that runs between $D$ and $D'$. We were going to do an isotopy that moves this 2-handle into a different position. When carrying out
such an isotopy we have to be careful that we do not pass through any other 2-handles. Therefore in order to do this in a correct manner we need to
keep track of the intersection points with the other planes. 

We are going to push the brown 2-handle in the following direction:

\centerline{\graphicspath{ {handle_cancellation/cancelling_B-B'/} }\includegraphics[width=5cm, height=4.5cm]{x-y-z_iso}}

This will cause the brown 2-handle (over time) to move in the following way:

\centerline{\graphicspath{ {handle_cancellation/cancelling_B-B'/} }\includegraphics[width=8cm, height=7cm]{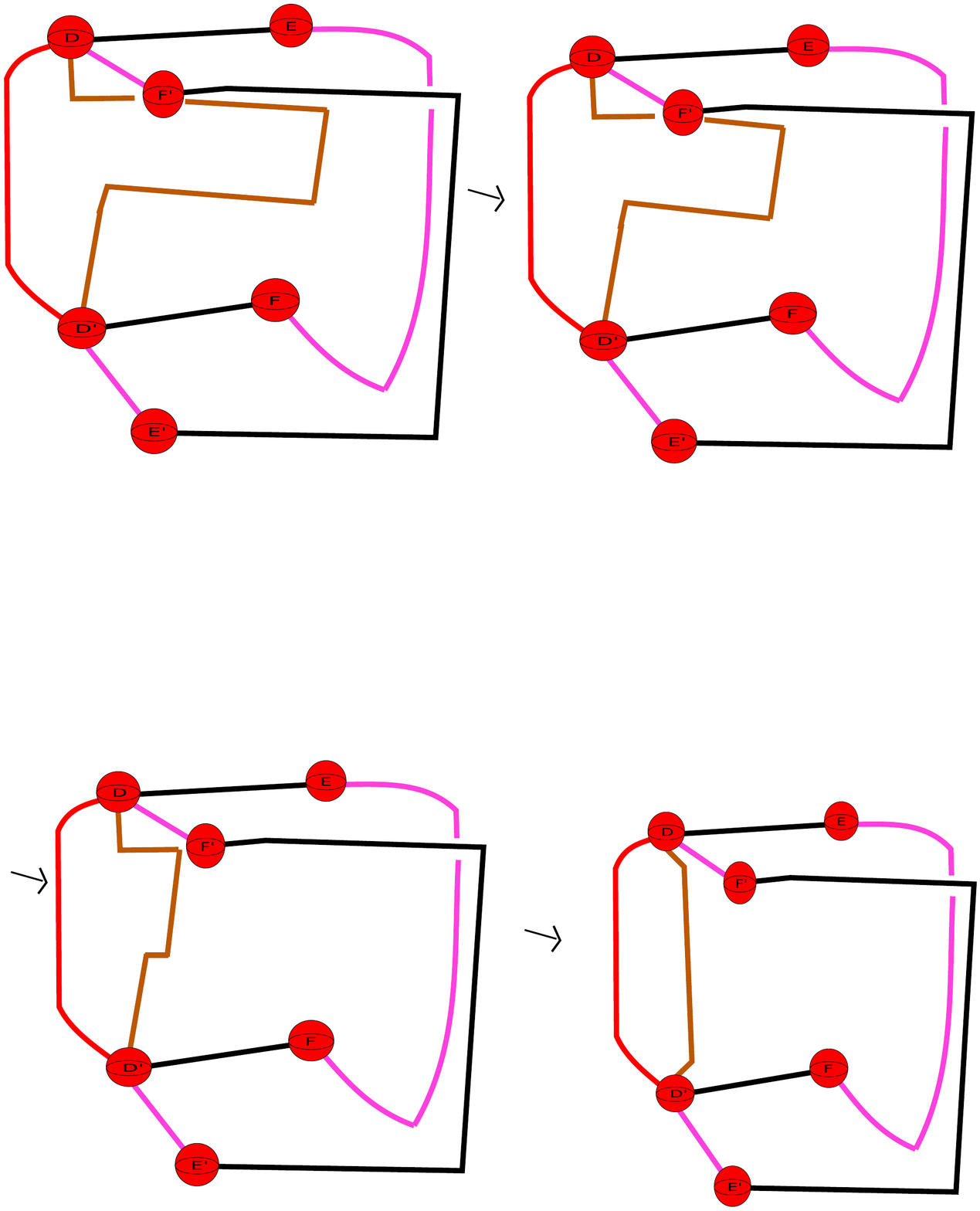}}

During this moving of the brown 2-handle the intersection point corresponding to this 2-handle in the
$x-y$ plane also moves, and the trajectory it takes looks like:

\centerline{\graphicspath{ {handle_cancellation/cancelling_B-B'/} }\includegraphics[width=6cm, height=5.5cm]{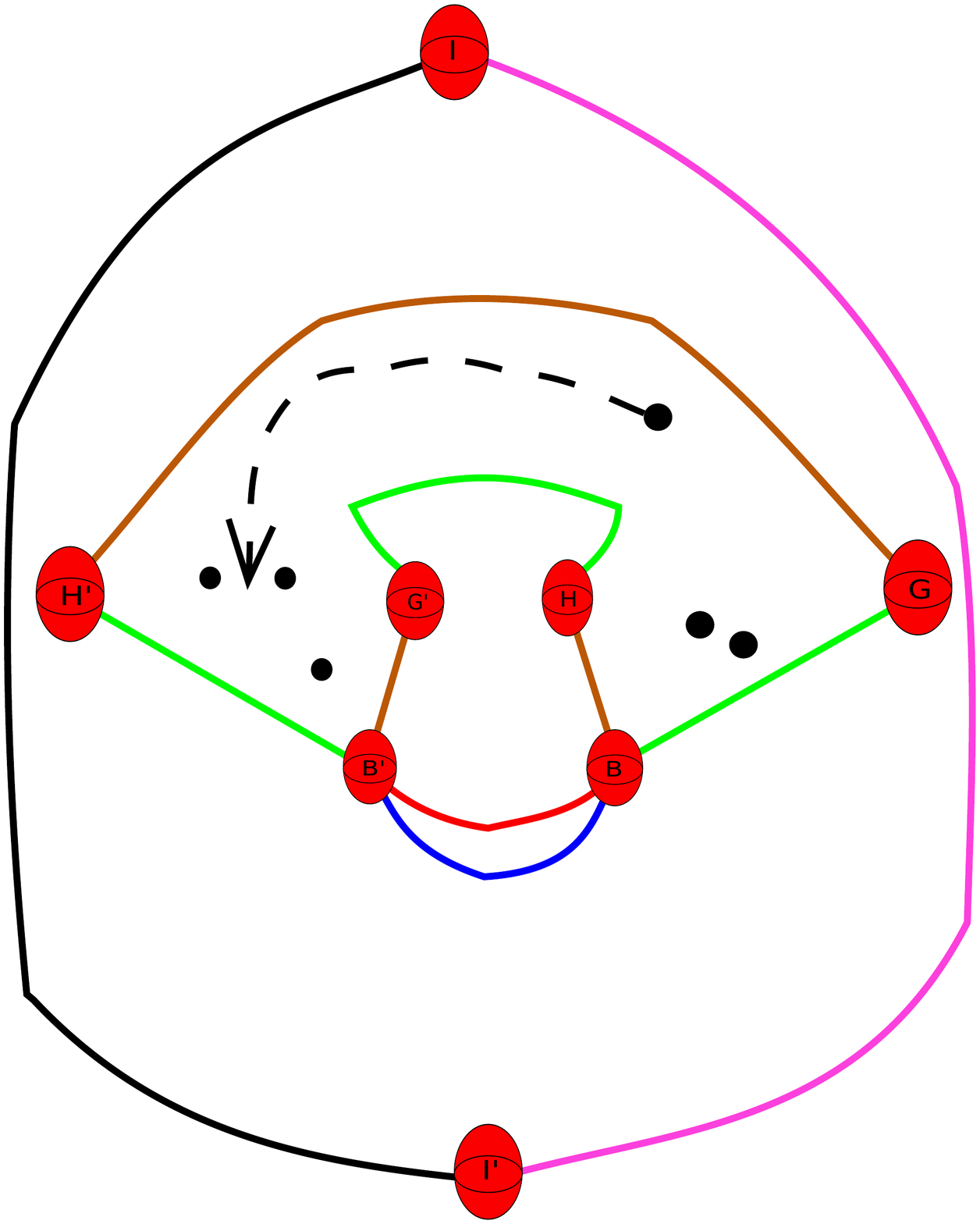}}

Observe that so far, in moving the brown 2-handle we have not passed through any of the other 2-handles, this is because in the $x-y$ plane we have kept track
of how the intersection point has moved. In the $x-z$ plane there is nothing to check as the brown 2-handle that we are moving does not intersect this plane. Finally, 
in the case of the $y-z$ plane we have that this 2-handle intersects it in two points, thus in carrying out this isotopy of the 2-handle in question there
will definitely be some change in the $y-z$ plane.
In order to understand what goes on observe that when we continue to move the 2-handle in the direction shown above we will come to a stage where a vertical arc of this 2-handle will lie in the $y-z$ plane. Thus when we look at the $y-z$ plane we will see an arc between the two points of intersection corresponding to this brown 2-handle. The following picture shows this arc in the $y-z$ plane.

\centerline{\graphicspath{ {handle_cancellation/cancelling_B-B'/} }\includegraphics[width=6cm, height=6cm]{y-z_iso}}

As we continue to move this 2-handle in the direction indicated we will see this arc in the $y-z$ plane disappear along with the original two points of intersection
that this 2-handle gave.

The final position of this 2-handle is shown in the following picture, with the bottom diagram showing the final position of the intersection point in the $x-y$ plane.

\centerline{\graphicspath{ {handle_cancellation/cancelling_B-B'/} }\includegraphics[width=7cm, height=7cm]{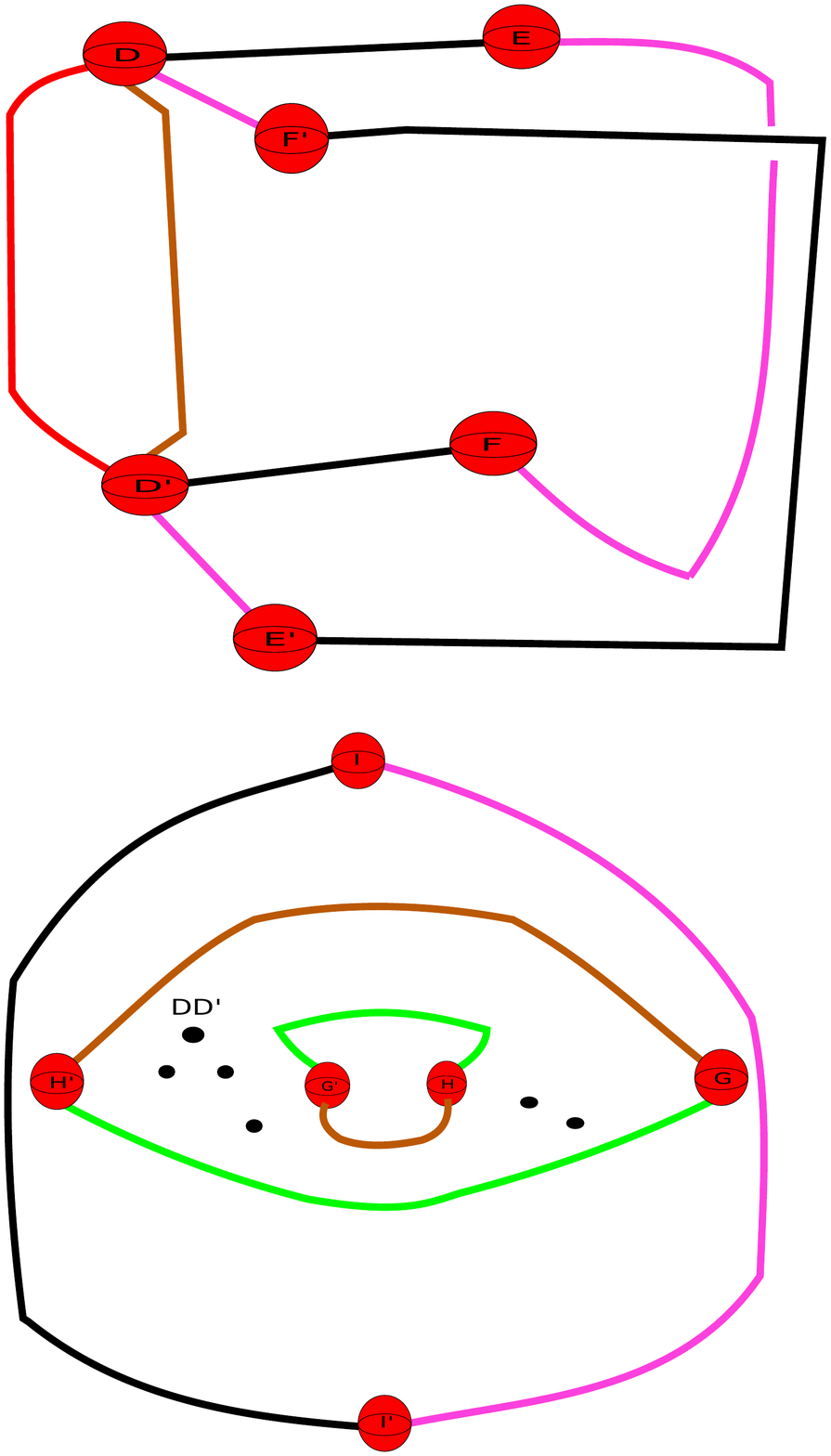}}

The following picture shows all diagrams so far, notice that two points of intersection in the $y-z$ plane have disappeared.

\centerline{\graphicspath{ {handle_cancellation/cancelling_B-B'/} }\includegraphics[width=10cm, height=10cm]{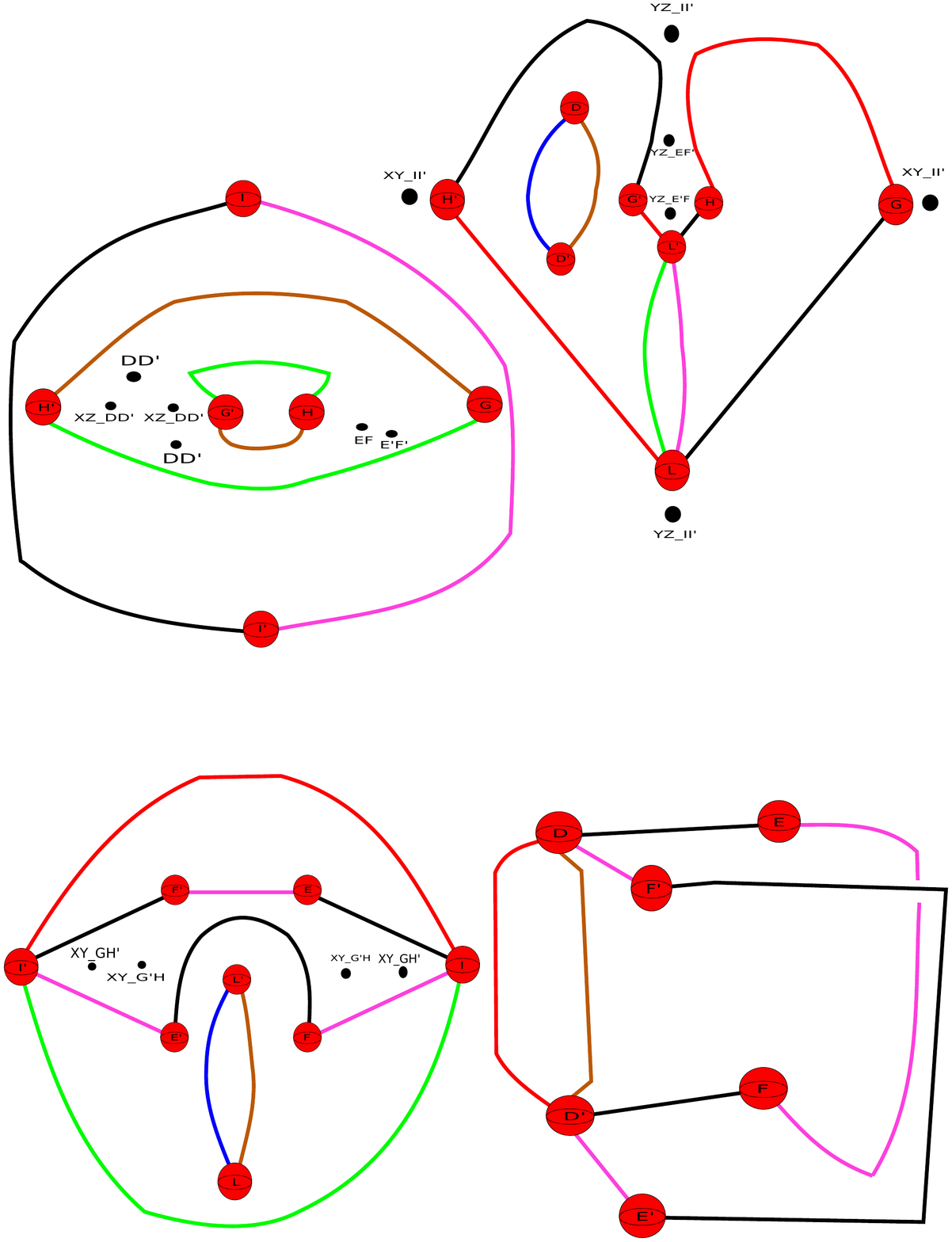}}

We are going to move on to cancelling $D-D'$. In this case we have a few options in regard to which 2-handle we use to carry out the cancellation. In the
$x-z$ plane we have the choice of the brown and blue 2-handles that pass between $D-D'$ once, and in the diagram that corresponded to the six 2-handles 
that did not all lie in a single 2-plane, we have the choice of the red and brown 2-handle that passes between $D-D'$ once. We will choose the red 2-handle
in the diagram that corresponds to the six 2-handles that did not all lie in a single 2-plane.

\centerline{\graphicspath{ {handle_cancellation/cancelling_D-D'/} }\includegraphics[width=9cm, height=8.5cm]{x-z_x-y-z}}

Observe that after cancelling $D-D'$ all the 2-handles that passed between $D-D'$ become loops with framing 0. We can see two in the $x-z$ plane, one
in blue and the other in brown, and we can see one brown one in the other diagram. As all these two handles have framing 0 we know that they form a cancellation
pair with a 3-handle, and hence we can erase them from our diagram. In doing so we will see the intersection points in the $x-y$ plane labelled 
\textbf{XZ\underline{\space}DD'} and \textbf{DD'} disappear, but we will get two new ones corresponding to the pink 2-handle running from $E'$ to $F'$ and
the black 2-handle running from $E$ to $F$, both in the diagram that corresponds to the six 2-handles that did not all lie in a single 2-plane.

\centerline{\graphicspath{ {handle_cancellation/cancelling_D-D'/} }\includegraphics[width=6cm, height=5cm]{x-y_int}}

All diagrams look like:

\centerline{\graphicspath{ {handle_cancellation/cancelling_D-D'/} }\includegraphics[width=8cm, height=8cm]{all}}

Let us remark that the cancellations and slides we have done so far have not interfered in any way with the two added 2-handles that run from $E$ to $G$ and
$E'$ to $G'$.

We move on to cancelling $L-L'$ using the green 2-handle in the $x-z$ plane. This cancellation only affects the diagrams in the $x-z$ and $y-z$ planes.

\centerline{\graphicspath{ {handle_cancellation/cancelling_L-L'/} }\includegraphics[width=6cm, height=7cm]{x-z_y-z}}

In undertaking this cancellation we obtain three 2-handles that become loops which carry framing number 0, these will then cancel a 3-handle and hence we can simply
erase them from our diagrams.

Note that as none of the 2-handles that ran between $L-L'$ intersected any of the other planes, the points of intersection in the various planes do not change.  
Furthermore it is easy to see that cancelling $L-L'$ does not cause the added 2-handles that run from $E$ to $G$ and $E'$ to $G'$ to change.

We can then cancel $I-I'$ using the red 2-handle that lies in the $y-z$ plane.

\centerline{\graphicspath{ {handle_cancellation/cancelling_I-I'/} }\includegraphics[width=8cm, height=6cm]{x-y_y-z}}

We obtain three 2-handles that are loops with framing 0, which we can simply erase from the diagrams. As for the intersection points, the only plane we have
to worry about is the $x-z$ plane. The intersection points labelled \textbf{XY\underline{\space}II'} disappear and nothing new takes their place. 
The points of intersection labelled \textbf{YZ\underline{\space}II'} also disappear, but we get two new points of intersection at the top of the $y-z$ plane
corresponding to the black 2-handle in the $y-z$ plane that runs between $E$ and $F'$ and the pink 2-handle in the $y-z$ plane that runs between 
$E'$ and $F$. The labelling of these points of intersection are \textbf{YZ\underline{\space}EF'} and \textbf{YZ\underline{\space}E'F} respectively.

The following picture shows the diagrams so far with labelled points of intersection. The first diagram represents the $x-y$ plane, the one below it
is the $x-z$ plane, and the next two diagrams are the $y-z$ plane and the diagram that corresponded to the six 2-handles that did not all lie in a single plane.

\centerline{\graphicspath{ {handle_cancellation/cancelling_I-I'/} }\includegraphics[width=9cm, height=8cm]{all}}

We can isotope the pink 2-handle in the $y-z$ plane joining $E'$ to $F$, so that it does not enclose the 1-handles $E$ and $F'$. The isotopy moves the 2-handle in
the direction of the $x$-axis. The following picture shows the direction in which we move this 2-handle.

\centerline{\graphicspath{ {handle_cancellation/cancelling_I-I'/} }\includegraphics[width=6cm, height=4cm]{y-z_iso}}

It is easy to see that in moving this 2-handle in this direction, we do not pass through any of the other 2-handles, hence it is a well defined isotopy. The
final position of this 2-handle can be seen in the following picture.

\centerline{\graphicspath{ {handle_cancellation/cancelling_I-I'/} }\includegraphics[width=5cm, height=4.5cm]{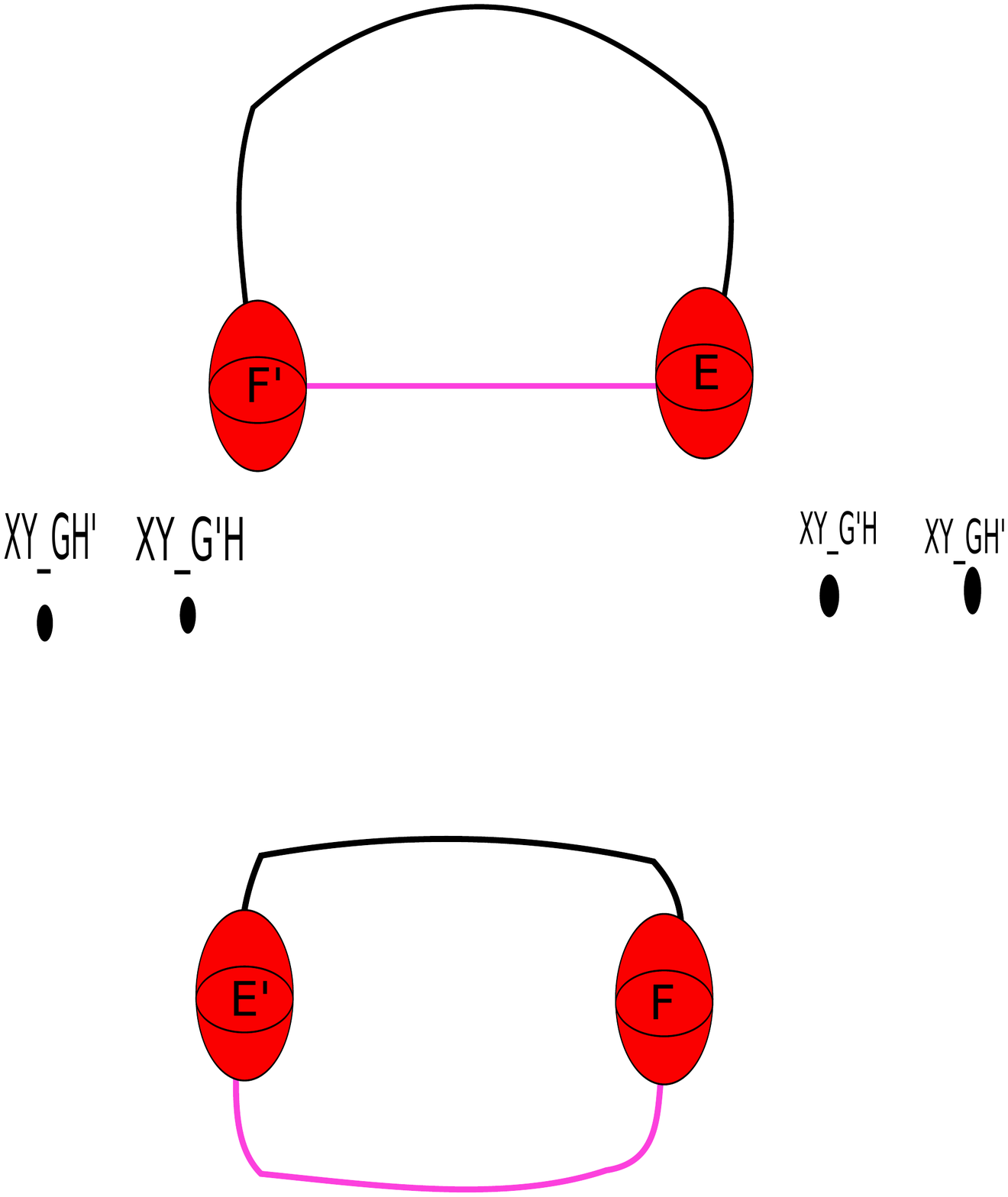}}

We can also isotope the 2-handles in the diagram that corresponded to the six 2-handles that did not all lie in a single plane to get:

\centerline{\graphicspath{ {handle_cancellation/cancelling_I-I'/} }\includegraphics[width=5cm, height=4cm]{x-y-z_iso}}

The 2-handles now lie in planes parallel to the $x-z$ plane.

We observe that this will cause the points of intersection in the $x-y$ plane to move into the following position:

\centerline{\graphicspath{ {handle_cancellation/cancelling_I-I'/} }\includegraphics[width=4cm, height=4cm]{x-y_int}}

We can now put everything together to give a picture in 3-space. Remember the co-ordinates of the 1-handles $H$, $H'$, $G$, $G'$ and $E$, $E'$, $F$, $F'$
are given as: \\

\begin{tabular}{|l|l|l||l|l|l|}
$E$ &  $S_{(0,+1,+1,0)}$ &  $(0, \frac{1}{\sqrt{2}} ,\frac{1}{\sqrt{2}})$ &  $E'$ & $S_{(0,-1,-1,0)}$ &  $(0, \frac{-1}{\sqrt{2}} ,\frac{-1}{\sqrt{2}})$ \\ 
 $F$ &  $S_{(0,+1,-1,0)}$ &   $(0, \frac{1}{\sqrt{2}} ,\frac{-1}{\sqrt{2}})$ &  $F'$ & $S_{(0,-1,+1,0)}$ &  $(0, \frac{-1}{\sqrt{2}} ,\frac{1}{\sqrt{2}})$ \\ 
 $G$ &  $S_{(+1,0,0,+1)}$ &  $(1 + \sqrt{2}, 0, 0)$ &  $G'$ & $S_{(-1,0,0,-1)}$ & $(1 - \sqrt{2}, 0, 0)$ \\ 
 $H$ &  $S_{(+1,0,0,-1)}$ &  $(-1 + \sqrt{2}, 0, 0)$ &  $H'$ & $S_{(-1,0,0,+1)}$ & $(-1 - \sqrt{2}, 0, 0)$ \\
\end{tabular} \\

The following picture shows a diagram of these 1-handles in 3-space along with the various 2-handles that run between them.

\centerline{\graphicspath{ {handle_cancellation/putting_together_2/} }\includegraphics[width=6cm, height=6cm]{x-y_y-z}}

Observe that in the top diagram the 2-handle denoted by $e$ passes over the 1-handle $H-H'$ once, hence they form a cancelling pair, and the cancellation
can be done in the $x-y$ plane without affecting any of the other 2-handles.

\centerline{\graphicspath{ {handle_cancellation/cancelling_H-H'/} }\includegraphics[width=7cm, height=6cm]{x-y}}

The 2-handles denoted by $f$ then cancels a 3-handle and can be erased from the diagram.

In our 3-space picture we then have:

\centerline{\graphicspath{ {handle_cancellation/cancelling_H-H'/} }\includegraphics[width=7cm, height=5cm]{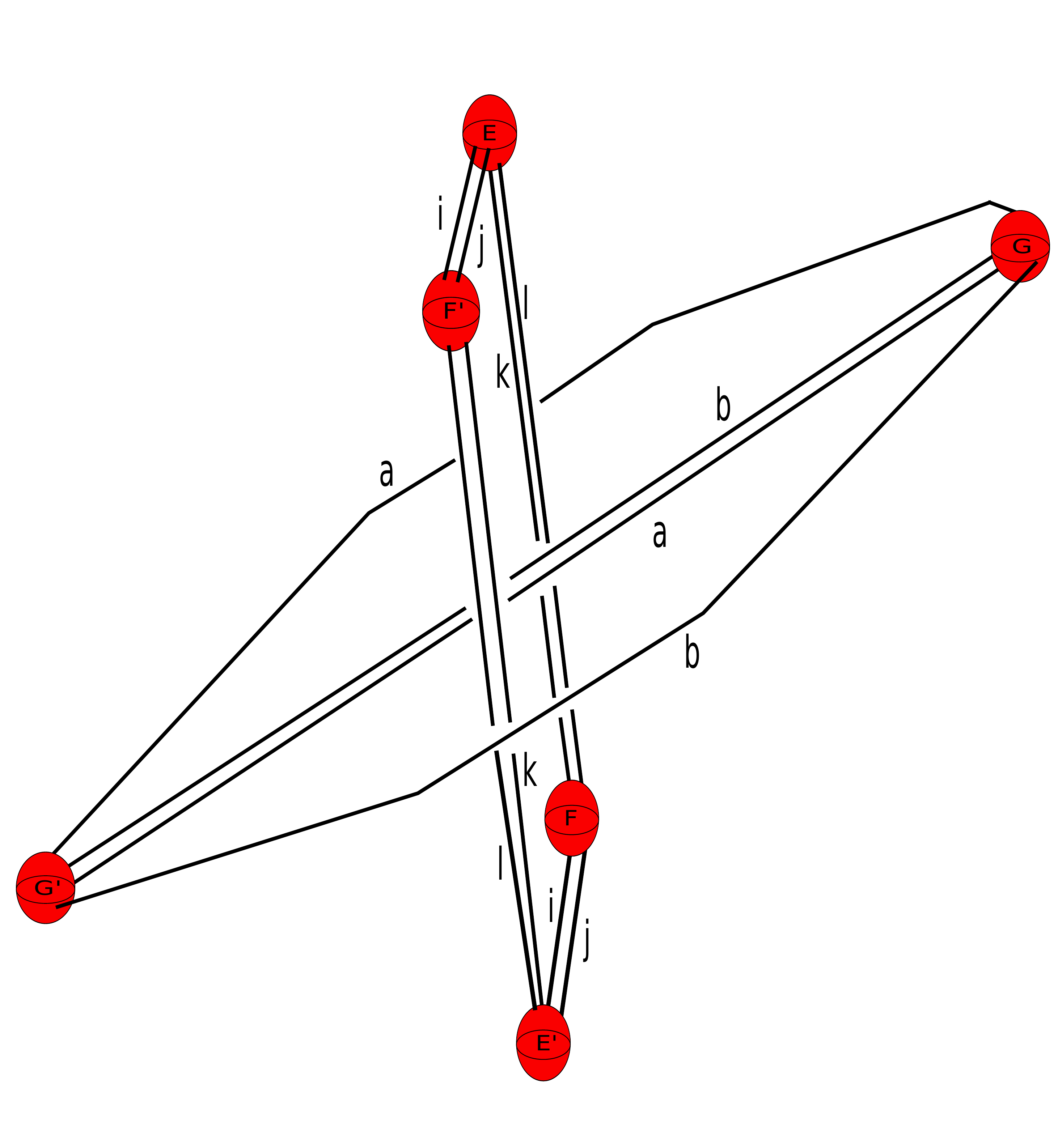}}

The 2-handle denoted $k$ runs over $F-F'$ once, and hence we can use it to cancel $F-F'$. This cancellation can be done in the $y-z$ plane without
affecting the other 2-handles.

\centerline{\graphicspath{ {handle_cancellation/cancelling_F-F'/} }\includegraphics[width=7cm, height=5cm]{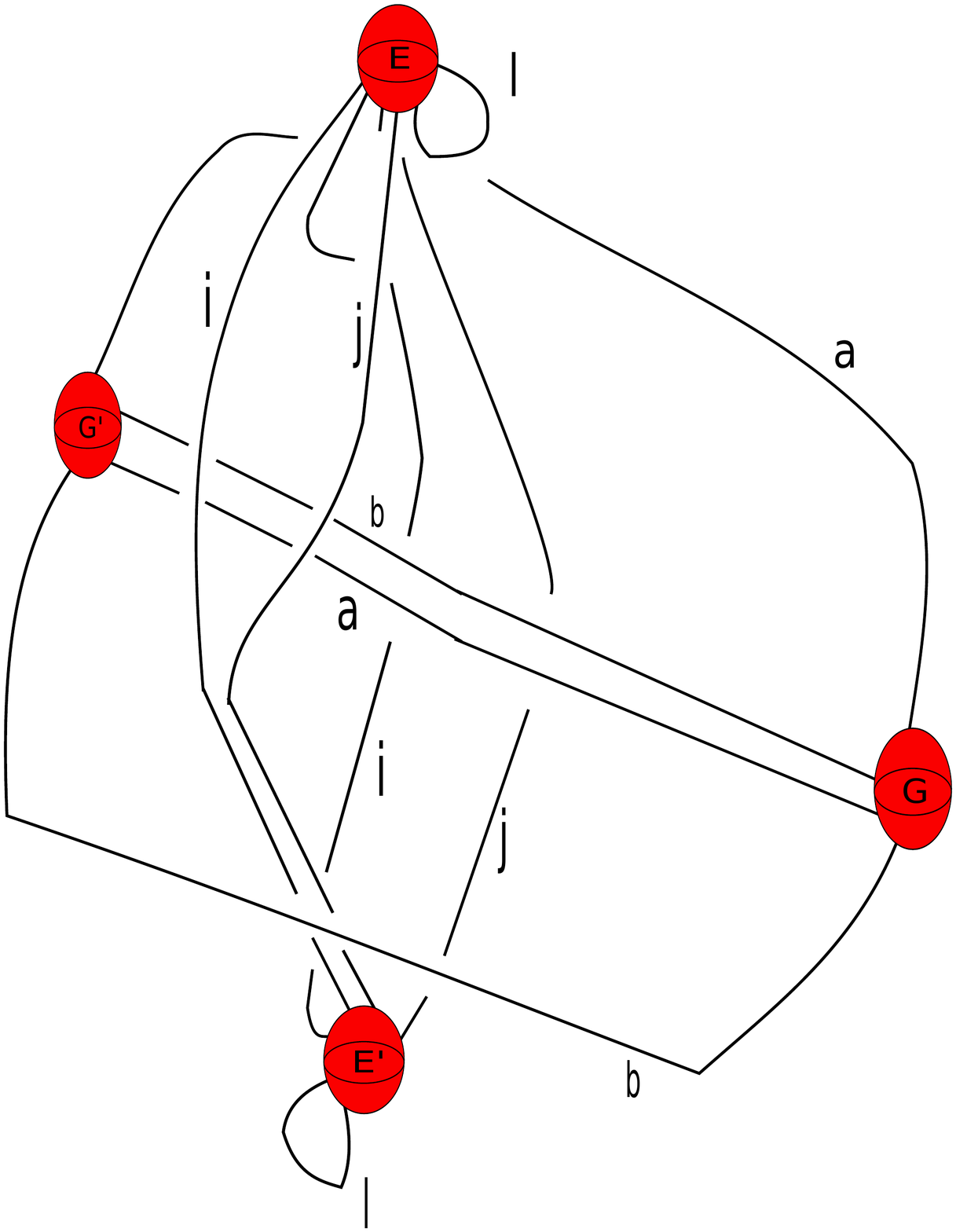}}

The 2-handle $l$ slides off $E-E'$ and cancels a 3-handle, we are then left with:

\centerline{\graphicspath{ {handle_cancellation/cancelling_F-F'/} }\includegraphics[width=7cm, height=6cm]{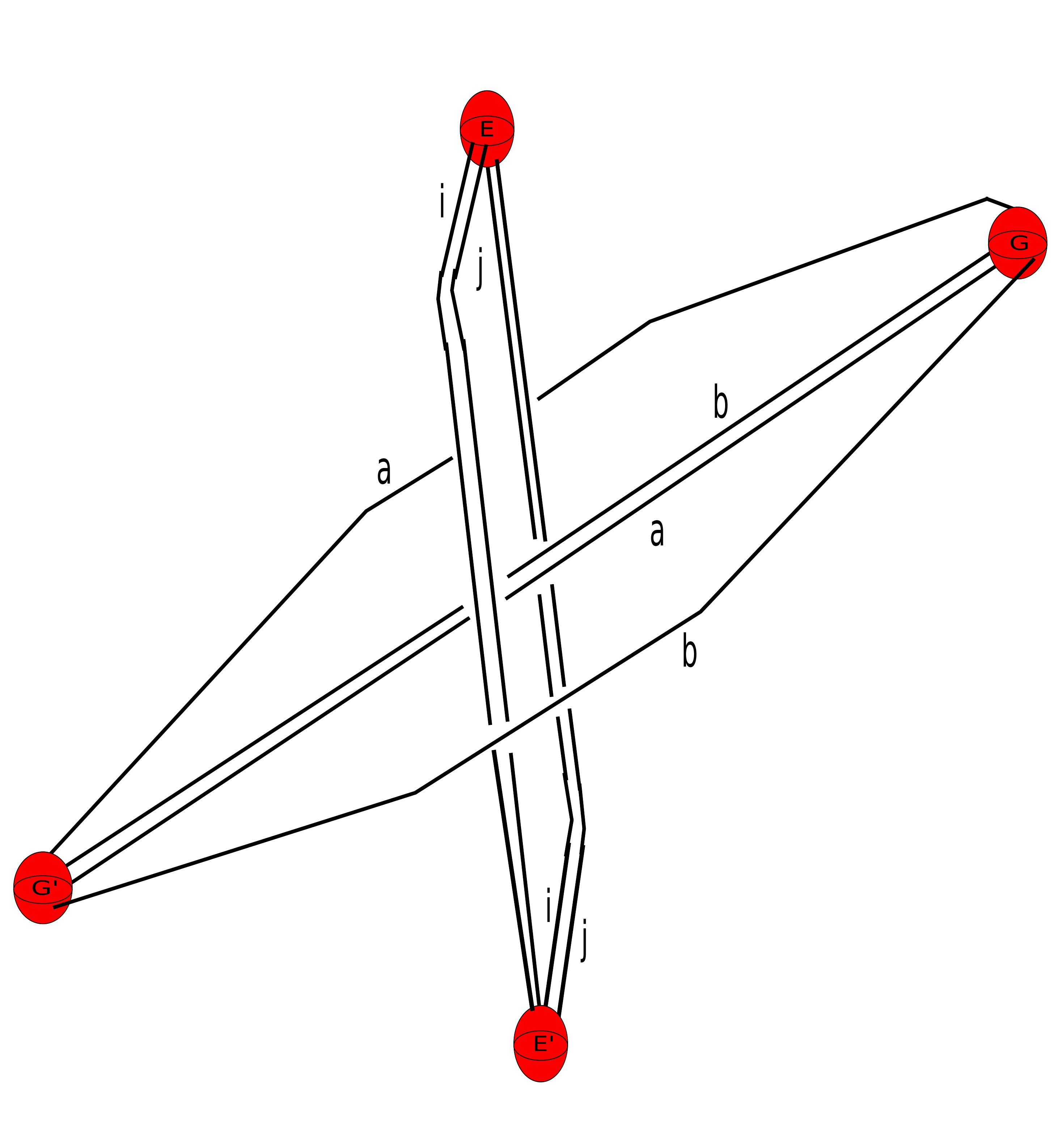}}

We can then slide the 2-handle denoted $j$ along $i$ to obtain a 2-handle that has one component looping back into $E$ and another into $E'$. We can push
either component through the attaching sphere it loops back into obtaining an unknot with framing zero. This then cancels a 3-handle and can be deleted
from the diagram.

Recall we still have not dealt with the 2-handle that corresponds to filling in the last boundary component. This is represented by
$e^{-1}g$, which corresponds to attaching a 2-handle component from $E$ to $G$, and then one from $E'$ to $G'$. So far the handle cancellations and slides
we have done have not affected these two components, and hence it suffices to put them in now.

\centerline{\graphicspath{ {handle_cancellation/adding_e=g/} }\includegraphics[width=7cm, height=5cm]{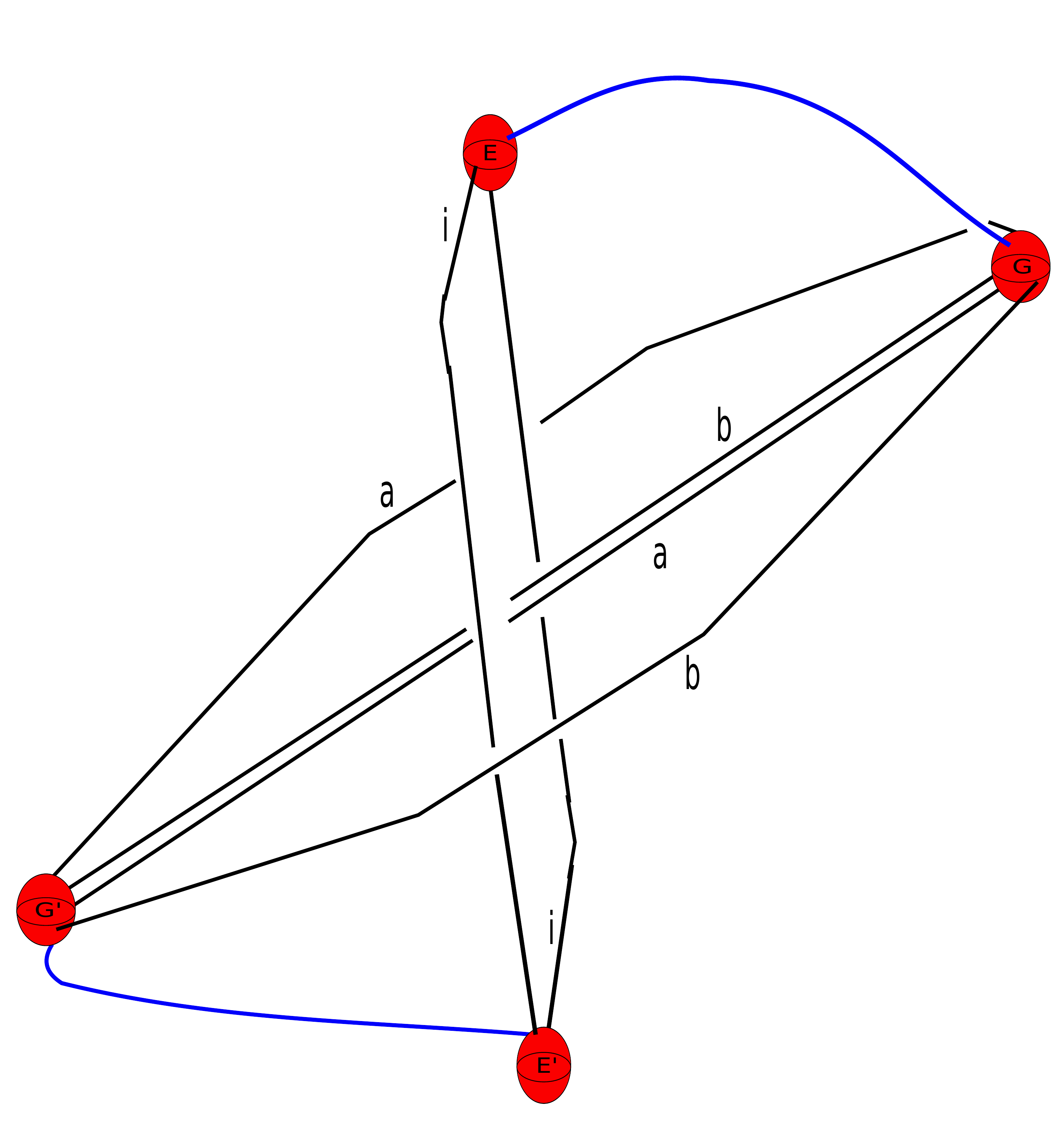}}

We can then use them to cancel $G-G'$:

\centerline{\graphicspath{ {handle_cancellation/cancelling_G-G'/} }\includegraphics[width=6cm, height=5cm]{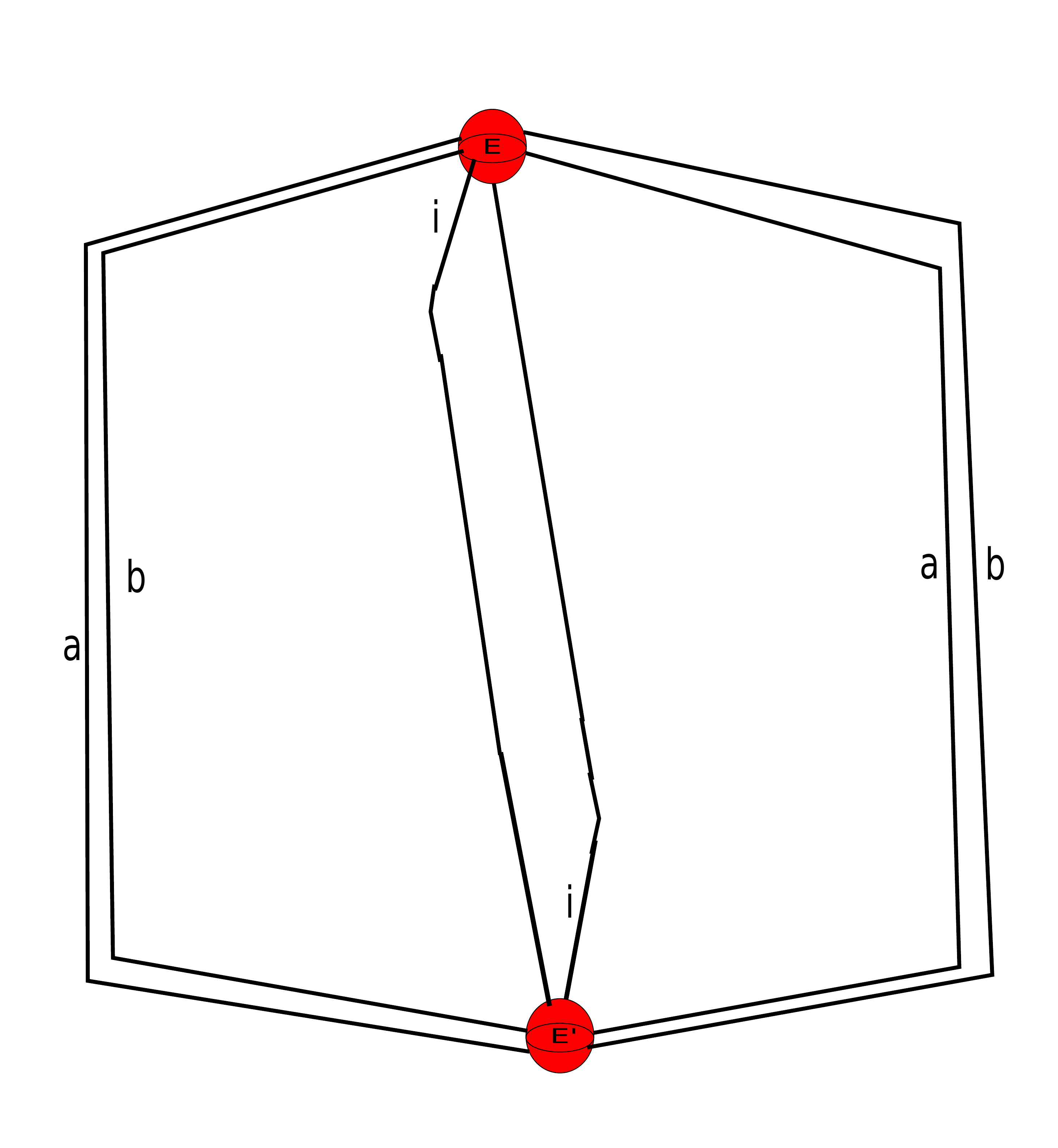}}

At this point one can perform a handle slide, one can slide the 2-handle labelled $b$ along $a$ the end result is that $b$ will have two components, one
that loops back into $E$ and one that loops back into $E'$. We can then push one component through the attaching sphere it loops back into obtaining an unknot
with zero framing. This cancels with a 3-handle and can be erased from the diagram.
We are then left with the same diagram as above with only the 2-handles
labelled $a$ and $i$ remaining.

The fundamental group of this manifold is clearly $\Z_2$. We want to take the double cover of the above 4-manifold to get a simply connected 4-manifold. Let us denote the above 4-manifold by $Y$, then we want to
form the double cover $X \rightarrow Y$. More precisely using the Kirby diagram of $Y$ obtained above, we want to obtain a Kirby diagram of $X$.
In order to do this we proceed as follows. First of all since $Y$ has only one 1-handle we can identity the 1-skeleton of $Y$ with $S^1 \times D^3$. As
$X$ is a cyclic double cover we have that the 1-skeleton of $Y$ is 2-fold covered in the obvious way. I.e. the 1-skeleton of $X$ is
$S^1 \times D^3$ and this double covers $S^1 \times D^3$ (corresponding to the 1-skeleton of $Y$) in the usual way. Each of the remaining handles of $Y$ lift
to two copies of that handle in $X$. There is a subtle issue in what we have said so far, namely we did not say how the data determining the trivialisation
of the normal bundle to each 2-handle lift to the double cover.
In our case this
is not a problem as all our 2-handles have a planar framing, and so when we lift them we still get 2-handles with a planar framing (i.e. parallel curves do not twist around). In the general case one only has
to observe that a choice of trivialisation on each normal bundle to each knot (representing a 2-handle) lifts to that trivialisation on each lifted knot, and changing the trivialisation to the normal bundle of a 2-handle in $Y$ by one twist
causes the trivialisation on each lifted knot to change by one twist, hence we know how to lift any choice of trivialisation.
Focusing on our case we make an important observation, the 2-handles in our manifold $Y$ do not twist around each other in any way, hence when we lift
them to $X$ we can focus on their lifts separately. In our situation we have two 2-handles denoted $a$ and $i$. When we lift each one, we will
get two 2-handles each passing the lifted 1-handle once. As none of them twist around each other we use any lifted one to cancel the lifted 1-handle, leaving
us with three unknot's with framing zero. Each of these then cancels a 3-handle, and we see that the manifold we are left with is $S^4$. The boundary
components of manifold 1011 where given by the non-orientable closed flat 3-manifold labelled $\textbf{G}$ ($\mathcal{B}_1$ in Wolf's notation). Appealing to the 
classification of closed
flat 3-manifolds we know that the orientable double cover of $\textbf{G}$ is the flat 3-manifold given by $\textbf{A}$ ($\mathcal{G}_1$ in Wolf's notation), which is
the 3-torus.

Thus we have proved the following theorem:

\begin{thm}\label{mainthm_1}
There exists a collection $L$ of five linked 2-tori embedded in a smooth 4-manifold $X$ such that $X$ is diffeomorphic to a standard $S^4$, and
$X-L$ admits a finite volume hyperbolic geometry.
\end{thm}

\end{document}